\documentclass[12pt]{amsart} 

\usepackage{lscape}
\usepackage{tabularx,mdwtab}

\usepackage{hyperref}

\makeindex

\usepackage{bbm}

\newcommand {\lra} {\longleftrightarrow}

\def\cal{\relax}
\usepackage{DLdef1}
\let\ssec\subsection

\renewcommand {\ssbegin}[2][*]
 {\refstepcounter{subsection}%
\if#1*
\addcontentsline{toc}{subsection}{\thesubsection.\hskip 1pc #2}%
\else
\addcontentsline{toc}{subsection}{\thesubsection.\hskip 1pc #2. #1}%
\fi
 \def \secno {\gdef \secno {}{\ssecfont
\thesubsection.\hskip 2ex}%
 }%
 \begin{#2}}

\renewcommand {\sssbegin}[2][*]
 {\refstepcounter{subsubsection}
\if#1*
\addcontentsline{toc}{subsubsection}{\thesubsubsection.\hskip 1pc #2}%
\else
\addcontentsline{toc}{subsubsection}{\thesubsubsection.\hskip 1pc #2. #1}
\fi
 \def \secno {\gdef \secno {}{\ssecfont \thesubsubsection.\hskip 2ex}%
 }%
 \begin{#2}}

\renewcommand {\parbegin}[2][*]
 {\refstepcounter{paragraph}
\if#1*
\addcontentsline{toc}{paragraph}{\theparagraph.\hskip 1pc #2}%
\else
\addcontentsline{toc}{paragraph}{\theparagraph.\hskip 1pc #2. #1}
\fi
 \def \secno {\gdef \secno {}{\ssecfont \theparagraph.\hskip 2ex}%
 }%
 \begin{#2}}

\newcommand {\fba}{{\mathfrak{ba}}}

 \def\mmat #1,#2,#3,#4,{\text{\small\arraycolsep=3pt $
\begin{pmatrix}#1&#2\\#3&#4\end{pmatrix}$}}

\newcommand{\Size}{\text{Size}}

\newcommand{\del}{\partial}
\newcommand{\LieSAlgs} {{\sf LieSAlgs}}

\makeatletter
\@addtoreset{table}{section}
\makeatother

\begin{document}

\title[Simple vectorial Lie superalgebras]{The
classification of simple complex Lie
superalgebras of~polynomial vector fields and their deformations}

\author{Dimitry Leites$^{a, *}$, Irina Shchepochkina$^{b}$}

\address{${}^{a}$Depart\-ment of Mathematics Stockholm University, Albanov\"agen 28, SE-114 19  Stockholm, Sweden; 
dimleites@gmail.com; ${}^{*}$The corresponding author\\
$^b$Independent University of Moscow, Bolshoj Vlasievsky per,
dom 11, RU-119 002 Moscow, Russia; irina@mccme.ru
}

\begin{abstract} We overview classifications of simple
infinite-dimensional complex $\mathbb{Z}$-graded Lie (super)algebras of polynomial growth,
and their deformations. A subset of such Lie (super)algebras consist of vectorial Lie (super)algebras whose elements are vector fields with polynomial, or
formal power series, or divided power coefficients. A given vectorial Lie (super)algebra with a (Weisfeiler) filtration
corresponding to a~maximal subalgebra of finite codimension is called W-filtered; the associated graded algebra is called
W-graded. Here, we correct several published results: (1) prove our old claim ``the superization of
 \'E. Cartan's problem (classify primitive Lie algebras) is wild", 
(2) solve a tame problem: classify simple W-graded and W-filtered vectorial Lie superalgebras, (3) describe the supervariety of deformation parameters for the serial W-graded simple vectorial superalgebras, (4) conjecture that the exceptional simple vectorial superalgebras are rigid.
We conjecture usefulness of our method in  classification of simple infinite-dimensional vectorial Lie (super)alge\-bras over fields of positive characteristic. 
\end{abstract}

\makeatletter
\@namedef{subjclassname@2020}{\textup{2020} Mathematics Subject Classification}
\makeatother
\subjclass[2020]{Primary 17B20, 17B66; Secondary 58A30}

\keywords{Lie superalgebra, Weisfeiler filtration, deformation of
Lie superalgebra, modular Lie superalgebra.}

\maketitle

\markboth{\itshape Dimitry Leites\textup{,} Irina
Shchepochkina}{{\itshape The classification of simple Lie
superalgebras of~polynomial vector fields over $\Cee$}}

\thispagestyle{empty}

\begin{quote}\rightline{In memory of A. N. Tyurin.}\end{quote}


\thispagestyle{empty}
\setcounter{tocdepth}{3}
\tableofcontents

\printindex

\section{Setting of the problem. Steps of our classification. Our results}\label{S:1}

Several times at least one of the classification problems we consider here was claimed to be solved, starting with \cite{K2}. 
We analyse these
claims, correcting them when needed, and complete the results and claims announced in 
 \cite{LSh2, LSh2$'$} and preprinted in \cite{LSh1, LSh5} with a~description of
 deformations, in particular, with odd parameters, see also \cite{Ld}. Observe that to shy away from some of odd parameters is contrary to the very
spirit of supersymmetry.\footnote{An interesting question: why all people
can easily accept extensions with odd center, such as the extension of $\fle(n)$ to $\fb(n)$, but some can not accept 
odd parameters of deformations or representations?... } When the deformation parameters are even, we describe \textit{deforms} --- the
\textit{results} of deformations.\footnote{A term coined by M.~Gerstenhaber, in the same manner as \textit{transform} is the result of
a~ transformation, or, \textit{e.g.}, history is an ideological \textit{construct} (as becomes more and more evident). We like these terms and suggest the term \textit{prolong} for the result of Cartan prolongation.}\index{Deform, the result (valuation) of a~deformation} Hereafter, we abbreviate ``Lie algebra or Lie superalgebra" to ``Lie (super)algebra".

\textbf{In this paper, the ground field is $\Cee$} unless otherwise specified, the term
``\textbf{vectorial} Lie (super)algebra"\index{Lie (super)algebra, vectorial} is short for ``Lie
(super)algebra of vector fields with polynomial coefficients", or with formal-power-series coefficients if a~topological completeness is needed. 
The adjective ``vectorial'' reflects the
fact that such Lie superalgebra $\cL$ is realized by
vector fields on the \textit{linear supermanifold}, \textit{i.e.}, a~ ringed space, corresponding (as described in formula~ \eqref{linSman}; for details, see \cite{KLLS, Ld}) to the 
\textit{superspace}, \textit{i.e.}, a~$\Zee/2$-graded space, $(\cL/\cL_{0})^*$, where $\cL_{0}\subset \cL$ is
a~maximal Lie subsuperalgebra of finite codimension. Let $\dim \cL=\infty$ unless otherwise stated.

\subsection{Weisfeiler filtrations and gradings}\label{WeisF} Let $\cL_{0}$ be a~maximal subalgebra of finite codimension, let $\cL_{-1}$ be a~minimal $\cL_{0}$-invariant subspace strictly
containing $\cL_{0}$; for $i\geq 1$, set:
\begin{equation}
\label{WeiGr}
\begin{array}{l}
\cL_{-i-1}=[\cL_{-1}, \cL_{-i}]+\cL_{-i},\\ 
\cL_i =\{D\in \cL_{i-1}\mid [D, \cL_{-1}]\subset\cL_{i-1}\}.
\end{array}
\end{equation}
We thus get a~filtration, called a~\textit{Weisfeiler filtration}, see \cite{W}:
\begin{equation}
\label{Wfilt} \cL= \cL_{-d}\supset \cL_{-d+1}\supset \dots \supset
\cL_{0}\supset \cL_{1}\supset \cdots
\end{equation}
The number $d$ in \eqref{Wfilt} is called the \textit{depth}\index{Depth of filtration or grading} of $\cL$, as well as
of the associated \textit{Weisfeiler-graded} Lie superalgebra
\be\label{grL}
 \text{$L=\mathop{\oplus}\limits_{-d\leq i}L_i$, where
$L_{i}=\cL_{i}/\cL_{i+1}$.}
\ee
Any Lie superalgebra $\cL$ \textit{together with a~Weisfeiler
filtration} will be briefly called \textit{W-filtered}; the graded algebra $L$
associated with $\cL$ will be called \textit{W-graded}.\index{Weisfeiler grading, W-grading}

For vectorial Lie
(super)algebras, the notion invariant with respect to automorphisms is that of filtration, not grading, but graded (super)algebras are easier to work with.

\ssec{A wider picture. Main Problem}\label{SS:1.2}\index{Problem, Main} \index{Growth, of a~
(super)algebra} The
\textit{growth} of a~
$\Zee$-graded algebra $A=\oplus A_i$ is defined, if $\dim A_i<\infty$ for all~ $i$, to be
\[
\text{gth}(A)=\mathop{\overline{\lim}}\limits_{n\tto\infty}\; \frac{\ln\
\dim \mathop{\oplus}\limits_{|i| \leq n}A_{i}}{\ln n}.
\]
A $\Zee$-graded algebra $A$ is said to be of \textit{polynomial growth} if
$\text{gth}(A)<\infty$, \textit{i.e.}, 
\be\label{PolGr}
 \text{$\dim \mathop{\oplus}\limits_{|i| \leq
n}A_{i}\sim n^{\text{gth}(A)}$ ~~as $n\tto\infty$.}
\ee

A \textit{filtered} Lie (super)algebra $\cdots \supset \cL_{i}\supset \cL_{i+1}\supset \cdots$ is said to be of \textit{polynomial growth} if the $\Zee$-graded Lie (super)algebra $L=\oplus L_i$, where $L_i:=\cL_i/\cL_{i+1}$, associated with $\cL$ is of polynomial growth. Vectorial Lie (super)algebras constitute a~subset in the set  of $\Zee$-graded (or filtered) Lie (super)algebras of polynomial growth (GLAPG for short), see eq.~\eqref{PolGr}, so to place vectorial Lie (super)algebras in a~natural wider picture we also review, briefly, the results and \textbf{Open problems}\index{Problem, open}  in the classification of simple Lie
(super)algebras of GLAPG type such as the following \textbf{partly classified} types 1)--4)  and their  ``\textit{relatives}"\index{Relative of a~Lie (super)algebra}, see Subsection~\ref{BrGo}:

1) Lie (super)algebras of vector fields with Laurent polynomials (or series when a~topological completeness is needed) as coefficients (a.k.a. \textit{stringy} Lie (super)algebras), their central extensions, see \cite{GLS1, KvdL}, and algebras of their derivations, see~\cite[Th.4.1]{Po2}; 

2) \textit{loop} (super)algebras of Laurent polynomials (or series) with values in a~simple Lie (super)algebra, their ``twisted" versions, see \cite{LSS, LS1, ShS}, and their non-simple ``relatives" --- \textit{double extensions}, see \cite{BLS}, especially the last two sections in it, and \cite{CCLL, HS}; 

3) Lie (super)algebras of ``matrices of complex size" type, see \cite{Lq};

4) Lie queerifications of the simple associative (super)algebras, see \cite{Lq}.

In the \textit{proof} of a~tame problem of classification of simple vectorial Lie superalgebras $\cL$,
we consider maximal subsuperalgebras $\cL_{0}$ of only \textit{finite}
codimension, and related \textit{Weisfeiler filtrations},\index{Weisfeiler filtration, W-filtration} see formulas
~\eqref{WeiGr}. However, filtrations with
$\dim(\cL/\cL_{0})=\infty$ are also useful, see Subsection~
\ref{SS:2.22}. 

We survey the current state of the classification of simple Lie
(super)algebras of polynomial growth, not only vectorial ones, and, briefly,
several related topics (continuous irreducible representations
of vectorial Lie (super)algebras (\cite{GLS3}), a~ successful attempt to see the exceptional finite-dimensional Lie algebras in a~series on Vogel's plane, see \cite{AM}, etc.).

Properties of Lie \textbf{super}algebras (even over fields of characteristic 0) resemble those of Lie algebras over fields $\Kee$ of positive characteristic, a.k.a. \textit{modular} Lie algebras.\index{Lie (super)algebra, modular} In the 1960s, Kostrikin and Shafarevich
suggested a~method for constructing simple finite-dimensional Lie
algebras over algebraically closed fields $\Kee$, and conjectured that their method produced all simple finite-dimensional Lie algebras for
$\Char \Kee:=p>7$, provided they are \textit{restricted} (for details and update, see \cite{BLLS1, BGLLS}). 

This KSh-conjecture is now proved, in a
generalized form, \textit{i.e.}, for any, not only restricted, simple Lie algebras, for $p>3$, see \cite{S, BGP, Sk}. The result was partly implicit until recently: the explicit formulas of the brackets of the deformed Hamiltonian Lie algebras can be found at the moment only in the Ph.D. thesis of S.~Kirillov in Russian, see \cite{Kir} based on \cite{Sk}. Here we provide with notions and correct formulas helping one to formulate conjectural lists of all
simple finite-dimensional Lie algebras over algebraically closed
fields of characteristic~ $p=2$ and $3$, and Lie superalgebras for $p>5$, 
see \cite{BLLS1, BGLLS} and references therein. 

\sssec{A broader goal}\label{BrGo} 
Starting with a~simple Lie
algebra of GLAPG type, list the results of iterated combinations
of the following operations and various gradings of these results:

a) $ \fg \longrightarrow \fe(\fg)$, the nontrivial central
extension of $\fg$. Examples: central extensions of loop algebras; the
Poisson algebra.

b) $ \fg \longrightarrow \fder~(\fg)$, the derivation algebra of
$\fg$, or a~double extension of $\fg$. Examples: affine Kac--Moody
algebras, see \cite{BLS}.

c) Deformations of $\fg$. (Although they sometimes lead out of the GLAPG type, such deformed Lie algebras are often no less important in
applications than the original algebras. Important examples: the result of the
deformation of the algebra of Hamiltonian vector fields induced by
the \textit{quantization} of the Poisson algebra. An example due to P.~Golod (Holod in modern
Ukrainian) \cite{GoHo1, GoHo2} evolved into the Krichever-Novikov theory of certain deformations
of loop algebras; for a~review, see \cite{Shei} and Open problem in Subsection~\ref{SSS:1.12.1}.)

d) All filtered completions of $\fg$ introducing topology to the
originally algebraic problem.

e) Forms of the Lie algebras of the above types a)-- d) over
non-closed fields.


We will call the algebras obtained from a~given simple Lie algebra
$\fg$ by the above procedures a)--e), most frequently, by just a) -- d), \textit{relatives} of each
other, and of $\fg$.\index{Relative of a~Lie (super)algebra}

\ssec{Problems we solve here over $\Cee$}\label{ssPWS} Our results (\cite{ALSh}, the key ingredients and ideas, see \cite{Sh5, Sh14}, and \cite{LSh3, LSh5} announced in \cite{LSh1,LSh2, LSh2$'$}), and
later papers \cite{Kga,Kbj,FaKa,CK,CK1,CK1a,CK2,CCK,CaKa1,
CaKa2,CaKa3} are devoted to exposition of
solutions of the following problems (cf. \cite{Klc,K2, K3}): 
\begin{equation}\label{2prob} 
\begin{minipage}[l]{14cm}
1) Classify simple W-graded vectorial Lie
superalgebras.\\
2) Classify deformations of simple W-graded vectorial Lie superalgebras.\end{minipage}
\end{equation}

There still remained points to be
clarified and corrected, \textit{e.g.},
see Section~\ref{Smb}. In particular, we recall the definition of Lie
superalgebras which shows how to interpret odd
parameters of deformation, see Subsection~ \ref{SS:2.2} and \cite{Ld}.

\textbf{Pivotal ideas of our classification of simple W-graded vectorial Lie superalgebras $\fg:=\oplus_{i\geq -d}\ \fg_i$ and their deformations}:
\begin{equation}\label{mainId} 
\begin{minipage}[l]{13.3cm}
\phantom{.}\hskip0.2cm 1) If $\dim \fg=\infty$, then among the operators in $\fg_0$ acting in $\fg_{-1}$, there is a~rank-$1$ operator with an even covector, except  for the $\Zee$-gradings of $\fg$ \textit{compatible}\footnote{Starting with \cite{K1.5},  instead of ``compatible" an inappropriate term ``consistent" is sometimes used.} with the $\Zee/2$-grading by parity, \textit{i.e.}, when $\fg_\ev=\oplus\, \fg_{2i}$ and $\fg_\od=\oplus\, \fg_{2i+1}$. \index{Grading, compatible}(The exceptional cases yield a \textit{conjectural} part of our classification of W-gradings.) \\
\phantom{.}\hskip0.2cm 2) Every simple vectorial Lie superalgebra has finitely many W-gradings.\\
\phantom{.}\hskip0.2cm 3) Deformations of any Lie superalgebra are parameterized by a~supervariety (degenerated to a~point if the algebra is rigid).
\end{minipage}
\end{equation}

\textbf{Conjecture}.\index{Conjecture} Our main method to obtain examples of simple vectorial Lie superalgebras --- \textbf{search for irreducible $\fg_0$ modules $\fg_{-1}$ with a~rank-1 operator} --- will be useful for obtaining new examples of simple W-graded \textbf{infinite-dimensional} vectorial Lie (super)algebras $\fg$ in characteristics~ $p=2$ and $3$; we conjecture that this classification will be easier to obtain than the classification of \textbf{finite-dimensional} simple Lie (super)algebras, cf.~\cite{BGLLS,CSS}.

\sssec{Classification of simple vectorial Lie algebras over $\Cee$. Their deformations} The well-known answer is as follows
\begin{equation}\label{summ} 
\begin{minipage}[l]{14cm}
\textbf{There are 4 series of simple W-graded vectorial Lie
algebras\footnote{This part is due to \'E.~Cartan, for an exposition, see
\cite{GQS}.}; they are rigid, see Section~\ref{S:11}, except for the Hamiltonian Lie
algebra\footnote{This part was a~kind of folklore, known since
1940s or 1950s; for details, see \cite{LSh3, KT}.} $\fh(2n)$ to which one can restrict the single non-trivial deformation\footnote{The physicists call this $\cQ$  \textit{quantization}, it leads both $\fpo(2n)$ and $\fh(2n)$ out of the class of W-graded or W-filtered Lie algebras: the deform $\cQ(\fpo(2n))$ is the Lie algebra of differential operators with polynomial coefficients and the commutators as products, see \cite{Lq}.} $\cQ$ of its central extension --- the Poisson Lie algebra $\fpo(2n)$.
}
\end{minipage}
\end{equation}

\ssec{On setting of the classification problems}\label{SS:1.1} Selection of Lie
(super)algebras with reasonably nice properties is a~matter of taste and is
influenced by the governing ``meta''-problem. One of the usual
choices is classification of \textbf{simple} Lie (super)algebras: they have
a~richer structure than other types of Lie (super)algebras, and therefore
are easier to study; they also illuminate important symmetries. Even
representatives of the ``complementary'' type --- solvable Lie (super)algebras
--- are often of significant interest only when 
considered as subalgebras of simple Lie (super)algebras.

$\bullet$ The \textbf{finite-dimensional} simple Lie algebras are naturally the first
to study. Note that this, ostensibly natural, choice is, however, unnatural in positive characteristic. \textbf{Conjecture}:\index{Conjecture} For $p>0$, to classify \textbf{infinite-dimensional} simple vectorial Lie (super)algebras will be easier than to classify finite-dimensional Lie (super)algebras and our pivotal idea 1) in \eqref{mainId}, will help.

\textbf{Over $\Cee$}, there are three ``classical" series
($\fsl$, $\fo$ and $\fsp$) and five exceptional (ostensibly not belonging to series) algebras; all these
simple Lie algebras are neatly encoded by (normalized) Cartan
matrices which are symbolically represented by very simple graphs --- Dynkin
diagrams, see \cite{Bbk, OV}. 

Killing and \'E.~Cartan, who discovered all 5 exceptional simple Lie algebras, described them as preserving a
non-integrable distribution, see \cite{Cart, Co}. (Actually, Killing and \'E.~Cartan considered serial simple Lie algebras in these terms as well. Later, with the discovery of the powerful method of studying the structure of simple Lie algebras and their representations in terms of roots and weights, the realization of simple Lie algebras as algebras preserving a~non-integrable distribution was abandoned and forgotten, although this approach can be fruitful in modern times as well, see~ \cite{GL5, BGLLS, KLS}.)

For an interesting interpretation of the exceptional Lie algebras, see
\cite{El3, BE}, where the picture is made clearer by superizing it and
passing to fields of positive characteristic
(Elduque's Super Magic Super Square). On breathtaking attempts, initiated by Vogel, to consider exceptional simple finite-dimensional complex Lie algebras as members of a~series, see a~ review~ \cite{AM}.

$\bullet$ Next on the agenda are simple infinite-dimensional
$\Zee$-\textit{graded Lie algebras of polynomial growth}; briefly: \textit{of GLAPG type}.\index{GLAPG} Recall that a~
$\Zee$-graded algebra $A=\oplus A_i$ is said to be \textit{simple
graded} algebra if it has no 
\textit{graded} ideals, \textit{i.e.}, no ideals $I=\oplus I_i$ such that
$I_i=I\cap A_i$. Note that a~simple $\Zee$-graded algebra can be nonsimple as an
abstract algebra: \textit{e.g.}, loop algebras with values in a~simple Lie
algebra $\fg$ are simple graded but not simple: indeed, the
collection of $\fg$-valued functions (loops)
that vanish at any non-empty set of points is an ideal. 

In \cite{K}, Kac classified infinite-dimensional 
\begin{equation}\label{Kaccl}
\begin{minipage}[c]{14.5cm}
simple $\Zee$-graded Lie algebras of GLAPG type\\ 
\textbf{provided they are generated by elements of
degree $\pm 1$}.
\end{minipage}
\end{equation}
Among the simple graded Lie algebras of polynomial growth, only $\fvect (1):=\fder~\Cee[x]$ and $\fwitt$ or $\fvect
^{L}(1):=\fder~\Cee[x^{-1}, x]$ are not generated by elements of
degree $\pm$1 with respect to any grading. Kac conjectured that
these two examples, plus his list of Lie algebras ~\eqref{Kaccl},
exhaust all the simple Lie algebras of GLAPG type (we recall this list in Subsection~
\ref{SS:1.3}).

Twenty years after \cite{K} has been published, O.~Mathieu proved this conjecture in a~string of
rather complicated papers culminating in \cite{M}.




\sssec{Certain Lie algebras of infinite growth are also useful} The reader should not think, as was customary circa 50
years ago, that all Lie (super)algebras of infinite growth are useless
or too difficult to study. For example, various versions of
$\fgl(\infty)$ and Lie (super)algebras studied by Borcherds,
Gritsenko and Nikulin, and Penkov with coauthors are, though huge,
of considerable interest; some of these (super)algebras are quite describable, together
with several types of their representations, see an approach in \cite{Eg} and \cite{CCLL, GN, PS, FPS, HP}.


\sssec{Types of simple Lie algebras of polynomial growth}\label{SS:1.3}
Let us qualitatively describe the simple Lie algebras of polynomial
growth to better visualize them. They break into the disjoint union
of the following types of which the types 1)--4) are $\Zee$-graded:

i) \textit{Finite-dimensional} algebras (growth 0). Each of them has a symmetrizable Cartan matrix.

ii) \textit{Loop algebras}, perhaps twisted (growth 1); more
important in applications are their nonsimple ``relatives'' called
\textit{affine Kac-Moody} algebras, cf. \cite{Kb3} and Subsection~
\ref{SSS:1.12.1}. Recently, these relatives were identified as \textit{double extensions};\index{Double extension}
for a~review and succinct description of this notion, see \cite{BLS}.

iii) \textit{Vectorial algebras},\index{Lie (super)algebra, vectorial} \textit{i.e.}, Lie algebras of vector
fields\footnote{These algebras are sometimes called 
``algebras of Cartan type", the term introduced in the 1960s when only 4 series of simple Lie algebras of vector fields were known and their deformations (results of quantization) were ignored. There are serial (pertaining to series) simple vectorial modular Lie (super)algebras (in any characteristic) and exceptions (from series) discovered from the 1970 till now, see \cite{BGL, BGLLS, KLS}. Note that $\fsvect_\lambda^L(1|2)$ is the only stringy  Lie superalgebra with Cartan matrix, see
\cite{GLS1}: a new super phenomenon.} with polynomial coefficients (growth is equal to the
number of indeterminates, which is \textit{assumed to be finite}) or their completions with formal power series as coefficients.

iv) The \textit{stringy algebra}\index{Lie (super)algebra, stringy} $\fvect^{L}(1) := \fder~\Cee[t^{-1},
t]$ (the superscript $L$ stands for Laurent). This class consists of one
algebra. This algebra and its completion $\overline{\fvect^{L}(1)} :=
\fder~\Cee[t^{-1}][[t]]$ are often called \textit{Witt algebra}
\index{Witt algebra} and denoted $\fwitt$ in
honor of Witt who considered its analog over fields of prime
characteristic; physicists call it the \textit{centerless Virasoro
algebra} because (the result of) its nontrivial central extension
(discovered by Fuchs and Gelfand) is called the
\textit{Virasoro algebra}, $\fvir$, in honor of Virasoro,\index{$\fvir$} who rediscovered the
corresponding central extension and was the first to demonstrate the importance of
$\fvir$ in physical models.

\begin{Remark} Strictly speaking, stringy algebras are also vectorial, but we
apply the term \textit{vectorial} only to (super)algebras whose coefficients are
\textit{polynomials} or \textit{formal series} or divided powers in the modular setting. The term ``stringy algebra'' comes from the lingo
of imaginative physicists who now play with the idea that the
elementary particle is not a~point, but rather a~slinky, a~ springy tiny
string, see \cite{GSW}. The 
term ``stringy algebra'' means
``pertaining to string theory", and also mirrors the structure of
stringy Lie superalgebras as collections (direct sums) of several strings --- the
modules over the Lie algebra $\fwitt$ every stringy Lie superalgebra
contains. In our papers, starting with \cite{FL}, ``stringy'' means either a~``simple
infinite-dimensional vectorial superalgebra on a~supermanifold whose
underlying manifold is a~circle, or a~``\textbf{relative}" (see
Subsection~ \ref{SS:1.2}) of such algebra". In some works, stringy Lie superalgebras (not necessarily simple) are called ``superconformal"; this is unjustified and inappropriate when applied to the central
extensions of simple stringy Lie superalgebras, as explained in
\cite{GLS1}.
\end{Remark} 

v) \textbf{Filtered Lie (super)algebras of type GLAPG} are not
classified. The graded Lie algebras associated to certain such filtered Lie
algebras are not simple and have a~huge radical; these filtered examples are
particular solutions to the Main Problem \ref{SS:1.2}, many of them cannot be
obtained as deformations of simple algebras of GLAPG type. Though there
are considerably many more simple filtered Lie algebras than simple
\textit{graded} ones, at the
moment we hope that it is possible to distill tame classification
problems. The examples known to us are:

(vA) The Poisson algebras of functions on the orbits in the coadjoint
representation of any simple finite-dimensional Lie algebra, the
multiparameter deformations of these Poisson algebras, and superizations of these constructions, as well as Lie (super)algebras
of matrices over rings of differential operators, cf. \cite{GL2,
LSe, DGS, Kon}.

(vB) The Lie superalgebras constructed by S.~Montgomery (\cite{My}), Vasiliev and Konstein-Tyutin, 
Pinczon and Ushirobira, see \cite{GL2, KV, PU, KT1} and references
therein; for the queerification of these Lie superalgebras, see \cite{Lq}.

(vC) The current algebras and vectorial algebras on varieties $M$,
cf. \cite{L7}; such examples are lately known for compact Riemann
surfaces $M$ as \textit{Krichever-Novikov algebras}, see \cite{KN, Shei}.

The set of examples of simple Lie \textbf{super}algebras of type v) is much
wider than the set of the non-super examples; there are several super 
versions of each of the non-super examples (vA)--(vC), see \cite{Lq}.

\paragraph{Additional features: Cartan matrix and NIS} Some of the Lie (super)algebras of the above types i) -- v) have certain properties which are sometimes overvalued and taken
for the definition of these Lie (super)algebras, as in
\cite{Dic}, complicating the nomenclature.

Each Lie algebra of type i) and each affine Kac-Moody algebra has
a~Cartan matrix; at least some of the algebras of type 5) also have
a~Cartan matrix albeit in a~generalized sense due to Saveliev and
Vershik \cite{SV}. For algebras of type i) and affine Kac-Moody
algebras, this matrix is always symmetrizable,y integer, and rather sparse so it
can be encoded by a~simple graph (Dynkin diagram). Lie
algebras with a~symmetrizable Cartan matrix and certain
algebras of type v) carry a~non-degenerate invariant symmetric
bilinear form
--- NIS\index{NIS} for short --- a~powerful tool for solving numerous problems (for examples, see the encyclopedia \cite{DSB}).

It was, perhaps, Vinberg who coined the term \textit{contragredient}\index{Lie (super)algebra, contragredient} for the $\Zee$-graded Lie algebras~ $\fg$ such that $\fg_i\simeq(\fg_{-i})^*$ with respect to the maximal torus in $\fg_0$; \textit{e.g.}, $\fg=\fvir$ or $\fq(n)$. Lately, this term is sometimes (thoughtlessly) 
used instead of ``Lie (super)algebra with Cartan matrix". 

\ssec{The super version of the Main Problem \ref{SS:1.2} is open at the moment}\label{SS:1.4} 
\textbf{Open is even its subproblem}:\index{Problem, Main}\index{Problem, open}
\begin{equation}\label{OpenP}
\begin{array}{l}
\textit{Classify simple $\Zee$-graded Lie superalgebras of
polynomial growth}\\ \textit{and their deformations}.
\end{array}
\end{equation}

Let us describe what part of problem ~\eqref{OpenP} is solved. The
above types i)--iv) of simple Lie algebras, and their properties,
become intermixed under superization as well as under passage to Lie
algebras over fields of positive characteristic. Nevertheless, there exists a~very simple stratification: There are only two types of
Lie (super)algebras
--- ``symmetric" and ``lopsided":

\underline{(SY) The {\sl symmetric Lie (super)algebras}} are those for which\index{Lie (super)algebra, symmetric}\index{Lie (super)algebra, lopsided}
\begin{equation}
\label{sy}
\begin{minipage}[l]{14cm}
there is a~root decomposition (related to a~maximal diagonalizing
subalgebra, called in what follows a~\textit{maximal torus})\index{Torus} such
that $\sdim \;\fg_{\alpha}=\sdim\; \fg_{-\alpha}$ for any root
 $\alpha$, where $\sdim$ is
\textit{superdimension}.
\end{minipage}
\end{equation}

\underline{(LO) The {\sl lopsided algebras}} are those for which the condition ~\eqref{sy}
does not hold.


$\bullet$ \textit{Finite-dimensional Lie superalgebras}. The simple finite-dimensional Lie superalgebras are classified in the ``naive" setting; for a~review together with the first results on the representation theory, see \cite{K2}. The deformations of simple finite-dimensional Lie superalgebras, in particular with odd parameters, were recently classified as well, see \cite{Ld}, thus completing the classification.

\underline{The Lie (super)algebras of type
SY} split into the following two subtypes:


\textbf{SYCM}, the class of symmetric algebras \textbf{with Cartan matrix};
\textit{a~posteriori} we see that the Cartan matrices of the finite-dimensional simple Lie (super)algebras are all symmetrizable, and
hence each of these Lie (super)algebras has an
even NIS, see \cite{BKLS, KLLS} or \cite{KLLS1}. 

\textbf{SYoCM}, the class of symmetric algebras \textbf{without Cartan matrix}
in the conventional sense (but, perhaps, having one in the sense of
Saveliev and Vershik, \cite{SV}) but with a~NIS that is either even or odd. 

Some of simple Lie (super)algebras without Cartan matrix have relatives
with a~Cartan matrix (\textit{e.g.}, the loop (super)algebras
and related to them affine Kac-Moody (super)algebras, cf. \cite{CCLL, BLS}).



$\bullet$ \textit{(Twisted) loop superalgebras}. Obviously,
any loop (super)algebra with values in any simple Lie (super)algebra is
simple as \textit{graded algebra} (but not as abstract algebras, see Subsection~\ref{SS:1.1}). 

For any (simple, or a~relative of a~simple) finite-dimensional Lie (super)algebra $\fg$ over $\Cee$, we consider spaces $\fg^{\ell(1)}:=\fg\otimes\Cee[t^{-1}, t]$\index{$\fg^{\ell(1)}:=\fg\otimes\Cee[t^{-1}, t]$} of ``loops", \textit{i.e.}, $\fg$-valued functions on the circle $S$ that admit expansions into Laurent polynomials (or their completions, Laurent series); here $t=\exp(i\varphi)$, where $\varphi$ is the angle parameter on $S$. We introduced the notation $\fg^{\ell(1)}$ to distinguish from $\fg^{(1)}:=[\fg, \fg]$\index{$\fg^{(1)}:=\fg':=[\fg, \fg]$} and to make the expression $(\fh^{(1)}(0|n))^{\ell(1)}$, and similar expressions, meaningful, see \cite{BLS}. In this paper, we consider loops only in this subsection, so we simply write $\fg':=[\fg, \fg]$.

Recall that if $\psi$ is an order-$m$ automorphism of $\fg$, and $\fg_{\bar k}$, where $0\leq \bar k;=k\pmod{m}\leq m-1$, are eigenspaces of $\psi$ with eigenvalue $\exp(\nfrac{2\pi k \sqrt{-1}}{m})$, then \index{$\fg_\psi^{\ell(m)}$}\index{Superalgebra, twisted loop} a~\lq\lq twisted" loop (super)algebra is
\[
\fg_\psi^{\ell(m)}:=\mathop{\oplus}\limits_{0\leq \bar k\leq m-1,\ j\in\Zee}\ \fg_{\bar k}t^{k+mj}. 
\]

Serganova classified non-isomorphic loop superalgebras $\fg^{\ell(1)}:=\fg\otimes\Cee[t^{-1}, t]$ for simple finite-dimensional Lie superalgebras $\fg$ over $\Cee$, and twisted loops $\fg_\psi^{\ell(m)}$ corresponding to order-$m$ automorphisms~ $\psi$ of $\fg$, see~\cite{Se2, Se3}; some ``relatives" of these simple algebras $\fg_\psi^{\ell(m)}$ have Cartan matrices classified in \cite{HS}. Cantarini rediscovered some of these examples, the ones with compatible grading, see \cite{Ca1}, obviously unaware, together with the referee of her paper, of Serganova's results (at least, Serganova was not cited).

Like the finite-dimensional superalgebras, loop superalgebras can
also be symmetric or lopsided. For details in the non-super case, where only symmetric loop algebras exist, see \cite{Kb3}. 

For some finite-dimensional Lie
superalgebras $\fg$ with a~Cartan matrix, there are twisted loop
superalgebras $\fg^{\ell(k)}_{\psi}$ no relative of which has a
Cartan matrix; an opposite phenomenon also occurs: for some simple
Lie superalgebras $\fg$ no relative of which has a~Cartan matrix,
there is a~relative of $\fg^{\ell(k)}_{\psi}$ with a~Cartan matrix.

The double extension of every simple (twisted) loop algebra has a~Cartan matrix corresponding to an extended Dynkin diagram. Unlike the non-super case, some of the loop superalgebras have no relative with Cartan matrix, and some have
a~non-symmetrizable Cartan matrix.

An intrinsic characterization of loop superalgebras without appeal
to Cartan matrices
together with an intrinsic characterization for stringy
superalgebras (both borrowed from \cite{M}) are given in \cite{GLS1}: both types of (super)algebras are $\Zee$-graded
$\fg=\oplus_{i=-d}^{\infty}\ \fg_{i}$ of infinite\index{Depth of filtration or grading}
depth $d=\infty$, but their elements act differently in the adjoint
representation:
\begin{equation}\label{loop}
\begin{minipage}[l]{13cm}
{\sl \underline{for the loop-type (super)algebras}, every root vector corresponding
to~any~real\footnote{Recall that a~root $\alpha$ is said to be
\textit{real} if $k\alpha$ is a~root only for finitely many distinct
values of the scalar $k$. Otherwise, the root $\alpha$ is called
\textit{imaginary}.}\index{Root, real}\index{Root, imaginary} root acts locally nilpotently
in the adjoint representation;\\ \underline{for~the~stringy (super)algebras}, this is
not the case.}
\end{minipage}
\end{equation}

Before O.~Mathieu gave the above intrinsic definition for the Lie
algebras (\cite{M}), Leites conjectured that, as for Lie
algebras, simple twisted loop superalgebras correspond to outer
automorphisms \hbox{$\varphi\in\Aut\fg/\Int\fg$}, where $\Aut\fg$
is the group of all automorphisms $\fg$ and $\Int\fg$ is the
subgroup of the inner automorphisms. Serganova listed these
automorphisms and amended the conjecture: to get non-isomorphic
algebras, one should factor $\Aut\fg$ by a~group 
larger than $\Int \fg$, namely, by the connected component of
the unit of $\Aut\fg$, see \cite{Se2,FLS}. 

An \textit{a~priori}
classification of Lie superalgebras of polynomial growth with
symmetrizable Cartan matrix is due to V.~Serganova and J.~van de Leur \cite{Se2,vdL}; these
classifications and Serganova's classification of Lie superalgebras
of polynomial growth with \textit{non-symmetrizable} Cartan matrix
cited in \cite{LS1,GLS1} (for a~published proof, see
\cite{HS}) prove the conjectured list of loop-type superalgebras
with indecomposable Cartan matrix given in \cite{FLS}.

Recall the conjecture in \cite{FLS}: the only other simple Lie superalgebras 
of loop type (see definition~\eqref{loop}) are (twisted) loops with values in simple finite-dimensional ``lopsided" or ``queer" Lie superalgebras. To prove this conjecture remains an \textbf{Open problem}, tackled in~\cite{Ca2}.\index{Problem, open}


$\bullet$ \textit{Stringy superalgebras}. These are Lie
superalgebras of vector fields on the super circles, and their
relatives. For the conjectured list of simple ones, see \cite{GLS1}.
Partly (with several extra conditions) this conjecture is proved in
\cite{KvdL}, where all non-trivial central extensions of simple stringy superlgebras are found; for a~proof which, inexplicably, ignores Ramond-type superalgebras, see \cite{Kcf, FaKa, Kcfe}. Observe
that some of the stringy superalgebras have Cartan matrices, though
non-symmetrizable ones \cite{GLS1}, and some of them have nontrivial
deformations, whereas their namesakes
with polynomial coefficients are rigid according to the modern
deformation theory that we will use in what follows, see Subsection~\ref{SSS:1.12.1}.

Some of the simple stringy superalgebras (\cite{FL, GLS1}) are
\textit{distinguished}: they have a~nontrivial central extension.
In \cite{GLS1}, there are also indicated exceptional stringy
superalgebras and occasional isomorphisms unnoticed in previous
papers and often ignored in later ones.

In this paper, we consider the remaining type of simple $\Zee$-graded
Lie superalgebras of polynomial growth:


$\bullet$ \textit{Infinite-dimensional vectorial Lie superalgebras}.
The main example to look at with the mind's eye is the Lie algebra
$L=\fder~\Cee[x]$, where $x:=(x_1,\dots, x_m)$, of polynomial vector fields with the grading and
filtration given by $\deg x_{i}=1$ for all~ $i$, and its
$(x)$-adic completion, the filtered Lie algebra
$\cL=\fder~\Cee[[x]]$ of formal vector fields.

\'E.~Cartan came to the classification problem of simple vectorial
algebras from geometrical problems in which \textit{primitive}, rather
than simple, Lie algebras naturally arise; let us recall their definition. Let a~Lie algebra $\cL$ contain \index{Lie (super)algebra, primitive}\index{Lie (super)algebra, transitive}
a~subalgebra $\cL_0$ which contains no non-zero ideal of $\cL$; then, $\cL$ is called \textit{transitive}.
A transitive Lie algebra $\cL$ is called \textit{primitive}\index{Lie (super)algebra, primitive} (relative a~maximal subalgebra $\cL_0$ of $\cL$), see \cite{O}, if 
\begin{equation}\label{prim}
\text{the space $\cL_1:=\{X\in \cL_0\mid [X,\cL] \in \cL_0\}$
is not $\{0\}$. }
\end{equation}

The set of infinite-dimensional
primitive Lie algebras of polynomial growth (resp., the set of finite-dimensional primitive Lie algebras) coincides with the set of Lie algebras of derivations
of simple vectorial Lie algebras (resp., the set of trivial central extensions of simple finite-dimensional Lie algebras, see \cite{O}). Since $\dim \fder~\fg/\fg\leq 1$ and is seldom equal to 1, the task of
classifying infinite-dimensional \textit{primitive} Lie algebras is
practically the same as that of classifying the simple vectorial Lie algebras. 

By contrast, the classification of primitive Lie
\textbf{super}algebras is a~wild problem, as shown already in \cite{ALSh}, see also
Subsection~ \ref{SS:10.2}. Having distinguished a~feasible
problem, namely,
\begin{equation}
\label{1.0} \textbf{``classify infinite-dimensional \textit{simple}
vectorial Lie superalgebras,"}
\end{equation} 
it is also natural to classify all relatives of these
algebras, see Subsection~\ref{BrGo}.

\ssec{Weisfeiler filtrations and gradings}\label{SS:1.5} Consider
infinite-dimensional filtered Lie (super)algebras $\cG$ with
decreasing filtration of the form
\begin{equation}
\label{1.1}
\cG= \cG_{-f}\supset \cG_{-f+1}\supset
 \cdots \supset \cG_{0}\supset
\cG_{1}\supset \cdots\ ,
\end{equation}
where \textit{depth}\index{Depth of filtration or grading} $f$ is finite and where
\begin{equation}\label{1.1cond}
\begin{array}{ll}
\text{a)}& \cG_{0}~~\text{is a~maximal subalgebra of finite
codimension;}\\
\text{b)}& \cG_{0}~~\text{does not contain ideals of the whole $\cG$.}
\end{array}
\end{equation}

The very term ``filtered algebra" implies that $[\cG_{i},
\cG_{j}]\subset \cG_{i+j}$. Set \hbox{$G_{i}:=\cG_{i}/\cG_{i+1}$}. We
assume that
\begin{equation}\label{1.1cond2}
\begin{array}{ll}
\text{a)}&\dim G_{i}<\infty\text{~~for all~ $i$;}\\
\text{b)}&\text{$\dim \oplus_{k\leq n} \, G_{k}$ grows as
a~polynomial in $n$.}
\end{array}
\end{equation}
We assume that these Lie (super)algebras $\cG$ are complete with
respect to a~natural topology whose basis of neighborhoods of zero
is formed by the spaces of finite codimension, \textit{e.g.}, by the~$\cG_{i}$.
(In the absence of odd indeterminates, this topology is the most
natural one: two vector fields are said to be \textit{$k$-close} if
their coefficients coincide up to terms of degree $\leq k$.) No
other topology will be encountered below, so all topological terms
(complete, open, closed) refer to this topology.\footnote{This
topology is naturally called \textit{projective limit
topology}, but the even clumsier term ``linearly compact\index{Lie (super)algebra, linearly compact}
topology'' is also used.}

Weisfeiler \cite{W} endowed every above-described filtered Lie
algebra $\cG$ with another, refined, filtration, called now the \textit{Weisfeiler filtration} \index{Weisfeiler filtration, W-filtration} (here $\cG=\cL$ as
abstract algebras and all inclusions are strict):
\[
\cL= \cL_{-d}\supset \cL_{-d+1}\supset \cdots \supset
\cL_{0}\supset \cL_{1}\supset \cdots ,
\]
where $\cL_{0}=\cG_{0}$ and $\cL_{-1}$ is the minimal
$\cL_{0}$-invariant subspace strictly containing $\cL_{0}$; the other
components are defined by formula \eqref{WeiGr}.

\sssec{On the importance of W-filtrations} Each simple infinite-dimensional vectorial Lie
algebra $\cL$ over $\Cee$ or $\Ree$ has only one (up to an automorphism) maximal
subalgebra $\cL_{0}$ of finite codimension; this $\cL_{0}$
determines the W-filtration; unlike W-grading, the W-filtration is invariant under automorphisms. The topology generated by the W-filtration
is precisely the topology in which $\cL$ can be considered as the Lie
(super)algebra of \textit{continuous} derivations of the algebra of
functions represented by formal power series, see~ \cite{BL3}.

It is often more natural to consider infinite-dimensional algebras,
and modules over them, as \textit{topological} algebras, rather than
as \textit{abstract} ones, and therefore vectorial algebras with
\textit{formal} coefficients and $(x)$-adic filtration (topology)
are more natural objects than vectorial algebras with
\textit{polynomial} coefficients. These ``more natural" algebras are
completions with respect to the $(x)$-adic filtration of the
vectorial algebras with polynomial coefficients, \textit{i.e.}, we complete
polynomials to series adding ``tails", \textit{i.e.}, elements of higher
degree. The complete topological algebras are natural objects, but it
is easier to classify first the associated graded algebras.

The W-filtration of a~given vectorial Lie (super)algebras $\cG$ is the filtration
corresponding to the $(x)$-adic topology on the space of functions,
or a~modification of this topology that takes into account a
non-integrable distribution preserved by $\cG$. Moreover, for the \textit{W-filtratered} Lie (super)algebra $\cL\simeq \cG$, the $L_{0}$-action on
$L_{-1}$ (recall definition \eqref{grL}) is irreducible.

For the \textbf{simple} filtered Lie (super)algebra~$\cL$, and the graded algebra $L$ associated
with it, the $L_{0}$-module $L_{-1}$ is faithful.
So, being primarily interested in \textbf{simple} Lie
superalgebras, \textbf{we assume that the $L_{0}$-module $L_{-1}$ is
faithful}.


\ssec{Classification problems: A to C}\label{ssProb}The following problems arise:

\textbf{A) Classify simple W-graded vectorial Lie superalgebras
as abstract ones, \textit{i.e.}, disregarding difference in grading.}

We immediately see that this problem is rather unnatural: quite
distinct algebras become equivalent. Indeed, as early as in
\cite{ALSh}, we observed that $\fvect(1|1)$, the Lie algebra of all
vector fields on $(1|1)$-dimensional superspace, is isomorphic as an
abstract algebra to $\fk(1|2)$, the Lie algebra of \textit{contact}
vector fields on $(1|2)$-dimensional superspace (and also to
$\fm(1)$, another, ``odd'' type of contact Lie superalgebra). Therefore, a~more
natural formulation is 

\textbf{A$'$)\label{ProbA} Classify simple W-graded vectorial Lie superalgebras.}

In some applications it suffices to confine ourselves to graded
algebras (\cite{L8}), but in most applications (\textit{e.g.}, in the study
of representations \cite{GLS4}), the Lie (super)algebras that are
\textit{complete} with respect to a~topology are more natural than
the graded ones associated with them. So the following problem might
look like a~reasonable goal:

On the road to solution of Problem B$'$ the following problem
\textit{seems} to be a~natural step:

\textbf{B)\label{B'} Classify simple W-filtered complete vectorial Lie
superalgebras as abstract ones, \textit{i.e.}, disregarding differences in filtrations.}

\textbf{B$'$)\label{ProbB} Classify simple W-filtered complete vectorial Lie
superalgebras.}

Natural (``straightforward'') \index{Completion, direct} examples of complete algebras are provided by the
algebras $\overline L$ obtained from the list of examples $L$ that
answer Problem A or A$'$ by taking formal series, rather than
polynomials, as coefficients. In addition to these examples, there
might occur nontrivial \textit{filtered deforms} of $\overline L$,
\textit{i.e.}, complete algebras $\widetilde {\cL}$ such that the graded
algebras associated with $\overline L$ and $\widetilde {\cL}$ are
both isomorphic to $L$, whereas $\overline L\not\simeq\widetilde {\cL}$. Examples of such filtered deforms will be given
later.

Observe that whereas the difference between Problems A$'$ and B$'$
is as big as that between Problems A and B, the difference between
Problems B and B$'$ is negligible. Actually, it is impossible to
solve Problem B$'$ and not solve Problem B. Indeed, various filtered
deforms are not \textit{a~priori} isomorphic as abstract algebras,
so, in order to solve Problem B$'$, we have to know, first, all the
W-filtrations or the associated W-gradings and then describe
filtered deforms for \textit{every} W-grading.

However, even Problem B$'$ is not the most natural one: there exist
deforms of algebras from class B$'$ that do not lie in class B$'$: \textit{e.g.},
the result of factorization modulo the center of the quantized Poisson superalgebra $\cQ(\fpo(2n|2m))$, cf.~\eqref{summ}; clearly, $\cQ(\fpo(2n|2m))$ is a~ familiar Lie superalgebra of differential operators with polynomial coefficients. So the true problem one should solve is

\textbf{C)\label{ProbC} Describe the simple W-gra\-ded vectorial Lie
superalgebras together with all their deformations (not only the filtered
ones).}

\subsubsection{Notation}\label{SS:1.5.1} We will
describe all the Lie superalgebras from the tables in our Main
Theorem in detail in due course, but to simplify grasping the
general picture from the displayed formulas of the following
theorem, let us immediately point out that the prime
example, the \textit{general vectorial algebra} $\fvect(m|n; r)$, is the Lie
superalgebra of vector fields whose coefficients are polynomials (or
formal power series, in which case it is more correct to denote it $\overline{\fvect(m|n; r)}$) in $m$ commuting and
$n$ anticommuting indeterminates with the grading (and filtration)
determined by equating the degrees of $r$, where $0\leq r\leq n$, of
odd indeterminates to 0, the degrees of all (even and odd) of the
remaining indeterminates being equal to 1. The regradings of other
series are determined similarly, but in a~more complicated way, described
below, see Subsection~\ref{secnon}. Usually, we do not indicate the parameter $r$ if $r=0$.

In the sequel, for any $\fg$, we write $\fcg:=\fg \oplus {\Cee}\cdot
z$ or $\fc(\fg)$\index{$\fcg$} to denote the trivial central extension with the
1-dimensional even center generated by $z$. 

Let $\fii\ltimes \fb$ or $\fb\rtimes \fii$\index{$\ltimes$ semi-direct sum} denote the semidirect sum
of modules (resp., algebras) in which $\fii$ is a~submodule (resp., an ideal).

Let $\fd(\fg):=\fg\ltimes\Kee D$\index{$\fd(\fg)$} be a~ \textbf{sub}algebra of the Lie (super)algebra $\fder~\ \fg$ of all derivations of the Lie (super)algebra $\fg$, where $D$ is usually (in examples below) the grading operator of $\fg$.\index{$\fd(\fg)$}

In proofs, we use standard notation of roots in the simple finite-dimensional Lie algebras, see \cite[Table 5]{OV} (or the respective tables in \cite{Bbk} or \cite{Hm}).

Let $\lambda=\frac{2a}{n(a-b)}\in\Cee\cup\{\infty\}$ for some
$a,b\in\Cee$. Although notation $\fb_{a, b}(n; r)$ introduced in
~\eqref{2.7.2} is more accurate than $\fb_{\lambda}(n; r)$, we
often use the shorter notation for brevity.

The notation $\fsvect_{a, b}(0|n)$ and $\fspe_{a, b}(n)$, though\index{$\fsvect_{a, b}(0\vert n)$}\index{$\fspe_{a, b}(n)$}
looks similar to $\fb_{a, b}(n; r)$, means something different: it stands for $\fg\ltimes \Cee(az+bD)$, where $a,b\in\Cee$ and
$\fg$ stands for $\fsvect(0|n)$ or $\fspe(n)$, respectively, $D$ is
the operator that determines the standard $\Zee$-grading of $\fg$,
and $z$ is a~central element generating the trivial center (\textit{i.e.}, a
direct summand).

The \textit{exterior algebra}\index{Algebra, exterior} of the (super)space $V$ is denoted by 
$E(V)$ or --- if $V=V_\od$ --- by $\Lambda(V)$ or $\Lambda(n)$, where $n=\dim V$, a.k.a. \textit{Grassmann algebra}\index{Algebra, Grassmann}\index{$\Lambda(n)$} $\Lambda(\xi)$ on $n$ generators $\xi_1,\dots, \xi_n$.

\paragraph{Overview of various properties of vectorial Lie (super)algebras}\label{over}{}~{}

$\bullet$ \textbf{On vectorial superalgebras, there are two analogs of trace}; for details, see Subsection \ref{traceDiv}:

1) \textit{Traces} (or \textit{super}traces, for emphasis) ---
various linear functionals that vanish on the first derived\index{Trace, supertrace}
subalgebra of $\fg$ or its (super)commutant $\fg':=[\fg,\fg]$ --- are traces in the classical sense, simpleminded analogs of the trace on the matrix Lie algebra. 

2) \textit{Divergences} --- ``Cartan prolongs'', as we will see, of traces on the linear subalgebras --- though not, actually, traces, are ``morally", as Manin might have dubbed it, closer to, and better resemble the traces on matrix Lie algebras than (super)traces in the classical sense 1).

Accordingly, the analog of the \textit{special} Lie (super)algebra --- the traceless subalgebra of the matrix Lie (super)algebra --- is the \textit{divergence-free}\index{Divergence}
subalgebra of a~vectorial algebra $\fg$; this analog is denoted by $\fs\fg$,
\textit{e.g.}, $\fsvect(n|m)$, $\fsm(n)$, etc. The even (resp., odd) part of the codimension of
$\fg'$ is equal to the number of even (resp., odd) supertraces on $\fg$.

$\bullet$ Lie superalgebras of series $\fk$ and $\fh$ are the straightforward
analogs of the Lie algebras of \textit{contact} and \textit{Hamiltonian} series,
respectively; the \textit{Poisson} Lie superalgebra $\fpo(2n|m)$ is the --- there is just one (class of) --- nontrivial
central extension of $\fh(2n|m)$.

$\bullet$ Lie superalgebras of series $\fm$ and $\fle$ are the \textit{``odd" analogs}
of the Lie algebras of contact and Hamiltonian series, respectively.
The bracket in the \textit{Buttin superalgebra} $\fb$, the central extension of $\fle$, is the
classical \textit{Schouten bracket}, currently more popular under the name of 
\textit{antibracket}. We say ``odd" in quotation marks: all these
``odd" analogs have both even and odd parts.

$\bullet$ The \textit{antibracket can be deformed}, the corresponding deformed Lie
superalgebras inside the GLAPG type are~$\fb_{\lambda}$, where
$\lambda\in\Cee \Pee^1$; the results of other deformations of $\fb_{\lambda}$, some of them with an odd parameter, lie outside of
the GLAPG class.

$\bullet$ The algebras $\fsm$, $\fs\fle$, and $\fs\fb$ are divergence-free
subalgebras in the respective algebras.

$\bullet$ The simple algebra $\widetilde{\fs\fb}$ is a~deformation (deform if the parameter of deformation is even) of a~nonsimple
Lie superalgebra $\fs\fb(n|n+1;n)$.

$\bullet$ We note \textbf{an outstanding property of the bracket in
$\fh_{\lambda}(2|2)\simeq\fb_{\lambda}(2; 2)$: it can be considered
as a~deform of the Poisson bracket as well as a~deform of the anti-bracket}.


$\bullet$ In tables ~\eqref{1.4}, ~\eqref{1.10} below: in parentheses after the
``family name" of the algebra there stands the superdimension of the
superspace spanned by the indeterminates and
--- after the semicolon --- an abbreviated description of the regrading.

The regradings of the series are governed by a~vector $\vec r$, often abbreviated to a~ number~ $r$ or a~ symbol, as
described in Subsection~\ref{secnon}. All regradings are given in
(hopefully) suggestive notation, \textit{e.g.}, $\fk\fas^\xi (1|6; 3\eta)$
means that, having taken $\fk\fas^\xi$ as the point of reference, we
set $\deg\eta=0$ for each of the three $\eta$'s (certainly, this
imposes some conditions on the degrees of the other indeterminates, see Subsection~\ref{secnon}). Passing from one regrading to another one, we take the realization with the minimal value of $\dim\cL/\cL_{0}$) as the
point of reference. For the exceptional Lie superalgebras, another
point of reference --- the \textit{compatible} grading is often more convenient; it is marked by~ $K$.

The exceptional grading $\Reg_{\fh}$ of $\fh_{\lambda}(2|2)$ is
described in passing in the list of occasional isomorphisms
~\eqref{1.7}; it is described in detail for another incarnation of
this algebra in Subsection~\ref{SS:1.3.1}.

$\bullet$ \textbf{On ``drop-outs''.} Several algebras are ``drop-outs'' from the
series.\index{Drop-outs from the
series} Examples: the algebras $\fsvect'(1|n; r)$ are ``drop-outs''
from the series $\fsvect(m|n; r)$ since the latter are not simple
for $m=1$, but contain the simple ideal $\fsvect'(1|n; r)$.

Similarly, $\fle(n|n; r)$, $\fb_{1}'(n|n+1; r)$ and
$\fb_{\infty}'(n|n+1; r)$ are ``drop-outs'' from the Lie superalgebras 
$\fb_{\lambda}(n|n+1; r)$ corresponding to $\lambda= 0$, 1 and
$\infty$, respectively, having either a~simple ideal of codimension
1, or a~center the quotient modulo which is simple. For $n=2$, there
are several more of exceptional values of $\lambda$.

The Lie superalgebra $\fsm(n|n+1; r)$ is a~drop-out from the series
$\fb_{\lambda}(n|n+1; r)$ for another reason: it is singled out
by its divergence-free property, hence deserves a~separate line.

The finite-dimensional Lie superalgebra $\fh(0|n)$ of Hamiltonian
vector fields is not simple: it contains a~simple ideal $\fh'(0|n)$. The Lie superalgebra $\fsl(m|n)$ for $m=n$ is not simple, and therefore the simple Lie superalgebra $\fpsl(n|n)$ for $n>1$ is a~ ``drop out" from the $\fsl$ series.

\textbf{One should carefully distinguish drop-outs} not only in the case of
characteristic~ $p$ or the super-case. People justly consider the
Lie algebra of Hamiltonian vector fields and its central extension, the
Poisson algebra 
as totally distinct algebras for any $p$ (or loop algebra and its double extension, \textit{i.e.}, affine Kac-Moody relative having a~ Cartan matrix). Likewise, the difference between
$\fsl$ and $\fp\fsl$ in characteristic~ $p$ and in the super-case is
enormous: \textit{e.g.}, the irreducible $\fp\fsl$-module of least dimension is the adjoint one.

Thus, all drop-outs, without exception, are quite separate items on
the list.

\ssec{The steps of our classification}\label{steps}
(From \cite[1988, pp.235--289]{LSh1}; compare with practically the same steps in
\cite[\S 3]{Kga}.) For Lie algebras, the following facts and ideas
are known, though, perhaps, not as well as they deserve, see
\cite{J,GQS}.

Namely, given a~Lie superalgebra $\cL$ with a~Weisfeiler filtration,
consider the associated graded Lie superalgebra
$\fg:= \oplus_{i\geq -d}\ \fg _{i}$, where $\fg
_{i}=\cL_{i}/\cL_{i+1}$. There are two distinct cases: that of
a~\textit{compatible} $\Zee$-grading and
when the grading is \textit{incompatible}.\index{Grading, compatible}

To the \underline{incompatible grading} certain arguments used in
classification of simple vectorial Lie algebras (\cite{J}) are
applicable: in particular, the proof of the facts that for W-graded Lie (super)algebras,
\begin{equation}\label{1-2-3}
\text{\begin{minipage}[c]{14cm} $1)$ $\fg _{-1}$ generates $\fg _{-}:=
\oplus_{i<0}\ \fg _{i}$;\\
$2)$ If $a\in \fg_{0}$ and $[\fg _{-1}, a]=0$, then $a=0$
(transitivity condition);\\
3) the $\fg _{0}$-module $\fg _{-1}$ is not just transitive, it is
irreducible.
\end{minipage}}
\end{equation}

The proof of the fact 
\begin{equation}
\label{3.1} \text{``a root of $\fg _{0}$ is the sum of positive
weights of $\fg _{-1}$"}
\end{equation}
is given in Section~\ref{S:3},
where ~\eqref{3.1} reads as follows:
\begin{equation}\label{rank1}
\begin{minipage}[l]{14cm}
{\sl the Cartan prolong of the pair $(\fg _{-}, \fg _{0})$ is of
infinite dimension if and only if the action of $\fg _{0}$ on $\fg
_{-1}$ has a~rank-$1$ operator --- the coroot corresponding to the root~
$\eqref{3.1}$ --- with an \textit{even} covector}.
\end{minipage}
\end{equation}

Having established this, we have to find the pairs $(\fg _{-1}, \fg
_{0})$ for which ~\eqref{3.1} holds. The search for such pairs, and the
proof of the completeness of their list occupy Sections~\ref{S:3}--\ref{S:6pc}.

\textbf{Step 1}. Formulation of the Pivotal Ideas of our proof.

(1) After Weisfeiler, we use the fact that every simple filtered
vectorial Lie superalgebra has a~W-filtration such that the
associated graded Lie superalgebra is simple, see \cite{W};

(2) We show that any simple $\Zee$-graded vectorial Lie superalgebra
$\oplus_{i\geq -d}\ \fg_{i}$, with
\textbf{in}compatible grading, is infinite-dimensional if and only
if $\fg_0$ acts on $\fg_{-1}$ so that the action has a~rank-1
operator with an even covector.

We have listed all simple W-graded vectorial Lie superalgebras of
depth 1 in the union of \cite{ALSh}, \cite{Sh3}, \cite{ShM},
\cite{Sh14}. This union contains (as
is now established) all the examples of simple vectorial
superalgebras of any depth $d$, 
considered as abstract ones, and all
their W-regradings, except for the W-regrading $CK$, see table \eqref{1.10}.

(3) The cases where $\fg_{-1}$ is a~purely odd superspace
(compatible $\Zee$-grading) or $d>1$. While we left these cases to be
considered later, Cantarini, Cheng and Kac investigated them
in \cite{K3,Kga, CK2} and \cite{CCK}.

For a~ proof of the classification, see \cite{LSh5}, a~
detailed version of \cite{LSh1,LSh3}. The case of compatible $\Zee$-grading cannot be considered by the above-described methods. This gap in our proof was filled in in the papers by Kac and his co-workers; no new examples were discovered as we predicted. 

\textbf{Step 2}. Criterion for infinite dimensionality of a~given Cartan
prolongation: existence of rank-1 operators with an even covector.

\textbf{Step 3}. Main reduction: the three cases. For $\fg_0$ and an
irreducible faithful representation $\rho$ of $\fg_0$ on $\fg_{-1}$
with a~ rank-1 operator, the following cases are possible:

I) $\rho$ is induced (from a~representation of a~Lie subsuperalgebra of
$\fg_0$ of codimension $0|k$).

II) $\fg_0$ is \textit{almost simple} (for this definition, see Subsection~\ref{SS:2.2.1}) or $\fg_0$ is a~central extension of an
almost simple Lie superalgebra;

III) $\sdim \fg_{-1}=2|2$ and $\fg_0$ is contained in a~particular Lie superalgebra.

\textbf{Step 4}. Description of prolongations of induced
representations (case I of step 3).

\textbf{Step 5}. Description of prolongations of almost simple Lie
superalgebras and their central extensions (case II of step 3).

\textbf{Step 6}. Case III of step 3. (The unique simple infinite-dimensional Lie superalgebra arising at this step is a~regrading of
$\fb'_{\infty}(2)$.)

\textbf{Step 7}. Description of W-gradings of the 
simple superalgebras found above.

\textbf{Step 8}. The case of compatible grading. 
Proof of our conjecture, see
\cite{K3,Kga, CK2,CCK}.

\begin{Theorem}[\textit{\cite{K3}}]
The only simple vectorial Lie superalgebras with compatible
Weisfeiler filtration are: $\fk(1|m)$ and the exceptional algebras
$\fk\fas$, $\fv\fle(4|3; K)$, $\fk\fle(9|6; K)$, and $\fm\fb(4|5;
K)$.
\end{Theorem}

\textbf{Step 9. Description of deformations and deforms}, see Sections~\ref{S:11} and \ref{S:14}.

\ssec{Warning: results of cohomology calculations might be confusing}\label{SSS:1.12.1} \index{Warning} There are several types of confusing phenomena: 

1) For Lie
(super)algebras, the deformation problem \textbf{ostensibly} has an
unequivocal answer: calculate the infinitesimal deformations --- the elements of $H^2(\fg;\fg)$, the tangent space to the
(super)variety of deformation parameters --- and that is it: if $H^2(\fg;\fg)=0$, the
algebra $\fg$ is \textit{rigid},\index{Lie (super)algebra, rigid} \textit{i.e.}, has no deformations. If not, construct the global deformations
whose tangent vectors correspond to the basis elements of
$H^2(\fg;\fg)$. To find isomorphisms among the members of
the parametric family of deforms still remains sometimes an \textbf{Open problem}, especially in the modular cases.\index{Problem, open}

There are, however, examples of simple
(finite-dimensional) Lie (super)algebras $\fg$ such that $H^2(\fg;\fg)\neq
0$, some or all the cocycles constituting a~basis of $H^2(\fg;\fg)$
can be integrated to global deformations, but the deforms corresponding to these cocycles are isomorphic to the initial Lie (super)algebra $\fg$, see
\cite{Ri,Car, BLW,BLLS}. 
We call the cocycles with this property, the
cohomology classes they represent, and the corresponding deformations
\textit{semi-trivial}.\index{Deformation, semi-trivial}\index{Deformation, trivial}\index{Deformation, true} We call the deforms that are neither trivial nor semi-trivial \textit{true} deforms. Fortunately, the reasons for a~deform to be
semi-trivial given in \cite{BLLS} do not exist over~ $\Cee$. Still,
for infinite-dimensional Lie (super)algebras, the knowledge of
$H^2(\fg;\fg)$ does not always help to describe the deforms:

2a) For example, ${H^2(\fg;\fg)=0}$ in terms of the modern cohomology theory (see
\cite{FF, LR}), if $\fg$ is any of the
simple loop algebras. Despite this, P.~Golod (Holod) gave the first example of a~deform of a~loop
algebra, see \cite{GoHo1, GoHo2}; Krichiver and Novikov initiated a~representation theory of such deforms,
for a~review, see \cite{Shei}.

2b) In \cite{FiFu}, an example is given of several non-isomorphic smooth deforms
with a~common tangent vector. In this case, to compute
$H^2(\fg;\fg)$ does not suffice to describe all deformations.

To figure out how to describe deformations that explain the above examples is an \textbf{Open problem}.\index{Problem, open}

\sssec{On filtered deformations}\label{SSS:1.12.3} For the
$\Zee$-graded Lie algebra $\fg$ of depth 1~and a~decreasing
filtration, it was well known (see, \textit{e.g.}, \cite{KfiD}; superization is streightforward) that\index{Deformation, filtered}
\begin{equation}
\label{filDefCK} \text{\begin{minipage}[l]{13cm} the filtered
deformations of $\fg$ correspond to ``certain elements" of the Spencer
cohomology\footnote{For examples of calculation and to compare the techniques and interpretations of the results, compare \cite{CK1, CK1a} with \cite{Po3} (by ``hand") and \cite{GL4,GLS2,BGLS, L8}.}
$H^{2,j}(\fg_{-1}; \fg):= \{c\in H^{2}(\fg_{-1}; \fg)\mid\deg c=j\}$.
\end{minipage}
}
\end{equation}

Let us describe which are the ``certain elements" alluded to in
~\eqref{filDefCK}, and how to treat the algebras of depth $d>0$.

The latter question was answered in \cite{GL4}, for more details,
see \cite{BGLS}: one has to pass from the Spencer cohomology to the Lie
(super)algebra cohomology: $H^{2}(\fg_{-1};
\fg)=\oplus_j H^{2,j}(\fg_{-1}; \fg)$. In this
formulation, it is obvious how to generalize the above-described passage to the algebras of depth $d>0$: let $\fg_{-}:=\oplus_{-d\leq j<0}\ \fg_j$, then deformations correspond to
``certain elements" of $H^{2}(\fg_{-}; \fg)$.
These ``certain elements" are the elements of $H^{2}(\fg; \fg)$ that are homogeneous of positive
degree. This
description allows the filtration of $\fg$ to begin with the component of degree
$\geq 0$, and hence answers the question Kontsevich posed to us in
2002: 
\[
\text{``How to describe filtered deformations of $\cL$ if its depth is
equal to 0?"}
\] 
However, for $d>0$, computing the cocycles of the
space $H^{2}(\fg; \fg)$ is a~much tougher task than computing even
the whole space $H^{2}(\fg_{-}; \fg)$.

In the paper \cite{CK1} by Cheng and Kac, the ``certain elements"
mentioned above are explicitly described (for the first time, as far as
we know): it is shown that filtered deformations of a~given graded Lie
superalgebra $\fg$ of non-zero depth and with a~decreasing filtration
correspond to \textbf{$\fg_0$-invariants of the superspace $H^2(\fg_-;
\fg)$}. 
(Note that Cheng and Kac did not believe in odd parameters
of deformations, and hence considered only even invariants, see \cite{CK1, CK1a}.)

A question arises: given two non-isomorphic $\Zee$-gradings $\fg$ and $\fh$ of the
same abstract algebra, when does an isomorphism
\begin{equation}
\label{1.5.5} (H^2(\fg_-; \fg))^{\fg_0}\simeq (H^2(\fh_-;
\fh))^{\fh_0}
\end{equation}
take place? The answer in the cases where all, not only filtered,
deformations are known clearly shows that the isomorphism ~\eqref{1.5.5}
takes place only if 
\[
\begin{array}{l}
\text{the cocycles $c_1\in(H^2(\fg_-; \fg))^{\fg_0}$
and $c_2\in(H^2(\fh_-; \fh))^{\fh_0}$}\\ 
\text{are compatible with both gradings.}
\end{array}
\]
For example, the unique cocycle in cases ~\eqref{1.5.4} is only
compatible with the W-gradings for which the product of an even number of odd
indeterminates is of degree 0. \textit{A posteriori} we know that for other $\Zee$-gradings, there are
no filtered deformations.

Therefore, \cite{CK1} does not contain a~complete list of filtered
deformations of all simple W-graded vectorial Lie superalgebras
$\cL$, as claimed. In fact, \cite{CK1} considered only $\frac1m$th of all possibilities for
each $n|m$-dimensional superspace $(\cL/\cL_0)^*$. Actually, the claim of \cite{CK1} is more vague and subject to interpretations: ``we determine all simple
filtered deformations of certain $\Zee$-graded Lie superalgebras
classified in \cite{K3}", but since in \cite{K3} the algebras are
classified only as abstract ones, it is impossible to say what
exactly is claimed in \cite{CK1}.

\ssec{Our results}\label{ssres} For the list of problems A to C, see Subsection~\ref{ssProb}.

1) \textbf{Simple vs primitive}. In \cite{ALSh} we have distilled a~tame subproblem ---
\textit{classify simple vectorial Lie superalgebras} --- from its
straightforward, but \textbf{wild in super setting}, see Subsection~ \ref{SS:10.2}, problem
\textit{classify primitive Lie superalgebras}, see \cite{K2, K3, Kga}.

2) \textbf{Solution of Problem A$'$}. The union of our
examples (\cite{ALSh}, \cite{LSh1}, \cite{Sh5, Sh14}) was our conjectural
list solving Problem A$'$, as announced in \cite{LSh2$'$}. The list we
announced is, indeed, complete. For compatible gradings, see \cite{K3, CK1, CK1a,CK2,CCK}, which prove our conjecture.

3) \textbf{Solution of Problem B$'$}. Combining solution of Problem A$'$ with the results on filtered deformations, see \cite{LSh3,GLS2}, we
obtain solution of Problem B$'$.

4) \textbf{Classification of W-gradings}. See Subsections~\ref{SS:7.2}, \ref{secnon}, \ref{SS:1.3.3}.

\sssbegin[Solution to Problems A and A${}'$, B and B${}'$]{Theorem}[Solution to Problems A and A${}'$, B and B${}'$]\label{mainth} 
The simple W-graded vectorial Lie superalgebras $\fg$ can be thriftily collected into $14$ families as abstract algebras~\eqref{1.4}, further divided into $34$ W-graded series, see Subsection~$\ref{SS:1.3.3}$, and $15$
W-graded exceptional algebras forming $5$ abstract algebras, see Subsection~$\ref{th1.10}$. They are pair-wise non-isomorphic,
as graded and filtered superalgebras, except for occasional isomorphisms
~\eqref{1.7}.  \end{Theorem}
\begin{equation}\label{1.4}\tiny 
\begin{array}{ll}
\renewcommand{\arraystretch}{1.4}
\begin{tabular}{|l|l|}
\hline $N$&algebra and conditions for its simplicity\cr 
\hline

1&$\fvect(m|n; r)$ for $m\geq 1$ and $0\leq r\leq n$ \cr 
\hline
\hline
FD &
$\fvect(0|n; r)$ for $n>1$ and $0\leq r\leq n$ \cr 
\hline

2&$\fsvect(m|n; r)$ for $m>1$, $0\leq r\leq n$ \cr 
\hline

FD &
$\fsvect(0|n; r)$ for $n>2$ and $0\leq r\leq n$ \cr 
\hline

3&$\fsvect'(1|n; r)$ for $n>1$, $0\leq r\leq n$\cr 
\hline

FD &$\widetilde{\fsvect}(0|n)$ for $n>2$ \cr \hline

\hline 
4&$\fk(2m+1|n; r)$ for $0\leq r\leq [\frac{n}{2}]$,
$(m|n)\neq (0|2k)$\cr 
\hline
&$\fk(1|2k; r)$ for $0\leq r\leq k$ except for
$r=k-1$\cr
 \hline
 \hline

5&$\fh(2m|n; r)$ for $m>0$ and $0\leq r\leq [\frac{n}{2}]$\cr
\hline

FD&$\fh'(0|n)$ for $n>3$ \cr 
\hline
\end{tabular}&\begin{tabular}{|l|l|}
\hline $N$&algebra and conditions for its simplicity\cr 
\hline

6&$\fh_{\lambda}(2|2; r)$ for 
$r=0$, $1$ and $\Reg_{\fh}$ \cr 
\hline


\hline
7&$\fm(n|n+1; r)$ for $0\leq r\leq n$, $r\neq n-1$\cr
 \hline
 
8&$\fsm(n|n+1; r)$ for $n>1$ but $n\neq 3$ and \cr
&$0\leq r\leq n$, $r\neq n-1$\cr \hline

9&$\fb_{\lambda}(n|n+1; r)$ for $n>1$, where $\lambda\neq 0, 1,
\infty$\cr
& and $0\leq r\leq n$, $r\neq n-1$;\cr
&$\fb_{\lambda}(2|3;
r)$, for 
$r=0$, $1$ and $\Reg_{\fb}$ \cr \hline

10&$\fb_{1}'(n|n+1; r)$ for $n>1$, $0\leq r\leq n$, $r\neq
n-1$\cr \hline

11&$\fb_{\infty}'(n|n+1; r)$ for $n>1$, $0\leq r\leq n$, $r\neq n-1$\cr \hline 

12&$\fle(n|n; r)$ for $n>1$, $0\leq r\leq n$, $r\neq n-1$\cr
\hline
\hline

13&$\fsle'(n|n; r)$ for $n>2$, $0\leq r\leq n$, $r\neq n-1$\cr
\hline

14&$\widetilde{\fs\fb}_{\mu}(2^{n-1}-1|2^{n-1})$ for $\mu\neq 0$ and
$n>2$ even\cr 
&$\widetilde{\fs\fb}_{\mu}(2^{n-1}|2^{n-1}-1)$ for $\mu\neq 0$ and
$n>1$ odd\cr\hline
\end{tabular}
\end{array}
 \end{equation}

Hereafter we abbreviate the algebras $\fm(n|n+1;
r)$, $\fb_{\lambda}(n|n+1; r)$, $\fle(n|n; r)$ and so on by $\fm(n; r)$, $\fb_{\lambda}(n; r)$, $\fle(n; r)$,
and so on, respectively.

\begin{proof}[Proof \nopoint] of this theorem occupies the main body of this article. \end{proof}

Observe that real forms (classified in an interesting paper \cite{CaKa1}) do not always survive regrading.

\sssec{W-gradings of serial simple vectorial Lie superalgebras}\label{remmain} Though $\fb_{\lambda}(2; r)$
and $\fh_{\lambda}(2|2; r)$ are isomorphic, we consider them
separately because they are deforms of very distinct structures: the superalgebras with the
odd and the even bracket, respectively. In Subsection~\ref{SS:2.22} we
show wherefrom the mysterious at this point isomorphism
$\fh_{\lambda}(2|2)\simeq \fh_{-1-\lambda}(2|2)$ comes.
\textbf{The family $\fh_{\lambda}(2|2)\simeq \fb_{\lambda}(2; 2)$
should be considered as an exceptional one.} The following table lists all W-gradings of serial simple vectorial Lie superalgebras, see Theorem~\ref{th7.3}.

\paragraph{Occasional isomorphisms}\label{Occasion} One can directly verify that\index{Isomorphisms, occasional}
\begin{equation}
\label{1.7} \footnotesize
\begin{array}{l}
\fvect(1|1)\simeq \fvect(1|1; 1);\\

\fsvect(2|1)\simeq \fle(2; 2); \; \fsvect(2|1; 1)\simeq \fle(2)\\

\fsm(n)\simeq \fb_{2/(n-1)}(n);\\

 \text{in particular,
$\fsm(2)\simeq \fb_{2}(2)$, and $\fsm(3)\simeq \fb_{1}(3)$; hence, $\fsm(3)$ is not
simple};\\

\fs\fle'(3)\simeq \fs\fle'(3; 3);\\

\fb_{1/2}(2; 2)\simeq \fh_{1/2}(2|2)=\fh(2|2);\quad
\fh_{\lambda}(2|2)\simeq \fb_{\lambda}(2;
2); \\

\fh_{\lambda}(2|2; 1)\simeq \fb_{\lambda}(2);\\

\text{$\fb_{1}'(2; {\rm Reg}_{\fb})\simeq \fle(2)$,}\\ 
\text{$\fle(2;{\rm
Reg}_{\fb})\simeq \fb_1'(2)$},\\
\text{$\fb_{\infty}'(2; {\rm
Reg}_{\fb}) \simeq \fb_{\infty}'(2)$};\\
\fb_{a, b}(2; \Reg_{\fb})\simeq \fb_{-b, -a}(2)\simeq 
\fb_{b, a}(2), \ \ \text{see eq.~\eqref{reg4}};\\
\fb_{\lambda}(2; {\rm Reg}_{\fb})\simeq \fb_{-1-\lambda}(2)\;\text{
and
}\; \fh_{\lambda}(2|2)\simeq \fh_{-1-\lambda}(2|2), \ \ \text{see eq.~\eqref{unmyst}; hence,}\\

\text{the fundamental domain is either } \{\lambda\in\Cee\mid \RE \lambda \geq -\frac12\} \text{ or
$\{\lambda\in\Cee\mid \RE \lambda \leq -\frac12\}$};\\

\widetilde{\fs\fb}_{\mu}(2^{n-1}-1|2^{n-1})\simeq 
\widetilde{\fs\fb}_{\nu}(2^{n-1}-1|2^{n-1})\text{ for $n>2$ even and $\mu\nu\neq 0$},\\
\widetilde{\fs\fb}_{\mu}(2^{n-1}|2^{n-1}-1)\simeq 
\widetilde{\fs\fb}_{\nu}(2^{n-1}|2^{n-1}-1)\text{ for $n>1$ odd}. 
\end{array}
\end{equation}

\textbf{Warning}.\index{Warning} Isomorphic abstract Lie superalgebras may be
quite distinct as filtered or graded: \textit{e.g.}, regradings provide us
with isomorphisms (\cite{ALSh}):
\begin{equation}
\label{1.9}  \fk(1|2;1)\simeq \fvect(1|1)\simeq \fm(1;1).
\end{equation}
Observe that only one of the above three non-isomorphic graded
algebras is W-graded over $\Cee$; the situation becomes more complicated over $\Ree$, cf. \cite{CaKa1}. 

The  standard grading of $\fk(1|2)$, though non-Weisfeiler, is often considered in physics
papers. 
Observe that some real forms of $\fk(1|2k; k-1)$ are W-graded.

\paragraph{Comments} For completeness, our table ~\eqref{1.4} includes occasional finite-dimensional
versions labeled by ``FD" instead of the number.

All the Lie superalgebras listed in Theorem~\ref{mainth} are the results of
\textit{Cartan prolongation}, see Subsection~$\ref{2.6}$ (or the derived algebras thereof),
and therefore are determined by the terms $\fg_{i}$
with $i\leq 0$ (or $i\leq 1$ in some cases), see Subsection~$\ref{SS:1.3.3}$. The ``shape" of these terms might drastically vary with $n$ and $r$, though they can be united in the
``families" ~\eqref{1.4} and ~\eqref{1.10}.

The Lie superalgebras $\widetilde{\fsvect}(0|n)$ and
$\widetilde{\fs\fb}_{\mu}(2^{n-1}|2^{n-1}-1)$ depend on an odd
parameter for $n$ odd.

The exceptional regrading ${\rm Reg}_{\fb}$ of $\fb_{\lambda}(2)$
and the isomorphism $\fh_{\lambda}(2|2)\simeq \fb_{\lambda}(2; 2)$
determine an exceptional grading ${\rm Reg}_{\fh}$ of
$\fh_{\lambda}(2|2)$. These regradings (${\rm Reg}_{\fb}$ and ${\rm
Reg}_{\fh}$, see eq.~\eqref{1.7}
and Subsection~\ref{SS:1.3.1}) are irrelevant as far as the
classification of simple Lie superalgebras is concerned, thanks to the isomorphism. On the
other hand, these regradings contribute to the group of
automorphisms of $\fh_{\lambda}(2|2)\simeq \fb_{\lambda}(2; 2)$.

5) \textbf{Solution of Problem C} is given by Theorem \ref{thdefb}, Theorem~\ref{ThRig} on deformations of serial simple vectorial Lie superalgebras, and Conjecture~ \ref{CjRig} on rigidity of the exceptional simple vectorial Lie superalgebras. The following propositions classify the filtered deformations.

\sssbegin[No filtered deformations of exceptional vectorial Lie
superalgebras]{Proposition}[No filtered deformations of exceptional vectorial Lie
superalgebras]\label{FilDefExc} 
For the
W-graded exceptional simple vectorial Lie superalgebras, the isomorphism ~\eqref{1.5.5} does not take place. The spaces $H^2(\fg_-; \fg)$ and
$H^2(\fh_-; \fh)$ are always distinct in these cases: they consist
of different $\fg_0$- and $\fh_0$-modules.
\end{Proposition}

\begin{proof}[Proof \nopoint] is given in \cite{GLS2}. \end{proof}

\sssbegin[Filtered deformations of simple serial vectorial Lie
superalgebras]{Proposition}[Filtered deformations of simple serial vectorial Lie
superalgebras]\label{FilDef} There are only $4$ series
$\widetilde \cL$ of filtered deformations of simple serial vectorial Lie
superalgebras $\cL$ corresponding to the following W-gradings: 
\begin{equation}
\label{3defor} \cL=\fsvect(0|n), \quad \fsle(n), \quad \fb_1(n), \quad \fb_1(n;n).
\end{equation}
The first three of these deformations 
can be described as follows.
Let $\fg\subset\fvect(n|m)$ be a~subalgebra with $\sdim\fg_-=n|m$,
where $\fg_-:=\oplus_{i<0}\ \fg_i$; let $\Theta$ be
the product of all $m$ odd indeterminates and let $t$ be a~parameter,
$p(t)=m\pmod2$. Then, $\widetilde \cL_t=(1+t\Theta)\cL$ and
\begin{equation}
\label{1.5.4} \text{$\widetilde \cL_t\simeq\begin{cases}\cL&\text{for $t$
even},\\
\cL\otimes\Lambda(t)&\text{for $t$ odd},
\end{cases}$~~is a~
subsuperalgebra of $\begin{cases}\fvect(n|m)&\text{for $t$
even},\\
\fvect(n|m)\otimes\Lambda(t)&\text{for $t$ odd}.
\end{cases}$}
\end{equation}

Moreover, $\widetilde \cL_{t_1}\simeq\widetilde \cL_{t_2}$ whenever $t_1t_2\neq 0$ for $t_1, t_2$ even, and for $t_1=at_2$ with $a\in\Cee$ and $a\neq 0$ for $t_1, t_2$ odd.

\end{Proposition}

\begin{proof}[Proof \nopoint] follows from the classification of deformations of the serial Lie superalgebras, see Section~\ref{S:11}. \end{proof}

\section{Background and basic examples}\label{S:2}

\ssec{Linear algebra in superspaces. Generalities}\label{SS:2.1} A
\textit{superspace} is a~$\Zee /2$-graded space 
$V=V_{\ev}\oplus V_{\od }$, where $\ev,\od\in\Zee/2$. A
\textit{supalgebra} is a~superspace with multiplication.

By $\Pi (V)$\index{$\Pi (-)$, the reversal of parity functor} we denote another copy of
the same superspace $V$, but with the inverted parity, \textit{i.e.}, $(\Pi(V))_{\bar
i}= V_{\bar i+\od }$. We call non-zero elements of $V_{\ev}$ \textit{even} (resp., those of $V_{\od}$ \textit{odd}); we denote the parity function by $p$ and write $p(v)=i$ if $v\in V_i$ and $v\neq 0$.

Let $a=\dim V_{\ev}$, $b=\dim V_{\od
}$. As is clear from Milnor's $K_0$-functor point of view, the \textit{superdimension}\index{Superdimension} of $V$ is $\sdim
V:=a+b\eps $, where $\eps ^2=1$. (Usually, the superdimension is written as a~pair $(a,b)$ or $a|b$;
this notation obscures the fact that $\sdim (V\otimes W)=\sdim V\cdot \sdim
W$.) For any superspace $V$, let $\Pty$\index{$\Pty$} be the parity operator,
\textit{i.e.}, $\Pty^2=\id_V$ and the eigenspaces of $\Pty$ corresponding to
the eigenvalues 1 and $-1$ are $V_{\ev}$ and $V_{\od}$,
respectively.

A superspace structure in $V$ induces a~superspace structure in
the space $\End (V)$. A \textit{basis of a~superspace} is always
a~basis consisting of vectors \textit{homogeneous} with respect to parity; let $\Size:=(p_1,
\dots, p_{\dim V})$ be an ordered collection of the parities of these vectors. We
call $\Size$\index{$\Size$} the \textit{format}\index{Format, of the basis, of the matrix} of (the basis of) $V$. A square
\textit{supermatrix} of format (size) $\Size$ is an $\sdim V\times
\sdim V$ matrix whose $i$th row and $i$th column are of the same
parity $p_i\in\Size$. We set $|\Size|=\dim V$. 

One usually considers one of the simplest formats
$\Size$, \textit{e.g.}, the format $\Size_{st}$ of the form $(\ev , \dots, \ev
; \od , \dots, \od)$ is called \textit{standard}. In this paper we
can do without nonstandard formats. But they are vital in the
classification of systems of simple roots of simple finite-dimensional 
Lie superalgebras, see \cite{Se4}. In particular, the reader
might be interested in them in connection with applications to
$q$-quantization or integrable systems. Systems of simple roots
corresponding to distinct nonstandard formats are related by
so-called odd reflections --- analogs of our nonstandard regradings, see \cite{CCLL, Leb}.


The matrix unit $E_{ij}$\index{$E_{ij}$, the matrix unit} is by default of parity $p_i+p_j$.

The formulas of Linear Algebra are defined via
the \textbf{Sign Rule}:\index{Sign Rule}
\begin{equation}\label{SignRule}
\begin{minipage}[l]{13cm}
\textit{if something of parity $\alpha$ moves past something of parity
$\beta$, the~sign~$(-1)^{\alpha\beta}$ accrues; the formulas defined only on
homogeneous elements are extended to arbitrary elements via linearity}.
\end{minipage}
\end{equation}

Examples of application of the Sign Rule: by setting $[X,
Y]:=XY-(-1)^{p(X)p(Y)}YX$ we get the notion of the \textit{super commutator}
and the ensuing notions of supercommutative and super anticommutative
superalgebras; a~\textit{Lie superalgebra} is a~ superalgebra $A$ that satisfies 
super anticommutativity and the super Jacobi identities:
\[
\begin{array}{l}
{}[x,y]=-(-1)^{p(x)p(y)}[y,x];\\
{}[x,[y,z]]=[[x,y], z]+ (-1)^{p(x)p(y)}[y, [x,z]]~~\text{for any $x, y, z\in A$}.
\end{array}
\]
A \textit{super derivation}
of a~superalgebra $A$ is a~linear mapping $D: A\tto A$ that satisfies
the super Leibniz rule
\[
D(ab)=D(a)b+(-1)^{p(D)p(a)}aD(b)~~\text{for any $a,b\in A$}.
\]

In particular, let $A=\Cee[x]$, where $x=(x_{1}, \dots , x_{n})$, be the free supercommutative
polynomial superalgebra in $x$; let the
superstructure in $A$ be determined by the parities of the indeterminates:
$p(x_{i})=p_{i}$. \textit{Partial derivatives} are defined (with the help of
the Sign and Leibniz Rules) by the formulas
$\pder{x_{j}}(x_{i})=\delta_{ij}$. Clearly, the collection $\fder~ \
\Cee[x]$ of all superderivations of $A$ is a~Lie superalgebra whose
elements are of the form $\sum f_i(x)\pder{x_{i}}$ with $f_i(x)\in
\Cee[x]$ for all~ $i$.

Note that applying the Sign Rule sometimes requires some dexterity.
For example, in non-super setting \textit{skew}- and \textit{anti-commutativity} are synonyms. We have to distinguish between (anti-)commutativity of the superalgebra $A$, and \textit{skew}- or \textit{anti-skew-commutativity}  commutativities of $\Pi(A)$; the Sign Rule applied for any $a,b\in A$ yields the four cases:
\[
\begin{array}{ll}
ba=(-1)^{p(b)p(a)}ab &(\text{super commutativity})\\
ba=-(-1)^{p(b)p(a)}ab &(\text{super anti-commutativity})\\
ba=(-1)^{(p(b)+1)(p(a)+1)}ab &(\text{super skew-commutativity})\\
ba=-(-1)^{(p(b)+1)(p(a)+1)}ab
&(\text{super antiskew-commutativity})\end{array}
\]
Thus, the
\textit{skew} formulas are those that can be ``straightened'' by the
reversal of parity of the space on which the multiplilcation is considered,
whereas the prefix \textit{anti} requires an overall minus sign regardless of
parity.

\textbf{Convention}.\index{Convention} In what follows, the supercommutative
superalgebra is always supposed to be associative with unit 1; morphisms of superalgebras are the even homomorphisms;
morphisms of supercommutative superalgebras send 1 into 1.

Let the indeterminates $x=(x_1, \ldots , x_{n+m})$, so that $n$ of them are
even and $m$ of them odd, generate the supercommutative superalgebra
$\cF$\index{$\cF$, superalgebra of functions} of ``functions'' (polynomials or series, etc.). Define the
supercommutative superalgebra $\Omega$\index{$\Omega$, superalgebra of differential forms} of differential forms as
the \textit{polynomial} algebra over $\cF$ in the $dx_i$, where $p(d)=\od$. 

Since
$dx_i$ is even for $x_i$ odd, we can go beyond polynomials
in $dx_i$: smooth or analytic functions in the $dx_i$ are called
\textit{pseudodifferential forms}\index{Form, pseudodifferential} on the supermanifold with
coordinates $x_i$, see \cite{BL2}. We will need such pseudo-forms to interpret
$\fb_{\lambda}(n)$. The exterior differential is defined on the
space of (pseudo) differential forms by the formulas (mind the Super
and Leibniz rules):
\[
d(f):=\sum dx_{i}\nfrac{\partial f}{\partial x_{i}}\text{~~and~~} d^2=0.
\]
The Lie derivative is defined (minding the same rules) by the formula
\[
L_{D}(df)=(-1)^{p(D)}d(D(f)).
\]
In particular,
\[
L_{D}\left
((df)^{\lambda}\right)=\lambda(-1)^{p(D)}d(D(f))(df)^{\lambda-1}\text{
for any $\lambda\in\Cee$ and $f$ odd}.
\]

\ssec{What the Lie superalgebra is}\label{SS:2.2} When dealing with
superalgebras it sometimes becomes useful to know their definition.
Lie superalgebras were distinguished in topology in the 1940s, so an offer of a~``better than usual'' definition of
a~notion which seemed to have been known since then might
look strange, to say the least. Nevertheless, the answer to the
question ``what is a~Lie superalgebra?'' is still not commonly
known.

So far, we operated naively: via the Sign Rule. However, the naive
definition suggested above (``apply the Sign Rule to the definition
of the Lie algebra'') is manifestly inadequate for considering 
supervarieties of deformations and for applications of
representation theory to mathematical physics, for example, in the
study of the coadjoint representation of a~ given Lie supergroup which
can act on a~supervariety, but never on a~vector superspace (an
object from another category). So, to deform Lie superalgebras and
apply group-theoretical methods in the ``super'' setting, we must be
able to recover a~supervariety (a~ringed space) from a~superspace (an object of the Linear Algebra), and vice versa.

A supervariety isomorphic to the ringed space whose base is a~variety $M$ and whose structure sheaf is the sheaf of sections of the Grassmann algebra $ \Lambda^{\bcdot}(E)$ of a~vector bundle $E$ over $M$, is called \textit{split}. Observe that every object in the category of smooth supervarieties (= supermanifolds) is split (for a~short proof of this fact, first established in \cite{Ga, Bat}, see \cite[Subsection 4.1.3]{MaG}), and therefore there is a~one--to--one correspondence between the set of objects in the category of vector bundles over manifolds $M$ and the set of objects in the category of smooth supermanifolds. The latter category has, however, many more morphisms than the former, see \cite[formulas ~(5) in Subsection~3.1]{Ld}.

A(n affine) \textit{superscheme}\index{Superscheme} $\Spec C$ --- a~purely algebraic version of the above over any field (or any commutative ring) of arbitrary characteristic --- is defined in \cite{L0} exactly as the affine scheme (see \cite{MaAG}): its points are prime ideals defined literally as in the commutative case, \textit{i.e.}, $\fp\subsetneq C$ is \textit{prime}\footnote{K.~Coulembier pointed out to us that the so far conventional definition in the non-commutative case is at variance with the common sense: at the moment, if condition~\eqref{primeID} holds, $\fp$ is called (say, in Wikipedia) \textit{completely prime} while it would be natural to retain the term \textit{prime}, as is done in \cite{L0}, since the definition is the same despite the fact that supercommutative rings are not commutative, whereas the term \textit{prime} refers (so far) to any ideal $P\subsetneq R$ in a~non-commutative ring $R$ with the property that
for any two ideals $A$ and $B$ in $R$ the following version of condition~\eqref{primeID} holds: 
\[ 
\text{if $AB\subset P$, then either $A\subset P$, or $B\subset P$}.
\]
} 
\be\label{primeID}
\text{if $a,b\in C$ and $ab\in \fp$, then either $a\in \fp$, or $b\in\fp$}.
\ee
The spectrum of the supercommutative superalgebra is a~ringed space, where \textbf{the left localization with respect to the multiplicative system $C\setminus \fp$ is equivalent\footnote{No proof of this fact was ever published to this day, so some authors try to evade discussing this fact and describe localizations only with respect to purely even multiplicative systems, thus giving a~wrong definition of the superscheme.} to the right localization, \textit{i.e.}, $cd^{-1}=d^{-1}c$ for any $c\in C$ and $d\in C\setminus \fp$}.

To every linear (a.k.a. vector) superspace $V$ there corresponds the \textit{linear supervariety} or supermanifold (in the $C^\infty$-category) which is a ringed space
\be\label{linSman}
\cV=(V_\ev, \ \ \cO_{V_\ev}\otimes \Lambda^{\bcdot}(V_\od^*)).
\ee
The \textit{supervariety of linear homomorphisms} $V\tto W$ is the linear supervariety corresponding to the superspace $\Hom(V,W)$. The \textit{morphisms} of linear superspaces constitute the space $\underline\Hom(V, W):=(\Hom(V,W))_\ev$.

\textbf{A proper definition of Lie superalgebras is as follows}. A
\textit{Lie superalgebra} is a~\textbf{Lie superalgebra in the category of supervarieties}. Instead of expressing the boldfaced requirement in terms of commutative diagrams (and how to check their commutativity?!) one uses the \textit{functor of points}, see \cite{Ld}. 
Namely, any Lie superalgebra $L$ in the sense of the naive definition represents a~functor $\ScommSalgs\leadsto \LieSAlgs$ from the category of supercommutative superalgebras to the category of Lie superalgebras (over the same ground field),
\[
\text{$C\longmapsto L(C):=L\otimes_{\Kee}C$ for any $C\in\ScommSalgs$},
\]
such that
\be\label{funktor}
\begin{minipage}[l]{14cm}
to any morphism of supercommutative superalgebras $C\tto C'$, there corresponds a~morphism of Lie superalgebras $L\otimes C\tto L\otimes C'$ so that a~composition of morphisms $C\tto C'\tto C''$ goes into the composition of morphisms $L\otimes C\tto L\otimes C'\tto L\otimes C''$;\ \
the identity map goes into the identity map.
\end{minipage}
\ee

In these terms, a~ \textit{homomorphism of Lie superalgebras} $\rho: L_1 \tto L_2$ is a~functor morphism, \textit{i.e.}, a~collection of Lie algebra
homomorphisms $\rho_C: L_1 (C)\tto L_2(C)$ such that any
homomorphism of supercommutative superalgebras $\varphi: C\tto C_1$
induces a~Lie algebra homomorphism ${\varphi: L(C)\tto L(C_1)}$, and
products of such homomorphisms are naturally compatible. In
particular, a~\textit{representation} of a~Lie superalgebra $L$ in
a~superspace $V$ is a~Lie superalgebra homomorphism $\rho: L\tto \fgl (V)$, \textit{i.e.},
a~collection of Lie superalgebra homomorphisms $\rho_C: L(C) \tto \fgl
(V )\otimes C$.

Observe that if the characteristic of the ground field is not equal to 3 or 2, one can use another functor of points, called ``Even rules" in \cite[Ch.1]{Del}, cf. \cite{KLLS}.

\sssec{Deformations with odd parameters. Deforms}{}~{} Which of the infinitesimal deformations can be extended to
a~global one is a~separate question, much tougher than the search of infinitesimal deformations, and it is usually solved
\textit{ad hoc}. For examples over fields of characteristics $3$ and
$2$, see \cite{BLW, BLLS, BGL1} and references therein. Deformations with odd parameters are always integrable.

1) \textbf{Deformations of representations}. Consider a~representation $\rho:\fg\tto\fgl(V)$.
The tangent space of the moduli superspace of deformations of $\rho$
is isomorphic to $H^1(\fg; V\otimes V^*)$. For example, if $\fg$ is
the $0|n$-dimensional (\textit{i.e.}, purely odd) Lie superalgebra (with the
only bracket possible: identically equal to zero), its only
irreducible representations are the 1-dimensional trivial one,
$\mathbbmss{1}$, and $\Pi(\mathbbmss{1})$. Clearly,
\[
\mathbbmss{1}\otimes \mathbbmss{1}^*\simeq
\Pi(\mathbbmss{1})\otimes \Pi(\mathbbmss{1})^*\simeq \mathbbmss{1},
\]
and, because the Lie superalgebra $\fg$ is commutative, the
differential in the cochain complex is zero. Therefore,
\[
H^1(\fg; \mathbbmss{1})\simeq\fg^*,
\]
so there are $\dim\, \fg$ of different odd parameters of deformations of the
trivial representation. If we consider~ $\fg$ ``naively'', we miss and lose all 
these deformations.

2) \textbf{Deformations of the brackets}.
Let $C$ be a~supercommutative superalgebra, let $\Spec C$ be the affine superscheme defined literally as the affine scheme (see \cite{MaAG}) of any commutative ring, see \cite{L0} later expounded in \cite{LSoS, MaG}.

Recall, see \cite{Ru}, where the non-super case is considered, that a~\textit{deformation} of 
a Lie superalgebra $\fg$ over $\Spec C$, is a~$C$-algebra $\fG$ such that $\fG\simeq\fg\otimes C$, as spaces. The deformation is \textit{trivial} if $\fG\simeq\fg\otimes C$, as Lie superalgebras, not just as spaces, and \textit{non-trivial} otherwise.

In particular, consider a~deformation with an odd parameter $\tau$. This is a~Lie superalgebra $\fG$ isomorphic to $\fg\otimes\Kee[\tau]$ as a~\textbf{super space}; if, moreover, $\fG$ is isomorphic to $\fg\otimes\Kee[\tau]$ as a~\textbf{Lie superalgebra}, \textit{i.e.},
\[
[a\otimes f, b\otimes g]=(-1)^{p(f)p(b)}[a,b]\otimes fg \text{~~for all~ $a,b\in \fg$ and $f,g\in\Kee[\tau]$},
\]
then the deformation is considered \textit{trivial} (and \textit{non-trivial} otherwise). Note that $\fg\otimes \tau$ is not an ideal of $\fG$: the ideal should be a~free $\Kee[\tau]$-module. 

\textbf{Comment}. In a~sense, the people who ignore odd parameters of deformations have a~point: they consider classification of simple Lie superalgebras over the ground field $\Kee$, right? This can be considered correct, but it is not right: the deformations parameterized by $(H^2(\fg;\fg))_\od$ are no less natural than the odd part of the deformed Lie superalgebra itself. To take these parameters into account, we have to consider everything not over $\Kee$, but over $\Kee[\tau]$. We do the same, actually, when $\tau$ is even and we consider formal deformations over $\Kee[[\tau]]$. If $\tau$ is even and the formal series converges in a~domain $D$, then we can evaluate the series and consider copies $\fg_\tau$ of $\fg$ for every~$\tau\in D$. 

When the parameter is formal (the series diverges) or odd, the only evaluation possible is $\tau=0$.

Examples that lucidly illustrate why one should always remember that
a~Lie superalgebra is not a~mere linear superspace, but a~linear
supervariety, are, \textit{e.g.}, deformations $\widetilde{\fsvect}(0|2n+1)$ and
$\widetilde{\fs\fb}_{\mu}(2^{2n}-1|2^{2n})$ with odd parameters. These deformations are simple Lie superalgebras in
the category of supervarieties.

\ssec{Examples of Lie superalgebras}\label{SS:2.3}{}~{}

\sssec{Algebras of operators} The \textit{general linear} Lie
superalgebra of all supermatrices of size $\Size$ is denoted by
$\fgl(\Size)$, where $\Size=(p_1, \dots, p_{|\Size|})$ is the ordered
collection of parities of the rows identical with the ordered
collection of parities of the columns;
usually, for the \textit{standard} (simplest in a~sense) format, $\fgl(\ev,
\dots, \ev, \od, \dots, \od)$ is abbreviated to $\fgl(\dim V_{\bar
0}|\dim V_{\bar 1})$. Any non-zero supermatrix from $\fgl(\Size)$ can
be uniquely expressed as the sum of its even and odd parts; in the
standard format, this is the following block expression with $p$ defined for $A,B,C,D$ non-zero: 
\[
\mmat A,B,C,D,=\mmat A,0,0,D,+\mmat 0,B,C,0,,\quad
 p\left(\mmat A,0,0,D,\right)=\ev, \ \ p\left(\mmat 0,B,C,0,\right)=\od.
\]

The \textit{supertrace} is the mapping 
\[
\text{$\fgl (\Size)\tto \Cee$,\index{Supertrace}\ \ 
$A=(A_{ij})\longmapsto \sum (-1)^{p_{i}(p(A)+1)}A_{ii}$,} 
\]
where $\Size=(p_{1},
\dots, p_{|\Size|})$. Since $\str [x, y]=0$, the subsuperspace of supertraceless
matrices constitutes the \textit{special linear} Lie subsuperalgebra
$\fsl(\Size)$.

There are, however, at least two super versions of $\fgl(n)$, not
one. Another version is the \textit{queer}
Lie superalgebra $\fq_J(n)$ defined as the Lie superalgebra that preserves the
complex structure given by an \textit{odd} operator $J$ in the $n|n$-dimensional superspace $V$, \textit{i.e.},
$\fq_J(n)$ is the centralizer $C(J)$ of $J$\index{$\fq_J(n)$, queer Lie superalgebra} in $\fgl(V)=\fgl(n|n)$:
\[
\fq_J(n)=C(J)=\{X\in\fgl(n|n)\mid [X, J]=0 \}, \text{ where $p(J)=\od$ and $
J^2=-\id$}.
\]
By a~change of basis one can reduce $J$ to the following normal shape
in the standard
format: $J=J_{2n}:=\mat {0&1_n\\-1_n&0}$, and then $\fq_J(n)$ takes the form
\begin{equation}\label{q}
\fq(n)=\left \{\mat {A&B\\B&A}\mid A, B\in\fgl(n)\right\}.
\end{equation}
The nonstandard formats of $\fq(\Size)$ for $|\Size|=2n$ are of the
form ~\eqref{q} with $A, B\in\fgl(\Size)$.

On $\fq(n)$, the \textit{queertrace} is defined by the formula
\[
\qtr\colon \mmat
A,B,B,A,\longmapsto \tr B, 
\]
first introduced in \cite{BL5}. Denote by $\fsq(n)$ the Lie superalgebra
of \textit{queertraceless} matrices.\index{Queer trace}

Note that the tautological representations of $\fq(V)$ and
$\fsq(V)$ in $V$, though irreducible in super sense (no invariant
sub\textbf{super}spaces), are not irreducible in the nongraded sense. Indeed, take
homogeneous (with respect to parity) and linearly independent
vectors $v_1, \dots , v_n$ from $V$; then, $\Span (v_1+J(v_1),
\dots , v_n+J(v_n))$ is an invariant subspace of $V$ which is not a
subsuperspace.

The representation of a~superalgebra is said to be \textit{irreducible}
\index{Representation, irreducible of type $G$ and $Q$}
\index{$G$-type irreducible representation} of
\textit{general type}, or just of \textit{$G$-type}, if there is no
invariant subspaces; the representation is called \textit{irreducible of
$Q$-type} \index{$Q$-type irreducible representation} if it has no invariant sub\textit{super}spaces, but
\textbf{has} an invariant subspace.

\paragraph{On normal shapes of supermatrices $J$ preserved by $\fq_J(V)$} Over an algebraically non-closed field $\Kee$, there are several types of queer Lie superalgebras $\fq_J(V)\subset\fgl(V)$, defined as centralizers of an odd operator $J$ such that $J^2=\lambda \id_V$ for a~$\lambda\in\Kee^\times$. Over algebraically closed fields, all these superalgebras $\fq_J(V)$ are isomorphic as is easy to see. Over $\Cee$, we can select 
$\lambda=-1$ and interpret $\fq_J(V)$ as the Lie subsuperalgebra of $\fgl(V)$ preserving the complex structure given by an odd operator~ $J$. For real forms of $\fq_J(V)$, see \cite{Se1} with details of the proof in \cite[Ch. D5]{B3}. Over algebraically closed fields of characteristic $2$, there is also one class of normal shapes, we take $\lambda=1$, see \cite{BGLLS}.

\subsubsection{Supermatrices of bilinear forms}\label{sssBilMatr} To a~given linear mapping $F: V\tto W$ of superspaces there corresponds the
dual mapping $F^*:W^*\tto V^*$ between the dual superspaces. Let a~basis of the superspace $V$ have format $\Size$ and
consist of vectors $v_{i}$; then, the formula
$F(v_{j})=\sum_{i}v_{i}X_{ij}$ assigns to the endomorphism $F$ the
supermatrix $X$. In the dual basis, the \textit{supertransposed}\index{Supertransposition}
matrix $X^{st}$ corresponds to the dual operator $F^*$:
\[
(X^{st})_{ij}=(-1)^{(p_{i}+p_{j})(p_{i}+p(X))}X_{ji}.
\]
Consider $X\in\fgl_\cC(\Size)$ over a~supercommutative superalgebra $\cC$, so the entries of the matrices $A,B,C,D$ belong to $\cC$ for any parity of $X$. In the standard format, this means that
\[
X=\mmat A,B,C,D,\longmapsto X^{st}:=\begin{cases}\begin{pmatrix}A^t&C^t\\-B^t&D^t\end{pmatrix},&\text{if $p(X)=\ev$},\\
\begin{pmatrix}A^t&-C^t\\B^t&D^t\end{pmatrix},&\text{if $p(X)=\od$}.\\
\end{cases}
\]

The supermatrices $X\in\fgl(\Size)$ such that
\[
X^{st}B+(-1)^{p(X)p(B)}BX=0\quad \text{for a~ homogeneous matrix
$B\in\fgl(\Size)$}
\]
constitute the Lie superalgebra $\faut (B)$ that preserves the
bilinear form $\cB$ on $V$ whose Gram matrix $B=(B_{ij})$ is given by the formula, see \cite[Ch.1]{LSoS},
\be\label{martBil}
B_{ij}:=(-1)^{p(B)p(v_i)}\cB(v_{i}, v_{j})\text{~~~~for the basis vectors $v_{i}\in V$}.
\ee
The formula \eqref{martBil} makes it possible to identify the bilinear form $\cB\in\Bil(V, W)$ with an element of $\Hom(V, W^*)$. 

Consider the \textit{upsetting} of bilinear forms\index{Upsetting of bilinear forms}
$u\colon\Bil (V, W)\tto\Bil(W, V)$, see \cite[Ch.1]{LSoS}, given by the formula 
\be\label{susyB}
u(\cB)(w, v):=(-1)^{p(v)p(w)}\cB(v,w)\text{~~for any $v \in V$ and $w\in W$.}
\ee
In terms of the Gram matrix $B$ of $\cB$, the form
$\cB$ is \textit{symmetric} if and only if 
\be\label{BilSy}
u(B)=B,\;\text{ where $u(B)=
\mmat R^{t},(-1)^{p(B)}T^{t},(-1)^{p(B)}S^{t},-U^{t},$ for $B=\mmat R,S,T,U,$.}
\ee
Similarly, \textit{anti-symmetry} of 
$\cB$ means that $u(B)=-B$.

Thus, we see that the \textit{upsetting} of bilinear forms corresponds to a~new operation on supermatrices, \textit{not} to the supertransposition.

Observe
that \textbf{the passage from $V$ to $\Pi (V)$ turns every symmetric (resp., anti-symmetric)
form $\cB$ on $V$ into an anti-symmetric (resp., symmetric) form $\cB^\Pi$ on $\Pi (V)$ by setting}
\[
\text{$\cB^\Pi(\Pi(x), \Pi(y)):=(-1)^{p(B)+p(x)+p(x)p(y)}\cB(x,y)$ for any $x,y\in V$}.
\]

Most popular normal shapes of the even non-degenerate symmetric
form are the following ones in the standard format:\index{$J_{2n}$}\index{$\Pi_{n}$}\index{$B_{ev}(2k\vert n)$, the normal shape of the even form}
\be\label{NorEv}
\begin{array}{l}B'_{ev}(m|2n)= \mmat 1_m,0,0,J_{2n},,\quad \text{where
$J_{2n}=\mmat 0,1_n,-1_n,0,$}\\
B_{ev}(2k|2n)= \mmat \Pi_{2k},0,0,J_{2n},, \quad \text{where
$\Pi_{2k}=\mmat 0,1_k,1_k,0,$},\\
B_{ev}(2k+1|2n)= \mmat
\Pi_{2k+1},0,0,J_{2n},, \quad \text{where
$\Pi_{2k+1}=\mat{0&0&1_k\\
0&1&0\\
1_k&0&0}$}.
\end{array}
\ee

The usual notation for the algebra $\faut (B_{ev}(m|2n))$ preserving a~non-degenerate symmetric even form is $\fosp(m|2n)$;
sometimes we write more precisely: $\fosp^{sy}(m|2n)$. 

The anti-symmetric non-degenerate even bilinear form is preserved by the
``symplectic-ortho\-go\-nal" Lie superalgebra, $\fsp\fo (2n|m)$ or $\fosp^{a}(m|2n)$, which is isomorphic to
$\fosp^{sy}(m|2n)$, but has a~different matrix realization.
Let $\fo(m):=\faut(\Pi_m)$ and $\fsp(2n):=\faut(J_{2n})$. In the standard format, we have:
\begin{gather}
\fosp (m|2n)=\left\{
 \mat{ A&B\\
C&D} \mid A\in\fo(m),\ \ 
D\in \fsp(2n),\ \ C=J_{2n}B^t\Pi_m\ \text{for any}\ B\right\}.
\end{gather}

The normal shapes of the odd non-degenerate bilinear\index{$B_{odd}(n\vert n)$, the normal shape of the odd form}
forms in the standard format (NB: the signs are correct) are 
\be\label{NorOd}
B_{odd}(n|n)=\begin{cases}J_{2n},&\text{if the form is symmetric},\\
\Pi_{2n},&\text{if the form is anti-symmetric}.\\
\end{cases}
\ee 
The present day usual notation for
$\faut (B_{odd}(\Size))$ is $\fpe(\Size)$.
The passage from $V$ to $\Pi (V)$ establishes an isomorphism
$\fpe^{sy}(\Size)\simeq \fpe^{a}(\Size)$. These isomorphic Lie
superalgebras are called, as A.~Weil suggested,\footnote{``Don't you
know the ancient Greek?!" he asked.} \textit{periplectic}.\index{$\fpe (n)$} In the
standard format, the matrix realizations of these superalgebras is:
\begin{equation}
\begin{split}
\fpe ^{sy}\ (n)&=\left\{\mat{ A&B\\
C&-A^t}\mid A\in\fgl(n),\, B=-B^t,\;
C=C^t\right\};\\
\fpe^{a}(n)&=\left\{\mat {A&B\\ C&-A^t} \mid A\in\fgl(n),\, B=B^t,\;
C=-C^t\right\}.
\end{split}
\end{equation}
Note that, although $\fosp^{sy} (m|2n)\simeq\fosp ^{a} (m|2n)$ and $\fpe ^{sy} (n)\simeq\fpe ^{a}
(n)$, certain properties of these incarnations are sometimes
very different, see Subsection~ \ref{2.6},
and
\cite{Sh5}.

The \textit{special periplectic} superalgebra is\index{$\fspe (n)$}
\[
\fspe(n):=\{X\in\fpe(n)\mid \str
X=0\}.
\]
Of particular interest to us will be the Lie superalgebras\index{$\fspe_{a, b}(n)$}
\begin{equation}
\label{spe}
 \fspe_{a, b}(n):=\fspe(n)\ltimes \Cee(aD+bz), \text{~~~
where $z=1_{2n}$, $D=\diag(1_{n}, -1_{n})$}
\end{equation}
and the nontrivial central extension $\fas$\index{$\fas$} of $\fspe(4)$ that
will be described after some preparation, see Subsection~\ref{ssas}.

\subsubsection{More notation} The tautological representation 
of a~matrix Lie superalgebra $\fg\subset \fgl(V)$ in $V$, and sometimes \textbf{the module $V$ itself}, are denoted by $\id$ or, for clarity,
$\id_{\fg}$. The context prevents confusion of this notation with
that of the identity (scalar) operator $\id_{V}$ on the space $V$,
as in the next paragraph:

For $\fg=\oplus_{i\in\Zee}\ \fg_{i}$, the trivial
representation of $\fg_0$ is denoted by $\Cee$ (if $\fg_0$ is
simple) whereas $\Cee[k]$ denotes the representation of $\fg_0$ that is 
trivial on its semi-simple part and such that $k$ is the value with which the
central element $z\in\fg_{0}$, where $z$ is chosen so that
$z|_{\fg_i}=i\cdot \id_{\fg_i}$, acts on $\Cee[k]$.

\paragraph{Projectivization}\label{SS:2.3.1} If $\fc$ is a~Lie algebra of scalar matrices, and $\fg\subset \fgl (n|n)$ is a~Lie subsuperalgebra containing $\fc$, then the \textit{projective}\index{Projectivization}
Lie superalgebra of type $\fg$ is $\fpg= \fg/\fc$. The Lie
superalgebras $\fg_1\bigodot \fg_2$ described in Subsection~\ref{ss4cases} are also
projective.

Projectivization sometimes leads to new Lie superalgebras, for
example: $\fpq (n)$, $\fpsq (n)$; $\fpgl (n|n)$, $\fpsl (n|n)$;
note that $\fpgl (p|q)\simeq \fsl (p|q)$ if $p\neq q$.

\ssec{Simple, almost simple and semi-simple Lie
superalgebras}\label{SS:2.2.1} Recall that any Lie superalgebra
$\fg$ without proper ideals and of dimension $>1$ is said to be
\textit{simple}. The list~\eqref{simple} contains all simple
finite-dimensional Lie superalgebras over $\Cee$, for various
subcases of the classification, see \cite{K1, Kapp, K1C, SNR, Dj1, Dj2, Dj3}; for a
summary ignoring odd parameters of deformations, see \cite{K2, Sch}; for a summary and classification of all deformations, see \cite{Ld}:
\begin{equation}\label{simple}
\begin{array}{l} \text{$\fsl(m|n)$ for $m> n\geq 1$, $\fpsl(n|n)$ for
$n>1$, $\fosp(m|2n)$ for $mn\neq 0$};\\
\text{$\fpsq(n)$ and $\fspe(n)$
for $n>2$;} \\
\text{a 1-parameter family of deformations $\fosp_a(4|2)$ for $a\neq 0, -1$};\\
\text{exceptional Lie superalgebras $\fa\fg(2)$ and
$\fa\fb(3)$;}\\
\text{finite-dimensional vectorial Lie superalgebras, namely}\\ 
\text{$\fvect(0|n)$ for $n>1$, $\fsvect(0|n)$ and
$\widetilde\fsvect(0|2n)$ for $n>2$, $\fh'(0|n)$ for $n>3$;}\\
\widetilde\fsvect(0|2n+1)\text{~~ for $n\geq 1$; the only (mind \eqref{ocIso}) deformation with an odd parameter.}\\
\end{array}
\end{equation}
For descriptions of these algebras, especially those with Cartan
matrices, see \cite{BGL,CCLL,GL1, GL3}, and \cite{K2,vdL}; for descriptions of the exceptional algebras, see \cite{DWN,Sud, El3, BE}.

\textbf{Occasional isomorphisms of finite-dimensional simple Lie superalgebras}:\index{Isomorphisms, occasional}
\be\label{ocIso}
\begin{array}{l}
\fsl(m|n)\simeq \fsl(n|m),\ \ \fvect(0|2)\simeq \fsl(1|2)\simeq\fosp(2|2),\\
\fspe(3)\simeq\fsvect(0|3),\ \ \fpsl(2|2)\simeq\fh'(0|4).\\
\end{array}
\ee
The isomorphisms $\fosp_{\alpha}(4|2)\simeq \fosp_{\alpha'}(4|2)$ are generated by the
transformations (see \cite{BGL}):
\begin{equation}\label{osp42symm}
\alpha\longmapsto \alpha':= -1-\alpha\ , \qquad \alpha\longmapsto \alpha':=
\nfrac{1}{\alpha}\ ,
\end{equation}
so the other isomorphisms are 
\begin{equation}\label{osp42symm1}
\alpha\longmapsto \alpha':=
-\nfrac{1+\alpha}{\alpha}\ ,\quad\alpha\longmapsto \alpha':=
-\nfrac{1}{\alpha+1}\ ,\quad\alpha\longmapsto \alpha':= -\nfrac{\alpha}{\alpha+1}.
\end{equation}


We say that a~Lie superalgebra $\fg$ is \textit{almost simple} if it can be sandwiched
(non-strictly) between a~simple Lie superalgebra $\fs$ and the Lie\index{Lie (super)algebra, simple}\index{Lie (super)algebra, almost simple}
superalgebra $\fder~\fs$ of derivations of $\fs$, \textit{i.e.},
\[
\fs\subset\fg\subset\fder~\fs.
\]

By definition, $\fg$ is \textit{semi-simple} if its \textit{radical} --- 
the maximal solvable ideal --- is $\{0\}$. \index{Lie (super)algebra, semi-simple}

Copying R.~Block's description of semi-simple Lie algebras over fields of prime
characteristic we describe
semi-simple Lie superalgebras as follows, cf. \cite{K2}.

Let $\fs_1, \dots , \fs_k$ be simple Lie superalgebras, let $n_1,
\dots , n_k$ be \textit{pairs} of nonnegative integers
$n_j=(n_j^\ev, n_j^\od)$, let $\cF(n_j)$ be the supercommutative
superalgebra of polynomials in $n_j^\ev$ even and $n_j^\od$ odd
indeterminates, and let $\fs=\oplus_{j}\left
(\fs_j\otimes\cF(n_j)\right)$. Then,
\begin{equation}
\label{2.1.3}
\fder~\fs=\mathop{\oplus}\limits_{j}\left((\fder~\fs_j)\otimes\cF(n_j)
\ltimes \id_{\fs_j}\otimes \fvect(n_j)\right). 
\end{equation}
{\sl Let $\fg$ be a~subalgebra of $\fder~\fs$ containing $\fs$. If
\begin{equation}
\label{proj} \text{the projection of $\fg$ to
$\id_{\fs_j}\otimes\fvect(n_j)_{-1}$ is \textbf{onto} for each $j$,}
\end{equation}
then $\fg$ is semi-simple}.

\sssec{Conjecture}\label{conj}\index{Conjecture} If the
(super)algebras $\fs_j$ are finite-dimensional or of polynomial growth, and
the (super)algebras of functions $\cF(n)$ depend on a~finite number of indeterminates, then
all semi-simple Lie superalgebras arise in the manner indicated,
\textit{i.e.}, as sums of subalgebras of the summands of ~\eqref{2.1.3}
satisfying condition~\eqref{proj}, cf. \cite{Che}.

\ssec{The irreducible $\fg$-modules of least
dimension $\neq 1$ over $\fg$ simple finite-dimen\-sio\-nal}\label{fact} Up to the~reversal of parity, such modules are: the adjoint one for
$\fg=\fpsq(n)$ and $\fpsl(n|n)$ for $n>2$, $\fag(2)$, $\fab(3)$ and $\fosp_a(4|2)$ for $a$ generic (\textit{i.e.}, $a\neq 1,2,3$ or obtained from such values of $a$ under the $S_3$-action, see formulas~ \eqref{osp42symm1}). 

For $a= 1,2,3$, the irreducible modules of least dimension are, up to the reversal of parity $\Pi$, of superdimension
$4|2$, $6|4$ and $8|6$, respectively, see \cite{GL3}. 
This fact was
first claimed in \cite{Kpr} without proof. The above claims are now deducible from results of Penkov and Serganova, and of Brundan. 


\ssec{Remarks on Lie (super)algebras} 1) Not all filtered Lie superalgebras $\cL$ or associated with them $\Zee$-graded Lie superalgebras $L:=\mathop{\oplus}_{i>-d}L_i$, where $L_i:=\cL_i/\cL_{i+1}$,
of finite depth are \textit{vectorial}, \textit{i.e.}, realizable by vector
fields on a~supermanifold of the same superdimension as that of
$\cL/\cL_{0}$; only those with faithful $L_{0}$-action on
$L_{-}:=\mathop{\oplus}_{i<0}\ L_i$ are.

\textbf{Unlike Lie algebras, simple infinite-dimensional vectorial \textit{super}algebras can have
\textit{several} non-isomor\-phic maximal subalgebras of finite
codimension}.

2) We denote the vectorial Lie (super)algebras as is customary in differential geometry: $\fvect (x)$ or $\fvect (V)$ if $V=\Span(x)$ for $x=(u_1, \dots
, u_n, \theta_1, \dots , \theta_m)$, where the $u_i$ are even
indeterminates and the $\theta_j$ are odd ones. We use a~similar
notation for the subalgebras of $\fvect$ introduced below. Some
algebraists sometimes abbreviate $\fvect (n):=\fvect (n|0)$ by
$W(n)$ in honor of Witt\footnote{Witt considered a~version of $\fvect
(n)$ over fields of characteristic~ $p>0$.}, and its subalgebra $\fs\fvect (n)$ of divergence-free fields
to $S(n)$, where $S$ stands
for ``special". Following Bourbaki, we mostly use Gothic font to name Lie (super)algebras.

3) Considering the subspaces $\cL_{i}$ as the basis of a~topology,
we can complete the graded or filtered Lie superalgebras $L$ or
$\cL$; the elements of the completion are the vector fields with
formal power series as coefficients. Though the structure of the
graded algebras is easier to describe, in applications one usually needs the completed
vectorial Lie superalgebras.


\sssec{Traces and divergences on vectorial Lie
superalgebras}\label{traceDiv} Recall Subsection \ref{over}. 

On any Lie (super)algebra $\fg$ over any field
$\Kee$, a~\textit{trace}\index{Trace, supertrace} is any linear mapping $\tr: \fg\tto \Kee$ such
that
\begin{equation}\label{deftr}
\tr (\fg')=0.
\end{equation}

Let $\Kee=\Cee$. Let $\cF$ be the superalgebra of ``functions". Let $\fg:=\fder~\,\cF$ be a~$\Zee$-graded vectorial Lie algebra with
$\fg_{-}:=\mathop{\oplus}\limits_{i<0}\fg_i$ generated by
$\fg_{-1}$, and let $\tr$ be a~trace on~ $\fg_0$. \textbf{Conjecturally}, any
$\Zee$-grading of a~given vectorial Lie algebra is given by the degrees
of the indeterminates with a~proof similar to that given in \cite{Sk91, Sk95}.\index{Conjecture} 

The \textit{divergence} $\Div:\fg\tto\cF$ is a~degree-preserving\index{Divergence}
$\ad_{\fg_{-1}}$-invariant prolongation of the trace on~ $\fg_0$ satisfying the following conditions,
so $\Div\in Z^1(\fg; \cF)$, \textit{i.e.}, is a~cocycle:
\[
\begin{array}{l}
X_i(\Div D)=\Div([X_i,D])\text{~~ for all elements
$X_i$ that span $\fg_{-1}$};\\
 \Div|_{\fg_0}=\tr;\\
 \Div|_{\fg_{-}}=0.\end{array}
 \]

Let $\Vol:=\cF^*$\index{$\Vol$} be the $\fvect(m |n)$-module of \textit{volume forms} dual to~ $\cF$. As an $\cF$-module, $\Vol$ is generated by the \textit{volume element} $\vvol_x=1^*$ for
fixed indeterminates (``coordinates") $x$ which we often do not
indicate. On the rank-1 $\cF$-module of \textit{weighted}\index{Density, weighted; $\lambda$-density}
$\lambda$-densities $\Vol^{\lambda}(m |n)$ with generator
$\vvol_x^{\lambda}$ over $\cF$, the $\fvect(m |n)$-action is given for any $f\in\cF$ and $D\in\fvect(m |n)$ by the
\textit{Lie derivative}\index{Lie derivative}
\begin{equation}
\label{LieDer}
L_D(f\vvol_x^{\lambda})=\left(D(f)+(-1)^{p(f)p(D)}\lambda
f\Div(D)\right)\vvol_x^{\lambda}.
\end{equation}

Often by a~\textit{volume form} we mean only elements $f\vvol_x$, where $f\in\cF$, such that $f(0)\neq 0$.\index{Volume form}

\ssec{Examples of vectorial Lie superalgebras}{}~{}

 \textbf{1) General vectorial Lie superalgebras}. Let $x=(u_1, \dots
, u_n, \theta_1, \dots , \theta_m)$. 
The \textit{the general vectorial
Lie superalgebra} is defined to be\index{$\fvect$, general vectorial Lie
superalgebra}\index{Lie (super)algebra, vectorial} $\fvect
(n|m):=\fder~\; \Cee[x]$.

\textbf{2) Special Lie superalgebras}. The
\textit{divergence}\index{Divergence} of the field
$D=\mathop{\sum}_if_i\pder{u_{i}} + \mathop{\sum}_j
g_j\pder{\theta_{j}}$ is the function 
\begin{equation}
\label{2.2.4} 
\Div D:=\mathop{\sum}\limits_i\pderf{f_{i}}{u_{i}}+
\mathop{\sum}\limits_j (-1)^{p(g_{j})}
\pderf{g_{i}}{\theta_{j}}.
\end{equation}

$\bullet$ The Lie superalgebra $\fsvect (n|m):=\{D \in \fvect
(n|m)\mid \Div D=0\}$ is called the \textit{special} (or
\textit{divergence-free}) \textit{vectorial superalgebra}.
\index{$\fsvect$, special vectorial Lie superalgebra}\index{Lie
(super)algebra, special vectorial}\index{Lie (super)algebra,
divergence-free}

It is clear that one can also describe $\fsvect(n|m)$ as
$\{ D\in \fvect (n|m)\mid L_D\vvol _x=0\}$. 

$\bullet$ The Lie superalgebra $\fsvect'(1|n):=[\fsvect(1|n),
\fsvect(1|n)]$ is called the \textit{traceless special vectorial
superalgebra}.\index{Lie (super)algebra, special vectorial, traceless}\index{$\fsvect'(1\vert n)$}

$\bullet$ The deformation of $\fsvect(0|m)$ is the Lie superalgebra
\[
\fsvect_{\lambda}(0|m)=\{D \in \fvect (0|m)\mid \Div
(1+\lambda\theta_1\cdot \ldots \cdot \theta_m)D=0\},
\]
where $p(\lambda)\equiv m\pmod 2$ and if $p(\lambda)=\od$, then $\fsvect_{\lambda}(0|m)$ is considered over $\Cee[\lambda]$. It is called the \textit{deformed
special} (or \textit{divergence-free}) \textit{vectorial
superalgebra}. Clearly, $\fsvect_{\lambda}(0|m)\simeq 
\fsvect_{\mu}(0|m)$ if $\lambda\neq 0$ and $\mu\neq 0$. Accordingly, we briefly denote
these deformations by $\widetilde{\fsvect}(0|m)$, see Subsection \ref{FilDef}.\index{$\widetilde{\fsvect}(0\vert m)$}


\textbf{3) Lie superalgebras that preserve Pfaff equations or, equivalently, non-integrable distributions, contact and
 symplectic forms}.
 
$\bullet$ Recall that a~\textit{symplectic form} is a~non-degenerate closed differential 2-form $\omega$.\index{Form differential, symplectic} In agreement with the normal shapes of the bilinear forms ~\eqref{NorEv} and ~\eqref{NorOd}, the descriptions of normal shapes of the symplectic forms (Darboux's theorem, when on manifolds) were formulated in \cite{L1}; this claim was proved in \cite{Shan}. The obstruction to reducing the non-degenerate differential 2-form $\omega$ to the normal shape is $d\omega$, see \cite{BGLS}. If $\omega$ is odd, it is called, as A.~Weil advised, \textit{periplectic}.\index{Form differential, periplectic}\index{Form differential, pericontact}

$\bullet$ The \textit{contact form} (resp., \textit{pericontact form}) is a~nowhere vanishing differential 1-form $\alpha$ on a~(super)manifold $M$ that singles out in the tangent bundle $TM$ a~\textbf{maximally} non-integrable distribution $\cD$ of codimension 1 (resp., $\eps=0|1$), \textit{i.e.}, at every point $m\in M$, there are given\index{Form differential, contact}
\be\label{contForm}
\begin{array}{l}
\text{a~subspace $\cD_m$ of codimension 1 (resp., $\eps=0|1$) in the tangent space $T_mM$, and}\\
\text{a~
differential 1-form $\alpha$ such that $d\alpha|_{\cD_m}$ is non-degenerate.}
\end{array}
\ee

Note that in Differential Geometry on odd-dimensional manifolds (over the ground field $\Kee$), the contact form is often defined as 
a~nowhere vanishing differential 1-form $\alpha$ such that
\be\label{notDef}
\text{$\alpha\wedge (d\alpha)^{\wedge n}=\lambda \vvol_x$ for some $\lambda\in\Kee$. }
\ee
Obviously, this definition cannot be superized since on supermanifolds there are no differential forms of highest degree. Condition \eqref{notDef} is a~\textit{property} of the contact form on a~manifold, whereas \eqref{contForm} is a~reasonable definition.

$\bullet$ Let $u=(t, p_1, \dots , p_n, q_1, \dots , q_n)$ be even indeterminates and $\Theta=(\theta_1,\dots , \theta_m)$ be odd indeterminates. The form\index{$\widetilde \alpha_1$}\index{$ \widetilde
\omega_0:=d\widetilde \alpha_1$}
\begin{equation}
\label{2.2.5} \widetilde \alpha_1 = dt +\mathop{\sum}\limits_{1\leq
i\leq n}(p_idq_i - q_idp_i)\ + \mathop{\sum}\limits_{1\leq j\leq
m}\theta_jd\theta_j\quad\text{and}\quad \widetilde
\omega_0:=d\widetilde
\alpha_1 
\end{equation}
is called \textit{contact}, and the form
$\widetilde \omega_0$ is called \textit{symplectic}.\index{Form differential, symplectic} Sometimes
it is more convenient to redenote the $\theta_j$ and set
$\Theta=\begin{cases}(\xi, \eta),& \text{ if }\ m=2k\\
(\xi, \eta,
\theta),&\text{ if }\ m=2k+1,\end{cases}$ where 
\begin{equation}
\label{zam}
\xi_j=\frac{1}{\sqrt{2}}(\theta_{j}-\sqrt{-1}\theta_{k+j}),\quad
\eta_j=\frac{1}{ \sqrt{2}}(\theta_{j}+\sqrt{-1}\theta_{k+j})\; \text{~
for~}\; j\leq k= \lfloor\nfrac{m}{2}\rfloor,\;
\; \theta =\theta_{2k+1},
\end{equation}
and in place of $\widetilde \omega_0$ or $\widetilde \alpha_1$ take\index{$\alpha_1$, form differential, contact}
\index{$\omega_0:=d \alpha_1$}
$\alpha_1$ and $\omega_0=d\alpha_1$, respectively, where
(we mostly consider the cases where $\ell=0$ or $1$)
\begin{equation}
\label{2.2.6} \alpha_1=dt+\mathop{\sum}\limits_{1\leq i\leq
n}(p_idq_i-q_idp_i)+ \mathop{\sum}\limits_{1\leq j\leq
k}(\xi_jd\eta_j+\eta_jd\xi_j)+
\begin{cases}0,&
\text{ if }\ m=2k\\
\mathop{\sum}\limits_{1\leq s\leq
\ell}\theta_sd\theta_s,&\text{ if }\ m=2k+\ell.\end{cases} 
\end{equation}

The Lie superalgebra that preserves the distribution given by the \textit{Pfaff equation} \index{Pfaff equation} 
\be\label{alpha1}
\text{$\tilde\alpha_1(X)=0$ for $X\in \fvect(2n+1|m)$}
\ee
(or by an equivalent equation $\alpha_1(X)=0$),
 \textit{i.e.}, the superalgebra
\begin{equation}
\label{k} \fk (2n+1|m):=\{ D\in \fvect (2n+1|m)\mid
L_D \alpha_1=f_D \alpha_1\text{ for some }f_D\in \Cee [t, p, q,
\theta]\},
\end{equation}
is called the \textit{contact superalgebra}.\index{$\fk$, contact
superalgebra} \index{Lie (super)algebra, contact} The divergence-free subalgebra of $\fk (2n+1|m)$ is the Lie superalgebra
\begin{equation}
\label{po}
\fpo (2n|m)=\{ D\in \fk (2n+1|m)\mid L_D\alpha_1=0\}.
\end{equation}
It is called the \textit{Poisson} superalgebra.\index{$\fpo$, Poisson
superalgebra} (The Poisson superalgebra can be interpreted as
the Lie superalgebra that, over a~symplectic supermanifold with the
symplectic form $d\alpha_1$, preserves the connection with form
$\alpha_1$ and curvature proportional to $d\alpha_1$ in a~line bundle.)

The above ``symmetric" expression of $\alpha_1$ is popular among
algebraists. However,
in analytical mechanics and differential geometry, the following
expression of the form $\alpha_1$ (without odd coordinates, of
course) naturally appears:
\begin{equation}
\label{alphachar2} \alpha_{1(2)}=dt-\mathop{\sum}\limits_{1\leq
i\leq n}p_idq_i + \mathop{\sum}\limits_{1\leq j\leq k}\xi_jd\eta_j+
\begin{cases}0,&
\text{ if }\ m=2k\\
\theta d\theta,&\text{ if }\ m=2k+1.\end{cases} 
\end{equation}
The form $\alpha_{1(2)}$ is the only reasonable expression of the
contact form over fields of characteristic~ 2, whereas the
symmetric expression of the contact forms is inapplicable:
$d\alpha_1=2\omega=0$.

$\bullet$ Similarly, let $u=q=(q_1, \dots , q_n)$ and $\theta=(\xi_1, \dots , \xi_n; \tau)$ be even 
and odd indeterminates, respectively. Set
\begin{equation}
\label{alp0}
\begin{array}{c}
\alpha_0=d\tau+\mathop{\sum}\limits_i(\xi_idq_i+q_id\xi_i),
\qquad\qquad \omega_1:=d\alpha_0
\end{array}
\end{equation}
and call the form $\omega_1$, as A.~Weil advised, \textit{periplectic}; accordingly, we call $\alpha_0$ \textit{pericontact}.\index{Form differential,
(peri)contact}\index{Form differential, periplectic}\index{$\alpha_0$, form differential, pericontact}
\index{$\omega_1:=d\alpha_0$, form differential, periplectic}


The Lie superalgebra that preserves the distribution given by the Pfaff equation
\be\label{alpha0}
\text{$\alpha_0(X)=0$ for $X\in \fvect(n|n+1)$,}
\ee
 \textit{i.e.}, the superalgebra
\begin{equation}
\label{m} \fm (n)=\{ D\in \fvect (n|n+1)\mid L_D\alpha_0=f_D\cdot
\alpha_0\text{ for some }\; f_D\in \Cee [q, \xi, \tau]\},
\end{equation}
is called the \textit{pericontact superalgebra}.\index{$\fm$,
pericontact superalgebra}\footnote{To call $\fm$ ``odd contact
superalgebra" is misleading, since the even part of $\fm$ is non-zero.
To say ``the ``odd" contact superalgebra'' is OK, though a~bit too long, and the pericontact form is even.}

\sssec{On normal shapes of contact and symplectic/periplectic forms} The
restriction of $\omega_0$ to $\ker \alpha_1$ is the orthosymplectic
form $B_{ev}(m|2n)$; the restriction of $\widetilde \omega_0$ to
$\ker \widetilde \alpha_1$ is $B'_{ev}(m|2n)$. Similarly, the
restriction of $\omega _1$ to $\ker \alpha_0$ is $B_{odd}(n|n)$. Clearly, $\ker \alpha_1= \ker \widetilde \alpha_1$. 

There is just one class of contact forms (both even and odd), whereas the obstruction to reducibility of the non-degenerate bilinear form $\omega$ (even or odd)
to the normal shape is $d\omega$, see \cite{BGLS}.


\subsubsection{Generating functions: contact series and their subalgebras}\label{SS:2.4.1} A laconic way to describe $\fk$,
$\fm$ and their subalgebras is via generating functions. \index{Function, generating} We express the contact fields of series $\fk$ (resp., $\fm$) in natural bases
$K_f$ (resp., $M_f$).

$\bullet$ \underline{Odd form $\alpha_1$}. For any $f\in\Cee [t, p,
q, \theta]$, set\index{$K_f$, contact vector field} \index{$H_f$,
Hamiltonian vector field}
\begin{equation}
\label{2.3.1} 
 K_f=(2-E)(f)\pder{t}-H_f + \pderf{f}{t} E,
\end{equation}
where $E=\mathop{\sum}\limits_i y_i \pder{y_{i}}$ (here the $y_{i}$
are all the coordinates except $t$) is the \textit{Euler operator},\index{$E$, the Euler operator}
and $H_f$ is the Hamiltonian field with Hamiltonian $f$ that
preserves $d\widetilde \alpha_1$:
\begin{equation}
\label{2.3.2} H_f=\mathop{\sum}\limits_{i\leq n}\left(\pderf{f}{p_i}
\pder{q_i}-\pderf{f}{q_i} \pder{p_i}\right )
-(-1)^{p(f)}\left(\mathop{\sum}\limits_{j\leq m}\pderf{ f}{\theta_j}
\pder{\theta_j}\right ) . 
\end{equation}
The choice of the form $\alpha_1$ instead of $\widetilde \alpha_1$
affects only the shape of $H_f$ which we give for $m=2k+\ell$:
\begin{equation}
\label{2.3.2'} H_f=\mathop{\sum}\limits_{i\leq n}\left
(\pderf{f}{p_i} \pder{q_i}-\pderf{f}{q_i} \pder{p_i}\right)
-(-1)^{p(f)}\left(\mathop{\sum}\limits_{j\leq
k}\left(\pderf{f}{\xi_j} \pder{\eta_j}+ \pderf{f}{\eta_j}
\pder{\xi_j}\right)+ \mathop{\sum}\limits_{j\leq \ell}\pderf{ f}{\theta_j}
\pder{\theta_j}\right).
\end{equation}
The expression of the contact field corresponding to the form
$\alpha_1$ or $\widetilde \alpha_1$ is 
\begin{equation}
\label{K_f(2)} K_f=(2-E)(f)\pder{t}-H_f + \pderf{f}{t} E,
\end{equation}
where $E=\mathop{\sum}\limits_i p_i
\pder{p_{i}}+\mathop{\sum}\limits_j \xi_j \pder{\xi_{j}}$, and $H_f$
is the Hamiltonian field with Hamiltonian $f$ that preserves
$d\widetilde \alpha_1$, see ~\eqref{2.3.2}.

 $\bullet$ \underline{Even form $\alpha_0$}. For any $f\in\Cee [q,
\xi, \tau]$, set
\begin{equation}
\label{2.3.3} M_f=(2-E)(f)\pder{\tau}- Le_f -(-1)^{p(f)}
\nfrac{\partial f}{\partial\tau} E, 
\end{equation}
where $E=\mathop{\sum}\limits_iy_i \pder{y_i}$ (here the $y_i$ are
all the coordinates except $\tau$), and\index{$E$, the Euler operator}
\begin{equation}
\label{2.3.4} Le_f=\mathop{\sum}\limits_{i\leq n}\left(
\pderf{f}{q_i}\ \pder{\xi_i}+(-1)^{p(f)} \pderf{f}{\xi_i}\
\pder{q_i}\right).
\end{equation}
\index{$M_f$, pericontact vector field} \index{$Le_f$, periplectic
vector field} Since
\begin{equation}
\label{2.3.5}
\renewcommand{\arraystretch}{1.4}
\begin{array}{l}
 L_{K_f}(\alpha_1)=2 \pderf{f}{t}\alpha_1=K_1(f)\alpha_1, \\
L_{M_f}(\alpha_0)=-(-1)^{p(f)}2 \pderf{
f}{\tau}\alpha_0=-(-1)^{p(f)}M_1(f)\alpha_0,
\end{array}
\end{equation}
it follows that $K_f\in \fk (2n+1|m)$ and $M_f\in \fm (n)$.
Note that
\[
p(Le_f)=p(M_f)=p(f)+\od.
\]

$\bullet$ To the (super)commutators $[K_f, K_g]$ and $[M_f, M_g]$
there correspond \textit{contact brackets}\index{Poisson
bracket}\index{Contact bracket}\index{$\{-, -\}_{k.b.}$, contact bracket}\index{$\{-, -\}_{m.b.}$, pericontact bracket} of the generating functions:
\[
\renewcommand{\arraystretch}{1.4}
\begin{array}{l}
{}[K_f, K_g]=K_{\{f, \; g\}_{k.b.}}\text{~~and~~}
 {}[M_f, M_g]=M_{\{f, \;
g\}_{m.b.}}.\end{array}
\]
To give the explicit formulas for the contact brackets, we first define the brackets on functions that do not depend on $t$
(resp., $\tau$).

The \textit{Poisson bracket} $\{-,-\}_{P.b.}$ in the
realization with the form $\widetilde{\omega}_0$ is given by the
formula\index{$\{-, -\}_{P.b.}$, Poisson bracket} \index{$\{-, -\}_{B.b.}$, Buttin, a.k.a. anti-bracket} 
\begin{equation}
\label{2.3.6} 
\renewcommand{\arraystretch}{1.4}
\begin{array}{ll}
\{f, g\}_{P.b.}&=\mathop{\sum}_{i\leq n}\
\bigg( \pderf{f}{p_i}\
 \pderf{g}{q_i}- \pderf{f}{q_i}
 \pderf{g}{p_i}\bigg)-\\
&(-1)^{p(f)}\mathop{\sum}_{j\leq m}\
 \pderf{f}{\theta_j}\
 \pderf{g}{\theta_j}\text{ for any }f, g\in
\Cee [p, q, \theta]\end{array}
\end{equation}
and in the realization with the form $\omega_0$ for $m=2k+1$ it is
given by the formula
\begin{equation}
\label{pb} 
\renewcommand{\arraystretch}{1.4}
\begin{array}{ll} 
 \{f, g\}_{P.b.}&=\mathop{\sum}\limits_{i\leq n}\ \bigg( \pderf{f}{p_i}\
\pderf{g}{q_i}-\ \pderf{f}{q_i}\
\pderf{g}{p_i}\bigg)-\\ 
&(-1)^{p(f)}\bigg(\mathop{\sum}\limits_{j\leq m}\left(
 \pderf{f}{\xi_j}\ \pderf{ g}{\eta_j}+\pderf{f}{\eta_j}\
\pderf{ g}{\xi_j}\right)+ \pderf{f}{\theta}\ \pderf{
g}{\theta}\bigg)\text{ for any }f, g\in \Cee [p, q, \xi, \eta,
\theta].
\end{array}
\end{equation}

The \textit{Buttin bracket} $\{-,-\}_{B.b.}$ \index{Antibracket$=$Buttin
bracket$=$Schouten bracket} \index{Buttin
bracket$=$Schouten bracket$=$Antibracket}\index{Schouten bracket$=$Antibracket$=$Buttin bracket} is given by the formula
\begin{equation}
\label{2.3.7} \{ f, g\}_{B.b.}=\mathop{\sum}\limits_{i\leq n}\
\bigg(\pderf{f}{q_i}\ \pderf{g}{\xi_i}+(-1)^{p(f)}\
\pderf{f}{\xi_i}\ \pderf{g}{q_i}\bigg)\text{ for any }f, g\in \Cee
[q, \xi].
\end{equation}

 In terms of the Poisson and Buttin brackets,
respectively, the (peri)contact brackets are\index{$\{ f, g\}_{k.b.}$}\index{$\{ f, g\}_{m.b.}$}
\begin{equation}
\label{2.3.8} \{ f, g\}_{k.b.}=(2-E) (f)\pderf{g}{t}-\pderf{f}
{t}(2-E) (g)-\{ f, g\}_{P.b.}
\end{equation}
and
\begin{equation}
\label{2.3.9} \{ f, g\}_{m.b.}=(2-E)
(f)\pderf{g}{\tau}+(-1)^{p(f)} \nfrac{\partial f}{\partial\tau}(2-E) (g)-\{ f,
g\}_{B.b.}. 
\end{equation}

\subsubsection{Divergence-free subalgebras}\label{SS:2.4.3} Since, as is easy to
calculate,
\begin{equation}
\label{div}
 \Div K_f =(2n+2-m)K_1(f),
\end{equation}
it follows that the divergence-free subalgebra of the contact Lie
superalgebra either coincides with the whole algebra (for $m=2n+2$)
or is the \textit{Poisson} superalgebra.\index{$\fpo$, Poisson
superalgebra}, cf. \eqref{po}
\be\label{po1}
\fpo(2n|m):=\{K_f\in\fk(2n+1 |m)\mid K_1(f)=0\}.
\ee

Since for $(n|m)\neq (0|2)$ the codimension of $(\fk(2n+1 |m)_0)'$ in $\fk(2n+1 |m)_0$ is equal to 1, there is exactly one (up to a~non-zero factor) trace on $\fk(2n+1 |m)_0$. Hence, there should be exactly one (up to a~non-zero factor) intrinsic divergence on $\fk(2n+1 |m)$, see \cite{BGLLS}. 
Setting $\tr(K_t)=1$, the divergence is given by $\partial_t$ for any $n$ and $m$. For $m\neq 2n+2$, this intrinsic divergence coincides with the restriction to $\fk$ of the divergence on $\fvect$.

\paragraph{The case of $\fk(1 |2)$} Although $\fk(1 |2)$ is not W-graded, let us consider it for completeness. This case is exceptional: the algebra $\fk(1 |2)_0$ of degree-0 elements is commutative and $\dim\fk(1 |2)_0=2$, so there are 2 traces on $\fk(1 |2)_0$, and hence 2 divergences on $\fk(1 |2)$. 
Let $\alpha_1:=dt+\xi d\eta+\eta d\xi$. To the operators
\[
\tilde K_\xi(f)=(-1)^{p(f)}(\partial_\eta-\xi\partial_t)(f), \ \tilde K_\eta(f)=(-1)^{p(f)}(\partial_\xi-\eta\partial_t)(f)\text{~~and~~}\tilde K_1=\partial_t
\]
there corresponds the second divergence
\begin{equation}\label{Div_2}
\Div_2:=\tilde K_\eta\tilde K_\xi-\tilde K_1.
\end{equation}

Let $\beta$ be the symbol of the class of a~pseudo-differential form\footnote{If the reader wishes to ``visualize'' this form, one can consider $\beta:=[d\xi(d\eta)^{-1}]$. Actually, the expression of~ $\beta$ in terms of $d\xi$ and $d\eta$ is irrelevant, we need to know only how the Lie derivative acts on it.} on which the Lie derivative acts as follows 
\begin{equation}\label{L_Dk12}
\begin{array}{l}
 L_{K_f}(\alpha_1^a\beta^b)=(2a\partial_t(f)+ (-1)^{p(f)} b \Div_2(K_f))(\alpha_1^a\beta^b)\\
 =((2a- (-1)^{p(f)} b )\partial_t(f)+ (-1)^{p(f)} b\tilde K_\eta\tilde K_\xi(f))(\alpha_1^a\beta^b).
\end{array}
\end{equation}

Let the space $\cF_{a,b}$ of weighted densities over $\fk(1|2)$ be defined as the rank-1 $\cF$-module generated by $\alpha^{a/2}_1\beta^{b}$. (For the modified formula \eqref{L_Dk12} when $\Char\Kee=2$, see~\cite{BGLLS}.)

\paragraph{The pericontact series} Here, the
situation is more interesting: the divergence-free subalgebra is
simple.

Since the restriction of the divergence from $\fvect$ to $\fm$ is equal to
\begin{equation}
\label{2.4.1} \Div M_f =(-1)^{p(f)}2\left ((1-E)\nfrac{\partial f}{\partial\tau} -
\mathop{\sum}\limits_{i\leq n}\frac{\partial^2 f}{\partial q_i
\partial\xi_i}\right ), 
\end{equation}
it follows that the divergence-free subalgebra of the pericontact\index{$\fsm(n)$}
superalgebra is
\begin{equation}
\label{sm2} \fsm (n) = \Span\left (M_f \in \fm (n)\mid
(1-E)\nfrac{\partial f}{\partial\tau} =\mathop{\sum}\limits_{i\leq n}\frac{\partial^2
f}{\partial q_i
\partial\xi_i}\right ).
\end{equation}
In particular,
\begin{equation}
\label{2.4.2}
 \Div Le_f = (-1)^{p(f)}2\mathop{\sum}\limits_{i\leq
n}\frac{\partial^2 f}{\partial q_i \partial\xi_i}.
\end{equation}
The odd analog of the Laplacian, namely, the operator\index{$\Delta=\mathop{\sum}_{i\leq n}\nfrac{\partial^2}{\partial q_i\partial\xi_i}$}
\begin{equation}
\label{2.4.3} \Delta:=\mathop{\sum}\limits_{i\leq n}\frac{\partial^2
}{\partial q_i
\partial\xi_i}
\end{equation}
on a~periplectic supermanifold appeared in physics under the name of
\textit{BRST operator}, cf. \cite{GPS}, or \textit{Batalin-Vilkovysky
operator}. Observe that $\Delta$ is the Fourier transform (with
respect to the ``ghost indeterminates'' $\check x$ (the odd ones, if
considered on manifolds) of the exterior differential~ $d$.

The divergence-free vector fields from $\fsle (n)$ are generated by
\textit{harmonic} functions, \textit{i.e.}, such that $\Delta(f)=0$.

The Lie superalgebra \index{$\fb(n)$, Buttin superalgebra}\index{Buttin superalgebra}
\begin{equation}
\label{but}
\begin{array}{c}
\fb (n)=\{ D\in \fm (n)\mid L_D\alpha_0=0\}\end{array}
\end{equation}
is called the \textit{Buttin} superalgebra in honor of C.~Buttin, see Section~\ref{rem2.3.7}. (A geometric interpretation: over a~\textit{periplectic} supermanifold, \textit{i.e.}, over
a~supermanifold with the periplectic form $d\alpha_0$, the Buttin
superalgebra is the Lie superalgebra that preserves the
connection with form $\alpha_0$ and curvature proportional to $d\alpha_0$ in a~ line bundle of superrank $\eps=(0|1)$.)

The Lie superalgebras\index{$\fsb(n)$}\index{$\fsm(n)$}
\begin{equation}
\label{sm1}
\begin{array}{c}
\fsm (n)=\{ D\in \fm (n)\mid \Div\ D=0\}, \\
 \fs\fb (n)=\{ D\in
\fb (n)\mid \Div\ D=0\}
\end{array}
\end{equation}
are called the \textit{divergence-free} (or \textit{special})
\textit{pericontact} and \textit{special Buttin} superalgebras,
respectively.

Lie superalgebras $\fsle (n)$, $\fs\fb (n)$ and $\fsvect (1|n)$
have\index{$\fsle'(n)$}\index{$\fs\fb'(n)$}\index{$\fsvect'(1\vert n)$} traceless ideals $\fsle '(n)$, $\fs\fb '(n)$ and $\fsvect
'(1|n)$ of codimension 1 (\textit{i.e.}, $\eps^n$, $\eps^{n+1}$ and $n$, respectively) defined from the exact sequences
\begin{equation}
\label{seqv}
\begin{array}{c}
0\tto \fsle '(n)\tto \fsle
(n)\tto \Cee\cdot \Le_{\xi_1\dots\xi_n} \tto 0, \\
0\tto \fs\fb '(n)\tto \fs\fb
(n)\tto \Cee\cdot M_{\xi_1\dots\xi_n} \tto 0, \\
 0\tto \fsvect '(1|n)\tto
\fsvect (1|n)\tto \Cee
\cdot\xi_1\dots\xi_n\partial_t\tto 0.\end{array}
\end{equation}

The Lie superalgebra of \textit{Hamiltonian
fields}\index{Hamiltonian vector fields} (also called \textit{Hamiltonian
superalgebra}) and its traceless subalgebra (defined only if $n=0$)
are\index{$\fh' (0\vert m)$}
\begin{equation}
\label{h}
\begin{array}{l}
\fh (2n|m):=\Span(D\in \fvect (2n|m)\mid L_D\omega_0=0),\\
\fh' (0|m):=\Span(H_f\in \fh (0|m)\mid \int f\vvol_{\theta}=0).
\end{array}
\end{equation}
The ``odd'' analogues of the Lie superalgebra of Hamiltonian
fields are the Lie superalgebra of vector fields $\Le_{f}$
introduced in \cite{L1}, and its special (divergence-free) subalgebra: \index{$\fsle (n)$}
\begin{equation}
\label{len}
\begin{array}{l}
\fle (n):=\Span(D\in \fvect (n|n)\mid L_D\omega_1=0),\\
\fsle (n):=\Span(D\in \fle (n)\mid \Div D=0).
\end{array}
\end{equation}


It is not difficult to prove the following isomorphisms (as superspaces):
\begin{equation}
\label{iso}
\begin{array}{rcl}
\fk (2n+1|m)&\simeq &\Span(K_f\mid f\in \Cee[t, p, q, \xi]),\\
\fle(n)&\simeq &\Span(Le_f\mid f\in \Cee [q, \xi]),\\
\fm (n)&\simeq &\Span(M_f\mid f\in \Cee [\tau, q, \xi]),\\
 \fh(2n|m)&\simeq &\Span(H_f\mid f\in \Cee [p, q, \xi]).
\end{array}
\end{equation}
We see that
\[
\begin{array}{l}
\fpo(0|m)\simeq \Span(K_f\in \fk (1|m)\mid K_1(f)=0),\\
\fpo' (0|m)\simeq \Span(K_f\in \fpo (0|m)\mid \int f\vvol_\theta=0),\\
\fh' (0|m)\simeq \fpo' (0|m)/\Cee\cdot K_1.
\end{array}\]

\sssec{On the history of the antibracket}\label{rem2.3.7} What we call here
the ``Buttin bracket'' was discovered in pre-super era by Schouten;
Buttin was the first to prove that this bracket satisfies the super anti-commutativity and Jacobi identities, see \cite{Bu}. The interpretations of the Buttin
superalgebra similar to that of the Poisson algebra and of the
elements of $\fle$ as analogs of Hamiltonian vector fields was given
in \cite{L1}. The Buttin bracket and ``odd mechanics'' introduced in
\cite{L1} were rediscovered by Batalin with Vilkovisky (and, even
earlier, by Zinn-Justin, but his papers went even more unnoticed, see
\cite{FLSf}, than \cite{L1}). The Schouten bracket gained a~great deal of currency under the name
\textit{antibracket}, cf. \cite{GPS}. The \textit{Schouten bracket$=$Antibracket$=$Buttin bracket}
was originally defined on the superspace of multi-vector fields on
a~manifold, \textit{i.e.}, on the superspace of sections of the exterior
algebra (over the algebra $\cF$ of functions) of the tangent bundle,
$\Gamma(\Lambda^{\bcdot}(T(M)))\simeq \Lambda^{\bcdot}_\cF(Vect(M))$.
The explicit formula for the Schouten bracket (hereafter, the hatted
slot should be ignored, as usual) is (mind the Sign Rule!)
\begin{equation}
\label{*2}
\renewcommand{\arraystretch}{1.4}
\begin{array}{c}
[X_1\wedge\dots \wedge\dots \wedge X_k, Y_1\wedge\dots \wedge Y_l]=\\
\mathop{\sum}\limits_{i, j}(-1)^{k-j+i-1}[X_i, Y_j]\wedge
X_1\wedge\dots\wedge \widehat X_i\wedge \dots\wedge X_k\wedge
Y_1\wedge\dots\wedge \widehat Y_j\wedge \dots \wedge Y_l.
\end{array}
\end{equation}
With the help of the Sign Rule we easily superize formula
~\eqref{*2}, \textit{i.e.}, replace any manifold $M$ by a~supermanifold $\cM$. Let $x$ and
$\xi$ be the even and odd coordinates on $\cM$. Setting
\begin{equation}
\label{**1} 
 \theta_i=\Pi\left(\pder{x_{i}}\right)=\check
x_{i},\quad q_j=\Pi\left(\pder{\xi_{j}}\right)= \check
\xi_{j}
\end{equation}
we get an identification of the Schouten bracket of multi-vector
fields on $\cM$ with the Buttin bracket of functions on the
supermanifold $\check\cM$ with coordinates $x, \xi$ and $\check x$,
$\check \xi$, and the transformation rule of the checked variables
induced by that of unchecked ones via ~\eqref{**1}. Grozman discovered a~deformation of the Schouten bracket, see \cite{G}, but his result was not understood at that time.


\subsubsection{The explicit form of the
pericontact bracket}\label{SS:2.4.2} Set $F_x:=\pder{x}F$. Let $M_{F}\alpha$ be a~short for $L_{M_{F}}\alpha$.
Since $M_{F}\alpha _0 = (-1)^{p(F)+1}2F_{\tau}\alpha_0$, it follows
that
\begin{equation}
\label{2.3.10}
\begin{array}{l}
(M_{F}M_{G} - (-1)^{(p(F)+1)(p(G)+1)}M_{G}M_{F})\alpha _0 =\\
\{(-1)^{p(G)+1}2[(\xi
F_{\tau}+F_{q})G_{\xi
\tau}+(-1)^{p(F)}(F_{\xi}-qF_{\tau})G_{q\tau}+
(-1)^{p(G)(p(F)+1)}2G_{\tau}F_{\tau}]\\
-(-1)^{(p(F)+1)p(G)}2[(\xi G_{\tau}+G_{q})F_{\xi \tau}+
(-1)^{p(G)}(G_{\xi}-qG_{\tau})F_{q\tau}+(-1)^{p(F)(p(G)+1)}
2F_{\tau}G_{\tau}]\} \alpha_0.
\end{array}
\end{equation}

On the other hand, we have
\begin{equation}
\label{2.3.11} ((2-E) (F)G_{\tau})_{\tau} = 2F_{\tau}G_{\tau}+\xi
F_{\tau\xi }G_{\tau}-qF_{\tau q}G_{\tau};\ 
\end{equation}
\begin{equation}
\label{2.3.12} (F_{\tau}(2-E) (G))_{\tau} =
(-1)^{p(F)+1}(2F_{\tau}G_{\tau}+F_{\tau}\xi G_{\tau\xi
}-F_{\tau}qG_{\tau q});\ 
\end{equation}
\begin{equation}
\label{2.3.13} (F_{q}G_{\xi }+(-1)^{p(F)}F_{\xi }G_{q})_{\tau} =
F_{\tau q}G_{q}+(-1)^{p(F)}F_{q}G_{\tau\xi }-F_{\xi }G_{\tau q}.
\end{equation}
Clearly,
\[
(-1)^{p(F)p(G)}~\eqref{2.3.11}+(-1)^{p(G)}~\eqref{2.3.12}+
(-1)^{p(F)+p(G)}~\eqref{2.3.13} = \{\cdots\},
\]
where $\{\cdots \}$ is the coefficient of $\alpha _0$ in the right-hand 
 side of ~\eqref{2.3.10}. Therefore,
\[
\{F, G\}_{m.b.} = (-1)^{p(F)+p(G)}\left((2-E) (F)G_{\tau}+
(-1)^{p(F)}F_{\tau}(2-E)
(G)+\{F, G\}_{\text{B.b.}}\right).
\]

The arguments for the contact bracket are similar (with $\alpha_0$
replaced with $\alpha_1$ or $\widetilde \alpha_1$, and the field $M_{f}$
with $K_{f}$, respectively).


\ssec{Cartan prolongations}\label{2.6}\index{Cartan prolongation}
Recall that the graded Lie superalgebra
$\fb=\oplus_{k \geqslant -d}\ \fb_k$
is said to be \textit{transitive}\index{Lie (super)algebra, transitive} if for all~ $k\geq 0$ we have
\begin{equation*}
\label{12}
\{x\in \fb_k \mid [x,\fb_{-}]=0\}=0, \text{~~where $\fb_{-}:=\oplus _{k <0}\,\fb_k$}.
\end{equation*}

Let
$\fg_-=\mathop{\oplus}_{-d\leq i\leq -1}\ \fg_i$ be a~nilpotent
$\Zee$-graded Lie (super)algebra and let $\fg_0$ be a~Lie sub(super)algebra of the Lie
(super)algebra $\mathfrak{der}_0(\fg_-)$ of degree-0 derivations of~ $\fg_-$. 
The maximal transitive $\Zee$-graded Lie (super)algebra $\fg=\oplus_{k \geqslant -d}\ \fg_k$ whose non-positive part is the given $\fg_-\oplus \fg_0$ is called the 
\textit{Cartan prolong}\index{Cartan prolongation} 
--- the result of Cartan prolongation --- of the pair $(\fg_-,\fg_0)$ and is denoted by $(\fg_-,\fg_0)_*$. 

We can realize $\fg_-$ by elements of negative degree of $\fvect(n|m; \vec r)$ and $\fg_0$ by degree-0 elements in $\fvect(n|m; \vec r)$ in a~non-standard (see Subsection~\ref{secnon}) grading $\vec r$ of $\fvect(n|m)$, where $n|m={\rm sdim} \fg_-$. Then, the degree-$k$ components of the Cartan prolongation $(\fg_-, \fg_0)_{*}:=\mathop{\oplus}_{k\geq -d}\ \fg_{k}$ 
 of the pair $(\fg_-, \fg_0)$ are obtained for any $k>0$ by setting
\[
\fg_k:=\{D\in \fvect(n|m; \vec r)_k\mid [D, \fg_i]\subset \fg_{k+i}\text{~~for any $i<0$}\}.
\]

In what follows ${\bcdot}$ in superscript denotes, as is now
customary, the collection of all degrees, while $*$ is reserved for
dualization; in the subscripts we retain the old-fashioned $*$ instead
of ${\bcdot}$ to avoid confusion with the punctuation
marks.


Suppose that the $\fg_0$-module $\fg_{-1}$ is \textit{faithful}.
Then, clearly,
\begin{equation}\label{87}
\begin{array}{ll}
(\fg_{-1}, \fg_{0})_*&\subset \fvect (n|m) = \fder~\Cee[x_1, \dots ,
x_n, \theta_1,\dots \theta_m],\; \text{ where }\; n|m = \sdim~\fg_{-1}\; \text{ and }\\ 
&\fg_i
= \{D\in \fvect(n|m)\mid \deg D=i, [D, X]\in\fg_{i-1}\text{ for any
} X\in\fg_{-1}\}.
\end{array}
\end{equation}
It is subject to an easy verification that the Lie algebra structure
on $\fvect (n|m)$ induces a~Lie algebra structure on $(\fg_{-1},
\fg_{0})_*$; actually, $(\fg_{-1}, \fg_{0})_*$ has a~natural Lie superalgebra
structure even if the $\fg_0$-module $\fg_{-1}$ is not faithful.

\subsubsection{Partial Cartan prolongation: involving
positive components}\label{CartProlPart}\index{Cartan prolongation, partial} Let $\fh_1\subset \fg_1$ be
a~proper $\fg_0$-sub\-module such that $[\fg_{-1}, \fh_1]=\fg_0$. If such
$\fh_1$ exists (usually, $[\fg_{-1}, \fh_1]\subsetneq\fg_0$), define
the 2nd prolongation of $(\mathop{\oplus}_{i\leq 0}\fg_i)\oplus
\fh_1$ to be
\begin{equation*}
\label{partprol} \fh_{2}:=\{D\in\fg_{2}\mid [D, \fg_{-1}]\subset
\fh_1\}.
\end{equation*}
The terms $\fh_{i}$, where $i>2$, are similarly defined by induction. Set
$\fh_i:=\fg_i$ for $i\leq 0$ and call $\fh_*:=\oplus\fh_i$ the \textit{partial Cartan prolong involving positive components}.

The Lie superalgebra $\fvect(1|n; n)$ is a~subalgebra of $\fk(1|2n;
n)$. The former is obtained as the Cartan prolong of the same
nonpositive part as $\fk(1|2n; n)$ and a~submodule of $\fk(1|2n;
n)_1$. The simple exceptional superalgebra $\fk\fas$ introduced in
\cite{Sh5, Sh14}, see also line $N=12$ in Table~\eqref{table3}, is
another example.

\subsubsection{Lie superalgebras of vector fields as Cartan
prolongs}\label{SS:2.6.2} The superization of the constructions from
Subsection~\ref{2.6} are straightforward, via the Sign Rule. We thus get Lie superalgebras
\begin{equation}
\label{carprol}
\begin{array}{l}
\fvect(m|n)=(\id, \fgl(m|n))_*,\\
\fsvect(m|n)=(\id, \fsl(m|n))_*, \\
\fh(2m|n)=(\id, \fosp^{a}(m|2n))_*, \\
\fle(n)=(\id, \fpe^{a}(n))_*,\\
\fs\fle(n)=(\id, \fspe^{a}(n))_*.\\
\end{array}
\end{equation}

For an algorithm for constructing prolongs, see \cite{Shch}; it is implemented in the code package \textit{SuperLie}, see \cite{Gr}.

Note that $(\id, \fosp^{sy}
(m|2n))_*$ and $(\id, \fpe ^{sy} (n))_*$ are finite-dimensional, unlike 
\[
\text{$(\Pi(\id), \fosp^{sy}
(m|2n))_*=(\id, \fosp^{a}(m|2n))_*$ ~~and~~ $(\Pi(\id), \fpe ^{sy} (n))_*=(\id, \fpe^{a}(n))_*$.}
\] 

There are \textbf{two} superizations of the contact
series: $\fk$ and $\fm$. Let us describe them again.

$\bullet$ Let the space of the Heisenberg Lie superalgebra $\fhei(2n|m)$ be the direct\index{$\fhei(2n\vert m)$, Heisenberg (super)algebra}
sum of a~$(2n, m)$-dimensional superspace $W$ endowed with
a~non-degenerate anti-symmetric \textbf{even} bilinear form $B$ and a~$(1,
0)$-dimensional space spanned by $z$. Define the bracket in $\fhei$ by the formula \begin{equation}
\label{2.5.3} \text{$z$ is in the center and $[v, w]=B(v, w)\cdot
z$ for any $v, w\in W$.}
\end{equation}

Clearly, 
\begin{equation}
\label{kont} \fk(2n+1|m)=(\fhei(2n|m), \fc(\fosp^{a}(m|2n)))_*\,.
\end{equation}
More generally, given $\fhei(2n|m)$ and a~subalgebra $\fg$ of
$\fc(\fosp^{a}(m|2n))$, we call $(\fhei(2n|m), \fg)_*$ the
\textit{$k$-prolong} of $(W, \fg)$, where $W$ is the tautological
$\fosp^{a}(m|2n)=\fosp(W)$-module.

$\bullet$ The ``odd'' analog of $\fk$ is associated with the
following ``odd'' analog of $\fhei(2n|m)$. Denote by $\fba(n)$ the\index{$\fba(n)$, anti-bracket superalgebra}
\textit{antibracket} Lie superalgebra ($\fba$ stands for Anti-Bracket read
backwards) whose space is $W\oplus \Cee\cdot z$, where $W$ is an
$n|n$-dimensional superspace endowed with a~non-degenerate
anti-symmetric \textbf{odd} bilinear form $B$ and $z$ is odd; the bracket in $\fba(n)$ is
given by the 
relations
\begin{equation}
\label{2.5.4} \text{$z$ lies in the center~~ and ~~$[v,
w]=B(v, w)\cdot z$ for any $v, w\in W$.}
\end{equation}

Clearly,
\begin{equation}
\label{mprol} \fm(n)=(\fba(n), \fc(\fpe^{a}(n)))_*\,.
\end{equation}
More generally, given $\fba(n)$ and a~subalgebra $\fg$ of $\fc\fpe^{a}(n)$, we
call $(\fba(n), \fg)_*$ the \textit{$m$-prolong} of $(W, \fg)$,
where $W$ is the tautological $\fpe^{a}(n)=\fpe(W)$-module.

Generally, given a~non-degenerate $\fg$-invariant form $B$ on a~superspace $W$ and a
Lie subsuperalgebra $\fg\subset\faut(B)$, we call any of the above
generalized prolongations \eqref{kont} and \eqref{mprol} the \textit{$mk$-prolongation} of the pair $(W, \fg)$ and denote it by 
$(W, \fg)_{*}^{mk}$.\index{$mk$-prolongation}\index{$(W, \fg)_{*}^{mk}$}

\paragraph{On the history of Cartan prolongations} The above-described procedure is sometimes called a~ \textit{generalized} prolongation because the initial Cartan prolongation was defined for $d=1$ only. Or rather it was customary to write so without reading works by Cartan and Killing who did consider the above-described generalized Cartan prolongation, see \cite{Cart, Ki, Co}. 

In the non-super setting, the ``generalized'' Cartan prolongation was introduced by Cartan (see \cite{Cart}) or, perhaps, Killing (see \cite{Ki}). The Cartan prolongation was rediscovered in two steps by Tanaka \cite{T}, see \cite{Y}.
The emphasis of these
papers is different from ours, this delayed recognition of the similarity of constructions. The
partial prolongation seems to be first described in \cite{ALSh, Sh14}. It also appears in the classification of simple Lie (super)algebras in characteristic $p>0$, see \cite{BGLLS}, where the correct definition of Cartan prolongation in characteristic $p>0$ was first given, after several approximations to the truth in earlier papers.


\ssec{Nonstandard realizations}\label{secnon}
The Step 7 of our proof (see Subsection \ref{steps}) implies that the nonstandard
gradings listed in table \eqref{nonstandgr} exhaust all the W-gradings of all the simple vectorial
Lie superalgebras. In particular, the gradings in the series
$\fvect$ induce the gradings in the series $\fsvect$, $\fsvect'$ and
the exceptional families $\fv\fle(4|3)$ and $\fv\fas(4|4)$; the
gradings in $\fm$ induce the gradings in $\fb_{\lambda}$, $\fle$,
$\fsle$, $\fsle'$, $\fb$, $\fs\fb$, $\fs\fb'$ and the exceptional
family $\fm\fb$; the gradings in $\fk$ induce the gradings in
$\fpo$, $\fh$, $\fh'$ and the exceptional families $\fk\fas$ and
$\fk\fle$.

The $\Zee$-gradings of vectorial Lie superalgebras are defined by
the vector $\vec r$ of degrees of the indeterminates. \textbf{For W-gradings}, this\index{$\vec r$ vector of degrees of the indeterminates}\index{Grading, non-standard}
vector can be shorthanded to a~symbol or the number $r$ of degree-0 indeterminates; we do not indicate $r$ if $r=0$. Let the indeterminates $t$ and $u_i$ be even, let $\tau$ and the $\theta_j$ be 
odd. In Table~\eqref{nonstandgr}, we consider $\fk (2n+1|m)$ preserving
the Pfaff equation with the 1-form $\alpha_1)$; here $(t; p,q)$ are the even indeterminates, the odd indeterminates being
 $\tau; \theta, \xi, \eta$, see ~eq.~\eqref{2.2.6}. For the series $\fm$, the
indeterminates are denoted as in eq.~\eqref{alp0}.

The standard realizations correspond to $r=0$, they are marked by an
$(*)$. Note that the codimension of ${\cal L}_0$ attains its
minimum in the standard realization.
\begin{equation}\label{nonstandgr}
\renewcommand{\arraystretch}{1.3}\footnotesize
\begin{tabular}{|c|c|}
\hline Lie superalgebra & its $\Zee$-grading \\ \hline

$\fvect (n|m; r)$, & $\deg u_i=\deg \xi_j=1$ for any $i, j$
\hfill $(*)$\\

\cline{2-2} $ 0\leq r\leq m$ & $\deg \xi_j=0$ for $1\leq j\leq r;$\\
&$\deg u_i=\deg \xi_{r+s}=1$ for any $i, s$
\\

\hline $\fm(n; r),$ & $\deg \tau=2$, $\deg q_i=\deg \xi_i=1$ for any
$i$ \hfill $(*)$\\ \cline{2-2} $\; 0\leq r< n-1$& $\deg \tau=\deg
q_i=2$, $\deg \xi_i=0$ for $1\leq i\leq r <n$;\\

& $\deg q_{r+j}=\deg \xi_{r+j}=1$ for any $j$\\
\hline

$\fm(n; n)$ & $\deg \tau=\deg q_i=1$, $\deg \xi_i=0$ for $1\leq
i\leq n$ \\ \hline

\cline{2-2} $\fk (2n+1|m; r)$, & $\deg t=2$, whereas, for any $i,
j, s$\hfill $(*)$\\
see eq.~\eqref{2.2.6};& $\deg
p_i=\deg q_i= \deg \xi_j=\deg \eta_j=\deg \theta_s=1$ \\

\cline{2-2} $0\leq r\leq [\frac{m}{2}]$& $\deg t=\deg \xi_i=2$,
$\deg \eta_{i}=0$ for $1\leq i\leq r\leq [\frac{m}{2}]$; \\
$r\neq k-1$ for $m=2k$ and $n=0$&$\deg p_i=\deg q_i=\deg
\theta_{j}=1$
for $j\geq 1$ and all $i$\\

\hline $\fk(1|2m; m)$ & $\deg t =\deg \xi_i=1$, $\deg \eta_{i}=0$
for $1\leq i\leq m$ \\ \hline
\end{tabular}
\end{equation}

In table \eqref{clsfeq10}, the five families of exceptional simple vectorial Lie superalgebras are
given in their W-gradings as Cartan prolongs $(\fg_{-1},
\fg_{0})_{*}$ or generalized Cartan prolongs $(\fg_{-},
\fg_{0})_{*}^{mk}$.
For \underline{depth $2$}, for
$\fg_{-}=\mathop{\oplus}_{-2\leq i\leq -1}\fg_{i}$, we write
$(\fg_{-2}, \fg_{-1}, \fg_{0})_{*}^{mk}$\index{$(\fg_{-2}, \fg_{-1}, \fg_{0})_{*}^{mk}$} instead of $(\fg_{-},
\fg_{0})_{*}^{mk}$. 
For $\fk\fa\fs^\xi$, the meaning of the notation $r=0$, or $1\xi$, or $3\xi$ is clear:
none, or one, or three of the $\xi$'s have degree 0 (and the
corresponding $\eta$'s acquire degree 2). Table~\eqref{clsfeq10} contains
exceptional superalgebras whose grading $r$ is described in table~\eqref{1.10}, indeterminates and their respective degrees
in the regrading $R(r)$ of the ambient superalgebra and
describes these regradings: 
\[
\text{(the degrees of the even
indeterminates $\mid$ the degrees of the odd indeterminates).}
\]
{\tiny%
\begin{equation}\label{clsfeq10}
\setcounter{MaxMatrixCols}{16}
\tabcolsep=3pt
\begin{tabular}{|l|l|}
\hlx{hv}
$\fv\fle(4|3)$\index{$\fv\fle(4\vert 3)$}& $ R(0)=\begin{pmatrix}
x_1&x_2&x_3&y&|&\xi_1&\xi_2&\xi_3\\ 1&1&1&1&|&1&1&1\\ \end{pmatrix}$\cr \hlx{vv}

$\fv\fle(5|4)$\index{$\fv\fle(5\vert 4)$}&$R(1)=\begin{pmatrix}
x_1&x_2&x_3&y&|&\xi_1&\xi_2&\xi_3\\ 2&1&1&0&|&0&1&1\\ \end{pmatrix}$\cr \hlx{vv}

$\fv\fle(3|6)$\index{$\fv\fle(3\vert 6)$}&$R(K)=\begin{pmatrix} 
x_1&x_2&x_3&y&|&\xi_1&\xi_2&\xi_3\\
2&2&2&0&|&1&1&1\\\end{pmatrix}$\cr \hlx{vhv}

$\fm\fb(4|5)$\index{$\fm\fb(4\vert 5)$}&$R(0)=\begin{pmatrix} 
x_0&x_1&x_2&x_3&|&\xi_0&\xi_1&\xi_2&\xi_3&\tau \\
1&1&1&1&|&1&1&1&1&2\\\end{pmatrix}$\cr
\hlx{vv}
$\fm\fb(5|6)$\index{$\fm\fb(5\vert 6)$}& $R(1)=\begin{pmatrix} 
x_0&x_1&x_2&x_3&|&\xi_0&\xi_1&\xi_2&\xi_3&\tau \\
0&2&1&1&|&2&0&1&1&2\\\end{pmatrix}$\cr
\hlx{vv}
$\fm\fb(3|8)$\index{$\fm\fb(3\vert 8)$}&$R(K)=\begin{pmatrix}
x_0&x_1&x_2&x_3&|&\xi_0&\xi_1&\xi_2&\xi_3&\tau \\
 0&2&2&2&|&3&1&1&1&3\\\end{pmatrix}$\cr

\hlx{vhv}
$\fk\fas(1|6)$\index{$\fk\fas(1\vert 6)$}&$R(0)=\begin{pmatrix}
t&|&\xi_1&\xi_2&\xi_3&\eta_1&\eta_2&\eta_3 \\
 2&|&1&1&1&1&1&1\\\end{pmatrix}$\cr
\hlx{vv}
$\fk\fas(5|5)$\index{$\fk\fas(5\vert 5)$}&$R(1\xi)=\begin{pmatrix}
t&|&\xi_1&\xi_2&\xi_3&\eta_1&\eta_2&\eta_3 \\
 2&|&0&1&1&2&1&1\\\end{pmatrix}$\cr
\hlx{vv}
$\fk\fas(4|4)$\index{$\fk\fas(4\vert 4)$}&$R(3\xi)=\begin{pmatrix} 
t&|&\xi_1&\xi_2&\xi_3&\eta_1&\eta_2&\eta_3 \\
1&|&0&0&0&1&1&1\\\end{pmatrix}$\cr
\hlx{vv}
$\fk\fas(4|3)$\index{$\fk\fas(4\vert 3)$}&$R(3\eta)=\begin{pmatrix}
t&|&\xi_1&\xi_2&\xi_3&\eta_1&\eta_2&\eta_3 \\
 1&|&1&1&1&0&0&0\\\end{pmatrix}$\cr
\hlx{vhv}
$\fk\fle(9|6)$\index{$\fk\fle(9\vert 6)$}&$R(0)=\begin{pmatrix} 
q_1&q_2&q_3&q_4&p_1&p_2&p_3&p_4&t&|&\xi_1&\xi_2&\xi_3&
\eta_1&\eta_2&\eta_3\\
1&1&1&1&1&1&1&1&2&|&1&1&1&1&1&1\\
\end{pmatrix}$\cr
\hlx{vv}
$\fk\fle(11|9)$\index{$\fk\fle(11\vert 9)$}&$R(2)=\begin{pmatrix} 
q_1&q_2&q_3&q_4&p_1&p_2&p_3&p_4&t&|&\xi_1&\xi_2&\xi_3&
\eta_1&\eta_2&\eta_3\\
1&1&2&2&1&1&0&0&2&|&0&1&1&2&1&1\\
\end{pmatrix}$\cr
\hlx{vv}
$\fk\fle(5|10)$\index{$\fk\fle(5\vert 10)$}&$R(K)=\begin{pmatrix} 
q_1&q_2&q_3&q_4&p_1&p_2&p_3&p_4&t&|&\xi_1&\xi_2&\xi_3&
\eta_1&\eta_2&\eta_3\\
2&2&2&2&1&1&1&1&2&|&1&1&1&1&1&1\\
\end{pmatrix}$\cr
\hlx{vv}
$\fk\fle(9|11)$\index{$\fk\fle(9\vert 11)$}&$R(CK)=\begin{pmatrix} 
q_1&q_2&q_3&q_4&p_1&p_2&p_3&p_4&t&|&\xi_1&\xi_2&\xi_3&
\eta_1&\eta_2&\eta_3\\
3&2&2&2&0&1&1&1&3&|&2&2&2&1&1&1\\
\end{pmatrix}$\cr
\hlx{vh}
\end{tabular}
\end{equation}}
In Table~\eqref{1.10}, indicated is also one of
the Lie superalgebras from the list of series ~\eqref{1.4} as an
ambient which contains the exceptional one as a~maximal subalgebra.
The W-graded superalgebras of \underline{depth $3$} appear as
regradings of the listed ones at certain values of $\vec r$; the
corresponding terms $\fg_{i}$ for $i\leq 0$ are given in
Subsection~\ref{SS:1.3.3}. Let
the label $CK$ commemorate an
exceptional grading found by Cheng and Kac, see \cite{CK2}, see Theorem~\ref{th1.10}. 

\begin{equation}
\label{1.10}\footnotesize
\renewcommand{\arraystretch}{1.4}
\begin{tabular}{|ll|l|}
\hline $\fv\fle(4|3; r)=(\Pi(\Lambda(3))/\Cee\cdot 1,
\fc\fvect(0|3))_{*}$& $\subset\fvect(4|3; R(r))$& $r= 0, 1, K$\cr

\hline $\fv\fas(4|4)=(\spin, \fas)_{*}$, see Subsection~\ref{ssas}& $\subset \fvect(4|4)$&\cr

\hline $\fk\fas^\xi(1|6; r)$& $\subset \fk(1|6; r)$&$r=0, 1\xi,
3\xi$\cr

$\fk\fas^\xi(1|6; 3\eta)=(\Vol_{0}(0|3), \fc(\fvect(0|3)))_{*}$&
$\subset \fsvect(4|3)$&$r=3\eta$\cr \hline

$\fm\fb(4|5; r)=(\fba(4), \fc\fvect(0|3))_{*}^m$& $\subset \fm(4|5;
R(r))$& $r=0, 1, K$\cr \hline

$\fk\fle(9|6; r)=(\fhei(8|6),
\fsvect(0|4)_{3, 4})_{*}^k$& $\subset \fk(9|6; r)$& $r=0, 2$, CK\cr

$\fk\fle(9|6; K)=(\id_{\fsl(5)}, \Lambda^2(\id^*_{\fsl(5)}),
\fsl(5))_{*}^k$& $\subset \fsvect(5|10; R(K))$&$r= K$\cr \hline
\end{tabular}
\end{equation}

\ssec{The modules of tensor fields}\label{SS:2.7} To advance
further, we have to recall the definition of the modules of tensor
fields over $\fvect(m|n)$ and its subalgebras, see \cite{BL4, L2, GLS3}.

Let $\fg=\fvect(m|n)$ and $\fg_{\geq}=\mathop{\oplus}_{i\geq
0}\ \fg_{i}$. For any other $\Zee$-graded vectorial Lie superalgebra
for whose component $\fg_0$ the notion of lowest weight can be
defined, the construction is identical.

Let $V$ be the $\fg_0$-module with the \textit{lowest} weight
$\lambda=\lwt(V)$ and \textit{even} lowest weight vector. Make $V$ into a~$\fg_{\geq}$-module by setting
$\fg_{+}\cdot V=0$ for $\fg_{+}=\mathop{\oplus}_{i>
0}\ \fg_{i}$. Let us realize $\fg=\fg_-\oplus \fg_0\oplus \fg_+$ by vector fields on the linear
supermanifold $\cC^{m|n}$ corresponding to the superspace $\fg_-:=\mathop{\oplus}_{i<
0}\ \fg_{i}$ with coordinates $x=(u, \xi)$. The
superspace $T(V):=\Hom_{U(\fg_{\geq})}(U(\fg), V)$\index{$T(V):=\Hom_{U(\fg_{\geq})}(U(\fg), V)$} is isomorphic, due
to the Poincar\'e--Birkhoff--Witt theorem, to ${\Cee}[[x]]\otimes
V$. Its elements admit a~natural interpretation as formal
\textit{tensor fields of type} $V$. When $\lambda=(a, \dots , a)$ we
will simply write $T(\vec a)$ instead of $T(\lambda)$. We will
usually need $\fg$-modules coinduced from \textit{irreducible}
$\fg_0$-modules.

\sssec{Examples}\label{tensF} As $\fvect(m|n)$- and $\fsvect(m|n)$-modules,
$\fvect(m|n)=T(\id)$.

$T(\vec 0)$ is the superspace of functions;\index{$T(\vec 0)$}

$\Vol(m|n):=T(\underbrace{1, \dots , 1}_m; \underbrace{-1, \dots , -1}_n)$, where the semicolon separates
the first $m$ (``even'') coordinates of the weight with respect to
the matrix units $E_{ii}$ of $\fgl(m|n)$, is the superspace of
\textit{volume forms}, see \cite[Ch.4]{LSoS} or \cite[Ch.1, \S1.11]{Del}.

We denote the generator of $\Vol(m|n)$, considered as a~$\Cee[[x]]$- or $\Cee[x]$-module,
corresponding to the ordered set of coordinates $x$ by $\vvol(x)$ or $\vvol_x$.
The space of $\lambda$-densities is \index{$\Vol^{\lambda}(m\vert n)$}
\index{$\Vol_0(0\vert n)$}
\[
\Vol^{\lambda}(m|n):=T(\underbrace{\lambda,\dots , \lambda}_m; \underbrace{-\lambda, \dots , -\lambda}_n). 
\]
In particular,\index{$T_0(\vec 0)$}
$\Vol^{\lambda}(m|0)=T(\vec \lambda)$, whereas 
$\Vol^{\lambda}(0|n)=T(\overrightarrow{-\lambda})$. In dimension $(0|n)$, we set
\begin{equation}
\begin{array}{l}
\label{vol} \Vol_0(0|n):=\{v\in \Vol(0|n)\mid\int v=0\},\\
T_0(\vec 0):=\Lambda(n)/\Cee\cdot 1.
\end{array}
\end{equation}

As $\fsvect(m|n)$-modules, $\Vol\simeq T(\vec 0)$. So, 
 we can define the $\fsvect(0|n)$-module \index{$T_0^0(\vec 0)$}
\begin{equation}
\label{T00} T^0_{0}(\vec
0):=\Vol_0(0|n)/\Cee\vvol_\xi.\end{equation}


\paragraph{Remark: Volume element as a~quotient} To view the volume element $\vvol(x)$ as ``$d^mud^n\xi$"
is wrong, but in some works this wrong notation is used to this day. Observe that the superdeterminant of the Jacobi matrix can never appear as a~factor
under the changes of ``variables'', \textit{i.e.}, indeterminates, in ``$d^mud^n\xi$". We can try to use the usual notation
of differentials provided \textit{all} the differentials
anticommute and $\deg du_i=1$ while $\deg d\xi_j=-1$ for all $i, j$. Then, at least the linear transformations that do not
intermix the even $u$'s with the odd $\xi$'s multiply the volume
element $\vvol_x$, viewed as the fraction $\frac{du_1\cdot ...\cdot
du_m}{d\xi_1\cdot ...\cdot d\xi_n}$, by the correct factor, the
Berezinian $\ber (A)$, see \cite{Lber}, of the transformation $A$. But how to justify this? Let
$x:=(u, \xi)$. If we consider the exterior differential
forms, then the $dX_{i}$'s super anticommute, hence, the $d\xi_i$
commute and so should $(d\xi_i)^{-1}$; whereas if we consider the \textit{symmetric} product of the differentials, as in metrics, then the $dx_{i}$'s
supercommute, hence, the $du_i$ commute, so neither the exterior, nor the 
symmetric product is what we need: all factors should anticommute.

Note that the $\pder{\xi_i}$ anticommute, and moreover, from the point
of view of linear transformations,
$\pder{\xi_i}=\frac{1}{d\xi_{i}}$. The expression $du_1\cdot ...\cdot
du_m\cdot\pder{\xi_1}\cdot \ldots \cdot \pder{\xi_n}$ is,
nevertheless, still wrong: many transformations $A: (u,
\xi)\longmapsto (v, \eta)$ send $du_1\cdot ....\cdot
du_m\cdot\pder{\xi_1}\cdot ...\cdot \pder{\xi_n}$ to the correct
element, $\ber (A)(dv_1\cdot ...\cdot dv_m\cdot\pder{\eta_1}\cdot
...\cdot \pder{\eta_n})$ plus extra terms. Indeed, the term
$dv_1\cdot ...\cdot dv_m\cdot\pder{\eta_1}\cdot ...\cdot
\pder{\eta_n}$ is the highest weight vector of an
\textit{indecomposable} $\fgl(m|n)$-module and $\vvol_x:=\vvol(u,\xi)$ is the
notation of the image of this vector in the 1-dimensional quotient
module by the invariant submodule that consists precisely of all
the extra terms. Deligne, Witten, and their followers denote the volume element --- the generator of the quotient module --- by $[dx/\partial_\xi]$, see \cite[Ch.1]{Del}.

\ssec{A.~Sergeev's central extension of $\fspe(4)$}\label{ssas}
In 1970's A.~Sergeev proved (cited in \cite{FuLe}, see also \cite{BGLL1}) that over $\Cee$ there is just one
nontrivial central extension of $\fspe(n)$ for $n>2$. It exists only
for $n=4$ and we denote it by Sergeev's initials: $\fas$. Let us represent an arbitrary
element $\widehat X\in\fas$ as a~pair $\widehat X=X+a\cdot 1_{4|4}$, where $X\in\fspe(4)$ and
$a\in{\Cee}$. 
The bracket in $\fas$
is
\begin{equation}
\label{2.1.4}
\left[\mat {A &B\\ C&-A^t} +a\cdot 1_{4|4},
\mat{A' & B' \cr C' & -A'{}^t} +a'\cdot
1_{4|4}\right]= \left[\mat{ A& B \cr C & -A^t},
\mat{ A' & B' \cr C' & -A'{}^t
}\right]+\tr(C\widetilde C')\cdot 1_{4|4},
\end{equation}
where $\ \widetilde {}\ $ is extended by linearity from matrices
$C_{ij}:=E_{ij}-E_{ji}$, on which $\widetilde C_{ij}=C_{kl}$ for any
even permutation $(1234)\longmapsto(ijkl)$.

The Lie superalgebra $\fas$\index{$\fas$} can also be described with the help of
the \textit{spinor representation}, see \cite{LSh4}.\index{Representation, spinor} Consider the particular case of $\fo(6)=(\fpo(0|6))_0$. The superspace of $\fpo(0|6)$ is the Grassmann superalgebra $\Lambda(\xi, \eta)$
generated by $\xi=(\xi_1, \xi_2, \xi_3)$ and $\eta=(\eta _1, \eta _2, \eta_3)$ and the Lie superalgebra structure given by the
Poisson bracket ~\eqref{2.3.6}.

Recall that $\fh(0|6)=\Span (H_f\mid f\in\Lambda (\xi, \eta))$. Now,
note that $\fspe(4)$ can be embedded into $\fh(0|6)$. Indeed,
setting $\deg \xi_i=\deg \eta _i=1$ for all~ $i$ we introduce
a~$\Zee$-grading on $\Lambda(\xi, \eta)$ which, in turn, induces
a~$\Zee$-grading on $\fh(0|6)$ of the form
$\fh(0|6)=\mathop{\oplus}_{i\geq -1}\ \fh(0|6)_i$. Since
$\fsl(4)\simeq \fo(6)$, we can identify $\fspe(4)_0$ with
$\fh(0|6)_0$.

It is not difficult to see that the elements of degree $-1$ in the
standard gradings of $\fspe(4)$ and $\fh(0|6)$ constitute isomorphic
$\fsl(4)\simeq \fo(6)$-modules; moreover, it is possible to embed $\fspe(4)_1$ into $\fh(0|6)_1$.

Sergeev's extension $\fas$ is the result of the restriction to
$\fspe(4)\subset\fh(0|6)$ of the cocycle that turns $\fh(0|6)$ into
$\fpo(0|6)$, see eq.~\eqref{HofSPE}. The quantization $\cQ$ deforms $\fpo(0|6)$ into
$\fgl(\Lambda(\xi))$; the through maps 
\[
T_\lambda:
\fas\tto\fpo(0|6)\stackrel{\cQ}{\tto}\fgl(\Lambda (\xi))
\]
are representations of
$\fas$ in the $4|4$-dimensional modules $\spin_\lambda$\index{$\spin_\lambda$} isomorphic
to each other for all~ $\lambda\neq 0$, so we assume $\lambda=1$ and do not indicate it. The explicit form of
$T_\lambda$ is 
\begin{equation}
\label{spinLambda} T_\lambda\colon \mat{A~& B \cr C & -A^t }+a\cdot
1_{4|4}\longmapsto \mat{A~& B-\lambda \widetilde C \cr C & -A^t
}+\lambda a\cdot 1_{4|4}, 
\end{equation}
where $\widetilde C$ is defined in
the line under eq.~\eqref{2.1.4}. Clearly, $T_\lambda$ is 
irreducible for any $\lambda$.

\ssec{Deformations of the Buttin superalgebra (after
\cite{Ko1, LSh3})}\label{SS:2.8} As is clear from the definition of the
Buttin bracket (and \cite{G} from where it went unnoticed), there is a~regrading (namely, $\fb (n; n)$ given by
$\deg\xi_i=0, \deg q_i=1$ for all~ $i$) under which $\fb(n)$,
initially of depth 2, takes the form
$\fg=\mathop{\oplus}_{i\geq -1}\ \fg_i$ with
$\fg_0=\fvect(0|n)$ and $\fg_{-1}\simeq \Pi(\Cee[\xi])$. Now let us replace 
the $\fvect(0|n)$-module $\fg_{-1}$ of functions (with inverted
parity) by the rank-1 module (over the algebra of functions) of
$\lambda$-densities (also with inverted parity), \textit{i.e.}, set
$\fg_{-1}\simeq \Pi(\Vol^\lambda(0|n))$, where the
$\fvect(0|n)$-action on the generator $\vvol_\xi^\lambda$ of the module $\Vol^\lambda(0|n)$ is given
by the formula
\begin{equation}
\label{2.7.1} L_D(\vvol_\xi^\lambda) =\lambda \Div D\cdot
\vvol_\xi^\lambda\; \text{ and }\;
p(\vvol_\xi^\lambda)=\od \text{~~when considered in $\fg_{-1}$}.
\end{equation}

Define $\fb_{\lambda}(n; n)$ to be the Cartan prolong
\begin{equation}
\label{fb} \fb_{\lambda}(n; n):=(\fg_{-1}, \fg_0)_*=(\Pi(\Vol^\lambda(0|n)),
\fvect(0|n))_*.
\end{equation}
Clearly, this is a~deform of $\fb(n; n)$. Kotchetkov called he collection of these prolongs 
$\fb_{\lambda}(n; n)$ for all~ $\lambda\in\Cee\Pee^1$ the \textit{main
deformation}, see \cite{Ko1}; the other deformations are called \textit{singular}, see Theorem~\ref{thdefb}.\index{Deformation of the Buttin
bracket, main}\index{Deformation of the Buttin
bracket, singular} The deformation of $\fb(n)$ with odd parameter is called \textit{quantization}: it is a~splitting image of the quantization of the Poisson (super)algebra.

The deform $\fb_{\lambda}(n)$ of $\fb(n)$ is a~regrading of
$\fb_{\lambda}(n; n)$ described as follows. Set\index{$\fb_{a, b}(n)$}\index{$\fb_{\lambda}(n)$}\index{$\Div_{\lambda}$}
\begin{equation}
\label{2.7.2} \fb_{a, b}(n) :=\{M_f\in \fm (n)\mid a\; \Div
M_f=(-1)^{p(f)}2(a-bn)\nfrac{\partial f}{\partial\tau}\}. 
\end{equation}

Using the explicit form of the divergence of $M_{f}$ we see that
\begin{equation}
\label{2.7.4}
\begin{array}{ll}
\fb_{a, b}(n) &=\{M_f\in \fm (n)\mid (bn-aE)\nfrac{\partial f}{\partial\tau} =
a\Delta
f\}=\\
&\{D\in\fvect(n|n+1) \mid L_{D}(\vvol_{q, \xi,
\tau}^a\alpha_{0}^{a-bn})=0\},\text{~~where $\Delta:=\mathop{\sum}_{i\leq
n}\frac{\partial^2 }{\partial q_i
\partial\xi_i}$
}.\end{array}
\end{equation}
It is subject to a~direct verification that $\fb_{a, b}(n)\simeq
\fb_\lambda(n)$ for $\lambda =\frac{2a}{n(a-b)}$. This isomorphism
shows that $\lambda$ actually runs over $\Cee \Pee^1$, not $\Cee$, as
one might hastily think.

For future use, we will denote the operator that singles out
$\fb_{\lambda}(n)$ in $\fm (n)$ as follows:
\begin{equation}
\label{2.7.3} 
 \Div_{\lambda}:=(bn-aE)\nfrac{\partial f}{\partial\tau}-a\Delta\;\text{ for
$\lambda :=\nfrac{2a}{n(a-b)}$ and $\Delta:=\mathop{\sum}\limits_{i\leq
n}\nfrac{\partial^2 }{\partial q_i
\partial\xi_i}$}.
\end{equation}

As follows from the description of $\fvect(m|n)$-modules
(\cite{BL4}) and the simplicity criteria for $\Zee$-graded Lie
superalgebras (\cite{K2}), the Lie superalgebras $\fb_\lambda(n)$
are simple for $n>1$ and $\lambda\neq 0$, 1, $\infty$, see eq.~\eqref{2.7.5}. It is also clear that the Lie superalgebras
$\fb_{\lambda}(n)$ are non-isomorphic for distinct values of $\lambda$, except for the 
isomorphisms already described in list~\eqref{1.7}.

The Lie superalgebra $\fb(n)=\fb_{0}(n)$ is not simple: it has an
$\eps=(0|1)$-dimensional, center. At
$\lambda=1$ and $\infty$, the Lie superalgebras $\fb_{\lambda}(n)$
are not simple either: each of them has an ideal of codimension
$\eps^{n+1}$ and $\eps^{n}$, respectively. The corresponding exact
sequences are\index{$\fb_{1} '(n)$}\index{$\fb_{\infty} '(n)$}
\begin{equation}
\label{2.7.5}
\begin{array}{c}
0\tto \Cee M_{1} \tto \fb(n)\tto \fle(n)\tto 0,\\
0\tto \fb_{1} '(n)\tto \fb_{1}(n)\tto \Cee\cdot
M_{\xi_1\dots\xi_n} \tto 0,\\
0\tto \fb_{\infty} '(n)\tto \fb_{\infty}(n)\tto \Cee\cdot
M_{\tau\xi_1\dots\xi_n} \tto 0.\\
\end{array}
\end{equation}
Clearly, at the exceptional values of $\lambda$, \textit{i.e.}, 0, 1, and
$\infty$, the deformations of $\fb_{\lambda}(n)$ should be
investigated extra carefully. As we will see immediately, this pays:
at each of the exceptional points we find additional deformations. 
The following is an \textbf{Open problem}\index{Problem, open}:

\textbf{What is the reason for the
exceptional deformation at $\lambda=-1$}? 

Other
exceptional values ($\lambda =\frac12$ and $-\frac32$) come from the
isomorphisms 
\[
\text{$\fb_{1/2}(2; 2)\simeq \fh_{1/2}(2|2)=\fh(2|2)$ and
$\fh_{\lambda}(2|2)\simeq \fh_{-1-\lambda}(2|2)$, ~~see list ~\eqref{1.7},}
\]
and the existence of an even non-degenerate bilinear form on $\Pi(\Vol^{1/2})$ which gives rise to a~central extension
\[
0\tto \Cee K_{1} \tto \fpo(2|2)\tto \fh(2|2)\tto 0.
\]
The restriction of the well-known deformation $\cQ$ --- quantization --- of $\fpo(2|2)$ to $\fh(2|2)$ explains the exceptional deformation at $\lambda=\frac12$. For the proof of the next theorem, see \S\ref{S:15}.

\sssbegin[Deformations of $\fb_{\lambda}(n)$]{Theorem}[Deformations of $\fb_{\lambda}(n)$, after \cite{Ko1,LSh3}]\label{thdefb}
For $\fg=\fb_{\lambda}(n)$, set $H=H^2(\fg; \fg)$ for brevity. Then,

\emph{1)} $\sdim ~H=(1|0)$ for $\fg=\fb_{\lambda}(n)$, except for 
$\lambda=0$, $-1$, $1$, $\infty$ for $n>2$. 

\emph{$1_2$)} For $n=2$, in addition
to the above, $\sdim ~H\neq (1|0)$ at $\lambda=\frac12$ and
$\lambda=-\frac32$.

\emph{2)} At the exceptional values of $\lambda$ listed in heading
$1)$, we have

$\sdim ~H=(2|0)$ at $\lambda=\pm 1$ and $n$ odd, or
$\lambda=\infty$ and $n$ even, or $n=2$ and $\lambda=\frac12$ or
$\lambda=-\frac32$.

$\sdim ~H=(1|1)$ at $\lambda=0$ or $\lambda=\infty$ and $n$ odd,
or $\lambda=\pm 1$ and $n$ even.

 Let $k=(k_{1}, \dots , k_{n})$; we set
$q^k=q_{1}^{k_{1}}\cdots q_{n}^{k_{n}}$ and $|k|=\sum k_{i}$. The corresponding cocycles $C$ are given by the following
non-zero values in terms of the generating functions $f$ and $g$,
where $d_{\od}(f)$ is the degree of $f$ with respect to odd
indeterminates only:
\[
\begin{tabular}{|c|c|c|}
 \hline
$\fb_{\lambda}(n)$&$p(C)$&$C(f, g)$\cr \hline \hline
$\fb_{0}(n)$&\textit{odd}&$(-1)^{p(f)}(d_{\od}(f)-1)(d_{\od}(g)-1)fg$\cr
\hline \hline $\fb_{-1}(n)$&$n+1\pmod 2$&$f=q^k, \ g=q^l \longmapsto
(4-|k|-|l|)q^{k+l}\xi_{1}\dots \xi_{n}+$\cr &&$\tau
\Delta(q^{k+l}\xi_{1}\dots \xi_{n})$\cr \hline
$\fb_{1}(n)$&$n+1\pmod 2$&$f=\xi_{1}\dots \xi_{n}, \ g\longmapsto
\begin{cases}(d_{\od}(g)-1)g,&\text{if $g\neq af$, $a\in\Cee$}\\
 2(n-1)f,&\text{if $g=f$ and $n$ is even}\end{cases}$
 \cr
\hline $\fb_{\infty}(n)$&$n\pmod 2$&$f=\tau\xi_{1}\dots \xi_{n}, \
g\longmapsto
\begin{cases}(d_{\od}(g)-1)g,&\text{if $g\neq af$, $a\in\Cee$}\\
2f,&\text{if $g=f$ and $n$ is odd}\end{cases}$\cr \hline \hline
$\fb_{1/2}(2)$&\textit{even}&\textit{induced by the quantization $\cQ$ of $\fpo$, see \cite{Fu,Vo,Tyut}}\cr
\hline
$\lambda$ generic& even&\textit{see} ~\eqref{cocycle}, ~\eqref{bbdef}\cr 
\hline
\end{tabular}
\]

On $\fb_{1/2}(2)\simeq\fb_{-3/2}(2)\simeq\fh(2|2; 1)$, the
cocycle $C$ is the one induced on $\fh(2|2)=\fpo(2|2)/\Cee\cdot K_1$
by the usual deformation (quantization) of $\fpo(2|2)$ with cocycle $Q$ applied to the nonstandard regrading $\fpo(2|2; 1)$ with subsequent passage to the quotient modulo the center
\emph{(spanned by constants)}.

\emph{3)} The space $H$ is diagonalizable with respect to the Cartan
subalgebra of $\fder~ ~\fg$. Let the cocycle $M$ corresponding to the
main deformation be one of the eigenvectors. Let $C$ be another
eigenvector in $H$, it determines a~singular deformation from
item $2)$. The only cocycles $kM+lC$, where $k, l\in\Cee$, that
can be extended to a~global deformation are those for $kl=0$, \textit{i.e.},
either $M$ or $C$.

All the singular deformations of the bracket $\{-, -\}_{old}$ in $\fb_{\lambda}(n)$, except the ones for
$\lambda=\frac12$ or $-\frac32$ and $n=2$, are linear in $\hbar$ even if $\hbar$ is even:
\begin{equation}
\label{cocycle} \{f, g\}_{\hbar}^{sing}=\{f, g\}_{old}+\hbar\cdot
C(f, g)\text{ for any } f, g\in \fb_{\lambda}(n).
\end{equation}
\end{Theorem}

\sssec{Remark: A relation with the Vey cocycle} C.~Roger observed that the singular
deformation (quantization) of $\fb_{0}(n)=\fb(n)$ is, up to a~sign,
the wedge product of two 1-cocycles, the derivations $f\longmapsto
(d_{\od}(f)-1)f$. He also noted that the cocycle on
$\fb_{1/2}(2)\simeq\fh(2|2; 1)$ induced by the quantization of
$\fpo(2|2)$ is a~straightforward superization of the well-known
Vey's cocycle~ \cite{V}.


Since the elements of $\fb_{\lambda}(n)$ are encoded by functions
(for us: polynomials or formal power series) in $\tau$, $q$ and
$\xi$ subject to one relation with an odd left-hand side which contains
$\tau$, it seems plausible that the bracket in
$\fb_{\lambda}(n)$ can be, at least for generic values of the parameter
$\lambda$, expressed solely in terms of $q$ and $\xi$. This is
indeed the case, and here is an explicit formula (in which $\{f,
g\}_{B.b.}$ is the antibracket and
$\Delta=\mathop{\sum}_{i\leq n}\ \frac{\partial^2 }{\partial
q_i\partial\xi_i}$):
\begin{equation}
\label{bbdef}
\renewcommand{\arraystretch}{1.4}
\begin{array}{l}
 \{f_1, f_2\}_{\lambda}^{main}=\{f_1,
f_2\}_{B.b.}+\lambda(c_{\lambda}(f_1, f_2)f_1\Delta f_2 +
(-1)^{p(f_1)}c_{\lambda}(f_2, f_1)(\Delta f_1)f_2),\\
 \text{where $ c_{\lambda}(f_1, f_2)=\frac{\deg
f_1-2}{2+\lambda(\deg f_2 -n)} $},
\end{array}
\end{equation}
and $\deg$ is computed with respect to the standard grading $\deg
q_{i}=\deg \xi_{i}=1$.

\sssec{The exceptional nonstandard regradings
$\Reg_{\fb}$ and $\Reg_{\fh}$}\label{SS:1.3.1}
This is a~regrading\index{$\Reg_{\fh}$}\index{$\Reg_{\fb}$}
of $\fb_{a, b}(2)$ given by the formulas:
\begin{equation}
\label{1.3.1}
 \deg\tau=0, \; \; \deg \xi_1=\deg\xi_2= -1, \; \;
\deg q_1=\deg q_2=1.
\end{equation}
We have the following two cases:
\begin{equation}
\label{1.3.2}
\renewcommand{\arraystretch}{1.4}
\begin{array}{ll}
1)&\text{\underline{$b=0$ or $a=b$}: ~$\fb_{a, 0}(2;
\Reg_{\fb})\simeq \fle(2)$ and
$\fb_{a, a}(2; \Reg_{\fb}) \simeq \fb_{\infty}'(2)$,}\\
&\text{in particular,
$\fg_{-2}=0$;}\\
2)&\text{\underline{$a$ and $b$ generic}: }
\fg_{-2}=\Span(\Le_{\xi_1\xi_2})\text{ and
}\fg_{-1}=\Span(\Le_{\xi_1}, \Le_{\xi_2}, \Le_{Q_1}, \Le_{Q_2}),
\end{array}
\end{equation}
where $Q_1=A\xi_1\xi_2 q_1+B\tau\xi_2$, $Q_2=A\xi_1\xi_2
q_2-B\tau\xi_1$ with some coefficients $A$ and $B$
determined by $a$ and $b$. The bracket on $\fg_{-1}$ is determined
by the odd form $\omega=c\sum dQ_i d\xi_i$, so~ $\fg_0$ must be
contained in $\fm(2)_0$. Direct calculations show that $\sdim
\fg_0=4|4$ and
\begin{equation}
\label{reg1} \fg_0=\fspe(2)\oplus \Cee X, \text{ where }
X=\Le_{a\tau+b\sum q_i\xi_i}.
\end{equation}
Indeed,
\begin{equation}
\label{reg2} \fspe(2)_0\simeq \fsl(2)=\Span(\Le_{q_1\xi_2},
\Le_{q_2\xi_1}, \Le_{q_1\xi_1-q_2\xi_2}), \quad
\fspe(2)_{-1}=\Cee\cdot \Le_{1},
\end{equation}
and
\begin{equation}
\label{reg3} \fspe(2)_1=\Cee \Le_{\alpha \xi_1\xi_2 P(q)+ \beta
\tau\Delta(\xi_1\xi_2P(q))},
\end{equation}
where $P(q)$ is a~degree-2 monomial, $\alpha,\beta\in \Cee$ and $\Delta=\sum \frac{\partial^2 }{\partial q_i
\partial\xi_i}$, see eq.~\eqref{2.4.3}.

The eigenvalues of $X$ on $\fg_{-1}$ are: $b-a$ on $(\fg_{-1})_\ev$
and $b+a$ an $(\fg_{-1})_\od$. Therefore,
\begin{equation}
\label{reg4} \fb_{a, b}(2; \Reg_{\fb})\simeq \fb_{-b, -a}(2)\simeq 
\fb_{b, a}(2).
\end{equation}

For $n>2$, as well as for $\fm(n)$ with $n>1$, similar regradings
are not of Weisfeiler type, as is not difficult to see.

The exceptional grading $\Reg_{\fb}$ of $\fb_{\lambda}(2)$ defines
the exceptional grading $\Reg_{\fh}$ of the isomorphic algebra
$\fh_{\lambda}(2|2)$, see eq.~\eqref{1.7}.

Thus, the exceptional regradings $\Reg_{\fb}$ or $\Reg_{\fh}$ do not
provide us with new W-graded vectorial algebras. Still, they are
needed to interpret the automorphisms ~\eqref{reg4}.

\ssec{The exceptional Lie subsuperalgebra $\fk\fas$ of $\fk
(1|6)$}\label{SS:1.3.4.1}
Like $\fvect(1|m; m)$, the Lie superalgebra $\fk\fas$ is not\index{$\fk\fas$}
determined by its non-positive part and requires a~closer study. The
Lie superalgebra $\fg=\fk(1|2n)$ is generated by the elements of
$\Cee[t, \xi_1, \dots, \xi_n, \eta_1, \dots, \eta_n]$. The standard
$\Zee$-grading of $\fg$ is induced by the $\Zee$-grading of $\Cee[t,
\xi, \eta]$ given by
\begin{equation}\label{a3}
\deg
K_f=\deg f-2, \text{ where }\deg t=2, \quad \deg
\xi_i=\deg\eta_i=1\ \ \text{for all $i$}.
\end{equation}
Clearly, in this grading, $\fg$ is of depth
2. Let us consider the functions that generate several first
homogeneous components of $\fg=\mathop{\oplus}_{i\geq
-2}\ \fg_i$:
\begin{equation}\label{a4}
\renewcommand{\arraystretch}{1.3}
\begin{array}{|c|c|c|c|c|}
\hline
\text{component}&\fg_{-2}&\fg_{-1}&\fg_{0}&\fg_{1}\\
\hline \text{its generators}&1&\Lambda^1(\xi, \eta)&\Lambda^2(\xi,
\eta)\oplus\Cee\cdot t&\Lambda^3(\xi,
\eta)\oplus t\Lambda^1(\xi, \eta)\\
\hline
\end{array}
\end{equation}
As one can prove directly, the component $\fg_1$ generates the whole
subalgebra $\fg_+=\mathop{\oplus}_{i>0}\ \fg_i$. 

The component
$\fg_1$ splits into two $\fg_0$-modules $\fg_{11}= \Lambda^3$ and
$\fg_{12}=t\Lambda^1$. It is obvious that $\fg_{12}$ is always
irreducible and the component $\fg_{11}$ is trivial for $n=1$.

The partial Cartan prolongs of $\fg_{11}$ and $\fg_{12}$ are
well known (here $\fd(\fpo(0|2n))$ is the Lie superalgebra $\fder~(\fpo(0|2n))=\Cee E\rtimes \fpo(0|2n)$, where $E=\sum \xi_i\partial_{ \xi_i} +\sum \eta_i\partial_{ \eta_i}$ is the grading operator):
\[
\begin{matrix}
(\fg_-\oplus\fg_0, \fg_{11})_*^{mk}\simeq \fpo(0|2n)\oplus\Cee\cdot
K_t\simeq \fd(\fpo(0|2n));\\
(\fg_-\oplus\fg_0,
\fg_{12})_*^{mk}=\fg_{-2}\oplus\fg_{-1}\oplus\fg_{0}
\oplus\fg_{12}\oplus\Cee\cdot
K_{t^{2}}\simeq \fosp(2n|2).\end{matrix}
\]

Observe a~remarkable property of $\fk(1|6)$: only for $2n=6$ the
component $\fg_{11}$ splits into 2 irreducible modules; we will
denote by $\fg_{11}^\xi$ the module containing
$\xi_{1}\xi_{2}\xi_{3}$, and by $\fg_{11}^\eta$ the module containing $\eta_{1}\eta_{2}\eta_{3}$.

Observe further, that $\fg_0=\fc\fo(6)\simeq \fgl(4)$. As
$\fgl(4)$-modules, $\fg_{11}^\xi$ and $\fg_{11}^\eta$ are the
symmetric squares $S^2(\id)$ and $S^2(\id^*)$ of the tautological
4-dimensional $\fgl(4)$-module $\id$ and its dual, respectively.

\sssbegin[The partial Cartan prolong $\fk\fas$]{Theorem}[The partial Cartan prolong $\fk\fas$, see \cite{Sh5},
\cite{Sh14}]\label{th1.13} 
The partial Cartan prolong $\fk\fas^\xi=(\fg_-\oplus\fg_0,
\fg_{11}^\xi\oplus\fg_{12})_*^{mk}$ is infinite-dimensional and
simple. It is isomorphic to $\fk\fas^\eta=(\fg_-\oplus\fg_0,
\fg_{11}^\eta\oplus\fg_{12})_*^{mk}$.
\end{Theorem}

\textbf{Convention}. When it does not matter which of the isomorphic algebras
$\fk\fas^\xi\simeq \fk\fas^\eta$ to take, we will simply write
$\fk\fas$, but by default we are always dealing with $\fk\fas^\xi$.

\ssec{The fifteen W-regradings of exceptional simple vectorial 
Lie superalgebras}\label{SS:1.3.2}\mbox{}

\sssbegin[The 15 W-regradings of exceptional 
Lie superalgebras]{Theorem}[The 15 W-regradings (\cite{Sh14}, \cite{CK2})]\label{th1.10}
All W-regradings of the exceptional simple vectorial Lie
superalgebras are the regradings of their ``standard'' ambients
listed in table~$\eqref{1.10}$.
\end{Theorem}

Thus, there are $5$ \textit{isomorphism classes} of abstract
Lie superalgebras consisting of $15$ W-filtered Lie superalgebras. For a~ proof, see Section~\ref{SS:7.2}. 

Table ~\eqref{table321} lists the terms $\fg_{i}$
for $-2\leq i\leq 0$ of the 15 exceptional W-graded algebras.

Observe that none of the simple W-graded vectorial Lie
superalgebras is of depth $>3$ and only two algebras are of depth 3:
$\fm\fb(4|5; K)$, for which 
\begin{equation}\label{a1}
\fm\fb(4|5; K)_{-3}\simeq \Pi(\id_{\fsl(2)}),
\end{equation}
and another one, $\fk\fle(9|6; CK)=\fk\fle(9|11)$, for which 
\begin{equation}\label{a2}
\fk\fle(9|11)_{-3}\simeq \Pi(\id_{\fsl(2)}\boxtimes\Cee[-3]),
\end{equation}
where $\Cee[-3]\simeq\Cee\,\One$ is the module over $\Cee\,\theta\partial_{\theta}\in\fvect(0|1)\subset\fk\fle(9|11)_0$, such that $(\theta\partial_{\theta})\One=-3\One$. Define the representation $T^{1/2}$ of $\fvect(0|n)$\index{$T^{1/2}$} by setting (cf.~\eqref{LieDer})
\begin{equation}\label{T12}
T^{1/2}_D(f\vvol^{1/2})=(D(f)+ (-1)^{p(f)p(D)}f\cdot \frac 12\Div D)\vvol^{1/2}
\end{equation}
for any $D\in\fvect(0|n)$. In table \eqref{table321}, the star ${\star}$ marks the Lie superalgebras $\fg$ of depth 3 whose
components $\fg_{-3}$ are described in eqs.~\eqref{a1} and ~\eqref{a2}.
 
\begin{equation}\label{table321}\footnotesize
\renewcommand{\arraystretch}{1.3}
\begin{tabular}{|c|c|c|c|c|}
\hline $\fg$&$\fg_{-2}$&$\fg_{-1}$&$\fg_0$&$\sdim\fg_{-}$\cr

\hline

\hline $\fv\fle(4|3)$&$-$&$\Pi(\Lambda(3)/\Cee
1)$&$\fc(\fvect(0|3))$&$4|3$\cr

\hline

$\fv\fle(4|3; 1)$&$\Cee[-2]$&$\id_{\fsl(2;\Lambda(2))}$&
$\fc(\fsl(2;\Lambda(2)))\ltimes T^{1/2}(\fvect(0|2))$&$5|4$\cr

\hline

$\fv\fle(4|3; K)$&$\id_{\fsl(3)}\boxtimes
\Cee[-2]$&$\Pi(\id^*_{\fsl(3)}\boxtimes \id_{\fsl(2)}\boxtimes
\Cee[-1]$)&$\fsl(3)\oplus\fsl(2)\oplus\Cee z$&$3|6$\cr

\hline \hline

$\fv\fas(4|4)$&$-$&$\spin$&$\fas$&$4|4$\cr \hline \hline

$\fk\fas$&$\Cee[-2]$&$\Pi(\id)$& $\fc\fo(6)$&$1|6$\cr

\hline

$\fk\fas(; 1\xi)$
&$\Lambda(1)$&$\id_{\fsl(2)}\boxtimes\id_{\fgl(2;\Lambda(1))}$&
$\fsl(2)\oplus(\fgl(2;\Lambda(1))\ltimes \fvect(0|1))$&$5|5$\cr

\hline

$\fk\fas(; 3\xi)$&$-$&
$\Lambda(3)$&$\Lambda(3)\oplus\fsl(1|3)$&$4|4$\cr \hline $\fk\fas(;
3\eta)$&$-$&$\Vol_{0}(0|3)$& $\fc(\fvect(0|3))$&$4|3$\cr

\hline \hline

$\fm\fb(4|5)$&$\Pi(\Cee[-2])$&$\Vol
(0|3)^{1/2}$&$\fc(\fvect(0|3))$&$4|5$\cr \hline

$\fm\fb(4|5; 1)$&$\Lambda(2)/\Cee 1$ &$\id_{\fsl(2;\Lambda(2))}$
&$\fc(\fsl(2;\Lambda(2))\ltimes T^{1/2}(\fvect(0|2))$&$5|6$\cr

\hline

$\fm\fb(4|5; K)^{\star}$&$\id_{\fsl(3)}\boxtimes
\Cee[-2]$&$\Pi(\id^*_{\fsl(3)}\boxtimes \id_{\fsl(2)}\boxtimes
\Cee[-1])$&$\fsl(3)\oplus\fsl(2)\oplus\Cee z$&$3|8$\cr

\hline \hline $\fk\fle(9|6)$&$\Cee[-2]$&$\Pi(T^0_{0}(\vec
0))$&$\fsvect(0|4)_{3, 4}$&$9|6$\cr

\hline

$\fk\fle(9|6; 2)$&$\Pi(\id_{\fsl(1|3)})$&$\id_{\fsl(2;\Lambda(3))}$&
$\fsl(2;\Lambda(3)) \ltimes \fsl(1|3)$&$11|9$\cr

\hline

$\fk\fle(9|6; K)$&$\id$&$\Pi(\Lambda^2(\id^*))$&$\fsl(5)$&$5|10$\cr

\hline

$\fk\fle(9|6;
CK)^{\star}$&$\id_{\fsl(3;\Lambda(1))}^*$&$\id_{\fsl(2)}\boxtimes
\id_{\fsl(3;\Lambda(1))}$&$\fsl(2)\oplus\left (\fsl(3;\Lambda(1))
\ltimes \fvect(0|1)\right)$&$9|11$\cr

\hline
\end{tabular}
\end{equation}

\ssec{W-gradings of serial simple vectorial Lie superalgebras}\mbox{}\label{SS:1.3.3}


\sssbegin[34 W-gradings]{Theorem}[34 W-gradings]\label{th7.3}
Any Weisfeiler regrading $\fh$ of a~given simple vectorial Lie
superalgebra is one of regradings from Tables~\eqref{table3},
~\eqref{table31}, \eqref{table32} and ~\eqref{table321}.
\end{Theorem}

For its proof, see Section~\ref{SS:7.2}. To facilitate
the comparison of various vectorial superalgebras, we offer the
following table~\eqref{table3}. The most interesting phenomena occur for extremal
values of 
$r$ and small values of the superdimension $m|n$.

The central element $z\in\fg_0$ is supposed to be chosen so that
it acts on $\fg_k$ as $k\cdot\id$.

Let $\Lambda (n):=\Cee[\xi_1, \dots, \xi_n]$ be the Grassmann
superalgebra generated by the elements $\xi_i$, each of degree~ 0. We denote (recall Examples
\ref{tensF})\index{$T(\vec 0)$}
\index{$T_0(\vec 0)$} 
\index{$T_0^0(\vec 0)$} 
\index{$\Vol_0(0\vert n)$}

\begin{equation}
\label{1.3.3}
\begin{array}{ll}
\fvect(0|n)\text{-modules}:&\Lambda (0):=\Cee, \quad T(\vec 0):=\Lambda (n),\quad
T_0(\vec 0):=\Lambda (n)/\Cee \cdot 1,\\
&\Vol _0(0|n):=\{v\in\Vol (0|n)\mid \int v=0\};\\
\fsvect(0|n)\text{-modules}:&T_0^0(\vec 0):=\Vol _0(0|n)/\Cee\cdot 1\ .
\end{array}
\end{equation}
Over $\fsvect(0|n)$, it is convenient to consider $\Vol _0(0|n)$
as a~submodule of $\Lambda(n)$.

For the range of the parameter $r$ in tables~\eqref{table3} and~\eqref{table31}, 
see table~\eqref{nonstandgr}. Certain values of~ $r$
(namely, $r=k-1$ for $\fk(1|2k)$, as well as $r=n-1$ for $\fm(n)$
and its subalgebras) are excluded because, for these values of~ $r$,
the corresponding grading is not a~W-grading: the $\fg_{0}$-module
$\fg_{-1}$ is reducible.

\underline{For $N=5', 3'$ and 9}: The $\fg_0$-modules $\fg_{-1}$
are described as $\fvect(0|2)$-modules.

\underline{For $N=16$}: we set $p(\mu)\equiv n\pmod 2$, so $\mu$ can
be an odd indeterminate. The Lie superalgebras $\widetilde{\fsvect}_\mu
(0|n)$ are isomorphic for non-zero values of $\mu$; hence, so are the
algebras $\widetilde{\fs\fb}_{\mu}(2^{n-1}-1|2^{n-1})$. Hence, for $n$
even, we can set $\mu=1$, whereas if $\mu$ is odd, we should
consider it as an additional indeterminate on which the coefficients
depend.

\underline{For $N=1$, 3, etc.}: the symbol $\id_{\fg}$ is short for $\id_{\fg(\sdim;\Lambda(r))}$.

\underline{For $N=11$, etc.}: the symbol
$\id_{\fosp}\boxtimes\Cee[-1]$ denotes the tensor product of $\id_{\fosp(m-2r|2n;\Lambda
(r))}$ by the $\fc$-module $\Cee[-1]$, where 
\[
\fc\subset \fc(\fosp(m-2r|2n))\subset \fc(\fosp(m-2r|2n;\Lambda(r))),
\]
with the natural action of $\fvect(0|r)$ in $\Lambda(r)$. Line $14$ should be similarly understood.

\underline{For $N=21$, 22}: the terms ``$\fg_{-i}$'' denote the
superspace isomorphic to the one in quotation marks, but with the
action given by the rules~\eqref{1.3.5} and ~\eqref{1.3.6}.

\underline{For $N=21-23$}: denote the tautological
$k|k$-dimensional $\fspe(k)$-module by $V$; let $D=\diag
(1_{k}, -1_{k})\in \fpe(k)$. Let $W=V\otimes \Lambda(r)$ and $X\in
\fvect(0|r)$. Let 
\be\label{top}
\Xi=\xi_1\cdots\xi_n\in \Lambda(\xi_1, \dots
,\xi_n).
\ee
Denote by $T^r$, where $r$ is integer, the representations of $\fvect(0|r)$ in
$W$ \index{$T^r$, the representations of $\fvect(0\vert r)$} in $\fspe(n-r)\otimes\Lambda(r)$ and $\fspe(n-r)\otimes
\Lambda(r)$ given by the formula
\begin{equation}
\label{1.3.4} T^r(X)=1\otimes X+D\otimes \nfrac{1}{n-r}\Div X.
\end{equation}
\underline{$\fg=\fsle'(n; r)_{0}$ for $r\neq n-2$}. For $\fg_{0}$,
we have:
\begin{equation}
\label{1.3.5}
\begin{array}{l}
\text{$\fvect(0|r)$ acts on the ideal $\fspe(n-r)\otimes
``\Lambda(r)''$ via
$T^r$, see eq.~\eqref{1.3.4}};\\
\text{any $X\otimes f\in\fspe(n-r)\otimes \Lambda(r)$ acts in
$\fg_{-1}$ as
$\id\otimes f$ and in $\fg_{-2}$ as $0$};\\
\text{any $X\in\fvect(0|r)$ acts in $\fg_{-1}$ via $T^r$ and in
$\fg_{-2}$ as $X$.}
\end{array}
\end{equation}
\underline{$\fg=\fsle'(n; n-2)$}. For $\fg_{0}$, we note that
\[
\fspe(2)\simeq \Cee(\Le_{q_1\xi_1-q_2\xi_2})\ltimes
\Cee\Le_{\xi_1\xi_2},
\]
whereas $\fg_{-2}$ and $\fg_{-1}$ are as above, for $r< n-2$. Set
$\fh=\Cee(\Le_{q_1\xi_1-q_2\xi_2})$. In this case
\begin{equation}
\label{1.3.6} \fg_0\simeq \underline{``\left (\fh\otimes
\Lambda(n-2)\ltimes \Cee\Le_{\xi_1\xi_2}\otimes
(\Lambda(n-2)\setminus \Cee\xi_3\cdots\xi_n)\right )''}\ltimes
T^1(\fvect(0|n-2)).
\end{equation}
The action of $\fvect(0|n-2)$, the quotient of $\fg_{0}$ modulo the
underlined ideal in quotation marks in \eqref{1.3.6} is performed via ~\eqref{1.3.5}. In the subspace
$\xi_1\xi_2\otimes\Lambda(n-2)\subset\fg_0$, this action is the same
as in the space of volume forms. So we can consider everything, except for
the terms proportional to $\Xi$, see eq.~\eqref{top}, or, speaking correctly, take
the irreducible submodule of functions with integral 0.

\underline{For $N=27$, 32, 33}: the terms ``$\fg_{i}$'' denote the
superspace isomorphic to the one in quotation marks, but with the
action given by eqs.~\eqref{1.3.5} and ~\eqref{1.3.7}. 

In the exceptional case $ar=bn$, \textit{i.e.},
$\lambda= \frac{2}{n-r}$, we see that the
$\fvect(0|r)$-action on the ideal
$(\fc\fspe(n-r)\otimes\Lambda(r))\ltimes \fvect(0|r)$ of $\fg_0$, and
on $\fg_-$, is the same as for $\fsle'$, see eq.~\eqref{1.3.5}.

In table~\eqref{table3}, to save space in line 18 we denote the $\fpe(n-r;\Lambda
(r))\ltimes \fvect(0|r)$-module $\id_{\fpe(n-r;\Lambda
(r))}$ by $\id_{\fpe}$, and we use similar abbreviations in other lines.
\index{$\fh_{\lambda}(2\vert 2)$}

Non-standard gradings in cases 1--23, with non-Weisfeiler gradings excluded,
are shown in tables~\eqref{table3}, \eqref{table31}.

{

\footnotesize

\renewcommand{\arraystretch}{1.3}
\begin{equation} \label{table3}
\begin{tabular}{|c|c|c|c|c|}
\hline $N$&$\fg$&$\fg_{-2}$&$\fg_{-1}$&$\fg_0$\cr \hline \hline

$1$&$\fvect(n|m; r)$&$-$&$\id_{\fgl}$&$\fgl(n|m-r;\Lambda
(r))\ltimes \fvect(0|r)$\cr \hline

$2$&$\fvect(1|m; m)$&$-$&$\Lambda (m)$& $\Lambda (m)\ltimes
\fvect(0|m)$\cr \hline \hline

$3$&$\fsvect(n|m; r),\ n\neq 1$&$-$&$\id_{\fsl}$&$\fsl(n|m-r;\Lambda
(r))\ltimes \fvect(0|r)$\cr \hline

$3'$&$\fsvect(2|1)$&$-$&$\Pi(T_0(\vec
0))$&$\fsl(2|1)\simeq \fvect(0|2)$\cr \hline

 \hline
 
$4$&$\fsvect'(1|m; r),\ r\neq m$&$-$&$\id_{\fsl}\boxtimes_{\Lambda
(r)}\Vol_0(0|r)$&$\fsl(n|m-r;\Lambda (r))\ltimes \fvect(0|r)$\cr
\hline

$5$&$\fsvect'(1|m; m)$&$-$&$\Vol_0(0|m)$&$\Lambda (m)\ltimes
\fsvect(0|m)$\cr \hline

$5'$&$\fsvect'(1|2)$&$-$&$T_0(\vec
0)$&$\fsl(1|2)\simeq \fvect(0|2)$\cr 
\hline
\hline

$6$&$\fh(2n|m)$&$-$&$\id_{\fosp}$&$\fosp(m|2n)$\cr \hline

$7$&$\fh(2n|m; r),\ 2r<m$&$T_0(\vec
0)$&$\id_{\fosp}$&$\fosp(m-2r|2n;\Lambda (r))\ltimes \fvect(0|r)$\cr
\hline

$8$&$\fh(2n|2r; r)$&$-$&$\id_{\fsp}$& $\fsp(2n;\Lambda (r))\ltimes
\fvect(0|r)$\cr \hline \hline

$9$&$\fh_{\lambda}(2|2)$&$-$&$\Pi(\Vol^{\lambda}(0|2))$&
$\fosp(2|2)\simeq \fvect(0|2)$\cr \hline

$10$&$\fh_{\lambda}(2|2; 1)$&$-$&
$\id_{\fsp}\boxtimes_{\Lambda
(1)}\Vol^{\lambda}(0|1)$&$\fsp(2; \Lambda
(1))\ltimes \fvect(0|1)$\cr \hline \hline
\end{tabular}
\end{equation}
\footnotesize

\renewcommand{\arraystretch}{1.3}
\begin{equation} \label{table31}
\begin{tabular}{|c|c|c|c|c|}
\hline $N$&$\fg$&$\fg_{-2}$&$\fg_{-1}$&$\fg_0$\cr \hline \hline

$11$&$\begin{matrix}\fk(2n+1|m;
r) \ \text{for $r\neq k-1$}\\ 
\text{if $m=2k$ and $n=0$}\end{matrix}$&$\Lambda(r)\boxtimes\Cee[-2]$&$\id_{\fosp}\boxtimes\Cee[-1]$&$\fc\fosp(m-2r|2n;\Lambda
(r))\ltimes \fvect(0|r)$\cr \hline

$12$&$\fk(1|2m; m)$&$-$&$\Lambda(m)$&$\Lambda (m) \ltimes
\fvect(0|m)$\cr \hline

$13$&$\fk(1|2m+1; m)$&$\Lambda(m)$&$\Pi(\Lambda (m))$&$\Lambda
(m)\ltimes \fvect(0|m)$\cr \hline
\hline

$14$&$\fm(n;
r),\ r\neq n$&$\Pi(\Lambda(r))\boxtimes\Cee[-2]$&$\id_{\fpe}$&
$\fc\fpe(n-r;\Lambda(r))\ltimes
\fvect(0|r)$\cr \hline

$15$&$\fm(n;
n)$&$-$&$\Pi(\Lambda(n))$&$ \Lambda(n)\ltimes \fvect(0|n)$\cr \hline
\hline

$16$&$\widetilde{\fs\fb}_{\mu}(2^{n-1}-1|2^{n-1})$&$-$
&$\frac{\Pi(\Vol(0|n))}{\Cee(1+\mu\xi_1\dots\xi_n)\vvol(\xi)}$&
$\widetilde{\fsvect}_\mu (0|n)$\cr \hline

\hline

$17$&$\fle(n)$&$-$&$\id_{\fpe}$&$\fpe(n)$\cr \hline

$18$&$\fle(n; r)$&$\Pi(T_0(\vec 0))$&$\id_{\fpe}$&$\fpe(n-r;\Lambda
(r))\ltimes \fvect(0|r)$\cr \hline

$19$&$\fle(n; n)$&$-$&$\Pi(T_{0}(\vec 0))$& $\fvect (0|n)$\cr \hline

\hline $20$&$\fsle'(n)$&$-$&$\id_{\fspe}$&$\fspe(n)$\cr \hline

$21$&$\fsle'(n; r),\ r\neq n-2, n$&``$\Pi(T_0(\vec 0))$'' & ``$\id\boxtimes\Lambda
(r)$'' &$(\fspe(n-r)\otimes \Lambda (r)')\ltimes
T^r(\fvect(0|r))$\cr \hline

$22$&$\fsle'(n; n-2)$&``$\Pi(T_0(\vec 0))$'' & ``$\id\boxtimes\Lambda
(r)$''&see eq.~\eqref{1.3.6}\cr \hline

$23$&$\fsle'(n; n)$&$-$&$\Pi(T^0_{0}(\vec 0))$&$\fsvect(0|n)$\cr
\hline
\end{tabular}
\end{equation}

}


\underline{For $N=24-34$}: 1) We assume that
$\lambda=\nfrac{2a}{n(a-b)}\neq 0$, 1, $\infty$; these three
exceptional cases (corresponding to the ``drop-outs'' $\fle(n)$,
$\fb_{1}'(n)$ and $\fb'_{\infty}(n)$, respectively) are considered
se\-pa\-rately.

2) The irreducibility condition of the $\fg_0$-module $\fg_{-1}$ for
$\fg=\fb'_{\infty}$ excludes $r=n-1$.

3) The case where $r=n-2$ is extra exceptional, so in table~\eqref{table3} we assume, unless otherwise specified, that
\begin{equation}
\label{*1} 0<r<n-2;\; \; \text{ additionally $a\neq b$ and $(a,
b)\neq \alpha(n, n-2)$ for any $\alpha\in\Cee$.} 
\end{equation}

\underline{For $N=24-27$}: we consider $\fb_{a,b}(n; r)$ for
$0<r<n-2$ and $ar-bn\neq 0$; in particular, this excludes
$\fb'_{\infty}(n; n)=\fb'_{a, a}(n; n)$ and $\fb'_{1}(n;
n-2)=\fb'_{n, n-2}(n; n-2)$.

If $z$ is the central element of $\fc\fspe(n-r)$ that acts on
$\fg_{-1}$ as $-\id$, then
\begin{equation}
\label{1.3.7} z\otimes \psi\text{ acts on $\fg_{-1}$ as
$-\id\otimes \psi$, and on $\fg_{-2}$ as $-2\id\otimes\psi$.} 
\end{equation}
Set
\begin{equation}
\label{c} c:=\nfrac{a}{ar-bn}.
\end{equation}

Let $a\cdot \str\otimes\id$ be the representation of
$\fpe(n-r)=\fs\fpe(n-r)\ltimes \Cee D$, which is $\id_{\fs\fpe(n-r)}$
on $\fs\fpe(n-r)$ and sends $D=\diag
(1_{n-r}, -1_{n-r})\in \fpe(n-r)$ to $2a\cdot\id$. 

Non-standard gradings in cases 24--34 are shown in table~\eqref{table32}.
{
\footnotesize
\renewcommand{\arraystretch}{1.3}

\begin{equation}\label{table32}
\begin{tabular}{|c|c|c|c|c|}\hline

 $24$&$\fb_{\lambda}(n)$&$\Pi(\Cee[-2])$&$\id$&
 $\fspe(n)\ltimes\Cee(az+bd)$\cr
 \hline

 $25$&
 $\arraycolsep=0pt
 \begin{array}{l}
 \fb_{\lambda}(n;r),\\ r<n-2
 \end{array}$
 & \tiny{$\Pi((-c)\str)\boxtimes
 \Vol^{2c}(0|r)$} &
 $\Big(\!\Big(\!\Big(\!\!-\frac{c}{2}\Big)\str\Big)\otimes\id\Big)
 \boxtimes\Vol^{c}(0|r)$&\tiny{
 $(\fpe(n{-}r)\otimes\Lambda(r))\ltimes \fvect(0|r)$}\\
 \hline

 $26$&$\fb_{\lambda}(n;n)$& --- &
 $\Pi(\Vol^{\lambda}(0|n))$&$\fvect(0|n)$\cr
 \hline
 \hline

 $27$ &
 $\arraycolsep=0pt
 \begin{array}{l}
 \fb_{2/(n-r)}(n; r),\cr
 r<n-2
 \end{array}$
 & $\Pi(\Cee)\boxtimes \Lambda
 (r)$&$\id\boxtimes ``\Lambda (r)$''& \tiny{$\fc\fpe(n-r)\otimes ``\Lambda
 (r)$'' $\ltimes T^r(\fvect(0|r))$}\cr
 \hline

 $28$&$\fb'_{\infty}(n)$&$\Pi(\Cee)$&$\id$& $\fspe(n)_{a, a}$\cr
 \hline

 $29$&
 $\arraycolsep=0pt
 \begin{array}{l}
 \fb'_{\infty}(n; r),\cr
 r<n-2
 \end{array}$ & $\Pi(\Cee)\boxtimes \Lambda(r)$
 &$\id\boxtimes\Lambda(r)$& \tiny{$((\fspe(n-r)_{a,a})\otimes\Lambda (r))\ltimes \fvect(0|r)$}\cr
 \hline

 $30$&$\fb'_{\infty}(n; n)$& --- &$\Pi(\Lambda(n))$&
 $(\Lambda(n)\setminus\Cee\cdot\Xi)\ltimes \fsvect(0|n)$\cr
\hline

 $31$&$\fb'_{1}(n)$&$\Pi(\Cee)$&$\id$&$\fspe(n)_{n, n-2}$\cr
 \hline

 $32$&
 $\arraycolsep=0pt
 \begin{array}{l}
 \fb'_{1}(n; r)\cr
 r<n-2
 \end{array}$
 &``$\Pi(\Vol_{0}(0|r))$''& $\id\boxtimes ``\Lambda(r)$''
 & \tiny{$((\fspe(n{-}r)_{n, n-2})\otimes ``\Lambda(r)'')\ltimes T^r(\fvect(0|r))$}\cr
 \hline

 $33$&$\fb'_{1}(n; n-2)$&``$\Pi(T_0(\vec 0))$'' & $\id\boxtimes``\Lambda (r)$'' &
 \parbox{53mm}{\tiny{\mbox{}\\ $((\fspe(n{-}r)_{n, n-2})
 \otimes ``\Lambda(r)$'') as in ${N=32}$
 with $\fc\fspe(2)$ instead of $\fspe(2)$}}
 \cr \hline

 $34$&$\fb'_{1}(n; n)$& --- &$\Pi(\Vol_{0}(0|n))$&$\fvect(0|n)$\cr
\hline\end{tabular}
\end{equation}%
}

\normalsize

\ssec{Grozman's theorems and a~description of $\fg$ as
$\fg_{\ev}\oplus\fg_{\od}$}\label{SS:2.19} It is convenient to
describe the Lie superalgebra $\fg$ of twisted multi-vector fields
as a~$\fvect$-module. Similarly, in \cite{CK2} the exceptional
algebras are described as $\fg=\fg_{\ev}\oplus\fg_{\od}$. For most
of the series, such description is of little value because each
homogeneous component $\fg_{\ev}$ and $\fg_{\od}$ has a~complicated
structure. For the exceptional simple Lie superalgebras,
the situation is totally different! 

Observe that apart from being
beautiful, such a~description is useful for constructing simple
\textit{Volichenko algebras}\index{Volichenko algebra} (inhomogeneous with respect to parity subalgebras of Lie superalgebras), cf. \cite{LS}, and simple Lie algebras in
characteristic~ 2.

In relation to this, we recall a~theorem of Grozman. He completely
described bilinear differential operators acting in spaces of
tensor fields on any $m$-dimensional manifold $M$ and invariant under all changes of coordinates.
Since the problem is local, we can assume that $M\simeq \Cee^m$, see \cite{BL4}. Miraculously, ``most" of the first-order invariant operators
determines a~Lie superalgebra on its domain of definition. Some of these
superalgebras turn out to be very close to simple. In the
constructions below we use some of these invariant operators.

Let $\rho$ be an irreducible representation of the group $GL(n)$ in
a~finite-dimensional vector space $V$ and let $\lambda=(\lambda_1 ,
\dots , \lambda_n)$ be its lowest weight. A~ \textit{tensor field of
type} $\rho$ (or of type $V$: depending on the emphasis) on an
$n$-dimensional connected manifold $M$ is any section $t$ of the
locally trivial vector bundle over $M$ with fiber $V$ such that,
under the change of coordinates, 
\[
t(y(x))=\rho\left(\pderf{y}{x}\right)t(x).
\] 
The space of tensor field of type $\rho$ or $V$ will be denoted
by $T(\rho)$, or $T(V)$, or even $T(\lambda)$. Examples:
\[
\begin{array}{l}
T(\mu, \dots , \mu)=\Vol^\mu,\text{~~the space of
$\mu$-densities};\\
\text{$\Omega^r:=T(0, \dots, 0, 1, \dots, 1) $, with\index{$\Omega^r$}\index{$L^r$}
$r$ symbols 1, is the space of
differential $r$-forms},\\
\text{in particular, $T(0)=\Omega^0$ is the space of
functions also denoted by $\cF$}; \\
\fvect(n)=T(-1, 0, \dots, 0); \\
L^r:=T(-1, \dots , -1, 0, \dots, 0)
\text{~~with
 symbols 1, is the space of $r$-vector fields,}\\
\text{the
$r$th exterior power (over the algebra of functions $\cF=\Omega^0$) of $\fvect(n)$}. \\
\end{array}
\]

\textit{The spaces of twisted $r$-forms and twisted $r$-vector
fields with twist} $\mu$ are defined to be, respectively,\index{$\Omega_{\mu}^r$}\index{$L_{\mu}^r$}
\[
\Omega^r_{\mu}=\Omega^r\otimes_{\Omega^0}\Vol^\mu\text{ and
}L^r_{\mu}=L^r\otimes_{\Omega^0}\Vol^\mu.
\]
Obviously, $L^r_{\mu}\simeq\Omega^{n-r}_{\mu-1}$ and
$\Vol^1=\Omega^n$.

The following statements are excerpts from Grozman's difficult
result \cite{G}. To describe one of the operators, $P_4$, we need
the \textit{Nijenhuis bracket}, defined by the formula
\begin{equation*}
\begin{array}{l}
 \omega^k\otimes \xi, \omega^l\otimes \eta\longmapsto
 (\omega^k\wedge \omega^l)\otimes [\xi, \eta]\ +\\
+\left(\omega^k\wedge L_{\xi}(\omega^l)+(-1)^k d\omega^k\wedge
\iota_\xi(\omega^l)\right )\otimes \eta+
\left(-L_{\eta}(\omega^k)\wedge\omega^l+(-1)^l
\iota_\eta(\omega^k)\wedge d\omega^l)\right )\otimes \xi,
\end{array}
\end{equation*}
where $\iota_{X}$ is the inner product by the field $X$ and $L_{X}$ is the Lie
derivative along $X$. The Nijenhuis bracket admits the
following interpretation which makes its invariance manifest: the
operator $D: (\Omega^{k} \otimes_{\Omega^0} \fvect(M),
\Omega^{\bcdot}) \tto \Omega^{\bcdot}$ given by the formula
\[
\begin{gathered}
D(\omega^k\otimes \xi, \omega):= 
d\left(\omega^k\wedge \iota_\xi(\omega) +(-1)^k\omega^k\wedge
\iota_\xi(d\omega)\right )= d\omega^k\wedge \iota_\xi(\omega)
+(-1)^k\omega^k\wedge L_{\xi}(\omega)
\end{gathered}
\]
is, for a~fixed $\omega^k\otimes \xi$, a~superderivation of the
supercommutative superalgebra $\Omega^{\bcdot}$, and the Nijenhuis
bracket is just the supercommutator of these superderivations. So we
can identify $\Omega^{\bcdot} \otimes_{\Omega^0} \fvect(M)$ with the
Lie subsuperalgebra $C(d)\subset \fvect(\widehat M)$, where
$\widehat M$ is the supermanifold $(M, \Omega^{\bcdot}(M))$, \textit{i.e.},
$C(d)$ is the centralizer of the exterior differential on $\Omega^{\bcdot}(M))$:
\[
C(d) = \{D\in\fvect(\widehat M)\mid [D, d] = 0\}.
\]

\textbf{Dualizations}. To any mapping $F:T(V)\tto T(W)$ there corresponds the dual mapping
$F^*:(T(W))^*\tto (T(V))^*$. If we consider tensors with
compact support, so integration can be performed, we can identify
\[
(T(V))^*\simeq T(V^*)\otimes_{\Omega^0} \Vol\simeq T(V^*\otimes
\str), 
\]
where $\str$ denotes the 1-dimensional $\fgl$-module given by the
trace (or supertrace, if $M$ is a~supermanifold). The
\textit{formal dual}\index{Dual, formal} is defined to be 
\[
(T(V))^*:=T(V^*\otimes \str).
\]

Given a~bilinear mapping $F:T(V_1)\otimes T(V_2)\tto T(W)$, it is possible to
dualize it with respect to each argument:
\[
\renewcommand{\arraystretch}{1.4}
\begin{array}{l}
F^{*1}:T(W^*\otimes \str)\otimes T(V_2)\tto T(V_1^*\otimes \str),\\
F^{*2}:T(V_1)\otimes T(W^*\otimes \str)\tto T(V_2^*\otimes \str).
\end{array}
\]

\sssbegin[Grozman's theorem: $1$st-order invariant operators]{Theorem}[\emph{\cite{G}}]\label{SS:2.20} The invariant
under arbitrary changes of variables irreducible differential
bilinear operators $D: T(\rho_1)\otimes T(\rho_2)\tto T(\rho_3)$ of
order $1$ are, up to dualizations and permutation of arguments, only
the following ones:
\[
{\bf P_{1}}: \Omega^r\otimes T(\rho_2) \tto T(\rho_3),\qquad (w,t)
\longmapsto Z(dw, t),
\]
where $Z$ is an
order-$0$ operator, the extension of the projection
$\rho_1\otimes \rho_2\tto \rho_3$ onto any of the irreducible
components;
\[
\renewcommand{\arraystretch}{1.4}
\begin{gathered}
{\bf P_2}: \fvect \otimes T(\rho) \tto T(\rho), \quad
\text{ \underline{the Lie derivative}}; \\
{\bf P_3}: T(S^p(\id^*))\otimes T(S^q(\id^*)) \tto T(S^{p + q
-1}(\id^*)), \quad \text{ \underline{the Poisson
bracket}};\end{gathered}
\]

${\bf P_4}$: On manifolds, the bracket in $C(d)$ is called the
\underline{ Nijenhuis bracket}. This bracket is a~linear combination
of operators $P_1$, $P_1^{*1}$, their composition with the
permutation operator $T(V)\otimes T(W) \tto T(W)\otimes T(V)$, and
a~new, irreducible, operator denoted $P_4$;
\[
{\bf P_5}: \Omega ^p\otimes \Omega ^q \tto \Omega ^{p + q + 1};
\quad \omega_1,\ \omega_2\longmapsto
(-1)^{p(\omega_1)}a(d\omega_1\cdot \omega_2) + b(\omega_1d\omega_2),
\text{ where } a, b\in\Cee;
\]

${\bf P_6} : \Omega ^p_{\mu}\otimes\Omega ^q_{\nu}\tto \Omega
^{p+q+1}_{\mu+\nu} $, where $|\mu|^2 + |\nu|^2 \neq 0$ and $p + q<n$
is given by the formula
\[
\omega_1\vvol^\mu,\ \omega_2\vvol^\nu\longmapsto \left (\nu
(-1)^{p(\omega_1)}d\omega_1\cdot \omega_2 -
\mu\omega_1d\omega_2\right)\vvol^{\mu+\nu};
\]
\[
{\bf P_7} : L^p\otimes L^q \tto L^{p + q - 1}\quad\text{ the
\underline{Schouten bracket}};
\]

${\bf P_8}: L^p_{\mu}\otimes L^q_{\nu} \tto L^{p+q-1}_{\mu+\nu}$, a
deformation of the Schouten bracket given by the next formula on
manifolds for $p + q \leq n$, and on supermanifolds of
superdimension $n|1$ for $p, q\in\Cee$:
{\footnotesize\begin{equation} \label{P8}
\renewcommand{\arraystretch}{1.4}
\begin{array}{l}
 X\vvol^\mu, \ Y\vvol^\nu \longmapsto\\
 \left ((\nu-1)(\mu + \nu -1)
 \Div X \cdot Y + 
 (-1)^{p(X)}(\mu-1)(\mu + \nu -1)X \Div Y
 - 
(\mu- 1)(\nu -1) \Div(XY)\right )\vvol^{\mu+ \nu},
\end{array}
\end{equation}
} where the \textit{divergence of a~multi-vector field} $f$ is
defined to be
\[
\Div(f)=\mathop{\sum}\limits_{i\leq n}\nfrac{\partial^2
f}{\partial_{x_i}
\partial_{\check x_{i}}}
\]
in local coordinates $(x, \check x)$ on the supermanifold $\check
M$ associated to the sheaf of sections of the exterior algebra of
the tangent bundle on any supermanifold $M$; here the checked
coordinates on $\check M$ are
\[
 \theta_i=\Pi\left(\pder{x_{i}}\right)=\check x_{i}, \emph{~~(see
Remark \ref{rem2.3.7})}.
\]
\end{Theorem}


\sssbegin[Lie and associative (super)algebras determined by invariant operators]{Theorem}[Lie and associative (super)algebras determined by invariant operators, \cite{G}]\label{SS:2.21} The following
natural invariant operators determine associative or Lie
superalgebras on their domains:

$P_3$, $P_4$, $P_7$, and $P_8$ define Lie
superalgebras. The vectorial Lie
superalgebra $C(d)$ is, however, not transitive.

$P_5$: For $ab=0$ and, of course, $|a|^2+|b|^2\neq 0$, it determines an associative superalgebra
structure.

For $a=b\neq 0$, it determines the structure of a~nilpotent Lie
superalgebra on $\Pi(\Omega^{\bcdot})$.

The bracket given by $d\omega_1\cdot \omega_2-\omega_1d\omega_2$
determines a~Lie algebra structure on the space $\Omega^{\bcdot}$.

$P_6$ multiplied by $\nfrac{\mu-\nu}{\mu\nu}$ determines structures
of nilpotent Lie superalgebras on the superspaces
\[
 \Pi\left(\mathop{\oplus}\limits_{\lambda\in\Cee}\Omega^{\bcdot}\otimes
\Vol^\lambda/d\Omega^{\bcdot}\right), \quad
\Omega^{\bcdot}_{+}=\Pi(\mathop{\oplus}\limits_{\lambda>0;\;
\lambda\in \Ree}\Omega^{\bcdot}\otimes \Vol^\lambda), \quad
\Omega^{\bcdot}_{-}=\Pi(\mathop{\oplus}\limits_{\lambda<0;\;
\lambda\in \Ree}\Omega^{\bcdot}\otimes \Vol^\lambda).
\]

$P_6$ multiplied by $\nfrac{\mu-\nu}{\mu+\nu}$ determines a~nilpotent Lie superalgebra structure on
\[
 \Pi\left(d\Omega^{\bcdot}\bigoplus \mathop{\oplus}\limits_{\lambda\neq
0;\; \lambda\in \Ree}\Omega^{\bcdot}\otimes \Vol^\lambda\right).
\]

$P_8$ multiplied by $\nfrac{1}{\mu\nu}$ determines a~nilpotent Lie
superalgebra structures on $\Omega^{\bcdot}_{+}$ and
$\Omega^{\bcdot}_{-}$.
\end{Theorem}

\paragraph{Remarks} The claims of this Theorem are breezily formulated in the short note \cite{G}, here they are formulated more explicitly. Their proof is done by a~ direct verification.

The operator $P_8$ is a~deformation of the
Schouten bracket we considered in $\fb_\lambda(n)$, so this simple Lie superalgebra could have been described 20 years earlier than in \cite{LSh3} if more attention was paid to Grozman's results. 

The operator $P_6$ is the ``Fourier transform" of the
Schouten bracket. 

The divergence of multi-vector fields is the ``Fourier transform" of the exterior derivation. 

Grozman proved that there are no invariant differential operators of degree $>3$, while the operators of degree 2 and 3 are compositions of degree-1 operators, with one exceptional degree-3 operator $Gz$ on the line now bearing Grozman's name. For a~version of Grozman's theorem over fields of characteristic $p>0$, see \cite{BoLe}. For the extension of Grozman's theorem to non-differential operators on the line, see \cite{IoMa}.

\ssec{Weisfeiler-type gradings with infinite-dimensional
components}\label{SS:2.22} It is sometimes convenient to consider
gradings of our vectorial Lie superalgebras which lead them out of
the GLAPG class, but still preserve the finiteness of the
depth and preserve the two characteristic features of the
W-gradings:
\begin{equation}\label{1-3}
\text{\begin{minipage}[c]{14cm} $1)$ for these grading
$\oplus\fg_i$, the subalgebra $\fg_{\geq
0}:=\mathop{\oplus}_{i\geq 0}\ \fg_i$ is maximal;\\
$2)$ the $\fg_0$-module $\fg_{-1}$ is irreducible. Namely, we set
$\deg x_{i}=0$ for $r$ even indeterminates $x$; we denote this
grading $\bar r$.
\end{minipage}}
\end{equation}
For the exceptional simple vectorial Lie superalgebras, and $\fb_{\lambda}(n; \bar
n)$, the multiplications are described as follows:

\underline{$\fg=\fk\fle(5|10)$}: here,
$\fg_\ev=\fsvect(5|0)\simeq d\Omega^{3}(5|0)$ and
$\fg_\od=\Pi(d\Omega^{1}(5|0))$ with the natural $\fg_\ev$-action on
$\fg_\od$ and the bracket of any two odd elements is given by their product. Since $\fsvect(5|0)$ preserves the volume element $\vvol$ and all its powers,
we can make, for any permutation
$(ijklm)$ of $(12345)$, the following identification
\[
dx_{i}\wedge dx_{j}\wedge dx_{k}\wedge dx_{l}:=dx_{i}\wedge dx_{j}\wedge dx_{k}\wedge dx_{l}\otimes
\vvol^{-1}=\sign(ijklm)\partial_{x_{m}}.
\]

\underline{$\fg=\fv\fas(4|4)$}: here, $\fg_\ev=\fvect(4|0)$
and $\fg_\od=\Omega^{1}(4|0)\otimes_{\Omega^0(4|0)} \Vol^{-1/2}(4|0)$ with the
natural $\fg_\ev$-action on $\fg_\od$ and the bracket of odd
elements is given by
\[
\left[\nfrac{\omega_{1}}{\sqrt{\vvol}},
\nfrac{\omega_{2}}{\sqrt{\vvol}}\right]=
\nfrac{d\omega_{1}\wedge\omega_{2}+ \omega_{1}\wedge d\omega_{2}}{
\vvol},
\]
where we identify
\[
\nfrac{dx_{i}\wedge dx_{j}\wedge dx_{k}}{
\vvol}=\sign(ijkl)\partial_{x_{l}}\text{ for any permutation $(ijkl)$
of $(1234)$}.
\]

\underline{$\fg=\fv\fle(3|6)$}: here,
$\fg_\ev=\fvect(3|0)\oplus \fsl(2)^{(1)}_{\geq 0}$, where
$\fg^{(1)}_{\geq 0}=\fg\otimes\Cee[x_{1}, x_2, x_3]$, and
\[
 \fg_\od=\left(\Omega^{1}(3|0)\otimes_{\Omega^0(3|0)} \Vol^{-1/2}(3|0)\right
)\boxtimes \id_{\fsl(2)^{(1)}_{\geq 0}}, \text{ with the natural
$\fg_\ev$-action on $\fg_\od$.}
\]

Recall that $\id_{\fsl(2)}$ denotes the irreducible $\fsl(2)$-module
$L^1$ with highest weight 1; its tensor square splits into the
symmetric square $L^2\simeq \fsl(2)$ and the exterior square
--- the trivial module $L^0$; accordingly, denote by $v_{1}\wedge
v_{2}$ and $v_{1}\bullet v_{2}$ the projections of $v_{1}\otimes
v_{2}\in L^1\otimes L^1$ onto the anti-symmetric and symmetric
components, respectively. For any 
$\omega_{1}, \omega_{2}\in\Omega ^1$ and $v_{1}, v_{2}\in L^1$, we
set
\[
\renewcommand{\arraystretch}{1.4}
\left[ \nfrac{\omega_{1}\otimes v_{1}}{\sqrt{\vvol}},
 \nfrac{\omega_{2}\otimes v_{2}}{\sqrt{\vvol}}\right]=
 \nfrac{(\omega_{1}\wedge \omega_{2})\otimes
(v_{1}\wedge v_{2})+(d\omega_{1}\wedge \omega_{2}+\omega_{1}\wedge
d\omega_{2})\otimes (v_{1}\bullet v_{2})}{\vvol},
\]
where we naturally identify $\Omega ^0$ with $\Omega
^3\otimes_{\Omega ^0} \Vol^{-1}$; we also naturally identify $\Omega ^2\otimes_{\Omega ^0}
\Vol^{-1}$ with $\fvect(3|0)$ by setting
\[
 \nfrac{dx_{i}\wedge dx_{j}}{ \vvol}=\sign(ijk)\partial_{x_{k}}\text{
for any permutation $(ijk)$ of $(123)$}.
\]

\underline{$\fg=\fb_{\lambda}(n; \bar n)$}: We have
\[
\fg_i=\left(\Pi(\Lambda ^{i+1}_{\Omega^0}(\fvect(n|0))\right)\otimes
\Vol^{i\lambda}
\] for $i=-1, 0, \dots , n-1$. The multiplication
is given by Grozman's operator $P_{8}$, see ~\eqref{P8}.

Consider the case where $n=2$ more closely. Clearly, one can
interchange $\fg_{-1}$ and $\fg_{1}$; this possibility explains the
isomorphism
\be\label{unmyst}
\fb_{\lambda}(2; \bar 2)\simeq \fb_{1-\lambda}(2; \bar 2), 
\ee
and
hence the mysterious isomorphism $\fh_{\lambda}(2|2)\simeq
\fh_{-1-\lambda}(2|2)$ mentioned in ~\eqref{1.7}.

In particular, we have additional outer automorphisms of
$\fb_{1/2}(2; \bar 2)$, whereas for $\lambda=\nfrac12$ and
$\lambda=-\nfrac32$ there are nontrivial central extensions of $\fh_{\lambda}(2|2)$.

Therefore, for $n=2$, to the exceptional values of $\lambda=0$, 1 and
$\infty$ existing for any $n>2$, we should add $\lambda=-1$ (for which $\fg_1\simeq\Vol$, and hence $\fg$ is not simple) and $\lambda=2$.

\underline{$\fg=\fm\fb(3|8)$}: See \S \ref{Smb}.

\ssec{On structures preserved}\label{2.10}
In the space
$\fg_{-1}=\Lambda(n)$, the Lie superalgebra $\fg_{0}=\fvect(0|n)$ preserves the tensor $\mult$ of valence $(2, 1)$ that
defines the multiplication. If $\fg_{-1}=\Lambda(n)/\Cee\cdot 1$,
then the structure preserved is the tensor $\mult'$ that defines
multiplication on the maximal ideal of $\Lambda(n)$.

The exceptional simple
vectorial Lie superalgebras $\fg$ preserve various distributions described in \cite{KLS}. The structures preserved by the serial simple Lie superalgebras are described by the structures preserved by $\fg_{0}$ on
$\fg_{-1}$. They are as follows:

(1) The tensor products $B\otimes\mult$ of a~bilinear (or a~volume)
form $B$ preserved (perhaps, conformally, up to multiplication by
a~scalar) in the fiber of the vector bundle over a~$0|r$-dimensional
supermanifold on which the structure governed by $\mult$ is
preserved,

(2) The tensor $\mult'$, or $\mult$ twisted by a~$\lambda$-density. Observe
that the volume element $B$ can be not only $\vvol_\xi$, but also
$(1+\alpha\xi_{1}\ldots\xi_{n})\vvol_\xi$.

The structures of other types, namely certain pseudodifferential
forms preserved by $\fb_{\lambda}(n)$, are already described.


\subsubsection{Remark: On Shmelev's interpretations of
$\fh_{\lambda}(2|2)\protect\simeq \fb_{\lambda}(2; 2)$
}\label{SS:2.10.1} It seemed that $\fb_{\lambda}(2)$ had one more
realization: G.~Shmelev interpreted $\fh_{\lambda}(2|2)\simeq 
\fb_{\lambda}(2; 2)$ as a~deformation $\fh_{\lambda}(2|2)$ of the
Lie superalgebra $\fh(2|2)$ of Hamiltonian vector fields as
preserving the pseudodifferential form
\[
d\eta^{\nfrac{2\lambda-1}{1-\lambda}}\left(\lambda dqdp+
(1-\lambda)d\xi d\eta\right);
\]
he also described $\fh_{\lambda}(2|2)$ as preserving the
pseudodifferential form
\[
d\eta^{\nfrac{1}{\lambda}-2}\left(dqdp+d\xi d\eta\right).
\]
These interpretations are cited (but not used) in \cite{ALSh, Po1}, but verification shows that Shmelev made a~mistake
calculating Lie derivatives, so these interpretations are wrong.


\section{Step 2. Criterion of infinite dimensionality of 
the Cartan prolong}\label{S:3}

In this section, we superize Wilson's proof
of the following Guillemin--Quillen--Sternberg theorem.

\ssbegin[Main technical theorem]{Theorem}[Main technical theorem]\label{th3.1}
Let $\fg=\mathop{\oplus}_{i\geq - d}\ \fg_i$ be an
infinite-dimensional transitive Lie algebra. Then, there exists an
element $g \in \fg_0$ such that $\rk (\ad_g) |_{\fg_{-1}}=1$.
\end{Theorem}

The converse statement is obvious: if such $g$ exists, then,
clearly, the Cartan prolong of $(\fg_{-}, \fg_0)$ is infinite-dimensional. Indeed, in some basis, $g$ takes the form
$x\partial_y$ and all the elements
$x^n\partial_y$ belong to the Cartan prolong.

\ssec{A reformulation of Theorem~\ref{th3.1}}
Let $\cA= \mathop{\oplus}_{i\leq 0}\ \cA_i$ be a~commutative
$\Zee$-graded algebra with unit generated by $A:=\cA_{-1}$ and let
$M= \mathop{\oplus}_{i\geq k}\ M_i$ be a~$\Zee$-graded
$\cA$-module. (In the case we are interested in, $\cA=
S^{\bcdot}(\fg_{-1})$ and $M=\fg$\,.) We will say that $M$ is a~
\textit{transitive module}\index{Module, transitive} if
\[
Am\neq 0~~\text{ for any non-zero }\ m\in M_i \text{ and any }i\geq
k+1.
\]
Notice that each element $m\in M$ determines a~linear mapping \[ m:
A\tto M, \quad a\longmapsto am. \] Clearly, $M$ is transitive if and
only if there exists an $i\geq k+1$ and a~non-zero $m\in M_i$ such
that $\ker m \neq A$.

Now it is clear that Theorem \ref{th3.1} can be reformulated in the
following equivalent form that is more convenient for us.

\sssbegin[Main technical theorem, a~version]{Theorem}[Main technical theorem, a~version]\label{th3.2}
Let $M= \mathop{\oplus}_{i\geq k}\ M_i$ be a~$\Zee$-graded
infinite-dimensional transitive $\cA$-module. Then, there exists
$m\in M_{k+1}$ such that ${\codim \ker m= 1}$; hence, $\rk~m=1$.
\end{Theorem}

\begin{proof}
Without loss of generality we can assume that $k=-1$, \textit{i.e.}, $M=
\mathop{\oplus}_{i\geq -1}\ M_i$ and $m\in M_0$.

$1)$ Since $\cA$ is commutative, $\ker am\supset \ker m$ for any
$a\in A$ and $m\in M$.

$2)$ Set $r_i=\max_{m\in M_i}\ \{\dim \ker m\}$. By the above,
$r_{i+1}\leq r_i$; hence, the sequence $r_i$ stabilizes:
$r_N=r_{N+1}=\dots :=r$ for some $N$. Due to transitivity, $r<\dim
A$, so it remains to prove that $r=\dim A-1$.

$3)$ We will say that $m$ is \textit{maximally degenerate}\index{Element, maximally degenerate} if $\dim
\ker m=r$.

By step $1)$, if $a\in A$ and $m\in M_n$ for $n\geq N+1$ is
maximally degenerate, then $am$ is maximally degenerate too and
$\ker am=\ker m$.

$4)$ Let $m_0\in M_{N+l}$ for $l\geq 0$ be a~maximally degenerate
element. For any $n\geq N$, set
\[
V_n(m_0):=\{m\in M_n\mid \ker m\supset
\ker m_0\}.
\]
It is clear that any non-zero element $m\in V_n(m_0)$ is
maximally degenerate and $\ker m=\ker m_0$.

$5)$ The following lemma is crucial in the next item. Actually, this
lemma is the key point of the whole proof of Theorem \ref{th3.2}.

\parbegin[A useful lemma]{Lemma}[A useful lemma]\label{SS:3.4} Let $V$ and $W$ be finite-dimensional vector spaces, and let ${U\subset \Hom(V, W)}$ be a~subspace such
that any non-zero homomorphism $u\in U$ is a~monomorphism. Then,
\begin{equation}
\label{lemma} \dim U\leq\dim W- \dim V+1. \end{equation}
\end{Lemma}

\begin{proof} (This elegant proof due to J.~Bernstein replaces an
original elementary but a~bit cumbersome proof.) Let $P(X)$ be the
projective space corresponding to a~vector space $X$. 

Since the mapping 
${U\otimes V\tto W}$ is a~monomorphism for any non-zero $u\in U$, there
exists a~mapping ${i: P(U)\times P(V)\tto P(W)}$ which for
each fixed $a\in P(U)$ and $b\in P(V)$ induces embeddings
${i_1: a\times P(V)\tto P(W)}$ and ${i_2: P(U)\times b\tto P(W)}$.

Let $h\in H^2(P(W))$, $h_{1}\in H^2(P(U))$ and $h_{2}\in
H^2(P(V))$ be the generators of the respective cohomology rings.
Since $i_1$ and $i_2$ are embeddings, it follows that $i^{*}(h)=h_{1}+h_{2}$.

On the other hand, the only relations $h$, $h_{1}$ and $h_{2}$ satisfy
are
\[
h^{\dim W+1}=0, \quad h_{1}^{\dim U+1}=0, \quad h_{2}^{\dim V+1}=0.
\]
Hence, replacing $h$ with $i^{*}(h)=h_{1}+h_{2}$ we see that the
first equation, \textit{i.e.}, $h^{\dim W+1}=0$, is a~new relation only if the inequality 
~\eqref{lemma} holds.
\end{proof}

$6)$ Each $m\in V_n(m_0)$ determines a~monomorphism $A/\ker m_0\tto
V_{n-1}(m_0)$. By Lemma \ref{SS:3.4} we have
\[
\renewcommand{\arraystretch}{1.4}
\begin{array}{l}
 \dim V_n(m_0)\leq \dim V_{n-1}(m_0) -(\dim A - r)+1\Longrightarrow\\
\dim V_{n-1}(m_0)\geq \dim V_n(m_0) +(\dim A - r)-1.\end{array}
\]

$7)$ Suppose $r<\dim A-1$. Then, by step $6)$, $\dim
V_{n-1}(m_0)\geq \dim V_n(m_0)+ 1$. Let us take $m_0\in M_{N+s}$,
where $s>0$ is arbitrary. Then, $V_{N+s}(m_0)\neq 0$; hence,
\[
\dim M_n\geq \dim V_{N}(m_0)\geq V_{N+1}(m_0)+1\geq\cdots\geq \dim
V_{N+s}(m_0)+s \geq 1+s.
\]
Since $\dim M_N<\infty$, we get a~contradiction. Thus, $r=\dim A-1$.

$8)$ So for some $ s\geq N$, we get an element $m_0\in M_s$ with
$\codim \ker m_0=1$. By transitivity there are $a_1,\ldots, a_s\in
A$ such that the element $m= a_1\cdots a_s m_0\in M_0$ is non-zero. Due to step $1)$,
$\ker m\supset \ker m_0$ and, by transitivity, $\ker m\neq A$. Then,
$\ker m =\ker m_0$ and $m$ is the sought-for element.
\end{proof}



\ssec{A superization of Wilson's theorem}\label{3.3}
Let $\cA=\mathop{\oplus}_{i\leq -1}\ \cA_i$ be the
supercommutative superalgebra generated by $A:=\cA_{-1}$. Let $\
M\mathop{\oplus}_{i\geq -1}\ M_i$ be an $\cA$-module. As above,
we will say that $M$ is \textit{transitive} if $\ker m\neq 0$ for
any non-zero $m\in M_i$, where $i\geq 0$.

\sssbegin[A superization of Wilson's theorem]{Theorem}[A superization of Wilson's theorem]\label{thWils}
In any infinite-dimensional transitive graded $\cA$-module $M$,
there exists an element $m\in M_0$ such that $\codim \ker m=1|0$.
\end{Theorem}

\begin{proof}
Let $\sdim A_{\bar 1}=0|n$.

$1)$ Due to supercommutativity of $\cA$, we have
\begin{equation}
\label{3.3.1}
a^2 m=0~~ \text{ for any $a\in A_{\bar 1}$ and $m\in M$}; 
\end{equation}
\begin{equation}
\label{3.3.2}
a_1\cdots a_i m=0~~ \text{ for any } m\in M \text{ and
} a_1, \ldots, a_i\in A_{\od} \text{ if }i>n. 
\end{equation}

$2)$ In view of condition ~\eqref{3.3.2} and the transitivity of $M$, we see that
$A_{\bar 0}m\neq 0$ for any $m\in M_i$ with $i\geq n$.

$3)$ Let $\cA_0\subset \cA$ be the subalgebra generated by
$A_{\ev}$. Let $\widetilde {M}=M/\oplus_{-1\leq i\leq n-1}\ M_i$. Let us
consider $\widetilde M$ as an $\cA_0$-module and apply
Theorem~\ref{th3.2} to $\widetilde M$. We get an element $m\in
M_{n+1}$ such that $\dim (\ker m\cap A_{\ev}) = \dim A_{\ev}-1$. If
$\ker m\cap A_{\od}\neq A_{\od}$, then there exists $a_1\in A_{\od}$ such
that $a_1\notin \ker m$. Then, $a_1m\neq 0$ and, thanks to condition ~\eqref{3.3.1}, we see that
$\ker a_1 m\supset\ker m\oplus \Cee a_1$.

Similarly, for $i\leq n$, there exist $a_1,\dots, a_i$ such that
\[
m_1=a_1\cdots a_i m\neq 0\text{~~ and $\ker m_1\supset \ker m +
A_{\od}$.}
\] 
Notice that $m_1\in M_s$, where $s=n+1-i>0$. Hence,
$\ker m_1\neq A$; thus, $\codim \ker m_1= 1|0$. Now we can complete
the proof in the same way as in step $8)$ of the
proof of Theorem~\ref{th3.2}. 
\end{proof}

\ssbegin[Existence of a~rank-1 operator]{Corollary}[Existence of a~rank-1 operator with even covector]\label{SS:3.7} Let
$\fg=\mathop{\oplus}_{i\geq -1}\ \fg_i$ be a~transitive
infinite-dimensional $\Zee$-graded Lie superalgebra with $\dim
\fg_i<\infty$ for all~ $i$. Then, there exists $v^*\otimes w\in
\fg_0$, where $v\in (\fg_{-1})_{\ev}$ and $w\in \fg_{-1}$.
\end{Corollary}

\section{Step 3. Main reduction: the three cases}\label{S:4}

\ssbegin[Corollaries of simplicity]{Lemma}[Corollaries of simplicity]\label{SS:4.1} If $\fg_*=\mathop{\oplus}_{-1\leq i\leq k}\ \fg_i$, where $k\geq 1$, is
simple, then

\emph{1)} $[\fg_{-1}, \fg_1]=\fg_0$;

\emph{2)} the representation of $\fg_0$ in $\fg_{-1}$ is
irreducible, \textit{i.e.}, $\fg_\geq=\mathop{\oplus}_{i\geq 0}\ \fg_i$
is a~maximal Lie subsuperalgebra;

\emph{3)} the representation of $\fg_0$ in $\fg_{-1}$ is faithful.
\end{Lemma}

\begin{proof} When one of the conditions
$1)-3)$ is violated, it is easy to find an ideal in $\fg_*$.

If 1) fails, take $ \fii=\fg_{-1}\oplus[\fg_{-1},
\fg_1]\oplus\left(\mathop{\oplus}_{i>0}\ \fg_i\right)$.

If 2) or 3) fails, take $\fii= \mathop{\oplus}_{n\geq
0}\ \fii_n$, where the $\fii_n$ are recursively defined (here $\rho$ is
the $\fg_0$-action on $\fg_{-1}$) by
\[
\renewcommand{\arraystretch}{1.4}\begin{array}{l}
\fii_{-1}=\begin{cases}\text{a $\fg_0$-invariant
subsuperspace},&\text{in case $2)$},\\
0,&\text{in case $3)$};\end{cases}\\
\fii_0=\begin{cases}\fg_0,&\text{in case $2)$},\\
\ker \rho,&\text{in case $3)$};
\end{cases}\\
\fii_n=\{g\in\fg_n\mid [\fg_{-1}, g]\subset \fii_{n-1}\}\quad \text{
for }\quad n\geq \begin{cases}0,&\text{in case
$2)$},\\
1,&\text{in case $3)$}. \end{cases}\end{array}
\]

The Jacobi identity implies the inclusion
\[
[\fg_{-1}, [\fg_k, \fii_n]]\subset[\fg_{k-1}, \fii_n]+[\fg_k, [\fg_{-1},
\fii_n]];
\]
the double induction on $k$ and $n$ shows that $\fii$ is an ideal in $\fg_*$.
\end{proof}

The conditions 2) and 3) of Lemma \ref{SS:4.1} allow us to think of $\fg_0$ as of
an irreducible linear Lie superalgebra, \textit{i.e.}, as a~subalgebra of
$\fgl(V)$, where $V=\fg_{-1}$, such that the tautological $\fgl(V)$-module $V$ is irreducible.

\subsubsection{Notation. The four cases}\label{ss4cases}
Consider the following four cases:

1) Let $V_1$ and $V_2$ be linear superspaces, $V=V_1\otimes V_2$.
The tensor product of the $\fgl (V_{i})$-actions on $V_{1}\otimes V_{2}$
gives us a~representation of $\fgl (V_{1})\oplus\fgl (V_{2})$:
\[
\rho (A, B)=A\otimes 1_{V_{2}}+1_{V_{1}}\otimes B \quad \text{ for
any }\; A\in \fgl (V_{1}), \ \ B\in \fgl (V_{2}).
\]
It is clear that $\rho$ has a~$1$-dimensional kernel
$\Cee(1_{V_{1}},\ -1_{V_{2}})$. Denote:\index{$\fg_1\bigodot\fg_2$}
\[
\fgl (V_{1})\bigodot\fgl
(V_{2}):=\rho(\fgl (V_{1})\oplus\fgl (V_{2})).
\]

2) The space of all operators in a~superspace $V$ can be endowed
with the following two structures: that of the Lie superalgebra $\fgl (V)$ and that of
the associative superalgebra $\Mat (V)$. We will carefully
distinguish them. (Recall that we endow with the superscript ${}^L$\index{$A^L$}\index{$A^{SL}$}
(resp., ${}^{SL}$) the Lie algebra (resp., Lie superalgebra)
constructed on the space of an associative algebra (resp., 
superalgebra) with the dot product replaced by the bracket (resp., 
superbracket). In particular, $\Mat (V)^L=\fgl(V)$. Various
\textit{super}structures on $V$ induce non-isomorphic, in general,
superstructures on both $\Mat(V)$ and the associated Lie superalgebra
$\fgl(V):=\Mat(V)^{SL}$, as is clear in the supermatrix format.)

Let $V_0$ be a~$(1, 1)$-dimensional superspace; define the two
associative subsuperalgebras:
\[
\begin{array}{l}
{\cA}\subset\Mat (1|1)=\Mat (V_0)=\Span\left\{\mat{a & b
\\ b & a}\mid a, b\in \Cee\right\},\\
C(\cA)\subset\Mat (1|1)=\Mat (V_0)=\Span\left\{\mat{ c & d
\\ - d& c}\mid c, d\in \Cee\right\}.
\end{array}
\]

Set $J:=\mat{0 & 1 \\
-1 & 0}$ and $\Pi:=\mat{0 & 1 \\
1 & 0}$; clearly, $1_2$ and $J$ generate ${\cA}$, while $1_2$ and
$\Pi$ generate $C({\cA})$. Note that
\[
J^2=-1, \quad \Pi^2=1; \quad J\Pi=\Pi J=0,
\] and therefore
${\cA}\simeq C({\cA})$ over $\Cee$ and $[{\cA}, C({\cA})]=0$; this
justifies our notation: $C({\cA})$ is the centralizer of $\cA$ in
$\Mat (1|1)$.

If $V_1$ and $V_2$ are linear superspaces of superdimension $m\eps$
and $n$, respectively, then the Lie superalgebras $(\Mat
(V_1)\otimes {\cA})^{SL}$ and $(C({\cA})\otimes\Mat (V_2))^{SL}$ of the
operators acting in $V_1\otimes V_0$ and $V_0\otimes V_2$,
respectively, can be identified with $\fq(m)$ and $\fq(n)$,
respectively.

Now, let $V:=V_{1}\otimes V_0\otimes V_{2}$. Then, a~natural
$\fq(n_{1})\oplus \fq (n_{2})$-action $\rho$ in $V$ is defined:
\begin{equation}
\label{qq}
\begin{split}
&\rho\left(\mat{ A & B \\
B & A}, \mat{ C & D \\ D & C }\right) :=(A\otimes 1_{V_{0}}+B\otimes
J)\otimes 1_{V_{2}}+
1_{V_{1}}\otimes(1_{V_{0}}\otimes C+\Pi\otimes D)\\
&\text{for any $\mat{ A & B \\ B &
A}\in\fq(m)$ and $\mat{ C
& D \\
D & C}\in\fq(n)$.}
\end{split}
\end{equation}
This representation has a~$1$-dimensional kernel.
Denote
\begin{equation}
\label{1Dker} 
 \fq (n_{1}) \bigodot \fq (n_{2}):=\rho(\fq
(n_{1})\oplus \fq (n_{2})).
\end{equation}

3) Let $V_{1}$ be a~linear superspace of superdimension $(r,s)$, let
$\fn= \fgl (V_{1})\otimes \Lambda (n) $, $\fk=\fvect (0|n)$ and let 
$\fg=\fn\ltimes \fk$ with the natural action of $\fk$ on the ideal
$\fn$. The Lie superalgebra $\fg$ has a~natural faithful
representation $\rho $ in the linear superspace $V=V_{1}\otimes
\Lambda (n) $ defined for any $A\otimes \varphi \in \fn$, $D\in
\fk$, and $v\otimes \psi\in V$ by the formulas
\[
\begin{array}{l}
\rho (A\otimes \varphi )(v\otimes
\psi)=(-1)^{p(\varphi)p(v)}Av\otimes\varphi\psi,\\
\rho (D)(v\otimes \psi )=-(-1)^{p(D)p(v)}v\otimes D(\psi).
\end{array}
\]
In what follows, we will always identify the elements of $\fg$ with
their images under $\rho$. Therefore, we will consider $\fg$ as embedded
in $\fgl (V)$.

4) Let $\fhei (0|2n)=\Span(\xi _{1}, \dots , \xi _{n}; \eta _{1},
\dots , \eta_{n}; z)$ be the Heisenberg superalgebra, with the
non-zero relations $[\xi _{i}, \eta _{j}]=\delta _{ij}z$. The
nontrivial irreducible finite-dimensional representation of $\fhei
(0|2n) $ that maps $z $ to the operator of multiplication by $1$ is
realized in $\Lambda (\xi) $ via \[\xi _{i}\longmapsto \xi
_{i},\quad \eta _{i}\longmapsto
\partial _{i}.\] The normalizer of $\fhei (0|2n) $ in $\fgl (\Lambda
(n)) $ is $\fg =\fhei (0|2n)\ltimes \fo (2n)$ or, in terms of
differential operators:
\[
\fhei (0|2n)\ltimes \fo (2n)=\Span(1, \xi_{i}, \partial_{i}, \xi
_{i}\partial_{j}, \partial_{i}\partial_{j}, \xi_{i}\xi _{j})_{i, j=
1}^n.
\]

\ssbegin[Irreducible linear Lie superalgebras]{Theorem}[Irreducible linear Lie superalgebras]\label{thShM}
Let $\fg\subset\fgl(V)$ be an irreducible linear Lie superalgebra
that is neither almost simple nor a~central extension of an almost simple Lie superalgebra. Then, $\fg$ is contained in one of the following four Lie
subsuperalgebras of $\fgl(V)$:

$1)$ $\fg\subset\fgl(V_1)\bigodot\fgl(V_2)\subset\fgl(V)$ for any
$V_1$, $V_2$ such that $V=V_1\otimes V_2$ and $\dim V_i>1$ for $i=1,
2$;

$2)$ $\fg\subset\fq(m)\bigodot\fq(n)\subset\fgl(V)$, where
$V=V_1\otimes V_0\otimes V_2$ for any $V_0$, $V_1$, $V_2$ such
that 
\[\text{$\sdim V_0=(1, 1)$, \ \ $\sdim V_1=m\eps $, \ \ $\sdim V_2= n$ and
$mn>1$;}
\]

$3)$ $\fg\subset\fgl(V_1)\otimes \Lambda (n)\ltimes
\fvect(0|n)\subset\fgl(V)$ for any $V_1$ and $n$ such that
$V=V_1\otimes \Lambda (n)$; the natural representation of $\fg$ in
$V$ is induced;

$4)$ $\fg\subset\fhei(0|2n)\ltimes \fo(2n)\subset\fgl(V)$ for $n>1$,
where $V\simeq \Lambda (n)$ or $\Pi(\Lambda (n))$.
\end{Theorem}

For proof, see Subsection~ \ref{SS:4.3.1p}.

\sssbegin[The 3 cases]{Corollary}[The 3 cases]\label{SS:4.3} If $\fg\subset\fgl(V)$ is an irreducible
linear superalgebra with a~rank-$1$ operator in $V$, then only one
of the following $3$ cases is possible:

$1)$ the tautological representation of $\fg$ in $V$ is induced;

$2)$ $\fg$ is almost simple or a~central extension of an almost simple
Lie superalgebra;

$3)$ $\sdim V=(2, 2)$ and $\fg\subset \fhei (0|4)\ltimes \fo (4)$.
\end{Corollary}

Indeed, none of the following matrix Lie superalgebras
\[
 \fgl(V_1)\bigodot \fgl(V_2)\;\text{ with $\dim V_i>1$ \qquad or
\qquad $\fq(n_1) \bigodot \fq(n_2)$ with $n_1 n_2>1$}
\]
have operators of rank 1; the only possibility for the Lie
superalgebra $\fg=\fhei(0|2n)\ltimes \fo(2n)$ to have a~
rank-1 operator occurs for $n= 2$.

\subsection{Proof of Theorem \ref{thShM} (\cite{ShM})}
\label{SS:4.3.1p} In this Subsection, $\fg\subset\fgl(V)$ is an irreducible linear Lie
superalgebra, and $\rho$ is its tautological representation in superspace
$V$. We will assume that $\dim V>1$, because the case where $\dim V=1$ is
trivial. By hypothesis, $\fg$ contains an ideal $\fk$.

Our proof of Theorem~\ref{thShM} is based on the following result
--- a~ superization of a~statement well-known from
Dixmier's book \cite[Subsection~5.4.1]{Di}.

\sssbegin[The 3 cases of irreducible linear Lie superalgebras: A, B, C]{Theorem}[The 3 cases of irreducible linear Lie superalgebras: A, B, C]\label{th3.3.1}
Let $\fg\subset\fgl(V)$ be an irreducible linear Lie superalgebra, $\rho$ its
tautological representation (in particular, $\rho$ is faithful). Let
$\fk $ be an ideal of $\fg $ and $\dim \fk>1$. Then, only one of the
following $3$ cases is possible:

\emph{A)} $\rho |_\fk $ is a~multiple of an irreducible $\fk
$-module $\tau $ and the multiplicity of $\tau $ is $>1$;

\emph{B)} there exists a~subalgebra $\fh\subset\fg$ such that
$\fk\subset \fh$ and $\rho\simeq \ind^{\fg }_{\fh }\sigma $ for an
irreducible $\fh$-module~ $\sigma $;

\emph{C)} $\rho |_\fk $ is irreducible.
\end{Theorem}

\begin{proof} The argument, see \cite{ShM}, follows the lines of \cite[Subsection~5.4.1]{Di} with one novel case:
irreducible modules of $Q$-type might occur. Fortunately, to
treat them is quite straightforward and we would rather save paper
by omitting the verification. \end{proof} 

We will say that the \textit{representation $\rho$ is of type A,
B, or C with respect to the ideal} $\fk$ if the corresponding
case of Theorem~\ref{th3.3.1} holds. Lemmas\index{Representation of type A,
B, or C with respect to the ideal}
\ref{L3.3.2}, \ref{L3.3.3}, and \ref{L3.3.5} deal with types A, B, and C,
respectively.

\parbegin[Case A]{Lemma}[Case A]\label{L3.3.2}
 If $\rho$ is of type \textit{A} with respect to the
ideal $\fk$, then for some $V_1$ and $V_2$, either $\fg\subset
\fgl(V_1)\bigodot\fgl(V_2)$ or $\fg\subset
\fq(V_1)\bigodot\fq(V_2)$.
\end{Lemma}

\begin{proof} Let $V$ be the superspace on which $\tau$ acts.
Then, $V=V_1\otimes U$ for some $U$ and
\[
\rho(h)(v\otimes u)=\tau(h)v\otimes u\; \text{ for any }\; h\in\fk, \; v\in
V_1, \; u\in U, \text{~~\textit{i.e.},}
\]
\begin{equation}
\label{3.3.2.1}
 \text{the operators $\rho (h)$ are of the form
$\tau(h)\otimes 1_U$.}
\end{equation}
On the other hand, $\fg$ is the stabilizer of $\tau$. So, for any
$g\in\fg$, there exists $s(g)\in \End V_1$ such that $\tau ([g,
h])=[s(g), \tau(h)]$ for any $h\in\fk$.

Let us consider the operator $T(g)=\rho(g)-s(g)\otimes 1_U$ in $V$.
For $h\in\fk$, we have
\[
[T(g), \rho(h)]=[\rho(g)-s(g)\otimes 1_U, \ \rho(h)]=\rho([g, h])-[s(g), \tau(h)]
\otimes 1_U=0\, ,
\]
\textit{i.e.}, $T(g)$ commutes with all the operators of the form $\rho(h)$, where
$h\in \fh$, and, since $\tau$ is irreducible, $T(g)$ is of the form
\[
T(g)=\begin{cases}
1_{V_1}\otimes A(g),& \text{ if }\; \tau\text{ is of $G$-type},\\
1_{V_1}\otimes A(g)+J\otimes B(g),& \text{ if $\tau$ is of $Q$-type},\ 
\end{cases}
\]
where $J$ is an odd generator of $C(\tau)$, the centralizer of $\tau$.

Therefore, $\rho(g)$ is of the form
\begin{equation}
\label{3.3.2.2} \rho(g)=\begin{cases}s(g)\otimes 1_U+1_{V_1}\otimes A(g),&
\text{if $\tau$ is of $G$-type},\\
s(g)\otimes 1_U+1_{V_1}\otimes A(g)+J\otimes B(g),&\text{if $\tau$
is of $Q$-type}.\end{cases}
\end{equation}
Formulas ~\eqref{3.3.2.1}--\eqref{3.3.2.2} show that $\fg$ is
contained in
\[
\begin{cases}\fgl(V_1)\bigodot\fgl(V_2)\subset\fgl(V), \text{ where $V_2=U$,}&\text{if $\tau$ is
of $G$-type},\\
\fq(V_1)\bigodot\fq(V_2)\subset\fgl(V), \text{ where $V_2=\Span(1,
J)\otimes U$,}&\text{if $\tau$ is of $Q$-type}. \qedhere
\end{cases}
\]
\noqed \end{proof}

\parbegin[Case B]{Lemma}[Case B]\label{L3.3.3} 
If $\rho$ is of
type \textit{B} with respect to the ideal $\fk$, then
\[
\fg\subset \fgl(V)\otimes \Lambda (n)\ltimes \fvect(0|n)\;\text{ for
some $V$ and $n$.}
\]
\end{Lemma}

\begin{proof}[Proof \nopoint] of this lemma follows from the definition of the induced
representation. \end{proof}

\subsubsection{Completion of the proof of Theorem \ref{thShM}}
\label{SS:4.3.5} It remains to consider the case where $\rho$ has the following
property: for any ideal $\fk\subset \fg$, either $\rho$ is of type
\textit{C} with respect to $\fk$ or $\rho\vert _\fk$ is the multiple
of a~character.

Let $\fr$ be the radical of the linear Lie superalgebra $\fg$.

\parbegin[On the radical $\fr$ in cases A and C]{Lemma}[On the radical $\fr$ in cases A and C]\label{L3.3.4}
If the radical $\fr$ is commutative and $\rho$ is of type A or C with respect to
$\fr$, then ${\dim\fr=0}$ or $1$.
\end{Lemma}

\begin{proof}[Proof \nopoint] follows from the facts that any irreducible module over
a~commutative Lie superalgebra is 1-dimensional (even or odd) and $\rho$ is faithful. \end{proof} 

\subsubsection{When the radical $\fr$ is not commutative}\label{SS:3.3.5}
In this case, consider the derived series of $\fr$:
\[
\fr\supset\fr_1\supset\dots\supset\fr_k\supset\fr_{k+1}=0,\;\text{
where $\fr_{i+1}=[\fr_i, \fr_i]$. }
\]
Clearly, each $\fr_i$ is an ideal in $\fg$ and the last non-zero
ideal, $\fr_k$, is commutative.

\parbegin[Case C]{Lemma}[Case C]\label{L3.3.5} 
If $\rho$ is of type
C with respect to
$\fr_{k-1}$ and $\rho| _{\fr_{k}}$ is scalar, then

either $\fr_{k-1}\simeq \fhei(0|2n)$ or $\fr_{k-1}\simeq \fhei(0|2n-1)$,
and

either $V\simeq \Lambda(n)$ or $V\simeq \Pi(\Lambda(n))$.

\noindent In these cases, $\fg\subset \fhei(0|2n)\ltimes \fo(2n)$.
\end{Lemma}

\begin{proof} We see that if the restrictions $\rho|_{\fr_{i}}$ are
irreducible for all~ $i$, then

1) $\dim\fr_k=1$ because $\rho$ is faithful;

2) $\fr_{k-1}=\fhei (0|m)$ for some $m$;

3) $\rho|_{\fr_{k-1}}$ is irreducible and faithful, so it can be
realized in the superspace of functions, $\Lambda (n)$, or in
$\Pi(\Lambda (n))$, where $n=\lfloor\nfrac m2\rfloor+1$;

4) $\fg\subset\fhei (0|2n)\ltimes \fo(2n)$, \textit{i.e.}, $\fg$ is contained
in the normalizer of $\fhei (0|2n)$ in the spinor representation of $\fo(n)$.
\end{proof}

\sssbegin[Ideal of type A or B]{Lemma}[Ideal of type A or B]\label{L3.3.6}
Let $\dim\fr\leq 1$, \textit{i.e.}, $\fg$ is either semi-simple or a~central
extension of a~semi-simple Lie superalgebra, but NOT almost simple
or a~central extension of an almost simple Lie superalgebra. Then, there exists an ideal $\fk$ such that $\rho$ is of type A or B
with respect to $\fk$.
\end{Lemma}

\begin{proof} Let $\dim\fr=0$, \textit{i.e.}, let $\fg$ be semi-simple. As
$\fg$ is not almost simple, then, see Subsection~\ref{SS:2.2.1}, the following alternative
arises:
\begin{equation}
\label{3.3.6}
\begin{array}{l}
\text{either $\fg$ contains an ideal $\fk$ of the form
$\fk=\fh\otimes \Lambda (n)$ with simple $\fh$ and
$n>0$},\\ 
\text{or $\fg=\mathop{\oplus}_{j\leq k}\ \fs_j$, where each
$\fs_j$ is almost simple and $k>1$.}\end{array}
\end{equation}

\sssbegin[Irreducible representation of type B]{Lemma}[Irreducible representation of type B]\label{L3.3.7}
If $\fg=\mathop{\oplus}_{j\leq k}\ \fs_j$ and $k\geq 2$, then
any irreducible faithful representation of $\fg$ is of type B with
respect to any its ideal $\fs_j$.
\end{Lemma}

\begin{proof} Since the stabilizer of any irreducible representation
of $\fs_j$ is the whole $\fg$, the type of any irreducible
representation of $\fg$ with respect to $\fs_j$ can be either B or
C. Due to faithfulness of the representation, type C is ruled
out.
\end{proof}

\sssbegin[No faithful irreducible representations]{Lemma}[No faithful irreducible representations]\label{L3.3.6.1}
Let $\fh$ be a~simple Lie superalgebra and let 
$\fk=\fh\otimes\Lambda(n)$. Then, $\fk$ has no faithful irreducible
finite-dimensional representations.
\end{Lemma}

For proof, see Subsection~\ref{Pf337a}; let us pass to corollaries. 

\sssbegin[No type A irreducible representations: conditions]{Corollary}[No type A irreducible representations: conditions]\label{cor3363}
If $\fg$ contains an ideal $\fk$ of the form described in eq.~\eqref{3.3.6}, then
$\fg$ cannot have any faithful irreducible finite-dimensional
representation of type A with respect to the ideal $\fk$.
\end{Corollary}

\sssbegin[Lemma
$\ref{L3.3.6}$ holds for semi-simple Lie superalgebras]{Corollary}[Lemma
$\ref{L3.3.6}$ holds for semi-simple Lie superalgebras]\label{cor3364}
Lemma $\ref{L3.3.6.1}$ and Corollary $\ref{cor3363}$ prove Lemma
$\ref{L3.3.6}$ for semi-simple Lie superalgebras.
\end{Corollary}

\sssbegin[A useful fact]{Lemma}[A useful fact]\label{l3365}
For any simple Lie
superalgebra $\fh$, we have
\[
[\fh_{\od}, \fh_{\od}]=\fh_{\ev}\quad \textit{and}\quad [\fh_{\ev},
\fh_{\od}]=\fh_{\od}\, .
\]
\end{Lemma}

\begin{proof}[Proof \nopoint] is obvious and so is left to the reader. (Or see \cite{K3}.)
\end{proof}

\paragraph{Proof of Lemma \ref{L3.3.6.1}}\label{Pf337a}
For $n> 0$, the Lie superalgebra $\fk$ contains a~supercommutative
ideal
\[
\fn_1=\fh\otimes\xi_1\cdots\xi_n.
\]
Moreover, if $n>2$, then $\fg$ contains a~supercommutative ideal
\[
\fn_2=\fh\otimes\Lambda^{n-1}(\xi)\oplus \fn_1 .
\]
Thanks to Theorem \ref{thShM}, any irreducible representation $\rho$
of $\fk$ in the superspace $V$ is obtained as follows: take an irreducible representation
(character $\chi$) of an ideal $\fn_i$, the stabilizer $\fst
_{\fk}(\chi)$, and an irreducible representation $\widetilde {\rho}$
of $\fst_{\fk} (\chi)$ whose restriction to $\fn_i$ is a~multiple
of $\chi$. Then,
\[
\rho=\ind^{\fk}_{\fst_{\fk}(\chi)} (\widetilde {\rho})\;\text{~~~ for
some $\widetilde \rho$.}
\]

Let us show that if $\dim V<\infty$, then
$\chi|_{\fh\otimes\xi_1\cdots \xi_n}=0$. Since
$\fh\otimes\xi_1\cdots \xi_n$ is an ideal in $\fk$, this would
contradict faithfulness of $\rho$. Indeed, $\dim V<\infty$ if
and only if $\fst_{\fk}(\chi)\supset\fk_{ \ev }$.

1) $n=2k+1$. For any $f\in \fh$ and $g=h\otimes 1\in
\fh_\ev\otimes 1$, we have
\[
0=\chi ([g, f\otimes\xi_1\cdots
\xi_n])=\chi ([h, f]\otimes\xi_1\cdots \xi_n),
\]
implying that
\[
\chi |_{[\fh_{\ev}, \fh_{\od}]\otimes\xi_1\cdots
\xi_n}\overset{\text{ (Lemma \ref{l3365})} }{=}
\chi|_{(\fh_{\od}\otimes\xi_1\cdots \xi_n)}=0.
\]
Since $\chi$ is a~character, it follows that
$\chi|_{\fh_{\ev}\otimes\xi_1\cdots \xi_n}=0$, hence
$\chi|_{\fh\otimes\xi_1\cdots \xi_n}=0$.

2) $n=2k>2$. For any pair
\[
\begin{array}{l}
g=h\otimes l_1(\xi)\in \fh_{\od}\otimes
\Lambda^1(\xi)\subset\fk _{\ev},\hspace{1.5cm}\\
g_1=h_1\otimes l_{n-1}(\xi)\in \fh_{\od}\otimes\Lambda^{n-1}(\xi)\subset
(\fn_2)_{\ev},
\end{array}
\]
we have
\[
0=\chi([g, g_1])=\pm \chi([h, h_1])\otimes l_1(\xi)l_{n-1}(\xi),
\]
\textit{i.e.},
\[
\chi|_{[\fh_{\od}, \fh_{\od}]\otimes\xi_1\cdots \xi_n}=\chi |
_{\fh_{\ev}\otimes \xi_1\cdots \xi_n}=0\ .
\]
But $\chi$ is a~character; therefore,
\[
\chi|_{\fh_{\od}\otimes\xi_1\cdots \xi_n}=0, \text{~~ and,
finally,~~} \chi|_{\fh\otimes\xi_1\cdots \xi_n=0} .
\]

3) $n=2$. Then, $\fk$ contains a~commutative ideal
$\fn=\fh\otimes\xi_1\xi_2$; let $\chi$ be a~character of $\fn$.

If $\chi$ generates a~finite-dimensional representation of $\fk$,
then $\fst_{\fk}(\chi)\supset\fk_{\ev}$, \textit{i.e.},
\begin{equation}\label{missed1}
\chi([f, g\otimes\xi_1\xi_2])=0~~ \text{~~~~for any~~} f\in \fh_{\ev},\;
g\in \fh .
\end{equation}

If $g\in \fh_{\od}$, then eq.~\eqref{missed1} holds automatically by
the parity considerations. Let $g\in \fh_{\ev}$. Then, nevertheless,
\[
\chi([f, g\otimes\xi_1\xi_2])=\chi([f, g]\otimes\xi_1\xi_2)=0,
\]
and therefore $\chi |_{[\fh_{\ev},
\fh_{\ev}]\otimes\xi_1\xi_2}=0$.

Set
\[
\fm:=\fh\otimes (\Lambda^1(\xi)\oplus\Lambda^2(\xi)).
\]

Clearly, $\fm$ is an ideal of $\fk$, and $\fst_{\fn}(\chi)=\fm$.
Therefore, the irreducible representation $\rho$ of $\fm$ given by
$\chi$ is such that $\rho|_{\fn}$ is a~multiple of $\chi$. Let
\[
\fn_1:=\fh_{\ev}\otimes\Lambda^1(\xi)\oplus([\fh_{\ev},
\fh_{\ev}]\oplus \fh_{\od})\otimes\Lambda^2(\xi).
\]
Clearly, $\fn_1$ is an ideal in $\fm$, and
\[
[\fn_1, \fm]\subset([\fh_{\ev}, \fh_{\ev}]\oplus
\fh_{\od})\otimes\Lambda^2 (\xi)\subset \ker\; \chi.
\]

Thus, $\fn_1\subset \ker\;\rho$ since the subspace
$V^{\rho(\fn_1)}=\{v\in V\mid \rho(\fn_1) v=0\}$ is
$\rho(\fn_1)$-invariant and non-zero. But $\fm/\fn_1$ is a~solvable, and by the super analog of \textit{Lie's theorem} (for its
formulation and proof over $\Cee$, see \cite{S0}) any of its irreducible finite-dimensional
representations is of the same form as $\chi$. Let it be $\chi$, for
definiteness sake.

This means that $\chi|_{[\fh_{\od},
\fh_{\od}]\otimes\Lambda^2(\xi)}=0$, \textit{i.e.},
$\chi|_{\fh\otimes\Lambda^2(\xi)}=0$, and since $\fh\otimes\xi_1
\xi_2$ is an ideal in $\fk$, the representation of $\fk$ given by
$\chi$ cannot be faithful. Lemma is proved. \qed

\ssec{Central extensions of semi-simple Lie
superalgebras}\label{SS:4.4} Let $\fg$ be a~Lie superalgebra, and let $\fr$ be 
its radical of dimension 1. Then, 
$\widetilde \fg=\fg/\fr$ is semi-simple by definition. 

First, consider the case where $\widetilde \fg=
\mathop{\oplus}_{i\leq k}\ \fs_i$ with each $\fs_i$ being almost
simple and $k>1$. Let $\pi:\fg\tto\widetilde \fg$ be the natural
projection. The Lie superalgebra $\fk:=\pi^{-1}(\fs_1)$ is an ideal
in $\fg$ and $\dim \fk>1$.

\sssbegin[No type C irreducible representations]{Lemma}[No type C irreducible representations]\label{L3371}
If $k>1$, then $\rho$ cannot be irreducible of type C with respect
to $\fk=\pi^{-1}(\fs_1)$.
\end{Lemma}

\begin{proof} Set $\fg_+=\mathop{\oplus}_{i\geq
1}\ (\fs_i)_{\ev}$. Clearly, $[\fs_1, \fg_+ ]=0$. Since $\rho(\fr)$
acts by scalar operators and $\rho$ is a~finite-dimensional
representation, it follows that $[\pi^{-1}(\fg_+), \fk]=0$, \textit{i.e.},
$\rho(\pi^{-1}(\fg_+))$ consists of intertwining operators of
$\rho|_{\fn}$, \textit{i.e.}, $\rho(\pi^{-1} (\fg_+))=\Cee\cdot 1$ and, since
$\rho$ is faithful, $\pi^{-1} (\fg_+)=\fr$ implying $\fg_+=0$.
\end{proof}

\sssec{The case where $\tilde \fg$ contains an ideal $\fk$ of
the form described in eq.~\eqref{3.3.6}}\label{SS:4.5} The central extension is
defined by a~cocycle $c\colon \widetilde \fg\times \widetilde
\fg\longrightarrow\Cee$. The cocycle condition reads
\[
c(f, [g, h])=c([f, g], h)+(-1)^{p(f)p(g)}c(g, [f, h]) ~~\text{~~for
any~~} f, g, h\in\widetilde \fg.
\]
As earlier, we assume that $\fg$ has a~faithful finite-dimensional
representation; so the restriction of $c$ to $\widetilde \fg_{\bar
0}\times \widetilde \fg_{\bar 0}$ is zero. Moreover, $c|_{\widetilde
\fg_{\bar 0}\times \widetilde \fg_{\bar 1}}=0$, by parity
considerations. Therefore,
\begin{equation}
\label{3.5} c|_{\widetilde \fg_{\bar
0}\times \widetilde \fg}=0.
\end{equation}

\parbegin[Technical]{Lemma}[Technical]\label{L4.16}
In the notation of Lemma~$\ref{L3.3.6.1}$, set $\fL^m:=\fh\otimes\Lambda^m(\xi_1, \dots , \xi_n)$. Then,
${c|_{\fL^{n}\times\fL^{n}}=0}$.
\end{Lemma}

\begin{proof} If $n=2k+1$, the condition ~\eqref{3.5}
means that
\[
c|_{\fh_{\bar 1}\otimes\Lambda^n(\xi)\times\fh_{\bar
1}\otimes\Lambda^n(\xi)}=0.
\]
Let $g_{\bar 1}, h_{\bar 1}\in\fh_{\bar 1}$ and $h_{\bar 0}\in\fh_{\bar
0}$. Then,
\[
\renewcommand{\arraystretch}{1.4}
\begin{array}{l}
c(f_{\bar0}\otimes\xi_1\cdots\xi_n, [g_{\bar 1}, h_{\bar
1}]\otimes\xi_1\cdots\xi_n)= 
-c(f_{\bar
0}\otimes\xi_1\cdots\xi_n, [g_{\bar
1}\otimes\xi_1\cdots\xi_n, h_{\bar
1}])=\\
-c([f_{\bar 0}\otimes\xi_1\cdots\xi_n, g_{\bar
1}\otimes\xi_1\cdots\xi_n], h_{\bar 1})- 
c(g_{\bar 1}\otimes\xi_1\cdots\xi_n, [f_{\bar
0}\otimes\xi_1\cdots\xi_n, h_{\bar 1}])=0.\end{array}
\]
Since $\fh_{\bar 0}=[\fh_{\bar 1}, \fh_{\bar 1}]$ (by Lemma
\ref{l3365}), it follows that $c|_{\fL^{n}\times\fL^{n}}=0$.

If $n=2k$, the condition ~\eqref{3.5} means that
\[
c|_{\fh_{\bar 0}\otimes\Lambda^n(\xi)\times\fh_{\bar
0}\otimes\Lambda^n(\xi)}=0.
\]
Let $f_{\bar 1}, h_{\bar 1}\in\fh_{\bar 1}$ and $g_{\bar 0}\in\fh_{\bar
0}$. Then,
\[
\begin{array}{l}
c(f_{\bar 1}\otimes\xi_1\cdots\xi_n, [g_{\bar 0}, h_{\bar
1}]\otimes\xi_1\cdots\xi_n)= 
-c(f_{\bar
1}\otimes\xi_1\cdots\xi_n, [g_{\bar
0}\otimes\xi_1\cdots\xi_n, h_{\bar
1}])=\\
-c([f_{\bar 1}\otimes\xi_1\cdots\xi_n, g_{\bar
0}\otimes\xi_1\cdots\xi_n], h_{\bar 1})- 
c(g_{\bar 0}\otimes\xi_1\cdots\xi_n, [f_{\bar
1}\otimes\xi_1\cdots\xi_n, h_{\bar 1}])=0.\end{array}
\]
Since $\fh_{\bar 1}=[\fh_{\bar 0}, \fh_{\bar 1}]$ (by Lemma
\ref{l3365}), it follows that $c|_{\fL^{n}\times\fL^{n}}=0$.
\end{proof}

\parbegin[Technical]{Lemma}[Technical]\label{L4.17}
For any $m=0, 1, \dots , n$
and $k=n-m$, $n-m+1, \dots, n$, we have $c|_{\fL^{k}\times\fL^{m}}=0$.
\end{Lemma}

\begin{proof} We will argue by reverse double induction on $k$ and
$m$. For $k=m=n$, Lem\-ma~\ref{L4.16} gives the statement.

Suppose the statement is true for all~ $k> k_0$ and $k_0+m\geq n$. Let us
show that it is true for $k=k_0$ as well. Observe that, see Subsection \ref{SS:2.2.1}, $\widetilde \fg$ contains $n$ elements $\eta _i$ such
that
\[
\ad_{\eta _i}|_{\fh}=\partial_{\xi_{i}}+D_i, \;\;\text{where
$D_i\in\fvect(0|n)$ and $\deg D_i>0$}
\]
in the standard grading
of $\fvect(0|n)$.

Let $\varphi\in \Lambda^{k_{0}+1}(\xi_1, \dots, \xi_n)$, $\psi\in
\Lambda^{m}(\xi_1, \dots, \xi_n)$; let $g, h\in\fh$ be such that
$p(g\otimes \psi)=\bar 0$ and $p(h\otimes \psi)=\bar 1$. Then, we
have
\[
\renewcommand{\arraystretch}{1.4}
\begin{gathered}
c([\eta_i, g\otimes\varphi],
h\otimes\psi)=(-1)^{p(g)}c(g\otimes\eta_i(\varphi), \ h\otimes\psi)=\\
\phantom{xxxxxxxx}\hskip 2cm
(-1)^{p(g)}\left(c(g\otimes\nfrac{\partial\varphi}{\partial\xi_{i}},\ 
h\otimes\psi)+c(g\otimes D_{i}\varphi,\ 
h\otimes\psi)\right).\end{gathered}
\]
The term $D_i\varphi\in\left (\mathop{\oplus}_{s\geq
k_{0}+1}\ \Lambda^s(\xi)\right )$ from the last summand vanishes by
the induction hypothesis.

On the other hand,
\begin{equation}
\label{coc}
 c([\eta_i, g\otimes\varphi], \ h\otimes\psi)=c(\eta_i,\ 
[g\otimes\varphi, h\otimes\psi])+c([\eta_i, \ h\otimes\psi],
g\otimes\varphi). 
\end{equation} 
Since
$\deg\varphi+\deg\psi=k_0+1+m>n$, we see that $[g\otimes\varphi,
h\otimes\psi]=0$. The second summand in the right-hand side of
~\eqref{coc} vanishes thanks to condition ~\eqref{3.5}. So
\[
c(g\otimes\nfrac{\partial\varphi}{\partial\xi_{i}}, \ h\otimes\psi)=0~~\text{~~
for any
$g\otimes\nfrac{\partial\varphi}{\partial\xi_{i}}\in\fL_{\bar
1}^{k_{0}}$ and $h\otimes \psi \in\fL_{\bar 1}^{m}$.} \hskip3cm \qed
\]
\noqed\end{proof}

\parbegin[Technical]{Lemma}[Technical]\label{L4.18}
The ideal $\fk$ has no 
faithful irreducible finite-dimensional representations.
\end{Lemma}

\begin{proof}[Proof \nopoint] reproduces the proof of Lemma~\ref{L3.3.6.1} word-for-word 
with the help of Lemmas~\ref{L4.16} and \ref{L4.17}. \end{proof}

\sssbegin[No type C representations]{Corollary}[No type C representations] The representation $\rho$ cannot be of type C with respect to the
ideal $\fk$.
\end{Corollary}

Lemmas A, B, C (\textit{i.e.}, \ref{L3.3.2}, \ref{L3.3.3}, \ref{L3.3.5}) combined with Lemmas \ref{L3.3.4} and \ref{L3.3.6} prove
Theorem~\ref{thShM}. \end{proof} 

\section{Step 4. Prolongation of the induced
representations}\label{Sstep4}

In this step we describe all the simple infinite-dimensional
$\Zee$-graded Lie superalgebras of depth 1 for which the irreducible
representation of $\fg_0$ on $\fg_{-1}$ is induced from
a~representation of a~subsuperalgebra $\fh\subset\fg_0$.

But, first of all, we describe certain universal Lie superalgebras.

\ssec{Universal property of $\fvect$, $\fk$ and $\fm$}\label{SS:5.1}\index{Property, universal}
Let $\fg$ be a~$\Zee$-graded transitive Lie superalgebra, and let $\fg_{-1}$ be 
an irreducible $\fg_0$-module. The following universal properties are immediate consequences of
the definition of Cartan prolongations.

1) If $\fg= \mathop{\oplus}_{-1\leq i}^{\infty}\ \fg_i$ and
$\sdim\fg_{-1}=(m|n)$, then $\fg\subset \fvect(m|n)$.

2) If $\fg= \mathop{\oplus}_{-2\leq i}^{\infty}\ \fg_i$ and
$\fg_{<0}=\fhei(2m|n)$, then $\fg\subset \fk(2m+1|n)$.

3) If $\fg= \mathop{\oplus}_{-2\leq i}^{\infty}\ \fg_i$ and
$\fg_{<0}=\fba(n)$, then $\fg\subset \fm(n)$.

\noindent The following are all the Lie superalgebras with ``big"
nonpositive part:

$1^\circ)$ Let $\fg= \mathop{\oplus}_{-1\leq i}^{\infty}\ \fg_i$ and
$\sdim\fg_{-1}=(m|n)$. Moreover, let $\fg_0\supset\fsl(m|n)$ and
$\dim \fg=\infty$. If $m>1$, then, $\fg = \fvect(m|n)$, or
$\fsvect(m|n)$, or $\fd(\fsvect(m|n)):=\fsvect(m|n)\ltimes \Cee E$, where $E:=\sum_{1\leq i\leq m+n}\ x_i\partial_{x_i}$.

If $m=1$, two more possibilities appear:
$\fg=\fsvect'(1|n)$ or $\fd(\fsvect'(1|n)):=\fsvect'(1|n)\ltimes \Cee E$.

$2^\circ)$ Let $\fg= \mathop{\oplus}_{-2\leq i}^{\infty}\ \fg_i
\subset \fk(2m+1|n)$ be such that $\fg_{\leq 0}=\fk(2m+1|n)_{\leq
0}$ and $\dim \fg=\infty$. Let either $m\neq 0$ or $m=0$ but $n\neq
6$; then, $\fg =\fk(2m+1|n)$ or $\fg=\fd(\fpo(2m|n)):=\fpo(2m|n)\ltimes \Cee E$, where $E:=\sum_{1\leq i\leq 2m+n}\ x_i\partial_{x_i}
$, where the $x_i$ are all indeterminates except $t$. If $m=0$
and $n=6$, then $\fg=\fk\fas$.

$3^\circ)$ If $\fg= \mathop{\oplus}_{-2\leq i}^{\infty}\ \fg_i
\subset \fm(n)$ is such that $\fg_{<0}=\fba(n)$, $\dim \fg=\infty$
and $\fg$ is simple, then $\fg =\fm(n)$, or $\fg=\fb_{\infty }(n)$, or
$\fg=\fb_{\lambda}(n)$ for $\lambda\neq 0$.

\ssbegin[Technical]{Lemma}[Technical]\label{L4.2}
Let $\fg= \mathop{\oplus}_{-1\leq i}^{\infty}\ \fg_i$ be
a~$\Zee$-graded transitive Lie superalgebra, and let ${W\subset \fg_{-1}}$ be 
a~subsuperspace. Let $F\in\fg_k$, where $k>0$, be such that
$[w, F]=0$ for any $w\in W$. Then, for any $u\in \fg_{-1}$, the
element $[u, F]$ commutes with any $w\in W$.
\end{Lemma}

\begin{proof} By the Jacobi identity
\[
[w, [u, F]]=(-1)^{p(w)p(u)}[u, [w, F]]+[[w, u], F]=0. 
\]
More explicitly, Lemma \ref{L4.2} says that if the expression
of $F$ does not involve certain indeterminates, then $[\partial_i,
F]$ does not involve those indeterminates either, for any $i$. \end{proof}

\ssbegin[On 3 prolongations: cases a), b), c)]{Theorem}[On 3 prolongations: cases a), b), c)]\label{th4.3}
Let $\xi=(\xi_1,\dots,\xi_k)$ and $\sdim V=m|n$, where $(m|n)\neq
(1|0) \text{ or } (0|1)$. Then,

\textup{a)} $(\Lambda(\xi), \Lambda(\xi)\ltimes\fvect(\xi))_* =
\fk(1|2k; k)$.

\textup{b)} $(\Pi(\Lambda(\xi)), \Lambda(\xi)\ltimes\fvect(\xi))_*
=\fm(k; k)$.

\textup{c)} $(V\otimes \Lambda(\xi), \fgl(V)\otimes
\Lambda(\xi)\ltimes \fvect(\xi))_*=\fvect(m|n+k; k)$.
\end{Theorem}

\begin{proof} It is clear that the Lie superalgebras in the right-hand sides of the isomorphisms a)--c) are subalgebras of the respective Lie superalgebras in the left-hand 
side. 

Let $\fg = \mathop{\oplus}_{-1\leq l}^{\infty}\ \fg_l$ be one
of~the superalgebras in the left-hand sides of the isomorphisms a)--c).
Let 
\be\label{xANDy}
\begin{array}{l}
y=-\id_V\otimes 1\in\fgl(V)\otimes \Lambda(\xi),\\
x=\Sigma \xi_i\partial_{\xi_i}\in\fg_0.\\
\end{array}
\ee 
Then, $\ad_ y$
determines the $\Zee$-grading of $\fg$ we consider in cases a)--c), while $\ad_x$
determines a~$\Zee$-grading of each component $\fg_l$:
\[
\fg_{l,i} = \{g\in\fg_l\mid \ad_x(g)=i\cdot g\}. 
\]

Note that $\fg_{-1}= \mathop{\oplus}_{i\geq 0}\ \fg_{-1,i}$
and $\fg_{0}= \mathop{\oplus}_{i\geq -1}\ \fg_{0,i}$. Set
\[
W:=\fg_{-1,0}\oplus \fg_{-1,1}=(V\otimes 1)\oplus (V\otimes
\Lambda^1(\xi)).
\]

We see that $\fg_0$ does not contain non-zero elements that commute
with $W$. Hence, by Lemma~\ref{L4.2}, none of the $\fg_l$ for $l>0$
contains non-zero elements that commute with $W$. But if
$\fg_l = \mathop{\oplus}_{i\geq -s}\ \fg_{l,i}$, then the
components $\fg_{l+1,t}$ with $t<-s-1$ should commute with $W$.
Hence, $\fg_{l+1} = \mathop{\oplus}_{i\geq -s-1}\ \fg_{l+1,i}$
and $\fg_l = \mathop{\oplus}_{i\geq -l-1}\ \fg_{l,i}$.

\subsubsection{Cases a) and b) of Theorem~\ref{th4.3}} Consider the operator $T=\ad_{2y+x}$, see eq.~\eqref{xANDy}.
Its eigenvalues are integers and, in each $\fg_l$, there are $\geq
l-1$ linearly independent eigenvectors. Set 
\[ 
\fg^T_s=\{g\in\fg\mid
T(g)=s\cdot g\}. 
\] 
Clearly, $\fg^T_s\subset
\mathop{\oplus}_{-1\leq l\leq s+1}\ \fg_l$ and $\dim \fg^T_s<\infty$.
Hence, we obtain a~new $\Zee$-grading: $\fg^T=
\mathop{\oplus}_{-2\leq s}^{\infty}\ \fg_s^T$ in $\fg$. The
negative part of this $\Zee$-grading is as follows, where $\pi$ is a~natural projection, see Subsection~ \ref{SS:4.4}: 
\[
\renewcommand{\arraystretch}{1.4}
\begin{array}{l}
\fg^T_{-2}=\fg_{-1,0}= \begin{cases}\Cee\cdot 1,&\text{in case
a)},\\
\Cee\cdot \pi(1),&\text{in case b)};\end{cases}\\
\fg^T_{-1}=\fg_{-1,1}\oplus \fg_{0,-1} =
\begin{cases}\Span(\xi_1,\dots,\xi_k,
\partial_{\xi_1} ,\dots,
\partial_{\xi_k}),&\text{in case a)}, \\
\Span(\pi(\xi_1),\dots,\pi(\xi_k),
\partial_{\xi_1} ,\dots, \partial_{\xi_k}),&
\text{in case b)}.\end{cases}
\end{array}
\]

Let us check that $\fg^T$ is transitive, \textit{i.e.},
\begin{equation}
\label{sth} \text{ if $g\in\fg^T_l$ with $l\geq 0$ and
$[\fg^T_{-1},g]=0$, then $g=0$.}
\end{equation}
 Indeed, for $g\in\fg$ satisfying ~\eqref{sth}, we have
$[W,g]=0$; hence, $g\in\fg_{-1}$. But only the subspace
$\fg^T_{-2}=\fg_{-1,0}\subset\fg_{-1}$ commutes with the
$\partial_{\xi_i}$ for all~ $i$.

Due to the universal property of the Lie superalgebras of series
$\fk$ and $\fm$ (see Subsection~\ref{SS:5.1}) this means that
\[
\fg^T\subset\begin{cases}\fk(1|k),\text{ and thus
$\fg^T=\fk(1|k)$},&\text{in
case a)},\\
\fm(k),\text{ and thus $\fg^T=\fm(k)$},&\text{in case b)}.\end{cases}
\]

\subsubsection{Case c) of Theorem~\ref{th4.3}}\label{SS:4.3.1s} Set
$W_0=\fg_{-1,0}=V\otimes 1$.

\parbegin[Technical]{Lemma}[Technical]\label{L4.1}
If $\sdim V\neq (1|0) \text{ or } (0|1)$ and $F\in \fg_{1}$ is such
that $[F,W_0]=0$, then $F=0$. In particular, $\fg_{1,-2}=0$.
\end{Lemma}

\begin{proof}
By Lemma \ref{L4.2}, the element $[F, v\otimes \varphi]\in \fg_0$
commutes with $W_0$ for any $v\in V$ and $\varphi\in \Lambda(\xi)$.
Hence, $[F, v\otimes \varphi]$ belongs to $\fvect(0|k)\subset\fg_0$.

Denote $D_{v, \varphi}:=[F, v\otimes \varphi]$. By transitivity,
$D_{v, \varphi}\psi \neq 0$ for some $v\in V$ and $\varphi, \psi\in
\Lambda(\xi)$. Take $u\in V$ linearly independent with $v$. Then
(the symbol $\pm$ stands for a~sign inessential to
us here),
\begin{equation}
\label{4.1}
\renewcommand{\arraystretch}{1.4}\begin{gathered} 
{}[[F, v\otimes \varphi],u\otimes\psi] \pm [v\otimes
\varphi, [F, u\otimes\psi]]= 
\pm u\otimes D_{v, \varphi}\psi \pm v\otimes D_{u, \psi}\varphi \neq
0.\end{gathered} 
\end{equation}
 By the Jacobi identity, the expression~\eqref{4.1} should vanish. This
contradiction completes the proof. \end{proof}

Lemmas \ref{L4.2} and \ref{L4.1} together mean that none of the components 
$\fg_l$ with $l>0$ contains elements that commute with $W_0$. Hence,
the least degree of $\Zee$-grading with respect to $\xi$ in
$\fg_{l+1}$ cannot be less than same in $\fg_l$ for $l\geq 0$. Thus,
$\fg_l= \mathop{\oplus}_{i\geq -1}\ \fg_{l,i}$.

Consider the operator $T=\ad_{y+x}$, see eq.~\eqref{xANDy}. Its eigenvalues are integers
and, in each $\fg_l$, there are $\geq l-1$ linearly independent
eigenvectors. So, as in cases a) and b) of Theorem~\ref{th4.3}, 
\[ 
\dim
\fg^T_s=\dim\{g\in\fg\mid T(g)=s\cdot g\}<\infty. 
\] 
Hence, in
$\fg$, we get a~new $\Zee$-grading $ \fg^T$.
The negative part of this $\Zee$-grading is
\[
\fg^T_{-1}=
\fg_{- 1,0}\oplus\fg_{0,-1} = V\otimes 1\oplus
\Span(\partial_{\xi_1},\dots, \partial_{\xi_k}).
\]

Moreover, by Lemma \ref{L4.1}, the elements of $\fg$ that commute
with $W_0=V\otimes 1\subset \fg^T_{-1}$ form the subspace
$\fg_{-1}\oplus \fvect(0|k)$. We can directly verify that the elements
of this subspace that commute with all the $\partial_{\xi_i}$ span
$\fg^T_{-1}$. Hence, $\fg^T$ is transitive. Due to the universality
property of the Lie superalgebra $\fvect(m|n+k)$, see
Subsection~\ref{SS:5.1}, this means that $\fg^T\subset\fvect(m|n+k)$;
thus, $\fg^T=\fvect(m|n+k)$.
\end{proof}

\ssbegin[Technical]{Corollary}[Technical]\label{cor4.3.2} Let 
$\xi=(\xi_1, \dots, \xi_k)$.

\textup{1)} If $\fg\subset\fk(1|2k; k)$ with $\fg_{-1} = \fk(1|2k;
k)_{- 1}=\Lambda(\xi)$ is such that $\fg_1$ does not
contain non-zero elements that commute with all elements of 
$W_0=\Cee\cdot 1\subset \fg_{-1}$, then $\fg\subset\fvect(1|k; k)$.

\textup{2)} If $\fg\subset\fm(k; k)$ with $\fg_{-1} = \fm(k;
k)_{-1}=\Pi(\Lambda(\xi))$ is such that $\fg_1$ contains no non-zero elements that commute with all elements of 
$W_0=\Cee\cdot \Pi(1)\subset \fg_{-1}$, then
$\fg\subset\fvect(0|k+1; k)$; hence, $\dim \fg<\infty$.
\end{Corollary}

\begin{proof}[Proof \nopoint] is the same as that in case c) of Theorem \ref{th4.3}. \end{proof}

\ssec{Back to the classification}\label{SS:5.6} Let, as above,
$\xi=(\xi_1, \dots, \xi_k)$ and let $V$ be an $m|n$-dimensional linear
superspace with coordinates $u=(u_1, \dots, u_{n+m})$.

\sssbegin[3 more cases]{Theorem}[3 more cases]\label{th5.5}
Let $\fg= \mathop{\oplus}_{-1\leq i}^{\infty}\ \fg_i$ be a~simple
infinite-dimensional transitive $\Zee$-graded Lie superalgebra with
$\fg_{-1}=V\otimes \Lambda(\xi)$ and $\fg_0\subset\fgl(V)\otimes
\Lambda(\xi)\ltimes \fvect(\xi)$.

\emph{1)} If $\dim V>1$, then $\fg$ is either $\fvect(1|n+k; k)$,
or $\fsvect(m|n+k; k)$ for $m>1$, or $\fsvect' (1|n+k; k)$;

\emph{2)} if $\sdim V=(1, 0)$, then $\fg$ is either $\fvect(1|k; k)$,
or $\fsvect'(1|k; k)$, or $\fk(1|2k; k)$, or $\fk\fas( ;3\xi)$.

\emph{3)} if $\sdim V=(0, 1)$, then $\fg$ is either $\fm(k; k)$, or
$\fb_{\lambda}(k; k)$ for $\lambda\in\Cee \Pee^1\setminus\{0\}$.
\end{Theorem}

\begin{proof} \underline{Case 1)}. By Theorem \ref{th4.3} we see that $\fg\subset\fvect(m|n+k; k)$.
Thus, we can use the realization of the elements of $\fg$ by vector
fields in indeterminates $u, \xi$. The $\Zee$-grading with respect to $\xi$ in each
degree-$l$ component $\fvect(m|n+k; k)_l$ determines a~filtration in
$\fg_l\subset \fg$.

Let $\gr(\fg)= \mathop{\oplus}_{-1\leq l}^{\infty}(\mathop{\oplus}_i\ \gr(\fg)_{l, i})$ be the corresponding bigraded
Lie superalgebra, \textit{i.e.},
\[
\gr(\fg)_{l,i}= \nfrac{\fg_l\cap \mathop{\oplus}_{j\geq i}\
\fvect(m|n+k; k)_{l,j}} {\fg_l\cap \mathop{\oplus}_{j> i}\
\fvect(m|n+k; k)_{l,j}}.
\]

By Lemma \ref{L4.1}, the $\fg_0$-module $\fg_{-1}$ is irreducible.
This gives the following conditions on $\gr(\fg)_{0, *}:=
\mathop{\oplus}_{i}\ \gr(\fg)_{0,i}$:
\begin{equation}
\label{4.2}
\begin{array}{l}\gr(\fg)_{0,-1} =
\Span(\partial_{\xi_1}, \dots,
\partial_{\xi_k})=\fvect(m|n+k; k)_{0,-1},\\ 

\gr(\fg)_{0,0} \subset \fvect(m|n+k; k)_{0,0}=\fgl(V)\oplus
\fgl(\xi)=\fgl(m|n)\oplus \fgl(0|k),\\

\gr(\fg)_{0,0} \text{ acts irreducibly in } V\otimes 1\subset
\fg_{-1}. \end{array}
\end{equation}

Without loss of generality we can assume that
\begin{equation}\label{woutloss} \text{the $\gr(\fg)_{0, 0}$-module
$\gr(\fg)_{0,-1} = \Span(\partial_{\xi_i}\mid i= 1, \dots, k)$ is
irreducible}.\end{equation} Indeed, if it is reducible and, \textit{e.g.}, $
\Span(
\partial_{\xi_1}, \dots, \partial_{\xi_s})$ is a~submodule, then
\[
\fg\subset \Big(\left(\fgl(V)\otimes \Lambda(\xi_1, \dots,
\xi_s)\ltimes \fvect(\xi_1, \dots, \xi_s)\right) \otimes
\Lambda(\xi_{s+1}, \dots, \xi_k) \Big)\ltimes \fvect(\xi_{s+1},
\dots, \xi_k) .
\]
So we set $\widetilde V=V\otimes \Lambda(\xi_1, \dots, \xi_s)$ and
$\widetilde \xi= (\xi_{s+1}, \cdots, \xi_k)$. Now, use the condition
$[\fg_{-1},\fg_1]= \fg_0$. As we saw in the proof of Theorem
\ref{th4.3},
\[
\fvect(m|n+k; k)_1=\mathop{\oplus}\limits_{i\geq -1} \fvect(m|n+k;
k)_{1, i},
\]
and $\gr(\fg)_{0, -1}$ can be obtained only by bracketing
\[
\gr(\fg)_{- 1, 0}=V\otimes 1=\Span(\partial_{u_1}, \dots,
\partial_{u_{m+n}})
\]
with 
\[
\fvect(m|n+k; k)_{1, -1}=
\Span(u_i\partial_{\xi_j}\mid \text{all possible }i, j).
\]
By condition ~\eqref{4.2}, for each $j=1, \dots, k$, there exist elements
$u_\alpha\partial_{\xi_j}\in \gr(\fg)_{1,-1}$ of the form
\[
u_\alpha
\partial_{\xi_j}\text{~~for any $\alpha=1, \dots, m+n$ and
$j=1, \dots, k$,}
\]
and $\gr(\fg)_{0,0}$ irreducibly acts on $V\otimes 1$, see assumption ~\eqref{woutloss}.

Hence,
$\gr(\fg)_{1,-1}=\fvect(m|n+k; k)_{1,-1}$.

Thus, for any $\alpha,\beta = 1, \dots, m+n$ and $i, j = 1, \dots,
k$, we have 
\[
[\xi_i\partial_{u_\beta}, u_\alpha
\partial_{\xi_j}]=\delta_{\alpha \beta}\xi_i\partial_{\xi_j} -
(-1)^{(p(u_\alpha )+1)(p(u_\beta)+1)}\delta_{ij}
u_\alpha\partial_{u_\beta} \in \gr(\fg)_{0, 0}.
\]
Consequently,
\[
\gr(\fg)_{0, 0}\supset (\fgl(V)\oplus
\fgl(\xi))\cap\fsvect(m|n+k).
\]
Now, consider the regrading $\fg^T$
of $\gr(\fg)$ into the standard $\Zee$-grading of $\fvect(m|n+k)$.
Then,
\[
\begin{array}{l}
\fg^T_{-1}= \gr(\fg)_{-1, 0}\oplus \gr(\fg)_{0, -1}=\fvect(m|n+k;
k)_{-1},\\

\fg^T_0= \gr(\fg)_{-1, 1}\oplus \gr(\fg)_{0, 0}\oplus \gr(\fg)_{1,
-1}\supset \fsl(m|n+k). \end{array}
\]
 Next, by the universality property $1^\circ)$, see Subsection~\ref{SS:5.1}, we
see that 
\[
\text{$\gr(\fg) = \fvect(m|n)$, or $\fsvect(m|n)$, or
$\fd(\fsvect(m|n))$, where $m>1$.}
\]

If $m=1$, then there are two more possibilities:
\[
\text{$\gr(\fg)=\fsvect'(1|n)$ and $\fd(\fsvect'(1|n)):=\fsvect'(1|n)\ltimes\Cee E$ for $E:=\sum x_i\partial_{x_i}$.}
\]

\underline{Case 2)}. If $\sdim V=(1|0)$, then
$\fg_{-1}=\Lambda(\xi)$ and $\fg_0\subset \Lambda(\xi)\ltimes
\fvect(0|k)$.

By Corollary \ref{cor4.3.2}, if the component $\fg_1$ does not
contain non-zero elements that commute with $\fg_{-1,0}=\Cee \cdot
1$, then $\fg\subset \fvect(1|k; k)$ and, by the above arguments, either $\fg =
\fvect(1|k; k)$ or $\fsvect'(1|k; k)$.

Suppose the component $\fg_1$ contains non-zero elements that commute with
$\fg_{- 1, 0}=\Cee \cdot 1$. Then, $\fg\subset \fk:=\fk (1|2k; k)$.

In terms of the standard indeterminates $t, \xi_1, \dots, \xi_k,
\eta_1, \dots, \eta_k$ for $\fk(1|2k)$, we have
\[
\begin{gathered}
\fk_{-1} = \fg_{-1}=\Lambda(\xi),\\
\Lambda(\xi)\otimes (t-\sum\xi_i\eta_i)\simeq 
\Lambda(k)\subset \fk_0,\\
\Lambda(\xi)\otimes\Lambda^1(\eta) \simeq 
\fvect(0|k)\subset \fk_0;\\
\fg_1= t\otimes\Lambda(\xi)\bigoplus t\otimes
\Lambda(\xi)\otimes\Lambda^1(\eta) \bigoplus
\Lambda(\xi)\otimes\Lambda^2(\eta).
\end{gathered}
\]
The subspace $W_1\subset \fk_1$ of elements that commute with $\fg_{-1,0}$
coincides with $\Lambda(\xi)\otimes\Lambda^2(\eta) $.

Consider $\gr(\fg)$. By conditions ~\eqref{4.2}, $\gr(\fg)_{0,-1}=\Span(\eta_1,
\dots, \eta_k)$, and $\eta_i$ acts in $W_1$ as $\partial_{\xi_i}$.

Since $\fg_1\cap W_1\neq 0$, we see that $\fg_1\cap \fk_{1, -2}\neq
0$, where $\fk_{1, - 2}=\Lambda^2(\eta)\otimes 1$. Hence,
$\gr(\fg)_{1, -2}\neq 0$. Set
\be\label{ha}
\fh:= \mathop{\oplus}\limits_{-1\leq i\leq k-1} \fh_i,
\text{~~where~~} \fh_i=\gr(\fg)_{0,
*}\cap \Lambda^1(\eta)\otimes \Lambda(\xi).
\ee
Clearly, $\fh\supset[\fg_{-1},
\gr(\fg)_{1,-2}]$.

\parbegin[Technical]{Lemma}[Technical]\label{SS:5.8.1} The component $\fh_{k-2}$, see eq.~\eqref{ha}, is not contained in
$\fsvect(0|k)_{k-2}$. In particular, $\fh_{k-2}\neq 0$.
\end{Lemma}

\begin{proof}
Let $X=\xi_1\cdots\xi_k\in\fg_{-1}$. Then,
\[
\renewcommand{\arraystretch}{1.4}
\begin{gathered}
{}\left[X, \
\mathop{\sum}\limits_{i<j}\alpha_{ij}\eta_i\eta_j\right]= (-
1)^{k+1}\mathop{\sum}\limits_{i<j}\alpha_{ij}\left(\nfrac{\partial X}
{\partial
\xi_i}\eta_j- \nfrac{\partial X}{\partial \xi_j}\eta_i\right)\longmapsto \\
 D:=(-1)^{k+1}\mathop{\sum}\limits_{i<j}\alpha_{ij}\left(\nfrac
{\partial X} {\partial \xi_i}\partial_{\xi_j}- \nfrac{\partial
X}{\partial \xi_j}\partial_{\xi_i}\right).
\end{gathered}
\]
But $\Div D=2 (- 1)^{k+1}\mathop{\sum}_{i<j}\ \alpha_{ij}
\nfrac{\partial ^2X}{\partial \xi_i \partial \xi_j} $. So $\Div D = 0
$ if and only if $\alpha_{ij}=0$ for all~ $i,j$. Since $
\gr(\fg)_{1,-2}\neq 0$, this proves the lemma.
\end{proof}

We have established that:

1) $\fh\subset \fvect(0|k)$,

2) $\fh_{-1}=\fvect(0|k)_{-1}$ and $\fh_{- 1}$ is an irreducible
$\fh_0$-module,

3) $\fh_{k-2}$ is not contained in $\fsvect(0|k)_{k-2}$.

\noindent Comparing these properties 1)--3) with the list of simple finite-dimensional Lie superalgebras, we conclude that $\fh=\fvect(0|k)$ for
$k>3$. If $k=3$, then $\fh$ can also be equal to
$\fsl(1|3)\subset \fvect(0|3)$.

Anyway, $\gr(\fg)_{0,*}\supset \Lambda^1(\eta)\otimes
\Lambda^1(\xi)\simeq \fgl(k)$. Thus, by the irreducibility of the
exterior square of the tautological $\fgl(k)$-module,
$\gr(\fg)_1\supset \Lambda^2(\eta)$.

The irreducibility of the $\fg_0$-module $\fg_{-1}$ means that
$\gr(\fg)_{0, *}$ should contain at least some elements of 
$\Lambda(\xi)\otimes (t- \sum\xi_i\eta_i)$. Hence, $t\in
\gr(\fg)_{0,*}$.

Now we are able to consider the regrading $\fg^T$ of $\gr(\fg)$ into
the standard $\Zee$-grading of $\fk(1|2k)$. By the above,
$\fg^T_{\leq 0}= \fk(1|2k)_{\leq 0}$. Hence, $\fg^T$ is either
$\fk(1|2k)$ or $\fk\fas\subset\fk(1|6)$, and therefore either
$\fg=\fk(1|2k; k)$ or $\fk\fas(;3\xi)$, respectively.

\underline{Case 3)}. If $\sdim V=0|1$, then
$\fg_{-1}=\Pi(\Lambda(\xi))$ and $\fg_0\subset\Lambda(\xi)\ltimes
\fvect(\xi)$, so $\fg\subset \fm(k; k)$. The rank-1 operators
with even covector acting in $\Lambda(\xi)\ltimes \fvect(\xi)$ are of the following form for $n$ even (resp., odd)
\[
g_i=\xi_1\cdots
\xi_n\partial_i,\text{ (resp., $h_i=\xi_1\dots \widehat{\xi_i}\dots
\xi_n\partial_i$) for }\ i=1, \dots, n.
\]
By Theorem \ref{th3.1}, the Lie superalgebra $\fg_0$ contains at
least one of them. Assumption~\eqref{woutloss} implies that $\fg_0$
contains either all the $g_i$ or all the $h_i$.

a) Let $\fg_0\supset\Span(g_1, \dots , g_n)$. Then, from
conditions~\eqref{4.2} we deduce that $\fg_0\supset \fvect(\xi_1, \ldots
\xi_n)$. Since the $\fg_0$-module $\fg_{-1}$ is irreducible, it follows that 
$\fg_0=\Lambda(\xi)\ltimes \fvect(\xi)$.

b) For the same reason as in case a), if $\fg_0\supset \Span(h_1,
\dots , h_n)$, then $\fg_0\supset\fsvect(\xi)$, and the following
cases are possible:
\[
\begin{array}{l}
1)\ \fg_0\simeq \fvect(\xi)\text{~~ realized as $\{v+\lambda\cdot
\Div\ (v)\mid v\in
\fvect(\xi)$ for $\lambda\neq 0\}$},\\
2)\ \fg_0\simeq \Vol_{0}(0|n)\ltimes\fsvect(\xi),\\
3)\ \fg_0\simeq \Lambda(\xi)\ltimes\fsvect(\xi).
\end{array}
\] In all these cases, $\fg_0$ is $\Zee$-graded by powers of
$\xi$ and $\fg_{0,0}\supset \fsl(0|k)\simeq \fsl(k)$. Due to
Corollary~\ref{cor4.3.2}, $\fg_{1, -2}\neq 0$; hence, $\fg_{1, -2}=
S^2(\id)$.

Consider the regrading $\fg^T$ of $\fg$ into the standard
$\Zee$-grading of $\fm(k)$. We get $\fg^T_{<0}=\fm_{<0}$ and
$\fg_0\supset \fspe(k)$. Then, by the universality property $3^\circ)$, see
Subsection~\ref{SS:5.1}, $\fg$ is either $\fm(k)$ or $\fb_\lambda(k)$, where
$\lambda\in\Cee \Pee^1\setminus\{0\}$.
\end{proof}

\ssbegin[Solvable $\fg_0$]{Corollary}[Solvable $\fg_0$]\label{SS:5.9} For the simple vectorial Lie
superalgebra $\fg$ with W-filtration, the Lie superalgebra $\fg_0$
is solvable only in the following cases:\footnote{Note that if $p=2$, there are examples of simple Lie superalgebras $\fg$ with solvable $\fg_\ev$, see \cite{BGL, BGLLS}.}
\[
\fg=\fvect (1|1), \quad \fsvect' (1|2; 1), \quad \fvect (1|2; 1)
\]
and the respective $\fg_0$ is isomorphic to
\[
\fgl(1|1), \quad \fsl(1|1)\otimes\Lambda(1)\ltimes \fvect(0|1),
\quad \fgl(1|1)\otimes\Lambda(1)\ltimes \fvect(0|1).\]
\end{Corollary}

\section{Step 5. Prolongations of almost simple Lie
superalgebras and their central extensions}\label{S:6}

In this step, we consider two essentially different cases: vectorial
and non-vectorial Lie superalgebras (Steps 5B and 5A, respectively).

Let $\fg$ be an almost simple Lie superalgebra or a~central
extension of an almost simple finite-dimensional
Lie superalgebra (they are listed in table ~\eqref{8.2}).
Let $T$ be a~faithful irreducible representation of $\fg$ in a
finite-dimensional superspace $V$. Our goal is to describe all the
collections $(\fg, T, V)$ for which there exists an element $g\in \fg$ such that 
$T(g)\in\fgl(V)$ is a~ rank-1 operator with an even covector.

Irreducible representations of simple vectorial Lie superalgebras
and their ``relatives'' --- central extensions of almost simple Lie
superalgebras connected to them --- are completely described, see Subsection \ref{SS:9.1}.
From this description we easily deduce
the list of admissible collections $(\fg, T, V)$.

For non-vectorial Lie superalgebras, there is a~complete description
of the typical (\textit{i.e.}, generic) representations (see \cite{S1}, where
the pioneering results of F.~Berezin and general ones due to V.~Kac
are presented with corrections and more lucidly; for a~modern
treatment of the problem, see \cite{SeVe}). But (bar a~couple of
exceptions) the rank-1 operators occur, if ever, in atypical
representations whose description was only conjectural when this
paper was being written (for reviews covering various aspects of the
problem, see \cite{L2}, \cite{KW}, \cite{Se5}).

Though we believed in the validity of the remarkable general
character formula due to Penkov and Serganova (at that time conjectural, now proved), our
proof in Step 5A did not appeal to the description of irreducible
representations. Inevitably then, the proof is rather lengthy, it occupies the rest of this Section.

Our method of description and of search for such collections $(\fg,
T, V)$ is based on pretty straightforward considerations which,
nevertheless, provide us with essential restrictions on the weight
diagram of $T$.

\ssbegin[List of rank-1 operators: Non-vectorial Lie
superalgebras]{Proposition}[List of rank-1 operators: Non-vectorial Lie
superalgebras]\label{prop5.0}
Let $T$ be a~faithful irreducible representation of $\fg$ in a
finite-dimensional superspace $V$. For the almost simple non-vectorial Lie superalgebra $\fg_{0}$ or any
central extension of such, the representations $(T, V)$ with a~rank-$1$ operator $T(g)$, where $g\in \fg_{0}$, with an even covector are
exactly those from Table ~\eqref{rank1T}, where the action of the
generator $z$ of the center of $\fg_{0}$ on $\fg_{-1}$ is normalized to
be equal to $-\id$.\end{Proposition}
\begin{equation}\label{rank1T}\footnotesize
\renewcommand{\arraystretch}{1.4}
\begin{tabular}{|c|c|c|}
\hline $V=\fg_{-1}$&$\fg_{0}$&$(\fg_{-1}, \fg_{0})_{*}$\cr 
\hline
$\id$ or $\id^{*}$&$\fsl(m|n)$ and $\fgl(m|n)$&$\fsvect(m|n)$ and
$\fvect(m|n)$\cr $\Pi(\id)$ or $\Pi(\id^{*})$&$\fsl(n|m)$ and
$\fgl(n|m)$&$\fsvect(n|m)$ and $\fvect(n|m)$\cr 
\hline
$\id$&$\fosp^{a}(m|2n)$ &$\fh(2n|m)$\cr &$\fc\fosp^{a}(m|2n)$
&$\fd(\fh(2n|m))$\cr 
\hline 
$\id$&$\fspe^{a}(n)$ and
$\fpe^{a}(n)$&$\fsle(n)$ and $\fle(n)$\cr 
\hline
$\id$&$\fc(\fspe^{a}(n))$ and $\fc(\fpe^{a}(n))$&$\fd(\fsle(n))$
and $\fd(\fle(n))$\cr 
\hline 
$\spin_{\lambda}$, $\lambda\neq
0$&$\fas$&$\fv\fle(4|4)$\cr
 \hline
\end{tabular}
 \end{equation}
 
\begin{proof}[Proof \nopoint] of Step 5A. We only have to test the simplicity of the prolongs $(\fg_{-1}, \fg_{0})_*$. 
\end{proof}

\sssbegin[On rank-1 operators]{Proposition}[On rank-1 operators]\label{prop5.1}
Let $T$ be a~faithful irreducible representation of $\fg$ in a
finite-dimensional superspace $V$. If $T(\fg)$ contains an odd operator of rank $1$, then either

$(1)$ $T(\fg)$ contains an even operator of rank $(1,0)$ and an even operator of rank $(0,
1)$, or

$(2)$ $T(\fg)$ contains an even operator of rank $(1, 1)$.
\end{Proposition}

\begin{proof} Let $\xi\in\fg_{\od}$ be such that
$\rk T(\xi)=1$, \textit{i.e.}, $T(\xi)v=\alpha(v)\del_\xi$, where $\alpha\in
V^*$. In other words,
$T(\xi)=\alpha\otimes v_\xi$. Let $\eta\in \fg_{\od}$, where $v_\xi\in V$ is a~fixed non-zero vector. Then,
\[
T([\eta, \xi])= T(\eta) T(\xi)+T(\xi) T(\eta).
\]
Consider the two operators 
\be\label{2ops}
\begin{array}{l}
A_\eta(v)= T(\eta) T(\xi)v=\alpha(v)\cdot T(\eta) v_\xi,\\
B_\eta(v)= T(\xi) T(\eta)v=\alpha(T(\eta v))\del_\xi .
\end{array}
\ee
Clearly, for each fixed $\eta$, both $A_\eta $ and $B_\eta $ are
of rank $\leq 1$; besides, they act in the spaces of opposite
parity.

\sssbegin[Technical]{Lemma}[Technical]\label{L5.2}
Let $T$ be a~faithful irreducible representation of $\fg$ in a
finite-dimensional superspace $V$. Let $W=\mathop{\cap}_{\eta\in\fg_{\od}}\ \ker T(\eta)$ and
$U=T(\fg_{\od})V$. Then, the subspaces $U$ and $W$ are
$\fg$-invariant.
\end{Lemma}

\begin{proof} Since $T$ is irreducible, it follows that $U=V$ and $W=0$. Now, the
$\fg_{\od}$-invariance is obvious. Let $g\in\fg_{\ev}$,
$\xi\in\fg_{\od}$, $v\in V$ and $u=\xi v\in U$. We claim that
$gu\in U$. Indeed,
\[
gu=g\xi v=\xi gv + [g, \xi]v\in T(\fg_{\od})V=V,
\]
because $gv\in V$ and $[g, \xi]\in\fg_{\od}$.

Let $w\in W$, $g\in \fg_{\ev}$. Let us show that
$gw\in W$. Indeed,
\[
\eta(gw)=g(\eta w)+[\eta, g]w=0
\]
because $\eta w=0$ and $[\eta, g]\in\fg_{\od}$. Hence, $gw\in W$.
\end{proof}

\parbegin[On rank-1 operators]{Corollary}[On rank-1 operators]\label{cor6.4}
There exist $\eta_1, \eta_2\in \fg_{\od}$ such that
$\rk A_{\eta_{1}}=1$ and $\rk B_{\eta_{2}}=1$, see \eqref{2ops}.
\end{Corollary}

To complete the proof of Proposition \ref{prop5.1}, note that if
$\rk B_{\eta_{1}}=1$ (resp., \hbox{$\rk A_{\eta_{2}}=1$}), then
\[
\rk T([\eta_{1}, \xi])=\rk
(A_{\eta_{1}}+B_{\eta_{1}})=(1, 1)\ \ \
\text{~~(resp., $\rk T([\eta_{2}, \xi])=(1, 1)$},
\]
whereas if $\rk B_{\eta_{1}}=\rk A_{\eta_{2}}=0$, then the operators
$T([\eta_{1}, \xi])$ and $T([\eta_{2}, \xi])$ are of rank 1 and act
on the spaces of opposite parity.
\end{proof}

\parbegin[On infinite prolongation of depth 1]{Corollary}[On infinite prolongation of depth 1] Let $T$ be a~faithful irreducible representation of $\fg_0$ in a
finite-dimensional superspace $\fg_{-1}$. If $\dim(\fg_{-1}, \fg_0)_*= \infty$, then there exists an element $g_0\in\fg_0$ such
that $\rk (T(g_0))|_{(\fg_{-1})_{\ev}}=1$ and $\rk (T(g_0)|_{(\fg_{-1})_{\od}}\leq
1$.\end{Corollary}

\ssbegin[On rank-1 operators in reductive algebras]{Proposition}[On rank-1 operators in reductive algebras]\label{prop5.3}
Let $\fG_j$ for $j=1, \dots, k$ be simple Lie algebras \emph{(\textbf{not
super}\-algebras)}, let $\fr$ be a~solvable Lie algebra and let $\fg=\fr\oplus
(\bigoplus \fG_j)$ be their direct sum. Let $T$ be a~\emph{(perhaps,
reducible)} representation of $\fg$ in a~space $V$ with a~rank-$1$
operator. Then, either

\emph{1)} there exists a~$1$-dimensional $\fg$-invariant subspace
of $V$, or

\emph{2)} up to renumbering, $\fG_1=\fsl(n)$ or $\fsp(2n)$ and
$V=W\oplus U$, where $U=\ker T(\fG_1)$ and $W$ is isomorphic, as a~
$\fg$-module, to either $\id\otimes\chi_0\otimes \lambda$ or
$\id^*\otimes\chi_0\otimes \lambda$, where $\id$ is the tautological
representation of $\fG_1$ and $\chi_0$ is the zero $1$-dimensional
representation of $\fG_2\oplus \dots \oplus \fG_k$, and $\lambda$ is
a~character of $\fr$.
\end{Proposition}

\begin{proof} For any
semi-simple Lie algebra ${\fg_{ss}:=\fG_1\oplus \dots \oplus \fG_k}$, any of its finite-dimensional representations is
completely reducible, and hence the sum of $\fg_{ss}$ with $\fr$ is 
direct; therefore, $V$ can be represented in the form
$V=\mathop{\oplus}_i\ V_i$, where $T|_{V_{i}}=\tau_1^i\otimes
\dots \otimes \tau_k^i\otimes \tau_\fr^i$ and $\tau_j^i$ is an
irreducible representation of $\fG_j$ and $\tau_\fr^i$ is an
irreducible representation of $\fr$.

Since $T$ contains a~rank-1 operator, so does the restriction
$T|_{V_{i_{0}}}$ for some $i_{0}$. Therefore, at most one of
the modules given by $\tau_1^{i_0}$, \dots,
$\tau_k^{i_0}$, $\tau_\fr^{i_0}$ can be of dimension $>1$.

If $\dim \tau_1^{i_0}= \dots =\dim \tau_k^{i_0}=\dim
\tau_\fr^{i_0}=1$, then $\dim V_{i_0}=1$ and possibility 1) of Proposition~
\ref{prop5.3} is realized.

If $\dim \tau_1^{i_0}= \dots =\dim \tau_k^{i_0}=1$, but $\dim
\tau_\fr^{i_0}>1$, then, since $\fr$ is solvable, $V_{i_0}$ contains a~$1$-dimensional $\fr$-submodule
$V_0$. Since $\fg_{ss}V_{i_0}=0$, it follows that $V_0$ is
a~$\fg$-submodule, \textit{i.e.}, once again the same possibility 1) of
Proposition~\ref{prop5.3} is realized.

Finally, let $\dim \tau_{j_0}^{i_0}>1$ for some $j_{0}$. Let us show
that in this case possibility 2) or Proposition~\ref{prop5.3} is realized.
Indeed, up to a~renumbering of the subalgebras $\fG_j$, we can assume
that $j_{0}=1$. Then, 
\[
\dim \tau_2^{i_0}=\dots =\dim
\tau_k^{i_0}=\dim \tau_\fr^{i_0}=1
\]
and the image
$T(\fg)|_{V_{i_{0}}}$ is isomorphic to either $\fG_1$ or to its trivial
central extension $\fc(\fG_1)$. Moreover,
$T(\fG_2\oplus\dots\oplus\fG_k)|_{V_{i_{0}}}=0$ and
$T(\fr)|_{V_{i_{0}}}=\lambda$ is a~character.

Since there are only three series of simple infinite-dimensional Lie
algebras with W-filtration of depth one, it follows that
a~\textit{faithful irreducible} representation $\tau$ of the Lie
algebra $\fh$ contains a~rank-1 operator if and only if for
a~character $\lambda$ of $\fh$ we have one of the following:
\[
\begin{array}{rl}
\fh=\fsl(n)&\text{and }\tau=\id \text{ or }
\id^*,\\
\fh=\fgl(n)&\text{and }\tau=\id\otimes \lambda
\text{ or } \id^*\otimes \lambda, \\
\fh=\fsp(2n)&\text{and }\tau=\id \simeq \id^*,\\
\fh=\fc\fsp(2n)&\text{and }\tau=\id \otimes
\lambda.
\end{array}
\]

In our case of the Lie algebra $\fg$, the role of $\fh$ is played by the image
$T(\fG_{1})|_{V_{i_{0}}}$ and we see that $\fG_1$ is isomorphic to
either $\fsl(n)$ or $\fsp(2n)$, while $T(\fG_{1})|_{V_{i_{0}}}$ is
either $\id$ or $\id^*$.

Thus, we take $V_{i_0}$ for the role of $W$. To complete the proof,
it suffices to observe that the rank of any operator $T(g)$ is
calculated by the formula
\[
\rk T(\fg)=\mathop{\sum}\limits_{1\leq i\leq k}\rk T(\fg)|_{V_{i}},
\]
and therefore $ T(\fg)|_{V_{i}}=0$ for $i\neq i_{0}$.
\end{proof}

\ssec{The cases of Proposition \ref{prop5.3}: two
regular and two exceptional subcases}
\label{SS:6.6} Now, let us return to the situation of interest to
us: $\fg$ is an almost simple Lie superalgebra or a~central
extension thereof, and $T$ its a~faithful representation of $\fg$ in
a~finite-dimensional superspace $V$ such that among the operators
$T(g)$ with $g\in\fg$, there is a~ rank-1 operator with an even
covector.

Consider the restriction $T|_{\fg_{\ev}}$. Thanks to Proposition~\ref{prop5.3}, the
representation of $\fg_{\ev}$ in the subspace $V_{\ev}$ must have
a~rank-1 operator, whereas the representation of $\fg_{\ev}$ in the
whole space $V$ has either a~rank $(1, 0)$ operator or a~rank $(1, 1)$
operator. Moreover, from Table ~\eqref{8.2} we see that $\fg_{\ev}$
is always of the form $\fG_1\oplus\dots\oplus\fG_k\oplus\fr$, where
the $\fG_i$ are simple and $\fr$ is a~solvable Lie algebra ($\fr$
can be zero and is almost always commutative). Hence, Proposition
\ref{prop5.3} is applicable to the $\fg_{\ev}$-action in the spaces
$V_{\ev}$ and $V_{\od}$.

Let, as in Proposition~\ref{prop5.3},
$\fg_{ss}=\fG_1\oplus\dots\oplus\fG_k$. Consider, separately, the
two cases corresponding to the possibilities of Proposition~\ref{prop5.3}:
\[
\begin{array}{l}
\text{1) among the operators $T(\fg_{ss})$, there are no operators
of rank $(1, 0)$ or $(1, 1)$;}\\

\text{2) among the operators $T(\fg_{ss})$, there is an operator 
of rank $(1, 0)$ or $(1, 1)$.}\\
\end{array}
\]

\sssec{Case 1) of Proposition \ref{prop5.3}} As is clear from the above discussion of the possible cases of this Proposition, if we
are in the setting of case 1), then the operator of rank
$(1, 0)$ or $(1, 1)$ is in this case of the form $T(h)$, where
$h\in\fr$. Moreover, since $T$ is irreducible, then
$h$ cannot belong to the center of the whole algebra~ $\fg$.

This possibility 1) is realized in the following regular subcases (as follows from
the classification of finite-dimensional simple Lie superalgebras and
their central extensions, see Table
~\eqref{8.2}):
\[
\fg=\fsl(m|n)\text{ for }m\neq n,\; \fgl(m|n), \;
\fosp(2|2n),\; \fc(\fosp(2|2n))
\]
and also in the \textit{exceptional subcases}, to be considered
separately:
\[
\begin{array}{l}
\text{a)~}\; \fpsl(2|2)\subset \fg\subset \fder~\fpsl(2|2)\text{
and }1\leq \dim \fg/\fpsl(2|2)\leq 2,\\
\text{b) the same as in a), but with
$\fc(\fpsl(2|2))$ instead of $\fpsl(2|2)$.}\\
\end{array}
\]

\paragraph{The regular subcases}\label{SS:5.4.1}
In these subcases, $\fg$ is of the form $\fg= \fg^{-1}\oplus\fg_{0}
\oplus\fg^{1}$, where $\fg_{0}=\fg_{\ev}$, and therefore
$\fg_{\od}=\fg^{-1}\oplus\fg^{1}$. Then, under a~proper normalization, the
operator $h$ acts as the grading operator:
\[ 
[h, g]=\pm g~~\text{~~ for any }g\in\fg^{\pm 1}. 
\] 
Let
$W_{\ev}=\Cee w_0$ be a~1-dimensional $\fg_{\ev}$-invariant
subspace of $V_{\ev}$ and $h\cdot w_0=\lambda w_0$. Then, \[
h(gw_0)=g(hw_0)+[h, g]w_0=\lambda(gw_0)\pm (gw_0)= (\lambda\pm
1)(g\cdot w_0), \] \textit{i.e.}, on each of the subspaces $W_{-1}=\fg^{-1}W$
and $W_{1}=\fg^{1}W$ the operator $h$ acts as multiplication by
a~scalar.

Consider the action of $\fg_{ss}$ on $W_{\pm 1}$. If
$g_0\in\fg_{ss}$ and $\xi\in \fg^{\pm 1}$, then
\[
g_0(\xi w_0)=\xi(g_0w_0)+[g_0, \xi]\cdot w_0=[g_0,\xi]w_0.
\]
Hence, the subspaces $W_{\pm 1}$ are $\fg_{\ev}$-invariant
and the representation of $\fg_{ss}$ in $W_{\pm 1}$ is isomorphic to
the quotient representation of $\fg_{ss}$ in $\fg^{\pm 1}$.

Since in all the cases considered the
representation of $\fg_{ss}$ in $\fg^{\pm 1}$ is
irreducible, we see that either $W_{\pm 1}=0$ or $\dim
W_{\pm 1}=\dim\fg^{\pm 1}>1$.

To ensure that $\rk T(h)=1$ or $(1, 1)$, we should have $\lambda=\pm
1$; moreover, if $\lambda+k\neq 0$, then $W_{k}=0$ for $k=\pm 1$.
Hence, $\fg_{\od}W_{\ev}=W_{-\lambda}$ and $h|_{W_{- \lambda}}=0$.
Then, $h|_{\fg^{\pm 1}\cdot W_{- \lambda}}=\pm \id$ and $\fg^{\pm
1}W_{-\lambda}\subset V_{\ev}$.

But since $\rk T(h)|_{V_{\ev}}=1$, we get 
\[
\fg_\lambda W_{-\lambda}=W_{\ev}, \quad \fg_{-
\lambda}W_{-\lambda}=0~~\text{ and }~~\fg_{\od}^2W_{\ev }=W_{\ev}.
\]
Therefore, $W_{\lambda}\oplus W_{\ev }$ is invariant under the whole $\fg$; hence, $V=W_{\lambda}\oplus W_{\ev }$ and
$\sdim V=(1, n)$, where $n=\dim \fg^1=\dim\fg^{-1}$.

Thus,
\[
\fg\subset \fgl(1|n) \text{ and }
\dim\fg_{\od}=2n=\dim\fgl(1|n)_{\od}.
\]
Therefore, either $\fg=\fsl(1|n)$ or $\fgl(1|n)$, and either $T=\id$
or $T=\id^*$ (recall how the action of the center is normalized, see
Proposition~\ref{prop5.0}).

\paragraph{The exceptional subcases}\label{SS:5.4.2}
These cases are
related to the simple Lie superalgebra $\fpsl(2|2)$ whose Lie
superalgebra of outer derivations is isomorphic to $\fsl(2)$.

Subcase a): $\fpsl(2|2)\subset \fg\subset \fder~\fpsl(2|2)$ and
$1\leq\dim \fg/\fpsl(2|2)\leq 2$.

Subcase b): $\fc\fpsl(2|2)\subset \fg\subset \fder~\fc\fpsl(2|2)$ and
$1\leq\dim \fg/\fc\fpsl(2|2)\leq 2$.
Let $z$ be the central element of $\fg$. 

Let $\fr_0\subset\fsl(2)$ be a~subalgebra of the Lie algebra of
outer derivations of $\fpsl(2|2)$ such that $1\leq\dim \fr_0\leq 2$.
Let
\[
\fr=\begin{cases}
\fr_0,&\text{ in subcase a)},\\
\fr_0\oplus\Cee z,&\text{ in subcase b).}
\end{cases}
\]
Denote by $R_0=\id_{\fsl(2)}|_{\fr_{0}}$ the restriction of the
tautological representation of $\fsl(2)$ to $\fr_0$. Then,
$\fg_{\ev}=\fsl_1(2)\oplus\fsl_2(2)\oplus\fr_0$ and
$\fg_{ss}=\fsl_1(2)\oplus\fsl_2(2)\simeq \fo(4)$, where the indices label copies of $\fsl(2)$ and --- in what follows shortly --- their tautological modules. Let $\chi_0$ be the
trivial character of $\Cee z$; then, as $\fg_{\ev}$-module,
$\fg_{\od}$ is isomorphic to $\id_1\boxtimes\id_2\boxtimes R$, where
\[
R= \begin{cases}
R_0,&\text{ in subcase a)},\\
R_0\otimes\chi_0,&\text{ in subcase b).}
\end{cases}
\]

Set $X=\mat{1&0\\0&- 1}\in\fsl(2)$ and $Y=\mat{0&1\\0&0}\in\fsl(2
)$. Clearly, up to the choice of a~basis in $\fg_{\od}$, there are
three essentially distinct possibilities for $\fr_0$:

1) $\fr_0=\Span(X, Y)$;

2) $\fr_0=\Cee X$;

3) $\fr_0=\Cee Y$.

Recall that we are looking for an irreducible representation $T$ of
$\fg$ in a~superspace $V$ such that $T(\fg_{\ev})$ contains an
operator of rank $(1, 0)$ or $(1, 1)$ whereas
$T(\fg_{ss})|_{V_{\ev}}$ has no such operators.

Then, the operator of rank $(1, 0)$ or $(1, 1)$ is $T(r)$, where
$r\in\fr_0$. Without loss of generality we can assume that either
$r=X$ or $r=Y$.

If $r=X$, then the arguments used in the regular case are applicable and we
conclude that $\fg$ has no such representations.

Let $r=Y$. Let, as in the general case, $W_{\ev}=\Cee w_0\subset
V_{\ev}$ be a~$1$-dimensional $\fg_{\ev}$-invariant subspace. Then,
$T(\fg_{ss})w_0=0$, and the representation $T$ determines
a~character $\chi_{\lambda , \alpha}$ of $\fr$ via the formulas
\[
\begin{array}{l}
T(z)w_0=\chi_{\lambda ,\alpha}(z)w_0=\lambda w_0\text{ for
}\lambda\neq
0,\\
T(Y)w_0=\chi_{\lambda,\alpha}(Y)w_0=\alpha w_0.
\end{array}
\]
(Clearly, if $\dim \fr_0=2$, then $\alpha =0$.)

Consider the $\fg_{\ev}$-submodule $\fg_{\od}W_{\ev}\subset V_{\od}$.
It is isomorphic to a~quotient of the $\fg_{\ev}$-module
\[
\fg_{\od}\otimes W_{\ev}=\id_1\boxtimes\id_2\boxtimes(R\otimes \chi_{\lambda,
\alpha}). 
\]
Therefore, the following three cases are possible:

1) $\fg_{\od}W_{\ev}=0$. Then, $W_{\ev}$ is a~$\fg$-invariant. Hence, $V=W_{\ev}$ and $T$ is not faithful.

2) $\fg_{\od}W_{\ev}\simeq \fg_{\od}\otimes W_{0}$. Then,
\[
\rk ~T(Y)|_{ \fg_{\od}\otimes W_{0}}=\begin{cases} 4,&\text{if }
\alpha
=0\\
8,&\text{if } \alpha \neq 0.\end{cases}
\] 
Both these
possibilities contradict our assumption on the rank.

3) $\fg_{\od}W_{\ev}\simeq _1\boxtimes\id_2\boxtimes\chi_{\lambda,
\alpha}$. If $\alpha\neq 0$, then $\rk ~T(Y)|_{\fg_{\od}W_{\ev}}=4$.
Hence, $\alpha=0$, in other words, $T(Y)|_{W_{0}}=T(Y)|_{
\fg_{\od}W_{\ev}}=0$. Set
\begin{equation}
\label{5.4.2.1}
\fg_Y:=\IM~\ad_Y\subset\fg_{\od}.
\end{equation}
Then,
\[
T(Y)W_{\ev}= T(Y)\fg_{\od}W_{\ev}- \fg_{\od}
T(Y)
W_{\ev}=0.
\]
Thus, $w_0$ is the highest weight vector in $T$.

In subcase a) this means that $T$ is $1$-dimensional, and therefore
not faithful.

In subcase b), let $c$ be the cocycle
that determines the central extension of $\fpsl(2|2)$ to $\fsl(2|2)$. As
is clear from Table~\eqref{8.2},
\begin{equation}
\label{5.4.2.2} c|_{\fg_{\od}\times \fg_Y}=0.
\end{equation}

Let $\fg_{-Y}\subset\fg_{\od}$ be a~$\fg_{ss}$-invariant subspace
complementary to $\fg_Y$. Then, $\fg_{\od}W_{\ev}=
\fg_{-Y}W_{\ev}$ by definition ~\eqref{5.4.2.1}. Take some non-zero $\xi\in\fg_Y$ and $\eta\in\fg_{-
Y}$. Then,
\[
\xi \eta w_0 -\eta\xi w_0=[\xi, \eta]w_0=c(\xi, \eta)\cdot
\lambda w_0.
\]
By condition ~\eqref{5.4.2.2}, both $Y$ and all the elements $\xi\in\fg_Y$
annihilate the space $\fg_{\od}W_{\ev}$. Hence, $\fg_{\od}W_{\ev}$
generates a~$\fg$-invariant subspace in $V$ which does not contain
$W_{\ev}$. This contradicts irreducibility of $T$.

Conclusion: in none of the exceptional subcases $\fg$ has irreducible
representations we are looking for.

\sssec{Case 2) of Proposition \ref{prop5.3}}\label{SS:5.5}
Due to Proposition~\ref{prop5.3}, $\fg_{ss}=\fG_1\oplus \fg_{comp}$, where
$\fG_1=\fsl(n)$ or $\fsp(2n)$ and the complementary Lie algebra
$\fg_{comp}$ is either semi-simple or zero. So $V=W\oplus U$ for
(decomposable) $\fG_{1}$--module $W$ and $\fg_{comp}$-module $U$,
where $T(\fG_1)|_{U}=0$, the $\fg_{ss}$-module $W_{\ev}$ is
isomorphic to either $\id\otimes\chi_0$ or $\id^*\otimes\chi_0$ and
the $\fg_{ss}$-module $W_{\od}$ is either zero or isomorphic
to either $\id\otimes\chi_0$ or $\id^*\otimes\chi_0$.

Consider the weight diagrams of $V_{\ev}$ and $V_{\od}$ with respect
to $\fg_{ss}$. Denote by $\eps_1, \dots , \eps_n$ the standard
vectors that span the space dual to the Cartan subalgebra of
$\fG_1$, see \cite{OV}, and by $\delta_1, \dots , \delta_N$ the
similar vectors for $\fg_{comp}$. (Clearly, if $\fg_{ss}$ is simple,
and therefore $\fg_{comp}=0$, these complementary vectors $\delta_i$
are absent.) Then, the collection of weights (with respect to
$\fg_{ss}$) in $W_{\ev}$ and $W_{\od}$ is as follows (with
multiplicity one):
\[
\fG_1=\begin{cases}\fsl(n):&\text{either $\eps_1$, \dots ,
$\eps_n$ or $-\eps_1, \dots , -\eps_n$;}\\
\fsp(2n):&\eps_1, \dots , \eps_n\; \text{ and $-\eps_1$, \dots ,
$-\eps_n$.}\end{cases}
\]
In the spaces $U_{\ev}$ and $U_{\od}$, all the weights are of the
form $\varphi(\delta_1, \dots , \delta_N)$, \textit{i.e.}, depend only on the
$\delta$'s. The weight diagram of such a~form will be referred to as
a~\textit{desirable} one. Let $\alpha$ be an odd root of $\fg$ and
$\xi_\alpha$ the corresponding root vector; let $\lambda$ be
a~weight of $V$ and $v_\lambda$ a~ corresponding weight vector.
Then, $\xi_\alpha v_\lambda$ is a~weight vector of weight
$\lambda+\alpha$.

\sssec{The scheme of our further considerations. Case-by-case study}\label{SS:6.7.1}

$1$) We start with a~$\fg_{\ev}$-invariant subspace $W_{\ev}$ whose
weights we know.

$2$) We apply $\fg_{\od}$ to this subspace and get a~subspace of
$V_{\od}$. By comparing the desirable weight diagram of $V_{\od}$
with weights $\lambda'$ of the form $\lambda+\alpha$ we get an
information about possible weight diagrams of $\fg_{\od}W_{\ev}$.

$3$) By applying item $2$) to a~$\fg_{\ev}$-invariant subspace of
$\fg_{\od}W_{\ev}$, we deduce that:


a) $\fg_{\od}^2W_{\ev}=W_{\ev}$, and therefore
$V=W_{\ev}\oplus\fg_{\od}W_{\ev}$;

b) a~ description of the weight diagram of $V$.


$4$) Since the dimension of $V$ is very low and its homogeneous
components are irreducible with respect to $\fg_{\ev}$, a~ direct
verification and dimension considerations quickly bring about the result.
If these considerations are not conclusive, we compare the
superspace $V$ obtained, the highest weights with
respect to $\fg_{\ev}$ of its homogeneous components we know, with the irreducible $\fg$-modules
with the same highest weights. Since all such modules are quotients
of the induced modules, it suffices to investigate only what Rudakov called
\textit{singular vectors} (see the review \cite{GLS3}) of two induced modules.


\paragraph{ $\fG_1=\fsp(2n)$}\label{SS:6.8.1} In this case, 
$\fg$ is either $\fosp(m|2n)$ or $\fosp(m|2n)\oplus\Cee z$. The weights
$\lambda$ of $W_{\ev}$ are $\pm \eps_1, \dots , \pm \eps_n$\,.

\underline{$m>2$}. For $\fosp(2m|2n)$, the odd roots are of the form
$\alpha=\pm\eps_i\pm\delta_j$. For
$\fosp(2m+1|2n)$, to the roots of this form we must add the
roots $\alpha=\pm\eps_i$.

Thus, the only sums $\lambda+\alpha$ that are desirable weights are
\[
\lambda'=\begin{cases}
\pm\delta_j,&\text{for 
$\fosp(2m|2n)$},\\
\pm\delta_j\text{ and $0$},&\text{for $\fosp(2m+1|2n)$. }\end{cases}
\]

By considering the sums $\lambda''$ of the form $\lambda'+\alpha$,
we see that, among them, the only desirable weights are those of the
form $\pm\eps_i$. Since they are of multiplicity $1$ in $V_{\ev}$, it follows that $\fg_{\od}^2W_{\ev}=W_{\ev}$; hence,
$V=W_{\ev}\oplus\fg_{\od}W_{\ev}$.

Consider the $\fg_{ss}$-module $\fg_{\od}W_{\ev}$ more closely.
Clearly, it is isomorphic to the quotient of $\fg_{\od}\otimes
W_{\ev}$. As $\fg_{ss}$-module,
$\fg_{\od}\simeq\id_{\fsp(2n)}\otimes\id_{\fo(m)}$ and
$W_{\ev}\simeq\id_{\fsp(2n)}\otimes\chi_0$. Thus, $\fg_{\od}\otimes
W_{\ev}$ is the quotient of
$\left(\id_{\fsp(2n)}\otimes\id_{\fsp(2n)}\right)\bigotimes\id_{\fo(m)}$
by the zero $\fsp(2n)$-action.

Now, observe that the tensor square of the tautological
$\fsp(2n)$-module is reducible:
\[
\id_{\fsp(2n)}\otimes\id_{\fsp(2n)}=
S^2(\id_{\fsp(2n)})\oplus\Lambda^2(\id_{\fsp(2n)}).
\]
Moreover, $S^2(\id_{\fsp(2n)}) \simeq \ad_{\fsp(2n)}$ (hence,
irreducible) and $\Lambda^2(\id_{\fsp(2n)})$ is the direct sum of
the trivial $1$-dimensional module $\mathbbmss{1}$ and the
 irreducible submodule complementary to it. Therefore, the subspace
of invariants in the tensor square
$\id_{\fsp(2n)}\otimes\id_{\fsp(2n)}$ is $1$-dimensional. Therefore,
the $\fg_{ss}$-module $\fg_{\od}\otimes W_{\ev}$ can be either
trivial or isomorphic to $\chi_0\otimes \id_{\fo(m)}$.

Thus, $\sdim V=(2n|1)$ or $(2n|m)$. 
Thanks to the faithfulness of the $\fo(m)$-action, dimension considerations show that only
the second possibility can be realized; the representation $T$ is,
actually, $\Pi(\id)$.

\underline{$m=2$}. Then, $\fg_{ss}=\fsp(2n)$. The weights $\lambda$
of $W_{\ev}$ are $\pm \eps_1, \dots, \pm \eps_n$ and the odd roots
are also $\pm \eps_1, \dots, \pm \eps_n$. Each root is of
multiplicity 2, and the center of $\fg_{\ev}$ separates the odd roots
into the negative and positive ones: $\fg_{\od}=\fg^{-1}\oplus
\fg^1$; the $\fg_{ss}$-action in $\fg^{-1}$ is $\id$, the $\fg_{ss}$-action on
$\fg^1$ is $\id^*$.

The desirable weight $\lambda'$ is of the form $\alpha+\lambda$ only for
$\lambda'=0$. The subspace $ \fg^{-1}W_{\ev} $ is the quotient of
$\fg^{-1}\otimes W_{\ev}\simeq \id^*\otimes \id$ and $\fg^{1}W_{\ev}$
is the quotient of $\fg^{1}\otimes W_{\ev}\simeq \id\otimes \id$. As
already noted, each of the modules $\id^*\otimes \id$ and
$\id\otimes \id$ contains only a~$1$-dimensional subspace of
invariants. Thus,
\[
\dim \fg_{\od}W_{\ev} = 0, 1 \text{ or } 2.
\]
The desirable weights of the form $\lambda'+\alpha$ are equal to
$\pm\eps_i$. Hence,
\[
\fg_{\od}^2 W_{\ev}=W_{\ev}\text{~~and~~}
V=W_{\ev}\oplus \fg_{\od}W_{\ev}.
\]
As above, dimension considerations imply that
$T(\fg_{\od})=\fgl(2n|1)_-$ in the standard format. Hence,
\[
\sdim V=(2n|2)\text{~~and~~}T=\Pi(\id).
\]

\underline{$m=1$}. Then, $\fg_{ss}=\fsp(2n)$ and $\fg_{\od}\simeq \id$,
as a~$\fg_{ss}$-module. The same arguments as above show that $\dim
\fg_{\od}W_{\ev}=1$ and $V=W_{\ev}\oplus \fg_{\od}W_{\ev}$. Thus,
$T=\Pi(\id)$.

\paragraph{$\fG_1=\fsl(n)$} This case splits into the following subcases A) --
F).\label{SS:6.8.2}

\underline{A) $\fg=\fspe(n)$, $\fspe(n)\oplus \Cee z$, $\fas$,
$\fpe(n)$, $\fpe(n)\oplus\Cee z$}. We will assume that $n>3$: as
$\fspe(3)\simeq \fsvect(0|3)$, this case and all it relatives are
considered in the study of vectorial Lie superalgebras.

For these Lie superalgebras, $\fg_{ss}=\fsl(n)$ and the weights of
$W_{\ev}$ are either $\lambda=\eps_1, \dots, \eps_n$ or
$\lambda=-\eps_1, \dots, -\eps_n$. The odd roots are
$\alpha=\pm(\eps_i +\eps_j)$ for $i\neq j$, and $2\eps_i$. Observe
that $ \fg_{\od}$ naturally splits into the sum of a~ positive and a~ negative
parts: $\fg_{\od}=\fg^{-1} \oplus\fg^1$.

The roots in $\fg^{-1}$ are $\alpha=-\eps_i- \eps_j$ while those in
$\fg^1$ are $\alpha=\eps_i+\eps_j, 2\eps_i$.

The weights $\lambda'$ of the form $\lambda+\alpha$ that occur in
the desirable diagram are as follows.

If $\lambda=\eps_1, \dots, \eps_n$, then $\lambda'=-\eps_1, \dots,-
\eps_n$ and $\fg^1W_{\ev}=0$.

If $\lambda=-\eps_1, \dots, -\eps_n$, then $\lambda'=\eps_1, \dots,
\eps_n$ and $\fg^{- 1}W_{\ev}=0$.

Recall that the multiplicity of the roots of the form $\pm\eps_i$ in
$V_{\od}$ is necessarily equal to 1. Therefore, $\dim
\fg_{\od}W_{\ev}=n$.

Similarly, the desirable weights $\lambda''$ of the form
$\lambda'+\alpha$ are as follows.

If $\lambda'=\eps_1, \dots, \eps_n$, then $\lambda'+\alpha=-\eps_1,
\dots, -\eps_n$ and $\fg^2_{\od}W_{\ev}=W_{\ev}$.

If $\lambda'=-\eps_1, \dots, -\eps_n$, then $\lambda'+\alpha=\eps_1,
\dots, \eps_n$ and $\fg^2_{\od}W_{\ev}=W_{\ev}$.

Thus, $V=W_{\ev}\oplus \fg_{\od} W_{\ev}$ and $\sdim V=(n|n)$.
Hence, if $\fg=\fspe(n)$ or $\fpe(n)$, then either $T=\id$ or
$T=\id^*$. Both these representations contain (odd) operators of
rank 1, but, for $T=\id^*$, the covector is odd.

For the Lie superalgebras $\fg=\fspe(n)\oplus \Cee z$ and
$\fg=\fpe(n)\oplus \Cee z$, we see that $T=\id\otimes\chi$, where
$\chi$ is an arbitrary character of the center, and, for $\fg=\fas$,
there is the one-parameter family $T_\lambda$ of $(4|4)$-dimensional
representations that differ from each other by the value on the
center.

\underline{B) $\fg=\fpsq(n)$, $\fsq(n)$, $\fpq(n)$, $\fq(n)$}. In
these cases, $\fg_{ss}=\fsl(n)$.

The odd roots are $\alpha=\eps_i-\eps_j$. Let, as above, $\lambda$
denote a~weight of $W_{\ev}$, let $\lambda'=\lambda+\alpha$ run over
the weights of desirable form from $V_{\od}$, let
$\lambda''=\lambda'+\alpha$ be the weights of desirable form in~
$V_{\ev}$. Then,
\begin{equation}
\label{bezN}
\begin{tabular}{|c|c|c|}
\hline
$\lambda$&$\lambda'$&$\lambda''$\\
\hline
$\eps_1,\dots,\eps_n$ &$\eps_1,\dots,\eps_n$
&$\eps_1,\dots,\eps_n$ \\
\hline
$-\eps_1,\dots,-\eps_n$ &$-\eps_1,\dots,-\eps_n$ &
$-\eps_1,\dots,-\eps_n$ \\
\hline
\end{tabular}
\end{equation}
Thus, $ \fg_{\od}^2 W_{\ev}=W_{\ev}$ and $\dim \fg_{\od} W_{\ev}=n$.
Hence, $\sdim V=(n|n)$. But, for $\fg=\fpsq(n)$ and $\fpq(n)$, there
are no such representations at all (as follows from Fact
\ref{fact}), while for $\fg=\fsq(n)$ and $\fq(n)$ this is the
tautological representation which has no rank-1 operators.

\underline{C) $\fg=\fsl(n|m)$, $\fgl(n|m)$, $\fpsl(n|n)$.} Observe
that, for $n=m=1$, the Lie superalgebra $\fgl(1|1)$ is solvable;
hence, all its irreducible representations are induced \cite{S0}; the
case considered in Step 4, see Section~\ref{Sstep4}. Due to the isomorphism $\fsl(2|1)\simeq 
\fvect(0|2)$, we can and will consider this case in Step~5B, see Subsection~\ref{SS:6.9}.

\underline{$m>1$}. Then, $\fg_{ss}=\fsl(n)\oplus \fsl(m)$ and
$\fg_{\od}=\fg^{-1}\oplus \fg^1$, where the roots of $\fg^{-1}$ are
of the form $\alpha_{-1}=- \eps_i+\delta_j$ and those of $\fg^{1}$ are of the form
$\alpha_{1}=\eps_i- \delta_j$.

First, suppose \underline{$n>2$}. Denote by $\lambda'_{\pm 1}$ the
desirable weights of the form $\lambda+\alpha_{\pm 1}$ and by
$\lambda''_{\pm 1}$ the desirable weights of the form $\lambda'_{\pm
1}+\alpha_{\pm 1}$. We get:
\begin{equation}
\label{5.6.2}
\renewcommand{\arraystretch}{1.4}
\begin{tabular}{|c|c|c|c|c|}
\hline $\lambda$&$\lambda'_{-1}$&$\lambda'_{1}$&$\lambda''_{-
1}$&$\lambda''_1$\\
\hline
$\eps_1,\dots,\eps_n$ &$\delta_1,\dots,
\delta_n$&$\emptyset$&$\emptyset$&$\eps_1, \dots, \eps_n$\\
\hline
$-\eps_1,\dots,-\eps_n$ &$\emptyset$&$\delta_1,\dots,
\delta_n$&$-\eps_1,\dots,-\eps_n$ &$\emptyset$\\
\hline
\end{tabular}
\end{equation}
Thus, we see that $ \fg_{\od}^2 W_{\ev}=W_{\ev}$. Let us consider
the $\fg_{ss}$-module $\fg_{\od} W_{\ev}$ in more detail. First, let
$W_{\ev}$, as $\fg_{ss}$-module, be $\id\otimes\chi_0$ (the weights
are $\eps_1, \dots, \eps_n$). Then, as is clear from table~\eqref{5.6.2},
$\fg^1 W_{\ev}=0$ and $ \fg_{\od} W_{\ev}=\fg^{-1} W_{\ev}$; hence,
$\fg_{\od} W_{\ev}$ is isomorphic to a~quotient of the tensor
product
\[
\fg^{-1}\otimes W_{\ev} \simeq \left( \id\otimes
\id^*\right)\bigotimes\left(\id\otimes\chi_0\right).
\]
But the $\fsl(n)$-module $\id\otimes \id^*$ is the sum of the
adjoint module and the trivial $1$-dimensional module. Since the
$\fsl(n)$-action on $\fg^{-1}\otimes W_{\ev}$ is trivial, we deduce
that $\sdim \fg^{-1}\otimes W_{\ev}=(0|m)$ and $\sdim V=(n|m)$.

Therefore,
\[
T=\begin{cases} \id,&\text{ for } \fg=\fsl(n|m)\text{ and }
\fgl(n|n),\\
\id\otimes \chi,\text{ where $\chi$ is an
arbitrary character,}&\text{ for } \fg=\fgl(n|m).
\end{cases}
\]
For $\fg=\fpsl(n|n)$, there are no such representations, as follows
from dimension considerations and Fact \ref{fact}.

The case where $W_{\ev}=\id^*\otimes \chi_0$ is
considered in the same manner. We obtain:
\[
T=\begin{cases} \id^*,&\text{ for } \fg=\fsl(n|m)\text{ and }
\fgl(n|n),\\
\id^*\otimes \chi, \text{ where $\chi$ is an arbitrary
character,}&\text{ for } \fg=\fgl(n|m). \end{cases}
\]

\underline{$n=2$}. In this case, by studying only the weight diagram
we cannot conclude that either $\fg^{-1}W_{\ev}=0$ or
$\fg^{1}W_{\ev}=0$. We will consider this case in
Lemma~\ref{L5.6.3}.

\underline{$m=1$ and $n>2$}. Then, $\fg_{ss}\simeq\fsl(n)$ and
$\fg_{\od}\simeq\fg^{-1}\oplus \fg^1$, where the roots of $\fg^{-1}$
are $\alpha_{-1}=-\eps_i $ and the roots of $\fg^{1}$ are
$\alpha_{1}=\eps_i $. The list of desirable weights in this case is
obtained from table~\eqref{5.6.2} by replacing $\delta_1,\dots,\delta_n$
with $0$. The remaining arguments are the same as in the above
subcase.

\sssbegin[Technical, on $\fsl(n)$-modules]{Lemma}[Technical, on $\fsl(n)$-modules]\label{L5.6.3}
Let the Lie superalgebra $\fh=\fh_{\ev}\oplus\fh_{\od}$ be such that
$\fh_\ev\simeq\fsl(n)$ and $\fh_\od\simeq(\id_{\fsl(n)})^*$, as
$\fsl(n)$-module. Let $V$ be an $\fh$-module and $W\subset V$ its
$\fsl(n)$-invariant subspace isomorphic to $\id_{\fsl(n)}$. Let,
moreover, $\fh_\od W$ be a~$1$-dimensional trivial $\fsl(n)$-module.
Then, $ \fh_{\od}^2 W=0$.
\end{Lemma}

\begin{proof} By hypothesis, $\fh_{\od}$ and $W$ can be considered
as dual linear spaces. Let $w_1, \dots, w_n$ be a~basis of $W$ and
$h_1, \dots, h_n$ the dual basis of $\fh_{\od}$. As already
noted, the $\fsl(n)$-module $\fh_{\od} W$ is the quotient of the
tensor product $\fh_{\od}\otimes W$ by the adjoint submodule.
The adjoint submodule is singled out as the kernel of the natural
bilinear form $B(h\otimes w)=h(w)$ on $\fh_{\od}\otimes W$. Then, the
generator $u$ of $\fh_{\od} W$ can be expressed as $u=h_1
w_1=\dots=h_n w_n$. For $i\neq j$, we get
\[
h_i u=h_i h_j w_j=-(-1)^{p(h_i)p(h_j)}h_j h_i w_j=0. \hskip3cm \qed
\]
\noqed\end{proof}

\parbegin[Particular cases of Lemma \ref{L5.6.3}]{Corollary}[Particular cases of Lemma \ref{L5.6.3}]\label{cor5.6.3}
\emph{1)} The statement of Lemma $\ref{L5.6.3}$ holds also if $\fh_{\od}=\id$ and
$W=\id^*$.

\emph{2)} For $n=2$, the modules $\id$ and $\id^*$ are isomorphic,
so Lemma $\ref{L5.6.3}$ holds also if ${\fh_{\od}\simeq W\simeq \id}$
and ${\fh_{\od}\simeq W\simeq \id^*}$.
\end{Corollary}

\ssec{Continuation of the story} Thus, let $\fg=\fsl(2|m)$,
$\fgl(2|m)$ for $m>2$ or let $\fg$ be any almost simple Lie superalgebra related
to $\fpsl(2|2)$, except for $\fder~ ~\fpsl(2|2)\simeq \fosp_0(4|2)$, see eqs.~\eqref{exSeqOsp4,2,a}. Then,
$\fg_{ss}=\fsl(2)\oplus\fsl(m)$, and it remains to consider the case
where $\fG_1=\fsl(2)$, \textit{i.e.}, the subspace $W_\ev\subset V$ is
isomorphic, as $\fg_{ss}$-module, to $\id\otimes \chi_0$.

Observe that the submodules $\fg^{-1}W_{\ev}$ and $\fg^1 W_{\ev}$
cannot vanish simultaneously. Let, for example,
$\fg^{-1}W_{\ev}\neq 0$. Then, $\fg^{-
1}W_{\ev}=\chi_0\otimes\id_{\fsl(m)}$ and, thanks to Lemma
\ref{L5.6.3}, $(\fg^{- 1})^2W_{\ev}= 0$. But then
$\fg^{1}(\fg^{-1}W_{\ev})\neq 0$ due to irreducibility of the
$\fg$-module $V$ and our usual arguments based on the weight
diagram show that $\fg^{1}(\fg^{-1}W_{\ev})=W_{\ev}$. Moreover,
$W_{\ev}\simeq \id_{\fsl(2)}\otimes \chi_0$ arises as the quotient
module of the tensor product
\[
\begin{array}{l}
\fg^{1}\otimes(\fg^{-1}W_{\ev})\simeq \left (\id_{\fsl(2)}\otimes
\id_{\fsl(m)}^*\right)\bigotimes\left(\chi_0\otimes
\id_{\fsl(m)}\right)= 
\id_{\fsl(2)}\bigotimes
\left(\id_{\fsl(m)}^*\otimes\id_{\fsl(m)}\right).
\end{array}
\]
Applying Lemma \ref{L5.6.3} once again we deduce that
$(\fg^{1})^2(\fg^{-1}W_\ev)= \fg^{1}W_\ev=0$, \textit{i.e.}, ${\sdim V=(2|
m)}$ and $T(\fg)=\fsl(2| m)$ or $\fgl(2|m)$.

Let $\fg=\fosp(4|2n)$ for $n>1$. Then,
\[
\begin{array}{l}
\fg_\ev=\fo(4)\oplus\fsp(2n)\simeq 
\fsl(2)\oplus\fsl(2)\oplus\fsp(2n),\\
\fg_\od\simeq \left(\id_{\fsl(2)}\otimes\id_{\fsl(2)}\right) \bigotimes\id_{\fsp(2n)}.
\end{array}
\]
Let $\fG_1=\fsl(2)$. Our routine arguments yield
\[
\fg_\od W_\ev\simeq \left(\chi_0(\fsl(2))\otimes\id_{\fsl(2)}\right)\bigotimes
\id_{\fsp(2n)}\text{ and } \fg_\od^2 W_\ev=W_\ev.
\]
This implies that $V=W_\ev\oplus \fg_\od W_\ev$ and $\sdim
V=(2|4n)$.

Let us show that there are no such $\fg$-modules. To this end, in the
space dual to the Cartan subalgebra of $\fo(4)$, select basis
elements $\eps_1, \eps_2$; for the roots of $\fG_1$ take $\pm
2\eps_1$ and let the roots of the complementary to $\fG_1$
subalgebra in $\fo(4)$ be $\pm 2\eps_2$; let
$\delta_1,\dots,\delta_n$ be the standard basis of the dual space to
the Cartan subalgebra of $\fsp(2n)$. Then, the weight basis of
$\fg_\od W_\ev$ consists of the vectors $v_{\lambda}$ of weight $
\lambda=\pm\eps_2\pm\delta_i$, whereas the weight basis of $W_\ev$
consists of $v_{\pm\eps_1}$. Let $ \xi_{\pm\eps_1\pm
\eps_2\pm\delta_i}$ be an odd root vector of $\fg$. Then,
\[
v_{\pm\eps_1}= \xi_{\pm\eps_1\pm \eps_2\pm\delta_i} v _{\mp
\eps_2\mp\delta_i}\text{ and }\xi_\alpha v_\lambda=0\text{ if
}\alpha+\lambda\neq\pm\eps_1.
\]
For $i\neq k$, we get
\[
\begin{array}{l}
\xi_{\eps_1\pm \eps_2\pm\delta_i} v_{\eps_1}=0,\\
\xi_{-\eps_1\pm \eps_2\pm\delta_i} v_{\eps_1}=
\xi_{-\eps_1\pm \eps_2\pm\delta_i}(\xi_{\eps_1\pm
\eps_2\pm\delta_k}v_{\mp\eps_2\mp\delta_k }) 
=- \xi_{\eps_1\pm \eps_2\pm\delta_k}\cdot \xi_{-\eps_1\pm
\eps_2\pm\delta_i}v_{\mp\eps_2\mp\delta_k }=0.
\end{array}
\]
The penultimate equality holds because $[\xi_{\eps_1\pm
\eps_2\pm\delta_k}, \ \xi_{-\eps_1\pm \eps_2\pm\delta_i}]=0$ 
since $\pm2\eps_2\pm(\delta_k+\delta_i)$ is not
a~root of $\fg$. Thus, $\fg_\od v_{\eps_1}=0$. Similarly, $\fg_\od
v_{-\eps_1}=0$, \textit{i.e.}, $\fg_\od W_\ev=0$ and $V=W_0$. This is a~contradiction.

\underline{D) $\fg=\fab(3)$}. Then,
$\fg_\ev=\fg_{ss}=\fsl(2)\oplus\fo(7)$. For this case, we modify our
notation somewhat. 
Let the roots of $\fsl(2)$
be $\pm 2\eps$.
Then, the weights of $W_\ev$ are $\lambda=\pm\eps$; the odd roots are
$\alpha=\pm\eps+ \nfrac12(\pm\delta_1\pm\delta_2 \pm\delta_3) $ and, as $\fg_\ev$-module,
$\fg_\od\simeq \id_{\fsl(2)}\otimes \spin_7$.

The desirable weights are
\[
\lambda '= \lambda+\alpha=\nfrac12(\pm\delta_1\pm\delta_2 \pm\delta_3)
\text{ and } \lambda ''=\lambda'+\alpha=\pm\eps.
\]
Therefore, $\fg_{\od}^2 W_\ev=W_\ev$ and $V=W_\ev\oplus \fg_\od
W_\ev$. Note that $\fg_\od W_\ev$ is the quotient of the
$\fg_\ev$-module
\[
\fg_\od \otimes W_\ev \simeq 
\left(\id_{\fsl(2)}\otimes\id_{\fsl(2)}\right)\bigotimes \left(\spin_7
\otimes\chi_0\right).
\]

The dimension of the trivial $\fsl(2)$-submodule in the $\fsl(2)$-module
$\id\otimes\id$ is equal to 1 and, since the $\fsl(2)$-action in
$\fg_\od W_\ev$ should be trivial, we deduce that $\dim \fg_\od
W_\ev= 1\cdot \dim(\spin_7)=8$, \textit{i.e.}, $\sdim V=(2|8)$, but $\fg=\fab(3)$ has no such irreducible
representations, see Fact
\ref{fact}.

\underline{E) $\fg=\fag(2)$}. Then,
$\fg_\ev=\fg_{ss}=\fsl(2)\oplus\fg(2)$ and, as $\fg_\ev$-module, $\fg_\od\simeq 
\id_{\fsl(2)}\otimes R(\pi_1)$, where
$R(\pi_1)$ is the irreducible representation of $\fg(2)$ with 
highest weight $\pi_1$, see \cite{OV}. Recall that $\dim
R(\pi_1)=7$. Arguments similar to those in the above case show
that the desirable representation should be of superdimension
$(2|7)$; but $\fag(2)$ has no such irreducible representations, see
Fact \ref{fact}.

\underline{F) $\fg=\fosp_a(4|2)$}. Then, $\fg_\ev=\fg_{ss}=
\fsl(2)\oplus \fsl(2)\oplus \fsl(2)$ and arguments similar to those
used in the above two subcases yield that $\sdim V= (2|4)$; there
are no such irreducible $\fg$-modules for $\alpha\neq 1$, except for
the equivalent to $a=1$ values $a=-\nfrac12$, $-2$, see formulas \eqref{osp42symm1}.

\ssec{Step 5B. Vectorial Lie superalgebras. The four cases}
\label{SS:6.9} Let $\fh_0$ be a~finite-dimensional Lie superalgebra.
Let $z$ be either its outer derivation or a~generator of the trivial center. Set
$\fg_0=\fh_0\ltimes \Cee\cdot z$ (if $z$ is central, this sum is
direct).

\sssbegin[Cartan prolongs with the same negative part]{Lemma}[Cartan prolongs with the same negative part]\label{L5.7.0} 
Let $\fh_{-1}=\fg_{-1}$ be an irreducible finite-dimensional
$\fg_0$-module. Consider the two Cartan prolongs:
\[
\fh=\mathop{\oplus}\limits_{i\geq -1}\fh_i=(\fh_{-1}, \fh_0)_*
\text{ and } \fg=\mathop{\oplus}\limits_{i\geq -1}\fg_i=(\fg_{-1},
\fg_0)_*.
\]
Either $\fg_1=\fh_1$, and then
$\fg=\fd(\fh)$, or $\fg_1=\fh_1\oplus \fg_{-1}^*$, direct sum of
vector spaces.
\end{Lemma}

\begin{proof}
For any $F\in\fg_1$ and $v\in\fg_{-1}$, the decomposition
\[
[F, v]=\alpha_F(v)z+h,
\]
where $h\in\fh_0$, determines a~linear form $\alpha_F
\in\fg_{-1}^*$. By the definition of Cartan prolongation,
$\alpha_F=0$ if and only if $F\in \fh_1$. Hence, as
a~$\fg_0$-module, the space $\fg_1/\fh_1$ is a~quotient of the
$\fg_0$-module $\fg_{-1}^*$ which is irreducible together with
$\fg_{-1}$. This completes the proof of the lemma.
\end{proof}

\paragraph{\underline{1) $\fh_0=\fvect(0|n)$, where $n\geq 3$, and
$\fg_0=\fcvect(0|n):=\fh_0\oplus \Cee \cdot z$}} (Recall, 
see Table
~\eqref{8.2}, that $\fh_0$ has no outer derivations or
nontrivial central extensions.) The rank-1 operator with even
covector exists only for the following $\fh_0$-modules $\fh_{-1}$:
\[
\renewcommand{\arraystretch}{1.4}
\begin{tabular}{|llll|}
\hline $1a)$&$\fh_{-1}=\Pi(\Lambda(n)/\Cee\cdot
1)$&$\Longrightarrow$&$\fh=(\fh_{-1}, \fh_0)_*=\fle(n; n)$;\\
\hline
$1b)$&$\fh_{-1}=\Vol_{0}(0|n)$&$\Longrightarrow$&$\fh=(\fh_{-1},
\fh_0)_*=\fsvect'(1|n; n)$; \\
\hline
$1c)$&$\fh_{-1}=\Pi(\Vol_{0}(0|n))$&$\Longrightarrow$&$
\fh=(\fh_{-1}, \fh_0)_*=\fb_1'(n; n)$;\\
\hline
$1d)$&$\fh_{-1}=\Pi(\Vol^\lambda(n))$&$\Longrightarrow$&$
\fh=(\fh_{-1} ,\fh_0)_*=\fb_\lambda (n; n)$ for $\lambda\neq 0,
1$.\\
\hline\end{tabular}
\]

The element $z$ acts in $\fg_{-1}$ as a~non-zero scalar. We normalize
$z$ so that $\ad_z|_{\fg_{-1}}= -\id$ and use for such
$\fcvect(0|n)$-modules the same notation as for
$\fvect(0|n)$-modules.

\parbegin[Prolongations with $\fg_0=\fcvect(0|n)$]{Proposition}[Prolongations with $\fg_0=\fcvect(0|n)$]\label{prop5.7.1}
Let $\fg_0=\fcvect(0|n)$, where $n\geq 3$, let $\fg_{-1}$ be
a~$\fg_0$-module of type \emph{1a) -- 1d)} below. Let $\fg=(\fg_{-1},
\fg_0)_*$ be the Cartan prolongation. Then, in the respective cases,
we have:

\emph{1a)} if $\fg_{-1}=\Pi(\Lambda(n)/\Cee\cdot 1)$, then
$\fg=\begin{cases}\fv\fle(4|3),&\text{for $n=3$,}\\
\fd(\fle(n; n)),&\text{for $n>3$;}\end{cases}$

\emph{1b)} if $\fg_{-1}=\Vol_{0}(0|n)$, then
$\fg=\begin{cases}\fk\fas(; 3\eta),&\text{for
$n=3$,}\\
\fg=\fd(\fsvect'(1|n; n)),&\text{for $n>3$;}\end{cases}$

\emph{1c)} if $\fg_{-1}=\Pi(\Vol_{0}(0|n))$, then
$\fg=\fd(\fb_1'(n; n))$;

\emph{1d)} if $\fg_{-1}=\Pi(\Vol^\lambda(n))$, then
$\fg=\fd(\fb_\lambda (n; n))$ for $\lambda\neq 0,1$.
\end{Proposition}

\begin{proof}
Note that $\fg$ is bigraded: each component $\fg_k$ is
$\Zee$-graded by powers of $\xi$ (call it \textit{$\xi$-grading})\index{$\xi$-grading} so
that $\fg_0=\mathop{\oplus}_{-1\leq i\leq n-1}\ \fg_{0, i}$ and
$z\in\fg_{0, 0}$; we have $\fg_{-1}=\mathop{\oplus}_{k\leq i\leq m}\ 
\fg_{-1, i}$, where
\[
k=\begin{cases}0,&\text{in cases 1b) -- 1d)},\\
1,&\text{in case 1a)};\end{cases}\quad\text{~~ and~~}m=\begin{cases}n,&\text{in cases 1a), 1d)},\\ 
n-1,&\text{in cases
1b), 1c).}\end{cases}
\]

By Lemma \ref{L5.7.0}, if $\fg\neq \fd(\fh)$, then the subspace
$\fg_{-1}^*\subset \fg_1$ contains the components of degree from $-m
$ to $-k$ with respect to the $\xi$-grading. In particular, $\fg_{1,
-m}\neq 0$.

But, in each case, we have a~subspace $W=\mathop{\oplus}
_{k\leq i\leq s}\ \fg_{-1, i}\subset\fg_{-1}$ such that $\fg_0$ does
not contain operators commuting with the whole $W$. By Lemma
\ref{L4.2}, this means that $\fg_{1,i }=0$ for all~ $i<-1-s$, \textit{i.e.},
the necessary condition for $\fg\neq \fd(\fh)$ is $-m\geq -1-s$, \textit{i.e.},
$m\leq 1+s$. Specifically,
\be\label{cases1}
\renewcommand{\arraystretch}{1.4}
\begin{tabular}{|llll|}
\hline $ 1a)$&$s=2$, $m=n$&$\Longrightarrow$&$n\leq 3$;\\
\hline $ 1b)$, $1c)$&$s=1$, $m=n-1$&$\Longrightarrow$&$n\leq 3$;\\
\hline $1c)$&$s=1, m=n$&$\Longrightarrow$&$n\leq 2$. \\
\hline
\end{tabular}
\ee
In cases 1a), 1b) for $n=3$, we get the exceptional simple Lie superalgebras
$\fv\fle$ and $\fk\fas$.

The case 1c) for $n=3$ calls for subtler arguments. We have
$\fg_{-1}=\mathop{\oplus}_{i=0, 1, 2}\ \fg_{-1, i}$, where
\[
\begin{array}{l}
\fg_{-1, 0}=\Cee\Pi(1),\; \\
\fg_{-1, 1}= \Pi(\Span(\xi_j\mid j=1, 2, 3)),\; \\
\fg_{-1, 2}= \Pi(\Span(\xi_j\xi_l\mid j, l=1, 2, 3)).\end{array}
\]
Suppose $\fg\neq \fd(\fh)$ and let
$F=\Pi(\xi_1\xi_2)^*\subset\fg_1$, \textit{i.e.},
\[
[F, \Pi(\xi_1\xi_2)]=z+h, \text{ where } h\in\fg_{0, 0}, \text{
\textit{i.e.}, } h=\mathop{\sum}\limits_{i, j}a_{ij}\xi_i\partial_{\xi_j}.
\]
Then, $F\in \fg_{1, -2}$, and thus, $[F, \Pi(1)]=0$. By the Jacobi
identity,
\[
0=[[F, \Pi(\xi_1\xi_2)], \Pi(1)]=(z+h)\Pi(1)=-\Pi(1)+h(\Pi(1)).
\]
Therefore, $h(\Pi(1))=\Pi(1)$ and
\begin{equation}
\label{5.7.1.1} 
 \Div h=-\mathop{\sum}\limits_{1\leq i\leq 3} a_{ii}=1. 
\end{equation}
Since $p( \Pi(\xi_1\xi_2) )=\od$, we have
\[
0=[[F, \Pi(\xi_1\xi_2)], \Pi(\xi_1\xi_2)]=(z+h)
\Pi(\xi_1\xi_2)=-
\Pi(\xi_1\xi_2)+h(\Pi(\xi_1\xi_2)).
\]
Hence, $h( \Pi(\xi_1\xi_2) )= \Pi(\xi_1\xi_2)$ and
$a_{11}+a_{22}+\Div h=1$ or
\begin{equation}
\label{5.7.1.2} a_{11}+a_{22}=0. 
\end{equation}

Clearly, $[F, \Pi(\xi_3)]\in\fg_{0,-1}$, \textit{i.e.},
$[F, \Pi(\xi_3)]=\sum b_j\partial_{\xi_j}$.
By the Jacobi identity
\[
\begin{array}{ll}
0=&[[F, \Pi(\xi_1\xi_2)], \Pi(\xi_3)]-[\Pi(\xi_1\xi_2), [F,
\Pi(\xi_3)]]= \\
&-\Pi(\xi_3)+h(\Pi(\xi_3))-(b_1\Pi(\xi_2)- b_2\Pi(\xi_1)).
\end{array}
\]
Hence, $a_{33}+\Div h=1$, \textit{i.e.}, $a_{33}=0$. The last equality
contradicts conditions~\eqref{5.7.1.1} and ~\eqref{5.7.1.2}.
\end{proof}

\paragraph{\underline{2) $\fh_0=\fsvect(0|n)$ for $ n\geq 3$}} A rank-1
operator with even covector exists only for the $\fh_0$-module 
$\fh_{-1}=\Pi(T^0_{0}(\vec 0))$, hence, $\fh=(\fh_{-1},
\fh_0)_*=\fsle'(n; n)$.

From Table ~\eqref{8.2} we know that $\fsvect(0|n)$ has no nontrivial
central extensions, but it has an outer derivation, $D=\sum
\xi_i\partial_{\xi_i}\in \fvect(0|n)$. So we have to consider the
following cases:
\be\label{cases2}
\begin{tabular}{|l|l|}
\hline 
$ 2a)$&$\fg_0=\fh_0\oplus \Cee z$, $\fg_{-1}=\fh_{-1}$ as 
spaces and $\fg_{-1}=\fh_{-1}\otimes (-1)$ as $\fg_0$-module;\\

$ 2b)$&$\fg_0=\fh_0\ltimes \Cee D$, $\fg_{-1}=\fh_{-1}$ as 
spaces, where $D$ acts as $\sum \xi_i\partial_{\xi_i}+\lambda\cdot
1$;\\

$ 2c)$&$\fg_0=(\fh_0\oplus \Cee z)\ltimes \Cee D$, $\fg_{-1}=\fh_{-1}$
as spaces,\\
& where $z$ acts as multiplication by $-1$ and $D$
as $\sum\xi_i\partial_{\xi_i}$;\\

$ 2d)$& $\fg_0=\fh_0\ltimes \Cee(az+b D)$ for $a, b\in\Cee$ and
$\fg_{-1}=\fh_{-1}$.\\
\hline 
\end{tabular}
\ee

\parbegin[Prolongations of relatives of $\fh_0=\fsvect(0|n)$]{Proposition}[Prolongations of relatives of $\fh_0=\fsvect(0|n)$]\label{prop5.7.2}
Let $(\fg_{-1}, \fg_0)$ be one of the pairs $2a)$ -- $2c)$ of cases $~\eqref{cases2}$ and
$\fg=(\fg_{-1}, \fg_0)_*$ its Cartan prolong. Then, in case $2b)$ for
$n=3$ and $\lambda=0$, we have $\fg=\fle(3)$. In all other cases,
$\fg=\fd(\fh)$.
\end{Proposition}

\begin{proof}
In all cases, $\fg_{-1}=\mathop{\oplus}_{1\leq i\leq n- 1}\ \fg_{-1,
i}$. By Lemma \ref{L5.7.0}, if $\fg\neq\fd(\fh)$, then $\fg_{1,
-n+1}\neq 0$. On the other hand,
$\fg_0=\mathop{\oplus}_{-1\leq i\leq n- 2}\ \fg_{-1, i}$ and $\fg_0$
contains no elements that commute with the subspace $W=\fg_{-1,
1}\oplus\fg_{-1, 2}$. Hence, by Lemma \ref{L4.2}, $\fg_{1, s}=0$ for
all $s$ such that $s+2<- 1$. Thus, $-n+1\geq -3$, \textit{i.e.}, $n\leq 4$.

(Note that $\fsvect(0|3)\simeq \fspe(3)$, and therefore this case of 
Proposition is a~particular case of a~case already considered.)

Let us examine the case where $n=4$ more carefully. We have
$\fg_{-1}=\mathop{\oplus}_{1\leq i\leq 3}\ \fg_{-1, i}$, where
\[
\begin{array}{l}
\fg_{-1, 1}=\Span(\Pi(\xi_i))_{i=1}^4,\\
\fg_{-1, 2}=\Span(\Pi(\xi_i\xi_j))_{i=1}^4,\\
\fg_{-1, 3}=\Span(\Pi(\xi_1\xi_2\xi_3)\}.
\end{array}
\]
Suppose $\fg_1\neq \fh_1$. In the cases 2a) -- 2d) of table \eqref{cases2}, there should exist
elements $F\in\fg_{1, -3}$ and $\alpha z+\beta D\in \fg_0$ such that
\[
g=[F, \Pi(\xi_1\xi_2\xi_3)]=-\alpha\cdot 1+\beta D+h,
\]
where $\alpha, \beta\in\Cee$, $\alpha\beta\neq 0$
and $h\in\fh_{0, 0}=\fsl(4)$.

By grading considerations, $[F, \fg_{-1, 1}]=0$. Further, by the Jacobi
identity,
\[
h|_ {\fg_{-1,1}}=\alpha\cdot 1-\beta D|_{\fg_{-1,1}}=(\alpha-\beta)\cdot 1.
\]
Since $h\in\fsl(4)$, it follows that $\alpha-\beta=0$ and $h=0$.
Without loss of generality we can assume that $\alpha=\beta=1$. Then,
$\ad_g|_{\fg_{-1,2}}=-1+D|_{\fg_{-1, 2}}= 1$. Hence, $\fg_{0,
0}=\fgl(4)\simeq \fc\fo(6)$ and $\fg_{-1,
2}\simeq \id_{\fo(6)}\otimes 1$ as $\fg_{0, 0}$-modules. Note that $[\fg_{-1, 2}, \fg_{1,
-2}]\subset\fg_{0, 0}$. Set
\[
V=\{F\in\fg_{1, -2}\mid [F, \fg_{-1, 2}]=0\}.
\]
Then, the quotient $\fg_{1, -2}/V$ naturally embeds into the first
component of Cartan prolong $(\Pi(\id), \fc\fo(6))_*=\fd(\fh(0|6))$.
Due to Lemma \ref{L5.7.0} this contradicts the assumption
$\fg_1\neq \fh_1$.
\end{proof}

\paragraph{\underline{3) $\fh_0=\widetilde{\fsvect}(0|n)$ for $n\geq 3$}} A
rank-1 operator with even covector exists only for the
$\fh_0$-module
$\fh_{-1}=\nfrac{\Pi(\Vol(0|n))}{\Cee(b+\alpha\xi_1\cdots\xi_n)\vvol(\xi)}$.
Hence, with $\mu=\nfrac{\alpha}{b}$, we have
\be\label{fh}
\fh=\begin{cases}(\fh_{-1},\fh_0)_*=\widetilde{\fs\fb}_{\mu}(2^{n-1}-1|2^{n-1})&\text{for $n$ even,}\\
(\fh_{-1},\fh_0)_*=\widetilde{\fs\fb}_{\mu}(2^{n-1}|2^{n-1}-1)&\text{for $n$ odd}.\\
\end{cases}
\ee

From Table ~\eqref{8.2} we know that $\fh_0$ has no
nontrivial central extensions or outer derivations. 

\parbegin[Prolongs of relatives of $\fh_0=\widetilde{\fsvect}(0|n)$]{Proposition}[Prolongs of relatives of $\fh_0=\widetilde{\fsvect}(0|n)$]\label{prop5.7.3}
For $\fh$ as in formula~\eqref{fh}, let
$\fg_0=\fh_0\oplus \Cee z$ and $\fg_{-1}=\fh_{-1}$, as
$\fh_0$-modules; let $z$ act on $\fg_{-1}$ as $-\id$. Then, $(\fg_{-1}, \fg_0)_*=\fd(\fh)$.
\end{Proposition}

\begin{proof}
The components $\fg_k$ are filtered by powers of $\xi$.
Consider the corresponding bigraded Lie superalgebra $\gr \fg$. If
$\fg_1\neq \fh_1$, then, by Lemma \ref{L5.7.0},
\[
(\gr\fg)_1=(\gr\fh)_1\oplus(\gr\fg)_{-1}^*.
\]
But $(\gr\fg)_{-1}=\mathop{\oplus}_{0\leq i\leq 2n-1}\ \fg_{-1, i}$.
Hence, $(\gr\fg)_{1, -2n+1}\neq 0$.

On the other hand,
$(\gr\fg)_0=\mathop{\oplus}_{-1\leq i\leq 2n-2}\ \fg_{0,i}$ and
$\gr\fg_0$ contains no elements that commute with the subspace
$W=\mathop{\oplus}_{i=0, 1}\ \fg_{-1, i}$. Hence, $(\gr\fg)_{1,
s}=0$ for all~ $s$ such that $s+1<-1$. Therefore, $-2n+1+1\geq -1$ or
$2n\leq 3$. This contradicts the fact that $n\geq 2$.
\end{proof}

\paragraph{\underline{4) $\fh_0=\fh'(0|n)$ for $n\geq 5$ and related
subcases}}
As follows from the description of the irreducible representations
of $\fh'(0|n)$ and $\fh(0|n)$, see \cite{Sha1,Sha2}, no induced
representation has a~rank-1 operator. 

\section{Step 6. A particular case: case 3) of Corollary
\ref{SS:4.3}} \label{S:6pc}

Let $\fg_{-1}=\Lambda(\xi_1, \xi_2)$ or $\Pi(\Lambda(\xi_1,\xi_2))$;
let $\fg_0\subset\fhei(0|4)\ltimes\fo(4)$, where $\fo(4)\simeq\fsl_1(2)\oplus\fsl_2(2)$ with the isomorphic copies of $\fsl(2)$ numbered for convenience:
\[
\begin{array}{rcl}
\fhei(0|4)&=&\Span(1,
\xi_1, \xi_2, \partial_1, \partial_2), \text{~~where $\partial_i:=\partial_{\xi_i}$},\\
\fsl_{1}(2)&=&\Span(\xi_1\partial_2, \xi_2\partial_1, \xi_1\partial_1-
\xi_2\partial_2),\\
\fsl_{2}(2)&=&\Span(\partial_1\partial_2, \xi_1\xi_2,
\xi_1\partial_1+\xi_2\partial_2-1).\end{array}
\]

We realize the elements of $\fg_{-1}\oplus\fg_0$ by vector fields,
\textit{i.e.}, we embed $\fg_{-1}\oplus\fg_0$ into $\fvect (2|2)$. Redenote
the respective vectors of the monomial basis of
$\fg_{-1}=\Pi(\Lambda(2))$:
\[
\pi(\xi_1)\longmapsto \partial _x, \; \pi(\xi_2)\longmapsto \partial
_y; \ \ \pi(1)\longmapsto \partial _\eta, \;
\pi(\xi_1\xi_2)\longmapsto \partial _\zeta.
\]
To the basis vectors from $\fg_0$ we assign the following vector
fields of degree 0 that act in the space of functions in the indeterminates $x, y, \eta,
\zeta$:
\[
\footnotesize
\begin{matrix}\begin{matrix}
\fhei(0|4): &-1&\longmapsto& x\partial_x
+y\partial_y+\eta\partial_\eta+\zeta\partial _\zeta,\\
&\xi_1&\longmapsto& \eta\partial_x -y\partial _\zeta, \\
&\xi_2& \longmapsto& \eta\partial _y+x\partial _\zeta,\\
 &\partial_1&\longmapsto& \zeta\partial _y- x\partial_\eta,\\
 &\partial_2&\longmapsto& -\zeta\partial _x-
 y\partial_\eta,\end{matrix}\\
\begin{matrix}
\fsl_{1}(2): &\xi_1\partial_2&\longmapsto& -y\partial _x,&
\quad\fsl_{2}(2):
&\partial_2\partial_1&\longmapsto& - \zeta\partial _\eta,\\
&\xi_2\partial_1&\longmapsto&-x\partial_y,&
 &\xi_1\xi_{2}&\longmapsto&-\eta\partial_\zeta,\\
 &\xi_1\partial_1- \xi_2\partial_2&\longmapsto&y\partial
 _y-x\partial_x,& &\xi_1\partial_1+\xi_2\partial_2-
 1&\longmapsto&-\eta\partial_\eta-
 \zeta\partial_\zeta.\end{matrix}\end{matrix}
\]

We see that this realization of $\fg_0$ is invariant under the
reversal of parity in $\fg_{- 1}$. So, both $\fg_{-1}=\Lambda(2)$ and
$\fg_{-1}=\Pi(\Lambda(2))$ bring about the same result (isomorphic
algebras $(\fg_{-1}, \fg_{0})_{*})$.

Consider the following subalgebra $\fh_0$ of $\fg_0$:
\[
\begin{array}{ll}
\fh_0&=\fhei(0|4)\ltimes\fsl_{1}(2)\\
&=\Span(1, \xi_1, \xi_2)\ltimes\Span(\partial_1, \partial_2,
\xi_1\partial_2, \xi_2\partial_1,
\xi_1\partial_1-\xi_2\partial_2)\\
&\simeq (\Lambda(2)\setminus\Cee\cdot\xi_1\xi_2)
\ltimes\fsvect(0|2).\end{array}
\]
Consequently, 
\[
(\fg_{-1}, \fh_0)_*\simeq \fb_{\infty}(2; 2).
\]

Direct calculations show that
\[
\text{$(\fg_{-1}, \fg_0)_*\simeq (\fg_{-1},
\fh_0)_*\oplus\fsl_{2}(2)$ and $(\fg_{-1}, \fhei(0|4))_*\simeq \fg_{-
1}\oplus\fhei(0|4)$. }
\]
So a~necessary condition for simplicity, \textit{i.e.}, $[\fg_{-1},
\fg_1]=\fg_0$, implies that 
\[
\text{$\fg_0=\fh_0$ and $\fg \simeq \fb_{a,
a}(2; 2)\simeq \fb_{\infty}(2; 2)$.}
\]

\section{Step 7. $\Zee$-gradings of vectorial Lie
superalgebras}\label{S:8}

\ssec{Weisfeiler regradings}\label{SS:8.1} Let
$\fg=\mathop{\oplus}_{i\geq -d}\ \fg_i$ be a~$\Zee$-graded
vectorial Lie superalgebra. By a~(\textit{Weisfeiler}) \textit{regrading} \index{Regrading} 
of $\fg$ we understand a~ transition from $\fg$ to the Lie
superalgebra $\fh=\mathop{\oplus}_{i\geq -D}\ \fh_i$ isomorphic
to $\fg$ as abstract Lie superalgebra and such that
\begin{equation}
\begin{array}{ll}
a)&\text{$\fh$ is transitive, \textit{i.e.}, for any non-zero $x \in\fh_k$ where $k \geq 0$, there is $y \in\fh_{-1}$}\\
&\text{such that $[x, y]\neq 0$;}\\
b)&\text{$\fh_{\geq 0}:=\mathop{\oplus}_{i\geq 0}\ \fh_i$ is
a~maximal
Lie subsuperalgebra of finite codimension;}\\
c)&\text{the $\fh_0$-module $\fh_{-1}$ is irreducible};\\
d)&\text{$D<\infty$, \textit{i.e.}, $\fh$ is of finite depth.}
\end{array}
\end{equation}

In what follows in this section, $\fh$ denotes a~Lie superalgebra
obtained from $\fg$ by a~regrading. Clearly, not every two gradings
on the same algebra determine a~bigrading; fortunately, for our purposes it suffices
 to consider only \textit{compatible gradings}\index{Grading, compatible} for which~\eqref{7.1.1}
is well defined thanks to the \textbf{hypothesis}: \textit{all $\Zee$-gradings of the Cartan prolongations are induced by the degrees of indeterminates} proved for $p>0$ by Skryabin, see \cite{Sk91,Sk95}.\index{Conjecture} So, we consider only gradings induced by the degrees of the indeterminates.

The proofs of the next three Lemmas were first published in \cite{Sh14}; since the article was not put in arXiv, let us reproduce the proofs here.

\sssbegin[On the number of
homogeneous components of prolongs]{Lemma}[On the number of
homogeneous components of prolongs]\label{L7.1.1}
Let $\fh=\mathop{\oplus}_{i\geq -D}\ \fh_{i}$ be a~W-regrading of
the Lie superalgebra $\fg=\mathop{\oplus}_{j\geq -d}\ \fg_{j}$ that
determines a~bigrading:
\begin{equation}
\label{7.1.1} \fh_{i}=\mathop{\oplus}\limits_{m(i)\leq j\leq M(i)}\fh_{i,
j},\text{ where $\fh_{i, j}=\fh_{i}\cap \fg_{j}$. }
\end{equation}
Let $\sdim (\fh_{0}\cap \fg_{-})=0|k$. Then, the number
$M:=M(-1)-m(-1)+1$ of the summands $\fh_{-1, j}$ in the decomposition
~\eqref{7.1.1} of $\fh_{-1}$ does not exceed the number of
homogeneous components with respect to the $\Zee$-grading of
$\Lambda(\fh_{0}\cap \fg_{-})$ induced by the $\Zee$-grading of
$\fg_{-}$.

In particular, if $\fh_{0}\cap \fg_{-}$ is homogeneous with respect
to the $\Zee$-grading of $\fg$, \textit{i.e.}, $\fh_{0}\cap
\fg_{-}\subset\fg_{j}$ for some $j$, then $M\leq k+1$.
\end{Lemma}

\begin{proof}
The component $V = \fh_{-1,M(-1)}$ of the maximal degree is invariant with 
respect to the subalgebra $\fh_{0,\geq 0}:= \oplus_{j\geq 0}\ \fh_{0,j}$. Since the $\fh_0$-module $\fh_{-1}$ is irreducible, it must be a~quotient of the induced module: as spaces, 
\[
\text{ind}^{\fh_0}_{\fh_{0,\geq 0}}V = \Lambda(\fh_0\cap \fg_-)\otimes V. \hskip 4cm\qed
\]
\noqed\end{proof}

\sssbegin[When $\fh_{-1}$ is homogeneous]{Lemma}[When $\fh_{-1}$ is homogeneous]\label{L7.1.2}
If $\fh_0\cap\fg_-=0$, then
$\fh_0\subset\fg_0$ and $\fh_{-1}$ is homogeneous.
\end{Lemma}

\begin{proof}
The Lie superalgebra $\fh_{0,0}$ transforms $\fh_{-1,j}$ into itself and the operators in $\fh_{0,k}$
send $\fh_{-1,j}$ into $\fh_{-1,j+k}$. Therefore, if the representation of $\fh_0$ on $\fh_{-1}$ is irreducible, then $\fh_{-1}$ is homogeneous with respect to the grading of $\fg$, \textit{i.e.}, $\fh_{-1} = \fh_{-1,j_0}$ for some $j_0$. But then $\fh_{0,k}$ for $k > 0$ sends $\fh_{-1}$ to $0$. Since $\fh$ is transitive, $\fh_{0,k} = 0$ for all~ $k > 0$, \textit{i.e.}, $\fh_0 = \fh_{0,0}$.\end{proof}

\sssbegin[When the gradings coincide]{Lemma}[When the gradings coincide]\label{L7.1.3}
If $\fh_0\cap\fg_-=0$ and
there exists a~non-zero $x\in\fg_{-1}\cap \fh_-$, then the gradings
of $\fh$ and $\fg$ coincide.
\end{Lemma}

\begin{proof} Since $\fh_{-1} =\fh_{-1,j_0}\subset\fg_{j_0}$, it follows that $\fh_{-k}\subset\fg_{k\cdot j_0}$. If $x\in\fh_{-k}$, then $-1=k\cdot j_0$, implying $k = 1$ and $j_0 = -1$, \textit{i.e.}, $\fh_{-1}\subset\fg_{-1}$, and therefore $\fh_{-1} = \fg_{-1}$ and $\fh_i = \fg_i$ for all~ $i$.
\end{proof}

\sssec{How to describe all the W-regradings}\label{SS:7.2} Given a~W-graded Lie superalgebra $\fg$, construct its W-regraded Lie superalgebra $\fh$; namely, let~$\fh$ be a~Lie (super)algebra obtained from~$\fg$ by a~regrading. First of all, we should list all the regradings that are not forbidden; our task is to select the ones that can be realized.

1) Determine $\fh_{0}\cap\fg_{-1}$.

2) Construct a~``minimal'' (better, though informally, saying: ``most
tightly compressed'') regrading with the given intersection, \textit{i.e.},
such that preserves in $\fh_{-1}$ all the elements of $\fg_{-1}$
except for those that are eliminated by the condition on the
intersection.

3) If the ``minimal'' regrading is of Weisfeiler type, then with the
help of Lemmas \ref{L7.1.2} and \ref{L7.1.3} we prove that any other
W-regrading with the given intersection coincides with the minimal
regrading.

4) If the ``minimal'' regrading is not of Weisfeiler type, then with
the help of Lemma \ref{L7.1.1} we prove that there are no
W-regradings with the given intersection.

\begin{proof}[Proof \nopoint] of Theorem \ref{th7.3}. Let $1\leq r\leq m$, let $\varphi(r)=\{i_1, \dots , i_r\}$ be
a~subset of distinct indices in $\{1,\dots, m\}$. Permuting, if necessary, the
indeterminates we can and will assume that
$\varphi(r)=\{1, \dots , r\}$. Set $\deg\xi_1=\dots =\deg\xi_r=0$
and define the remaining degrees as in Tables~\eqref{clsfeq10} and
~\eqref{nonstandgr}. In Tables~\eqref{table3},~\eqref{table31}, \eqref{table32} and
~\eqref{table321}, the terms $\fh_{-1}$ and $\fh_0$ are indicated for
the regraded algebras~ $\fh$ for the serial and the exceptional
algebras. Now, consider various cases.

1) Let either $\fg=\fvect(n|m)=\fvect(x,\xi)$ or
$\fg:=\fsvect(n|m)=\fsvect(x,\xi)$. Then,
\[
\fg_{-}=\fg_{-1}=\Span(\partial_{x_{1}}, \dots , \partial_{x_{n}},
\partial_{\xi_{1}}, \dots ,\partial_{\xi_{n}}). 
\]
Since the depth of
$\fh$ if finite, $\deg x_i>0$ for all~ $i$. Hence,
$\partial_{x_{i}}\in\fh_{-}\cap\fg_{-1}$ for all~ $i$.

If $\deg\xi_j\neq 0$ for all~ $j$, then $\fh_0\cap\fg_-=0$ and $\fh=\fg$ by
Lemma \ref{L7.1.3}.

Let $\deg \xi_i=0$ for $i\leq r$ and $\deg \xi_i\not=0$ for
$i>r\geq0$. Denote the thus regraded superalgebras by $\widetilde \fg=\fvect(n|m;
\widetilde {\vec r})$ or $\widetilde \fg=\fsvect(n|m; \widetilde {\vec r})$. 
Then, \[\widetilde \fg_{-1}=
\Span(\partial_{x_{1}}, \dots , \partial_{x_{n}};
\partial_{\xi_{j}}\mid j>r)\otimes\Lambda (\xi_i\mid i\leq r).
\]
Let us show that
$\fh\simeq \widetilde \fg$.

Since $\partial_{x_{i}}\in \widetilde \fg_{-1}$ for all~ $i$, then by Lemma \ref{L7.1.3} it
suffices to check that $\widetilde \fg_{-1}\cap \fh_0=0$. By definition of $\fh$,
\[
\left(\Span(\partial_{x_{i}}, \partial_{\xi_{j}})\otimes 1\right)\cap\fh_0=0
\text{ for } j>r.
\]
For $i<r$, we have $\partial_{\xi_{i}}\in\fh_0$. Since $\fh_0$ is
a~subalgebra, the following implication holds:
\begin{equation}\label{implics}
\begin{minipage}[c]{12cm}
{\sl if some element of $\widetilde \fg_{-1}$ belongs to $\fh_0$, then all the
elements of $\Span(\partial_{x_{i}}, \partial_{\xi_{j}}\mid j>r)$
belong to $\fh_0$.}
\end{minipage}
\end{equation}
Therefore, $\widetilde \fg_{-1}\cap\fh_0=0$.

2) For the other series, take into account that the corresponding preserved (perhaps, conformally) form
($\omega_0$, $\omega_1$, $\alpha_0$, $\alpha_1$, $\widetilde
{\alpha_1}$, and so on) should be homogeneous.

For the series $\fh$ and $\fk$, it is convenient to split the
indeterminates into dual pairs: the even ones, $p_1, \dots , p_m$
and $q_1, \dots , q_m$ and the odd ones: $\xi_1, \dots , \xi_k$ and
$\eta_1, \dots, \eta_k$; the remaining odd coordinates are the self-dual
$\theta$'s, see eq.~\eqref{2.3.2'}. 
The corresponding modifications for contact series are obvious.

Observe that the cases of $\fk(1|2n; n)$ and $\fm(n; n)$ are obtained
from the general considerations given above by ``compressing" the
grading to eliminate the zero components. 

The case of the exceptional grading of $\fb_{\lambda}(2)$ was
considered in Subsection~\ref{SS:1.3.1}.
\end{proof}

\sssec{W-regradings of $\fk\fle(9|6)$}\label{ssWgradksle}
The W-regrading $\fk\fle(9|11)$ is the only one overlooked in \cite{Sh14}. Though the description of all possible regradings of $\fk\fle(9|6)$ in
\cite{CK2} is correct, only an outline of the idea and the final result are
given and one W-regrading of $\fk\fas$ was overlooked. Let us give a~detailed proof of the completeness of the
list of regradings of $\fk\fle(9|6)$. 

There are two ``least-codimension realizations'' for $\fk\fle(9|6)$.
It is more convenient to describe regradings starting with the
realization where $\fg_{\ev}=\fsvect(5)$ and $\fg_{\od}=\Pi(d\Omega^1(5))$, cf.
Subsection~\ref{SS:2.22}.

Set $\fg_{-2}=\Span(\partial_{1}, \dots, \partial_{5})$ and
$\fg_{-1}=\Span(dx_{i}dx_{j}\mid 1\leq i<j\leq 5)$.

Let $\fh$ be a~regrading of $\fg$. Then, clearly, $\fg_{-2}\subset \fh_{-}$
and $\fh_{0}\cap \fg_{-}\subset \fg_{-1}$.

Since, for any permutation $\sigma\in S_{5}$, we have
\[
[dx_{\sigma(1)}dx_{\sigma(2)},\
dx_{\sigma(3)}dx_{\sigma(4)}]=(-1)^{\sign(\sigma)}\partial_{x_{\sigma(5)}},
\]
it follows that $dx_{\sigma(1)}dx_{\sigma(2)}$ and
$dx_{\sigma(3)}dx_{\sigma(4)}$ cannot simultaneously appear in
$\fh_{0}$ after the regrading. Hence, $\dim \fh_{0}\cap \fg_{-}\leq 4$.

So we have to consider the following cases:

1) $\dim \fh_{0}\cap \fg_{-}=1$. In this case $\fh_{0}\cap
\fg_{-}\simeq\Span(dx_{1}dx_{2})$.

2) $\dim \fh_{0}\cap \fg_{-}=2$. In this case $\fh_{0}\cap
\fg_{-}\simeq\Span(dx_{1}dx_{2}, \ dx_{1}dx_{3})$.

3) $\dim \fh_{0}\cap \fg_{-}=3$. There are the two subcases:
\begin{equation}
\begin{array}{ll}
\text{3a)}&\fh_{0}\cap \fg_{-}\simeq\Span(dx_{1}dx_{2},\ dx_{1}dx_{3},\ 
dx_{1}dx_{4}),\\
\text{3b)}&\fh_{0}\cap \fg_{-}\simeq\Span(dx_{1}dx_{2},\ dx_{1}dx_{3},\ 
dx_{2}dx_{3}).
\end{array}
\end{equation}

4) $\dim \fh_{0}\cap \fg_{-}=4$. In this case, $\fh_{0}\cap
\fg_{-}\simeq\Span(dx_{1}dx_{2},\ dx_{1}dx_{3},\ dx_{1}dx_{4},\ 
dx_{1}dx_{5})$.

The standard arguments (either based on Lemma 7.1.3 from \cite{Sh14},
or see Subsection~\ref{SS:7.2}) show that in the cases 2) and 3a) $\fh$ is
not a~W-regrading, whereas in case 1) we obtain $\fk\fle(9|11)$,
in case 3b) we obtain $\fk\fle(11|9)$, and in case 4) we obtain
$\fk\fle(9|6)$. The corresponding W-gradings listed in \cite{CK2}
are:

1) $\deg x_{1}=\deg x_{2}=3$, \;\; $\deg x_{3}=\deg x_{4}=\deg
x_{5}=2$; \;\; $\deg d=-3$;

3b) $\deg x_{1}=\deg x_{2}=1$,\;\; $\deg x_{3}=\deg x_{4}=\deg
x_{5}=2$; \;\; $\deg d=-2$;

4) $\deg x_{i}=1$ for $i=1, \dots , 4$, \;\; $\deg x_{5}=2$; \;\;
$\deg d=-\frac32$.


\section{The deformations}\label{S:11}

For a~(more or less) conventional theory, see \cite{BLW, BLLS}; for the subtleties that arise over fields of positive characteristic, see \cite{BGL1}.\footnote{For Lie
\textbf{super}algebras, there is a~new type of (co)homology
--- the divided power one, depending on a~shearing vector ---
introduced in \cite{BGLL}.

If $p=2$, there are 2 types of (co)homology theories for Lie
algebras: with the usual (alternating) cocycles and of new type:
with non-alternating cocycles; accordingly, there are 6 types
of (co)homology theories for Lie superalgebras with non-zero odd part,
see \cite{LZ}. See also \cite{BGL1} with partial answers for $p=2$ and \textbf{open problems}\index{Problem, open} for $p=3$.}

\ssec{The tools to compute (co)homology}\label{SS:11.1} There
are not many
general theorems which help one to calculate Lie algebra
(co)homology. 
There
are less than a~half dozen general theorems plus the
Hochschild-Serre spectral sequence, see \cite{Fu}. Superization of
these theorems is immediate.

The simplest and most powerful theorem states that
\begin{equation}\label{powerTh}
\begin{minipage}[c]{12cm}
{\sl ``Any Lie superalgebra $\fg$ trivially acts on its own cohomology
$H^i(\fg; M)$ and homology $H_i(\fg; M)$ for any $i$.''}
\end{minipage}
\end{equation}

This theorem is particularly useful when $\fg$ is $\Zee$-graded and contains a~grading
operator. Most of the results on Lie algebra (co)homology are
obtained by scrutinizing the ``zero mode space" of the grading
operator in the space of cycles of the (co)homology to compute.

Observe also that 
\begin{equation}\label{3facts}
\begin{array}{l}
H^{\bcdot}(\fgl(1))=\Lambda^{\bcdot}(c_1)\text{~~for the element $c_1:=1^*\in(\fgl(1))^*$};\\
H^{n}(\fg_1\oplus\fg_2)\simeq
\mathop{\oplus}_{a+b=n}\ H^{a}(\fg_1)\otimes H^{b}(\fg_2),\\
H^{n}(\fg; M_1\oplus M_2)\simeq H^{n}(\fg; M_1)\oplus H^{n}(\fg;
M_2).
\end{array}
\end{equation}

If $\fg$ is a~finite-dimensional reductive Lie algebra, then $H^{n}(\fg; M)=0$ for any
finite-dimensional irreducible module $M$ and any $n$; the cohomology $H^{\bcdot}(\fg)$ with coefficients in the trivial module, are also known. The published
proofs known to us reduce the computation to compact forms of the
complex Lie groups with the given Lie algebra, but using the 
Hochschild-Serre spectral sequence and induction on rank one can
easily obtain the answer algebraically; the base of induction being $\fgl(1)$ and
$\fsl(2)$; for them, the computations can be easily performed by hand.

\sssec{(Co)homology with coefficients in the (co)induced
modules}\label{SS:11.2} Let $\fg$ be a~Lie superalgebra, let $\fh$ be its subalgebra, let $M$ be an $\fh$-module. We construct the induced
and coinduced $\fg$-modules from $M$ by setting
\begin{equation}
\label{h1} \Ind^\fg_\fh(M):=U(\fg)\otimes_{U(\fh)}M, \quad\quad
\Coind^\fg_\fh(M):=\Hom_{U(\fh)} (U(\fg), M).
\end{equation}
Clearly, $\Ind
^\fg_\fh(M)$ and $\Coind ^\fg_\fh(M)$ are $\fg$-modules.

For any $\fh$-module $M$, we have (\cite{Fu})
\begin{equation}
\label{h2} H^q(\fg; \Coind^\fg_\fh(M))\simeq H^q(\fh;
M),\quad\quad H_q(\fg; \Ind^\fg_\fh(M))\simeq H_q(\fh; M).
\end{equation}

In particular, for any vectorial Lie superalgebra $\cL$ considered with a
Weisfeiler filtration, the $L_0$-module $M$ such that $\cL_{1}M=0$,
and the $\cL$-module $T(M)$ of tensor fields (see Subsection~
\ref{SS:2.7}), we have
\begin{equation}\label{(2)} H^{\bcdot}(\cL; T(M))\stackrel{~\eqref{h2}}{\simeq}
H^{\bcdot}(\cL_0; M)\stackrel{\text{[Fu]}}{\simeq} H^{\bcdot}(L_0)\otimes
(H^{\bcdot}(\cL_1)\otimes M)^{L_0},
\end{equation}
where $N^{L_0}$ denotes the space of ${L_0}$-invariants of the ${L_0}$-module $N$.

For the proof of the second isomorphism in ~\eqref{(2)}, see Theorem~2.2.8
(due to Losik) in Fuchs' book \cite{Fu}. It is given there for
$\cL=\fvect(n|0)$ only, but it can be directly translated to any Lie
(super)algebra $\cL$ of vector fields. 
We use the
Hochshield-Serre spectral sequence for the algebra $\cL_0$ and its
ideal $\cL_1$; according to \cite{Fu}, $E_\infty=E_2$.

Anyway, we always obtain an estimate of $H^{i}(\cL_0; M)$ from above by 
just looking at
\begin{equation}\label{(4)}
\begin{array}{l}
E_2^{p,q}:=H^p(\cL_0/\cL_1; H^q(\cL_1; M))=H^p(L_0; (H^q(\cL_1;
M))^{L_0})\stackrel{\text{since $\cL_1M=0$}}{=} \\
\phantom{XX}=H^{p}(L_0)\otimes
(H^{q}(\cL_1)\otimes M)^{L_0}\text{~~for $p+q=i$}.
\end{array}\end{equation}

Thus, the computation of $H^2(\cL; \Coind^\cL_{L_0}(M))$ is reduced to computation (description) of the following spaces: 
\[
\begin{array}{ll}
(1)&\text{$H^p(L_0)$ for $p=0, 1, 2$,}\\
(2) &\text{$L_0$-modules $H^{q}(\cL_1)$ for $q=2, 1, 0$, see \cite{GLP},}\\
(3)&\text{$L_0$-invariants in the space $H^{q}(\cL_1)\otimes M$ for $q=2, 1, 0$.}\\
\end{array}
\]

Computing cohomology of Lie \textbf{super}algebras $\fg$ is a~more
involved business than computing cohomology of Lie algebras, even if $\fg$ is simple, because (except for $\fosp(1|2n)$ whose properties startlingly resemble those of $\fo(2n+1)$, as was first observed by Scheunert and expounded in \cite{DLZ}) there
are many nontrivial irreducible $\fg$-modules $M$ for which
$H^{n}(\fg; M)\neq 0$ for various values of $n$. Therefore, if possible, we have
to reduce the computation to studying cohomology of the reductive
part of $\fg_\ev$.

\subsection{Cohomology of finite-dimensional classical Lie
superalgebras}\label{S:14}
For completeness and references, we recall the results of calculation of the cohomology of the finite-dimensional simple Lie superalgebras and their
``close'' relatives. The known results with trivial
coefficients are announced in \cite{FuLe} (for proofs in certain cases,
see \cite{Fu} and \cite{Kz}; for a~summary with corrections of earlier claims, see \cite{BoKN}). Nothing is known about cohomology with
nontrivial coefficients, except for the case where the module of coefficients is
``typical'' and $\fg$ admits an invariant non-degenerate
symmetric bilinear form (in which case the arguments from Lie
algebra theory are applicable).
\begin{equation}\label{FuLe}
\renewcommand{\arraystretch}{1.4}
\begin{array}{l}
H^{\bcdot}(\fgl(m|n))\simeq \begin{cases}
H^{\bcdot}(\fgl(m))\simeq\Cee[\xi_1, \xi_3, \dots, \xi_{2m+1}],&\text{if $m\geq n$},\cr
H^{\bcdot}(\fgl(n))\simeq\Cee[\xi_1, \xi_3, \dots, \xi_{2n+1}],&\text{otherwise.}\end{cases} \\
H^{\bcdot}(\fp\fgl(n|n))\simeq H^{\bcdot}(\fgl(n))\otimes\Cee[h_2],\; \text{ where
$h_2$ represents the extension}\\
0\tto\Cee\tto\fgl(n|n)\tto \fp\fgl(n|n)\tto 0.\\
H^{\bcdot}(\fsl(m|n))\simeq \begin{cases} H^{\bcdot}(\fsl(m))\simeq\Cee[\xi_3, \dots, \xi_{2m+1}],&\text{if
$m\geq n$},\cr
H^{\bcdot}(\fsl(n))\simeq\Cee[\xi_3, \dots, \xi_{2n+1}],&\text{otherwise.}\end{cases} \\
H^{\bcdot}(\fp\fsl(n|n))\simeq H^{\bcdot}(\fsl(n))\otimes\Cee[h_2]
\text{ for $n>2$, where $h_2$ represents the extension}\\
0\tto\Cee\tto\fsl(n|n)\tto \fp\fsl(n|n)\tto 0.\\
H^{\bcdot}(\fp\fsl(2|2))\simeq \Cee[y_2, h_2, x_2, \xi_3]/(x_2y_2-h_2^2).
\end{array}
\end{equation}
Observe that the odd cocycle $\xi_3$ in the bottom line of ~\eqref{FuLe}
was missed in the original calculations of \cite{FuLe} and the answer was
rectified when, at our request, A.~Shapovalov used Grozman's
\textit{SuperLie} (\cite{Gr}), see also \cite{Kor}. The
2-cocycles representing the central extensions in the 1st line of~\eqref{FDcoc} are:
\[
\renewcommand{\arraystretch}{1.4} 
\mat{ A_1&B_1\cr C_1&D_1}, 
\mat{ A_2&B_2\cr C_2&D_2}
\longmapsto \begin{cases}
y_2: &\tr C_1C_2,\cr 
h_2:& 
\tr(B_1C_2+C_1B_2),\cr
x_2:& \tr B_1B_2.\cr 
 \end{cases}
\]
Clearly, the \textit{degrees} of the cocycles are:
\[
 \deg y_2=-2,\quad \deg h_2=0, \quad\deg x_2=2,
\]
and therefore they are precisely the weights of the basis elements of $\fsl(2)$.

Having defined an $\fsl(2)$-structure on the space of these
cocycles, we cook from it, together with $\fp\fsl(2|2)$, the superalgebra $\fosp_a(4|2)$ for $a=0$, see eq.~\eqref{exSeqOsp4,2,a} and
\cite{BGLL1}. The other results of \cite{FuLe} we might need are as
follows: 
\begin{equation}\label{HofSPE}
\begin{array}{l}
H^{\bcdot}(\fsl(1|1))\simeq \Cee[\xi_1, \eta_1]/(\xi_1\eta_1).\\
H^{\bcdot}(\fpe(n))\simeq \Cee[\xi_{1}, \xi_{5}, \dots , \xi_{4k+1\leq 2n+1};
\xi_{2n-1}]\text{~~for $n>1$}.\\
H^{\bcdot}(\fspe(2n+3))\simeq H^{\bcdot}(\fpe(2n+3))\text{~~for $n\geq 0$}.\\
H^{\bcdot}(\fspe(2n+4))\simeq H^{\bcdot}(\fpe(2n+4))\otimes\begin{cases}
\Kee[x_{n+2}],&\text{if $n\equiv 0\pmod 4$},\cr
\Kee[\xi_{n+2}],&\text{if
$n\equiv 2\pmod 4$},\end{cases}\text{~~for $n\geq 0$}.\\
H^{\bcdot}(\fspe(2))\simeq \Kee[x_1, y_2]/(x_1^2-y_2).
\end{array}
\end{equation}
In particular, only $\fspe(4)$ has a~non-trivial central extension, see Subsection~\ref{ssas} and \cite{BGLL1}.
\begin{equation}\label{HofOSP}
H^{\bcdot}(\fosp(m|2n))\simeq\begin{cases}
H^{\bcdot}(\fo(m)),&\text{for $m>2n$},\cr
H^{\bcdot}(\fsp(2n)),&\text{otherwise}.\end{cases}
\end{equation}
For cohomology of series $\fq$ and its
relatives, see \cite{BoKN}; we will not need it in this paper.

The cohomology of exceptional simple Lie superalgebras with trivial
coefficients are computed by a~remarkably ingenious method due to
Gruson \cite{Gr1}--\cite{Gr3}:
\[
\begin{array}{l}
H^{\bcdot}(\fa\fg(2))\simeq H^{\bcdot}(\fg(2)),\\
H^{\bcdot}(\fa\fb(3))\simeq H^{\bcdot}(\fo(7)),\\
H^{\bcdot}(\fosp_a(4|2))\simeq H^{\bcdot}(\fo(4))\simeq
H^{\bcdot}(\fsl(2)\oplus\fsl(2)).\end{array}
\]
Gruson's method has, however, limited applicability (at least, at the moment): sometimes we need various \textit{nontrivial} coefficients in which case Gruson's method does not work. To calculate (co)homology of the exceptional superalgebras by
hands is hardly possible, but fortunately, for every given algebra
and module, we have Grozman's \textit{SuperLie} code package, see
\cite{Gr}.

\ssec{Preliminary calculations}\label{SS:11.5} Let $\fg$ be a~Lie
superalgebra and let
\begin{equation}
\label{h3} 0\tto A\buildrel {\partial _0}\over \tto C\buildrel
{\partial _1}\over \tto B\tto 0, \; \; \text{~~where~~}\; \; \;
p(\partial _0)=\bar{0} \; \; \text{ and }\; \; \; \partial _1 \;
\; \text{ is either even or odd}, 
\end{equation}
be a~short exact sequence of $\fg$-modules. Let $d$ be the
differential in the standard cochain complex of $\fg$, cf.
\cite{Fu}.

Consider the following long cohomology sequence induced by the sequence ~\eqref{h3}:
\begin{equation}
\label{h4} \dots\buildrel {\partial }\over \tto H^i(\fg ;
A)\buildrel {\partial _0}\over \tto H^i(\fg ; C)\buildrel
{\partial _1}\over \tto H^i(\fg ; B)\buildrel {\partial }\over
\tto H^{i+1}(\fg ; A)
\buildrel {\partial _0}\over \tto \dots 
\end{equation}
where the $\partial _i $ are the differentials induced by the
namesake differentials in ~\eqref{h3}, and ${\partial =d\circ
\partial _1^{-1}}$. 

Since $ \partial _0$ and $\partial _1 $ commute
with $d$, the sequence ~\eqref{h3} is well-defined and the same
arguments as for Lie algebras (\cite{Fu}) demonstrate that the long
cohomology sequence ~\eqref{h4} induced by ~\eqref{h3} is exact
(\cite{LPS}).


 Let $\mathbbmss{1}$ be the 1-dimensional
trivial $\fg_0$-module and let $\mathbbmss{1}[k]$ be the $\fg_0$-module
trivial on the semi-simple part while the
distinguished central element acts on it by the scalar $k$. Then,
\begin{equation}
\label{h5}
\renewcommand{\arraystretch}{1.4}
\begin{array}{ll}
\fvect(m|n)=&\begin{cases}T(\id_{\fgl(m|n)}))&\text{ as $\fvect(m|n)$-module},\\
T(\id_{\fsl(m|n)}))&\text{ as $\fsvect(m|n)$-module};\end{cases}\\
\fk(2m+1|n)=&\begin{cases}T(\mathbbmss{1}[-2])&\text{ as $\fk(2m+1|n)$-module},\\
T(\mathbbmss{1})&\text{ as $\fpo(2m|n)$-module};\end{cases}\\
\fm(n)=&T(\Pi(\mathbbmss{1}[-2]))\text{ as $\fm(n)$\defis{} and $\fb_{a, b}(n)$-module};\\
\fpo(2m|n)=&T(\mathbbmss{1})\text{ as an $\fh(2m|n)$-module;}\\
\fb(n)=&T(\Pi(\mathbbmss{1})) \text{ as a~$\fle(n)$-, $\fsle(n)$-
and $\fsle'(n)$-modules.}
\end{array}
\end{equation}
We also denote by $\cF:=T(\mathbbmss{1})=T(\vec 0)$ the space of
functions. 

The isomorphism ~\eqref{(2)} yields
\be\label{LongF}
H^2(\cL; T(M))=H^2(L_0)\bigoplus(H^1(L_0)\otimes(H^1(\cL_1)\otimes M)^{L_0})\bigoplus
(H^0(L_0)\otimes(H^2(\cL_1)\otimes M)^{L_0}).
\ee

Denote for brevity: $\fv:= \fvect(m|n)$, $\fs:=
\fsvect(m|n)$.
Thus, we see that
\begin{equation}
\label{h6} H^i(\fs; \fv)\simeq H^i(\fsl(m|n);
\id_{\fsl(m|n)})=0 \text{ for $i=1, 2$ and $(m|n)\neq (1|1)$};
\end{equation}
\begin{equation}
\label{h7} H^i(\fsvect(m|n); \cF)\simeq H^i(\fsl(m|n))=0 \text{ for
}
i=1, 2.
\end{equation}

So, from Subsection~\ref{SS:11.1}, the above formulas
~\eqref{h5}---\eqref{h7}, and Section~\ref{S:14} we get (for $(m|n)$ such
that $\fh(2m|n)$ and $\fk(2m+1|n)$ are simple)
\begin{equation}
\label{4.4} H^i(\fh(2m|n); \fpo(2m|n))\simeq H^i(\fosp(n|2m);
\id_{\fosp(n|2m)})=0~\text{for
$i=1, 2, 3$}; 
\end{equation}
\begin{equation}
\label{4.5} H^2(\fle(n); \fb(n))\simeq H^2(\fpe(n);
\id_{\fpe(n)})=0\text{~~for
$n>1$};
\end{equation}
\begin{equation}
\label{4.6} H^i(\fsle(n); \fsb(n))\simeq H^i(\fspe(n);
\id_{\fspe(n)})=0~\text{for $i=1, 2$ and $n>1$}; 
\end{equation}
\begin{equation}
\label{4.7} H^2(\fpo(2m|n); \fk(2m+1|n))\simeq H^2(\fosp(n|2m);
\id_{\fosp(n|2m)})=0~\text{for $i=1, 2$}; 
\end{equation}
\begin{equation}
\label{4.8} H^i(\fb_{a, b}(n); \fm(n))\simeq H^i(\fspe_{a, b}(n);
\id_{\fspe_{a, b}(n)})=0\text{~~for
$n>1$ and $i=1,2$}. 
\end{equation}

Recall that for any (finite-dimensional) Lie algebra $\fg$ and any irreducible $\fg$-module $M$, we have
\begin{equation}
\label{H_0} 
H^0(\fg; M)\simeq \begin{cases}0&\text{if
$M$ is non-trivial},\\\Cee&\text{otherwise}.
\end{cases}
\end{equation}

Recall also that for any simple Lie superalgebra $\fg$, we have
\begin{equation}
\label{H_simple} 
\begin{array}{l}
H^0(\fg)= \Cee,\quad \quad
H^1(\fg)\simeq (\fg/[\fg, \fg])^*=0.
\end{array}
\end{equation}

Now, let us apply the above preliminary calculations and Subsection \ref{S:14}.

\ssbegin[Rigidity of simple vectorial Lie superalgebras]{Theorem}[Rigidity of simple vectorial Lie superalgebras]\label{ThRig} The simple infinite-dimensional serial vectorial Lie
superalgebras are rigid, except for $\fh:=\fh(2m|n)$ for $(m|n)\neq (1|2)$ with $H^2(\fh; \fh)\simeq \Cee$, $\fle:=\fle(n)$ for $n\neq 2$ with $H^2(\fle; \fle)\simeq \Pi\Cee$, and
$\fb_{a, b}(n)$ whose deformations, together with those of $\fh(2|2)$ and $\fle(2)$, are listed in Theorem~$\ref{thdefb}$ and Section~$\ref{S:15}$. The exceptional simple vectorial Lie superalgebras have no filtered deformations, see \cite{GLS2}.
\end{Theorem}

It is instructive to compare the above theorem with the description of deformations of Lie superalgebras and their desuperisations over fields of characteristic 2, see \cite{BGLLS}.

\begin{proof} \underline{$\fg=\fvect(m|n)$, or $
\fk(2m+1|n)$, or $\fm(n)$}. In these cases, we deduce from eq.~\eqref{LongF} with the help of \cite{GLP}, Subsection \ref{S:14}, and \cite[Table 5]{OV} that $H^2(\fg;\fg)=0$. Shmelev \cite{Sm2} (resp., Kotchetkov \cite{Ko1}) had already proved that $\fk(2n+1|m)$ (resp., $\fm(n)$) are rigid by painstaking calculations (repeated
application of the Hochshield-Serre spectral sequence) at the time the results of \cite{GLP} were not obtained yet.

\underline{$\fg=\fsvect(m|n)$}. Recall the abbreviations $\fv:= \fvect(m|n)$, $\fs:=
\fsvect(m|n)$.

$1\degree$: $m\neq 0, 1$. The short exact sequence of
$\fsvect(m|n)$-modules
\begin{equation}
\label{4.9} 0\tto \fs\tto \fv\stackrel {\Div}{\tto}\cF\tto 0
\end{equation}
gives rise to the long exact sequence
\begin{equation}
\label{4.10}
\renewcommand{\arraystretch}{1.4}
\begin{array}{l}
 \ldots\tto H^1(\fs; \fs)\tto H^1(\fs; \fv)\tto H^1(\fs; \cF)\tto
H^2(\fs; \fs)\tto H^2(\fs; \fv)\tto \ldots
\end{array}
\end{equation}
Note that, for any $(m|n)\neq (1|1), \ m+n> 2$, we have
\begin{equation}
\label{4.11} 
\begin{array}{l}
H^1(\fs; \cF)= H^1(\fsl(m|n))=0;\\
H^2(\fs; \fv)= H^2(\fsl(m|n); \id_{\fsl(m|n)})=0.
\end{array}
\end{equation}
This and ~\eqref{h7} imply that $H^2(\fs; \fs)=0$.

$2\degree$: $m=0$. In this case (considered in \cite{Ld}) the divergence is not a~surjective
mapping onto~ $\cF$, its image is the space $\Vol_0(0\vert n)$ \index{$\Vol_0(0\vert n)$} of functions with
integral 0:
\begin{equation}
\label{SinV} 0\tto \fs\tto \fv\stackrel {\Div}{\tto}\Vol_0\tto 0.
\end{equation}
In the associated long exact sequence
\begin{equation}
\label{SinVCoh} \ldots\tto H^1(\fs; \fs)\tto H^1(\fs; \fv)\tto
H^1(\fs; \Vol_0)\tto H^2(\fs; \fs)\tto H^2(\fs; \fv)\tto
\ldots
\end{equation}
we have $ H^1(\fs; \fv)= 0$ and $H^2(\fs; \fv)=0$ by ~\eqref{4.11}, so $
H^2(\fs; \fs)\simeq H^1(\fs; \Vol_0)$. To compute it, consider the short
exact sequence
\begin{equation}
\label{4.12} 0\tto \Vol_0\tto \cF\stackrel {\int}{\tto}\Cee\tto 0 
\end{equation}
which gives rise to the long exact sequence
\begin{equation}
\label{SinVCoh2}
\renewcommand{\arraystretch}{1.4}
\begin{array}{l}
0\tto H^0(\fs; \Vol_0)\tto H^0(\fs; \cF)\tto H^0(\fs)
\tto H^1(\fs; \Vol_0)\tto H^1(\fs; \cF)\tto \ldots
\end{array}
\end{equation}
and since
\[
\renewcommand{\arraystretch}{1.4}
\begin{array}{l}
H^0(\fs; \Vol_0)\simeq \Cee,\quad
 H^0(\fs; \cF)\simeq
H^0(\fsl(n))\simeq\Cee;\quad
H^1(\fs; \cF)\simeq H^1(\fsl(n))= 0, 
\end{array}
\]
it follows that $H^1(\fs; \Vol_0)\simeq H^0(\fs)\simeq\Cee$. Since we
already know one global deformation, there are no other deformations.

$3\degree$: $m=1$, \textit{i.e.}, $\fg=\fsvect'(1|n)$ for $n>1$. From the exact sequence
\[
0\tto \fsvect '(1|n)\tto
\fsvect (1|n)\tto \Cee
\cdot\xi_1\cdots\xi_n\partial_t\tto 0
\]
we derive the exact sequence (here $\fs'=\fsvect'(1|n)$ and
$\fs=\fsvect(1|n)$)
\begin{equation}
\label{SinVCoh3}
\begin{array}{l}
 \ldots\tto H^1(\fs'; \fs')\tto H^1(\fs'; \fs)\tto
H^1(\fs')\tto
H^2(\fs'; \fs')\tto \\
\tto H^2(\fs'; \fs)\tto H^2(\fs')\tto\ldots
\end{array}
\end{equation}
Since $H^1(\fs')=0$ by \eqref{H_simple} and $H^2(\fs')=0$, see Subsection~\ref{bezN2}, we have ${H^2(\fs'; \fs')\simeq H^2(\fs'; \fs)}$. To compute $H^2(\fs'; \fs)$ we use the exact sequence \eqref{SinV} of $\fs'$-modules leading to the exact sequence
\begin{equation}
\label{SprimInS}
0\tto H^0(\fs'; \fs)\tto H^0(\fs'; \fv)\tto H^0(\fs'; \Vol_0)
\tto H^1(\fs'; \fs)\tto H^1(\fs'; \fv)\tto \ldots
\end{equation}
Since $H^0(\fs'; \Vol_0)=0$ by \eqref{H_0} and $H^1(\fs'; \fv)\simeq H^1(\fsl(n); \id_{\fsl(n)})=0$, it follows that 
\[
H^2(\fs'; \fs')\simeq H^2(\fs'; \fs)=0.
\]


\underline{$\fg=\fb_{\lambda}(n)$}.

$1\degree$: $\lambda$ generic. From the short exact sequence of
$\fb_{\lambda}(n)$-modules
\begin{equation}
\label{4.11'} 0\tto \fb_{\lambda}(n)\tto \fm(n)\stackrel
{\Div_{\lambda}}{\tto}T(\vec 0)\tto 0, 
\end{equation}
where $\Div_{\lambda}$ is given by formula ~\eqref{2.7.3}, we derive
the long exact sequence
\begin{equation}
\label{4.14}
\renewcommand{\arraystretch}{1.4}
\begin{array}{l}
 \ldots\tto H^1(\fb_{\lambda}(n);
\fb_{\lambda}(n))\tto H^1(\fb_{\lambda}(n); \fm(n))\tto
H^1(\fb_{\lambda}(n); T(\vec 0))\tto\\
\tto H^2(\fb_{\lambda}(n);
\fb_{\lambda}(n))\tto H^2(\fb_{\lambda}(n); \fm(n))\tto
\ldots
\end{array}
\end{equation}
Observe that for $\lambda=\frac{2a}{n(a-b)}\in\Cee\cup\{\infty\}$,
we have
\begin{equation}
\label{4.15} H^1(\fb_{\lambda}(n); T(\vec 0))\simeq H^1(\fspe_{a,
b}(n))\simeq \Cee.
\end{equation}
Together with ~\eqref{4.8}, \textit{i.e.}, the fact that $H^1(\fb_{\lambda}(n); \fm(n))=0$ and $H^2(\fb_{\lambda}(n); \fm(n))=0$, this implies that $H^2(\fb_{\lambda}(n);
\fb_{\lambda}(n))\simeq \Cee$.

$2\degree$: the exceptional $\lambda=0$, 1, $\infty$ (and extra exceptional values $\lambda=\frac12$, $-\frac32$ if
$n=2$), see \S\ref{S:15}.


\underline{$\fg=\fh(2n|m)$}. Set $\fpo:=\fpo(2n|m)$, $\cF:=\cF(2n|m)$, $\fh:=\fh(2n|m)$, and $\fk:=\fk(2n+1|m)$.

\underline{$1\degree$: Generic case: $(n|m)\neq(1,2)$, $(n|m)\neq(0,m)$ for
$m>4$}. From the short exact sequence of $\fg$-modules
\begin{equation}
\label{4.15'} 0\tto \fpo(2n|m)\tto \fk(2n+1|m)\stackrel
{\partial_t}{\tto}\cF(2n|m)\tto 0 
\end{equation}
we derive the long exact sequence
\begin{equation}
\label{4.16}
\renewcommand{\arraystretch}{1.4}
\begin{array}{l}
 \ldots\tto H^1(\fpo; \fpo)\tto H^1(\fpo; \fk)\tto H^1(\fpo;
\cF)\tto
H^2(\fpo; \fpo)\tto\\
\tto H^2(\fpo; \fk)\tto H^2(\fpo; \cF)\tto
\ldots\end{array}
\end{equation}
From ~\eqref{4.7} we have $H^i(\fpo; \fk)=0$ for $i=1, 2$ and
$H^1(\fpo; \cF)\simeq H^1(\fpo; \fpo)\simeq \Cee$. Hence, $H^2(\fpo;
\fpo)\simeq \Cee$.

\sssbegin[Different cohomology for different function spaces]{Remark}[Different cohomology for different function spaces] The above is the shortest known to us proof of
the uniqueness of quantization of the Poisson Lie superalgebra
realized on the space of polynomials or formal series. Warning: for
other types of functions, say for Laurent polynomials or functions with compact support, the answer
is different, cf. \cite{KT1,KT}.
\end{Remark}


Let $\fh:=\fh(2n|m)$ and $\fpo:=\fpo(2n|m)$. From the exact sequence
\begin{equation}
\label{4.17} 0\tto \Cee\tto \fpo \tto\fh\tto 0
\end{equation}
we derive the long exact sequence
\begin{equation}
\label{4.18}
\renewcommand{\arraystretch}{1.4}
\begin{array}{l}\ldots\tto H^1(\fh)\tto H^1(\fh; \fpo)\tto H^1(\fh;
\fh)\tto
H^2(\fh)\tto H^2(\fh; \fpo)\tto \\
\tto H^2(\fh; \fh)\tto H^3(\fh)\tto
H^3(\fh; \fpo)\tto H^3(\fh; \fh)\tto
H^4(\fh)\tto\ldots\end{array}
\end{equation}
Let $\omega$ be the form preserved by $\fh$, let $k=\lfloor\nfrac{m}{2}\rfloor$ and 
\[
1_{n+k}\in\fgl(n|k)=\fgl(V)\subset\fh_0=S^2(V)\oplus \fgl(V)\oplus S^2(V^*)~~\text{for $m=2k$} 
\]
and similarly for $m=2k+1$. Note that $H^3(\fh)\simeq \Cee=\Cee[\omega_0^*\wedge 1_{n+k}^*]$, where 
\[
\omega_0^*\wedge 1_{n+k}^*=(\sum p_i^*\wedge q_i^* +\sum \xi_j^* \wedge \eta_j^*)\wedge (\sum p_iq_i +\sum \xi_j\eta_j)^*. 
\]
Then, from ~\eqref{4.4} we have $H^i(\fh; \fpo)=H^i(\fosp; \id_{\fosp})=0$ for
$i=1, 2, 3$, so $H^2(\fh; \fh)\simeq \Cee$. Clearly, this deformation of $\fh$ is the restriction of the quantization of~$\fpo$.

\underline{$2\degree$: $(2n|m)=(2|2)$}. See Section~ \ref{S:15} where this case is considered separately: the above argument fails, 
because $\fosp(2|2)\simeq\fvect(0|2)$:
\begin{equation}\label{4.18a}
H^3(\fh(2|2); \fpo(2|2))\simeq H^3(\fosp(2|2))=
H^3(\fvect(0|2))\simeq H^3(\fgl(2))=\Cee.
\end{equation}

\underline{$3\degree$: $(n|m)=(0|m)$ for
$m>4$}, see \cite{Tyut}. There is only one class of deformations. Clearly, this deformation of $\fh$ is the restriction of the quantization of~$\fpo$.

\underline{$\fg=\fle(n)$, $n>2$}. Let $\fb:=\fb(n)$ and $\fle:=\fle(n)$. From the exact sequence
\begin{equation}
\label{4.117} 0\tto \Cee\tto \fb(n) \tto\fle(n)\tto 0
\end{equation}
we derive (similarly to the case of Hamiltonian Lie
(super)algebras) the long exact sequence
\begin{equation}
\label{4.118}
\begin{array}{l}
\ldots\tto H^1(\fle)\tto H^1(\fle; \fb)\tto H^1(\fle;
\fle)\tto
H^2(\fle)\tto H^2(\fle; \fb)\tto \\
\tto H^2(\fle; \fle)\tto H^3(\fle)\tto H^3(\fle; \fb)\tto H^3(\fle;
\fle)\tto
H^4(\fle)\tto\ldots\end{array}
\end{equation}
From ~\eqref{4.5} we have $H^i(\fle; \fb)=H^i(\fpe; \id_{\fpe})=0$ for $i=1,
2, 3$. Since $H^3(\fle)\simeq \Pi\Cee=\Cee[\omega_1^*\wedge 1_{n}^*]$ by the same arguments as in the case where $\fg=\fh$, it follows that $H^2(\fle; \fle)\simeq \Pi\Cee$. Clearly, this deformation of $\fle$ is the restriction of the quantization of $\fb$, see Subsection~ \ref{SS:2.8}.


\underline{$\fg=\fsle'(n)$, $n>2$}. The answer is checked via computer for $n$ small: $\fsle'(n)$ is rigid. \end{proof}

\ssbegin[Rigidity of exceptional simple vectorial Lie superalgebras]{Conjecture}[Rigidity of exceptional simple vectorial Lie superalgebras]\label{CjRig} The exceptional simple vectorial Lie
superalgebras are rigid.\index{Conjecture}
\end{Conjecture}

\section{Deformations of $\fb_{\lambda}(n;n)$
(proof of claims in \cite{LSh3})}\label{S:15}

In the journal where \cite{LSh3} was published, the technical mathematical details were out of place. Let us give them here. Recall that $\fb_{\lambda}(n;n)=(\Pi(\Vol^{\lambda}), \fvect(0|n))_*$ and the $\fvect(0|n)$-action \eqref{LieDer} on 
the space of $\lambda$-densities. 
For $\lambda=\frac12$, the action preserves the
non-degenerate symmetric bilinear form
\[
(a\sqrt{\vvol}, b\sqrt{\vvol})=\int ab\ \vvol=(-1)^{p(a)p(b)}\int ba\ \vvol.
\]
Hence, on $\Pi(\sqrt{\Vol})$, the form 
\[
(\Pi(a\sqrt{\vvol}), \Pi(b\sqrt{\vvol}))=\int ab\ \vvol
\]
is non-degenerate and \textbf{anti}-symmetric.

\ssec{If $n=2$} Due to an isomorphism $\fosp(2|2)\simeq \fvect(0|2)$, we deduce that $\fb_{1/2}(2;2)\simeq \fh(2|2)$. The elements
of $ \fh(2|2)$ are given in terms of generating functions as $H_f$,
and our first step in the description of the main deformation of
$\fb_{\lambda }$ is a~realization of the elements of $\fb_{\lambda
}$ as deformations of the fields $H_f$, \textit{i.e.}, we construct a~map
\[
f\longmapsto H_f+\hbar(\lambda)W_f\in\fb_{\lambda }
\]
such that $\hbar(\lambda)\tto 0$ as $\lambda\tto \frac12$ and
$\fb_{\lambda }=\Span(H_f+\hbar(\lambda)W_f \mid f\in\cF)$.

Here is the description of the map
$\fvect(0|2)\tto\fvect(2|2)_0\simeq \fgl(2|2)$. 
Let a~basis in $\Pi(\Vol^{\lambda})$ be
\begin{equation}
\label{7.2} \vvol^{\lambda}, \quad \alpha \vvol^{\lambda}, \quad
\beta \vvol^{\lambda}, \quad \alpha\beta
\vvol^{\lambda},
\end{equation}
\textit{i.e.}, as space, $\Pi(\Vol^{\lambda})$ is isomorphic to the Grassmann algebra generated
by $\alpha$ and $\beta$.

Let $\fvect(2|2)=\fder~\Cee[p, q, \xi, \eta]$. Let us identify
$\Pi(\Vol^{\lambda})$ with $\fvect(2|2)_{-1}$ by setting
\begin{equation}
\label{7.3} \partial_{\xi}:=\vvol^{\lambda}, \quad
\partial_{p}:=\alpha \vvol^{\lambda}, \quad\partial_{q}:= \beta
\vvol^{\lambda}, \quad
\partial_{\eta}:=\alpha \beta \vvol^{\lambda}.
\end{equation}
Here is the representation of $\fvect(0|2)=\fder~\Cee[\alpha,
\beta]$ in terms of $p, q, \xi, \eta$:
\begin{equation}
\label{7.4}
\renewcommand{\arraystretch}{1.4}\begin{array}{ll}
\pder{\alpha}&=-p\pder{\xi}+\eta\pder{q},\\
\pder{\beta}&=-q\pder{\xi}-\eta\pder{p},\\
\alpha \pder{\beta}&=-q\pder{p},\\
\beta \pder{\alpha}&=-p \pder{q},\\
\alpha\pder{\alpha}&=\lambda\xi\pder{\xi}+(\lambda-1)p\pder{p}+\lambda
q\pder{q} + (\lambda-1) \eta\pder{\eta},\\
\beta \pder{\beta}&=\lambda\xi\pder{\xi}+\lambda p\pder{p}+
(\lambda-1)q\pder{q} + (\lambda-1)\eta\pder{\eta},\\
\alpha\beta\pder{\alpha}&=\lambda \xi \pder{q} +
(\lambda-1)p\pder{\eta},\\
\alpha\beta \pder{\beta}&=-\lambda \xi \pder{p}-
(\lambda-1)q\pder{\eta}.
\end{array}
\end{equation}

Looking at the explicit basis of $(\fb_{\lambda})_0$, we see that the
\textit{linear} operator 
\[
D=P\pder{p}+Q\pder{q}+X\pder{\xi}+Y\pder{\eta}\in\fvect(2|2)
\]
belongs to $(\fb_{\lambda})_0$ if it satisfies the
following system of differential equations (here $\lambda=\frac12+
h$)
\be\label{1-8}
\left\{\renewcommand{\arraystretch}{1.4}\begin{array}{ll}
\pderf{X}{\eta}=0&(1)\\
\pderf{Y}{\xi}=0&(2)\\
\pderf{X}{p}+(-1)^{p(D)}\pderf{Q}{\eta}=0&(3)\\
\pderf{X}{q}-(-1)^{p(D)}\pderf{P}{\eta}=0&(4)\\
\pderf{P}{p}+\pderf{Q}{q}=(-1)^{p(D)}(\pderf{X}{\xi}+\pderf{Y}{\eta})&(5)\\
(1+2h)\pderf{Y}{p}=(2h-1)(-1)^{p(D)}\pderf{Q}{\xi}&(6)\\
(1+2h)\pderf{Y}{q}=-(2h-1)(-1)^{p(D)}\pderf{P}{\xi}&(7)\\
(1+2h)(\pderf{P}{p}+\pderf{Q}{q})=4h (-1)^{p(D)}\pderf{X}{\xi}&(8)
\end{array}\right.
\ee
By dimensions considerations, any \textit{linear} operator $D$ satisfying conditions
(1)--(8) of \eqref{1-8} belongs to $(\fb_{\lambda})_0$. By the definition of the Cartan
prolongation, \textit{any} $D$ satisfying conditions (1)--(8) belongs to
$\fb_{\lambda }$.

The explicit solution of the system (1)--(8) of \eqref{1-8} (here we omit a~page of obvious
transformations) is given in terms of one function $f\in\Cee[p, q,
\xi, \eta]$:
\[
\left\{\renewcommand{\arraystretch}{1.4}\begin{array}{l}
P=-\pderf{f}{q}+P_0(p, q),\\
Q=\pderf{f}{p},\\
X=-(-1)^{p(f)}\pderf{f}{\eta},\\
Y=(-1)^{p(f)}\frac{2h-1}{1+2h}\pderf{f}{\xi},\\
\end{array}\right.
\]
where
\begin{equation}
\label{7.5}
 \pderf{P_0}{p}=-\frac{4h}{1+2h}\frac{\partial^2f}{\xi\eta}.
\end{equation}
The relation ~\eqref{7.5} determines $P_0$ up to a~term that 
depends only on $q$. It remains to observe that adding to $f$ such
a~term does not affecty the values of $Q$, $X$, $Y$. Therefore, we
can set
\begin{equation}
\label{7.6}
 P_0=-\frac{4h}{1+2h}\int^p_0\frac{\partial^2f}{\xi\eta}dp=
\frac{4h}{1+2h}\int^p_0\frac{\partial^2f}{\eta\xi }dp.
\end{equation}
Now let
\[
 \hbar=\frac{4h}{1+2h}=\frac{2\lambda -1}{\lambda }\tto 0
\iff \lambda \tto \frac12.
\]
To summarize:

\sssbegin[Deformation of the Buttin bracket]{Theorem}[Deformation of the Buttin bracket]\label{SS:15.1} Set
$\hbar(\lambda)=\nfrac{2\lambda-1}{\lambda}$. 

Then, $D\in\fvect(2|2)$
belongs to $\fb_{\lambda}(2; 2)$ if and only if
\begin{equation}
\label{7.7} 
 D=D_{f}:=H_{f}+\hbar(\lambda)W_{f},\text{ where
$W_{f}=\left (\int_{0}^p \pder{\eta}\pderf{f}{\xi}dp\right
)\partial_{p}+
(-1)^{p(f)}\pderf{f}{\xi}\partial_{\eta}$}
\end{equation}
for some $f\in \Cee[p, q, \xi, \eta]$. Then,
\begin{equation}
\label{7.8} [D_{f}, D_{g}]=D_{\{f, g\}}+\hbar(\lambda)D_{c(f, g)}
\text{~~for any $f, g\in
\Cee[p, q, \xi, \eta]$},
\end{equation}
where
\begin{equation}
\label{7.9}\footnotesize
\renewcommand{\arraystretch}{1.4}
\begin{array}{l}
 c(f, g)=-\pderf{f}{p}\int_{0}^p\pder{\eta}\pderf{g}{\xi}
dp+ \pderf{g}{p}\int_{0}^p\pder{\eta}\pderf{f}{\xi} dp+
 \xi\pder{\xi}\left(
(-1)^{p(f)}\pderf{f}{\xi}\pderf{g}{\eta}+
\pderf{f}{\eta}\pderf{g}{\xi}\right)\Big|_{p=0, q=0}+\\
+ \pder{\eta}\left(\int_{(0, q)}^{(p, q)}\left(
(-1)^{p(f)}\pderf{f}{p}\pderf{g}{\xi}-\pderf{f}{\xi}\pderf{g}{p}\right)dp+
\int_{(0, 0)}^{(0, q)}\left( (-1)^{p(f)}\pderf{f}{q}\pderf{g}{\xi}-
\pderf{f}{\xi}\pderf{g}{q}\right)\Big|_{p=0}dq\right).
\end{array}
\end{equation}
\end{Theorem}

\sssec{$\fh(2|2)$ has more deformations than $\fpo(2|2)$}\label{HvsPo} Note that the formula
\[
[H_f, H_g]_{new}=H_{\{f, g\}_{P.B.}}+\hbar (\lambda)\cdot H_{c(f, g)}
\]
defines a~deformation of $\fh(2|2)$ (which is the main
deformation of $\fb_{1/2}(2)$), but (and this agrees with an
unpublished preprint by Batalin and Tyutin) the formula
\begin{equation}
\label{7.10} \{f, g\}_{new}=\{f, g\}_{P.B.}+\hbar (\lambda)\cdot
c(f, g)
\end{equation}
does not define a~deformation of $\fpo(2|2)$ because the expression~\eqref{7.10}
does not satisfy the Jacobi identity.

\ssec{Deformations of $\fg=\fb_{1/2}(n; n)$}\label{SS:15.2}
Clearly, $\fg_{-1}$ is isomorphic to $\Pi(\sqrt{\Vol})$. Therefore,
there is an embedding
\begin{equation}
\label{7.11} \fb_{1/2}(n; n)\subset\begin{cases}
\fh:=\fh(2^{2k-1}|2^{2k-1}),&\text{
for $n=2k$},\\
\fle:=\fle(2^{2k}),&\text{ for $n=2k+1$.}\end{cases}
\end{equation}
It is tempting to determine quantizations of $\fg$ in addition to
those considered by Kotchetkov, as the composition of the embedding
~\eqref{7.11} and subsequent quantization.

For $n=2$, when ~\eqref{7.11} is not just an embedding, but an
isomorphism, this idea certainly works and we get the following extra
quantization of the antibracket described in Theorem~\ref{thdefb}:
we first deform the antibracket to the point $\lambda =\frac12$
along the main deformation, and then quantize it as the quotient of
the Poisson superalgebra. Let us show that this scheme fails to give new algebras for
$n=2k>2$.

\sssbegin[Deformations of $\fb_{1/2}(n; n)$]{Theorem}[Deformations of $\fb_{1/2}(n; n)$]\label{bezN3} For $n=2k>2$, the image of
$\fb_{1/2}(n; n)$ under the embedding ~\eqref{7.11} is rigid under the
quantization of the ambient.
\end{Theorem}

\begin{proof} Consider the Lie superalgebra $\fg:=\fb_{1/2}(n;
n)\subset\fm(n; n)$. Recall that, in this grading, the depth of $\fg$ is
equal to 1, and
\[
\fg_{-1}=\Pi(\sqrt{\Vol}),\qquad\fg_0=\fvect(0|n).
\]

Consider the case of embedding into $\fh$, see ~\eqref{7.11}; the
case of embedding into $\fle$ is similar. Denote the indeterminates
corresponding to $\Pi(\Lambda^{\bcdot}\sqrt{\vvol})\simeq \Pi(\oplus
\Lambda^i)$ as follows:
\[
\begin{array}{l}
\text{the ones corresponding to $i\leq k-2$ by $q=(q_1, \dots , q_N)$,}\\
\text{the ones corresponding to $i=k-1$ by $Q=(Q_1, \dots , Q_M)$,}\\
\text{the ones corresponding to $i=k$ by $x=(x_1, \dots , x_L)$,}\\
\text{the ones corresponding to $i=k+1$ by $P=(P_1, \dots , P_M)$,}\\
\text{the ones corresponding to $i>k+1$ by $p=(p_1, \dots , p_N)$}
.\end{array}
\]

The symplectic form determines the pairing (of the differentials) of
$p$ with $q$, $P$ with $Q$ and the $x$'s with themselves splitting
them into either two dual groups for $k=2l$ (call them $y$ and~ $z$),
or into three groups (two dual ones and one, consisting of one element $x_0$,
self-dual) for $k=2l+1$.

If all types of indeterminates are present, then $n\geq 4$.

Let $V=\mathop{\oplus}_{i\geq -1}\ V_i$ be the projection of
$\fg$ to the space of vector fields generated by polynomials in $p,
P, x, Q, q$. Of course, $V$ must not be a~subalgebra in $\fh$, but
it is invariant with respect to the partial derivatives $\pder{p},
\pder{P}, \pder{x}, \pder{Q}, \pder{q}$ because
$V_{-1}=\fg_{-1}=\fh_{-1}$.

Consider $V_0$. The $\fg_0$-action on $\fg_{-1}$ is such that
$\Lambda^s$ goes into $\mathop{\oplus}_{i\geq s-1}\ \Lambda^i$, so we deduce the following 1)--6):

1) The indeterminates $(p, P)$ can go only into $(p, P, x)$, but
never into $(Q, q)$. In turn, this means that $V_0$ does not contain
monomials of the form $q_iq_j$, $q_iQ_j$ and $Q_iQ_j$.

Therefore, the subspace $V_0$ contains only the monomials of degree
$\leq 1$ in $q, Q$.

2) The indeterminates $x$ can go only into $(p, P)$, into $x$, and
$Q$, but never into $q$. In turn, this means that $V_0$ can contain
monomials of the form $x_ix_j$ and $x_iQ_j$ but not $x_iq_j$.

Consider a~copy $\widetilde \fg$ of $\fg$, with the $\Zee$-grading
given by 
\[
\deg \xi_i=\deg \tau=1,\ \ \deg u_i=0 \text{~~for the remaining indeterminates $u_i$}.
\]
Consider the bigrading of $\fg$ given by
\[
\fg:=\mathop{\oplus}\limits_{i,
j\geq -1}\fg_{i, j}, \text{ where } \fg_{i, j}=\fg_{i}\cap
\widetilde \fg_{j}.
\]
The monomials of the form $x_{i}x_j$ enter the decomposion of elements
from $\fg_{0,0}$, whereas the monomials of the form $x_{i}Q_j$ enter
the decomposion of elements of $\fg_{0,-1}$.

Does $V_0$ contain degree-3 monomials in $x$? If yes,
they correspond to the elements of $\fg_{1, t}$, where $k-1+t=0$,
\textit{i.e.}, $k=1-t$, and hence $k\leq 2$ since $t\geq -1$.

Therefore, if $k>2$, and hence $n>4$, then the subspace $V_0$ contains only the monomials of degree $\leq 2$ in $x$.

3) Similarly, for $k>1$, and hence $n>2$, the subspace $V_1$ does not
contain monomials of the form $x_{i}x_jQ_s$.

Therefore, if $n>4$, then the monomials belonging to $V$ can be only of
degree $\leq 2$ in $x, Q, q$, and monomials of degree 2 in $x, Q, q$
should be of the form $f(p, P)x_{i}x_j$ or $f(p, P) x_{i} Q_j$, and
\[
\renewcommand{\arraystretch}{1.4}
\begin{array}{l}
V\subset \Cee[p, P]
\bigoplus (\mathop{\oplus}\limits_i\Cee[p, P] x_{i})\bigoplus
(\mathop{\oplus}\limits_i\Cee[p, P]Q_{i})\bigoplus \\
(\mathop{\oplus}\limits_i\Cee[p, P]q_{i})\bigoplus
(\mathop{\oplus}\limits_i\Cee[p, P] x_{i} x_{j})\bigoplus
(\mathop{\oplus}\limits_i\Cee[p, P]
x_{i}Q_{j}).\end{array}
\]
Denote the summands above by $W$, $W_x$, $W_Q$, $W_q$, $W_{xx}$ and
$W_{xQ}$, respectively.

Consider the quantization of $\fh$:
\[
\begin{array}{l}
p_{i}\longmapsto p_{i}, \quad P_{i}\longmapsto P_{i}, \quad
y_{i}\longmapsto
y_{i}, \\
q_{i}\longmapsto -(-1)^{p(p_{i})}\pder{p_i}, \quad Q_{i}\longmapsto
-(-1)^{p(P_{i})}\pder{P_i}, \quad z_{i}\longmapsto
-(-1)^{p(y_{i})}\pder{y_i}.
\end{array}
\]
We see that $W$ goes into operators of multiplication by functions,
whereas $W_x$, $W_Q$ and $W_q$ turn into vector fields. Therefore, the
commutators of the elements of $W\oplus W_x\oplus W_Q\oplus W_q$
do not change under the quantization.

Since the space $W\oplus W_Q\oplus W_q$ does not depend on $x$,
the commutators of its elements with the elements of $W_{xx}\oplus
W_{xQ}$ do not change, either.

Thus, all deformations (``quantizations") that change the
commutators concern only the brackets of $ W_{xQ}$ with elements of
$W_x\oplus W_{xx}\oplus W_{xQ}$.

4) Since the monomials $x_{i}x_{j}$ go to $\fg_{0, 0}$ and $Q_{j}\in\fg_{-1,
k-2}$ and since $\{P_s x_{i}x_{j}, Q_{s}\}_{P.b.}=\pm x_{i}x_{j}$,
it follows that the monomials $P_s x_{i}x_{j}$ go to $\fg_{1, l}$, where
$l+k-2=0$ and $l\geq -1$. Hence, $k=2-l$ and for $k>3$ (\textit{i.e.}, for
$n>6$) the subspace $V_1$ does not contain elements of the form $P_s
x_{i}x_{j}$, and hence $W_{xx}$ is of the form
$\mathop{\oplus}_{i, j}\ \Cee[p] x_{i}x_{j}$.

5) Since the monomials $x_{i}Q_{j}$ go to $\fg_{0, -1}$ and since $\{P_s
x_{i}Q_{j}, Q_{s}\}_{P.b.}=\pm x_{i}Q_{j}$, it follows that the monomials $P_s
x_{i}Q_{j}$ go to $\fg_{1, l}$, where $l+k-2=-1$ and $l\geq -1$.
Hence, $k=1-l$ and for $k>2$ (\textit{i.e.}, for $n>4$) the subspace $V_1$ does
not contain elements of the form $P_s x_{i}Q_{j}$, and hence $W_{xQ}$
is of the form $\mathop{\oplus}_{i, j}\ \Cee[p] x_{i}Q_{j}$.

6) Since the monomials $x_{i}P_{j}$ go to $\fg_{0, 1}$ and since
$\{P_sP_{i}x_{j}, Q_{s}\}_{P.b.}=\pm P_{i}x_{j}$, it follows that the monomials 
$P_s P_{i}x_{j}$ go to $\fg_{1, l}$, where $l+k-2=1$ and $l\geq -1$.
Hence, $k=3-l$ and for $k>4$ (\textit{i.e.}, for $n>8$) the subspace $V_1$ does
not contain elements of the form $P_s P_{i}x_{j}$, and hence $W_{x}$
is of the form $\mathop{\oplus}_{i}\,\left(\Cee[p]\bigoplus (\mathop{\oplus}_{j}\, \Cee[p]
P_{j}) x_{i}\right)$.

Thus, for $n>8$, the space $W_x\oplus W_{xx}\oplus W_{xQ}$ has only
elements of degree $\leq 2$ in $P, x, Q$:
\[
\renewcommand{\arraystretch}{1.4}\begin{array}{l}
W_x\oplus W_{xx}\oplus W_{xQ}=\\
\Cee[p]\otimes\left((\mathop{\oplus}\limits_{i} \Cee x_{i})\bigoplus (\mathop{\oplus}\limits_{I, j}
\Cee \Pee_jx_{i})\bigoplus (\mathop{\oplus}\limits_{i} \Cee x_{i} x_{j})\bigoplus
(\mathop{\oplus}\limits_{i} \Cee x_{i} Q_{j}) \right).
\end{array}
\]
The quantization preserves the brackets of embedded vector fields.

Let $n=8$. Then, $W_{xx}=\mathop{\oplus}_{i, j}\ \Cee[p] x_{i} x_{j}$ and
$W_{xQ}=\mathop{\oplus}_{i, j}\ \Cee[p] x_{i} Q_{j}$, so commutators of the
elements of these spaces do not change under quantization.

Consider brackets of elements of $W_{xQ}$ and $W_{x}$. The effect of
quantization can appear only in brackets of order-2
differential operators with functions, but in the case under consideration
the coefficients of order-2 differential operators are constant
with respect to $x, P$, so
\[
\renewcommand{\arraystretch}{1.4}\begin{array}{l}
\left[\frac{\partial ^2}{\partial x_j\partial P_k}, f(P)x_j\right]=
(-1)^{p(x_j)}\pderf{f(P)}{P_k}x_j\partial_{x_j}+f(P) \partial_{P_k}+
\pderf{f(P)}{P_k}=\\
\widehat{\left \{ x_jQ_k, f(P)x_j\right\}}\stackrel{\text{under quantization}}{=}\\
\frac12\Div \left((-1)^{p(x_j)}\pderf{f(P)}{P_k}x_j\partial_{x_j} +
f(P) \partial_{P_k} \right).
\end{array}
\]
Here, it is important that $W_{xQ}$ does not depend on $P$, \textit{i.e.},
the coefficients of order-2 differential operators are constant with
respect to $x, P$ and \textbf{the degree of functions is immaterial}, so
for $n=6$ the brackets of elements of $W_{xQ}$ and $W_{x}$ do not change
under quantization, either.

Let now $n=6$. Then,
\[
W_{xx}=\{ \mathop{\oplus}\limits_{i, j} \Cee[p, P] x_{i}
x_{j}\mid\text{ the degree with respect to $P$ is $\leq 1\}$. }
\]
Indeed, the elements $P_s x_{i} x_{j}$ go to $\fg_{1, -1}$ and
$Q_t\in\fg_{-1, 1}$, and since $\{P_tP_s x_{i} x_{j},
Q_t\}_{P.b.}=\pm P_s x_{i} x_{j}$, it follows that the elements $P_tP_s x_{i}
x_{j}$ go to $\fg_{2, l}$, where $l+1=-1$, which is a~contradiction.

The elements of the form $ \Cee[p] x_{i} x_{j}$ and $\Cee[p] x_{i}
x_{j}$ lie in $\widetilde V_{-1}$. \end{proof}

\section{The exceptional simple vectorial Lie superalgebra
$\fm\fb(4|5)$}\label{Smb}

In this section, already used in \cite{BGLLS}, 
we illustrate the algorithm presented in detail in
\cite{Shch} and rectify one formula from \cite{CCK}. This algorithm
allows one to describe vectorial Lie superalgebras by means of
differential equations. In \cite{Sh5, Sh14} the algorithm was used
to describe the exceptional simple vectorial Lie superalgebras of
depth~ 1.

The Lie superalgebra $\fm\fb(4|5)$ has three realizations as a
transitive and primitive (\textit{i.e.}, not preserving invariant foliations
on the space it is realized by means of vector fields) vectorial Lie
superalgebra. Speaking algebraically, the requirement to be a~
transitive and primitive vectorial Lie superalgebra is formulated as
having a~W-filtration, so $\fm\fb(4|5)$ has three such
filtrations.

Two of these W-filtrations are of depth 2, and one is of depth 3.
In this section, we consider the latter case (the grading $K$), \textit{i.e.},
we explicitly solve the differential equations that single out our Lie
superalgebra. We thus explicitly obtain the expressions for the
elements of $\fm\fb(4|5; K)$. In particular, we verify and
rectify the formulas for the brackets given in \cite{CCK}.

The Lie superalgebra $\fm\fb(4|5)$ first appeared in \cite{Sh3}. The
construction was performed in two steps:

1) Take $\fg_0=\fvect(0|3)$ with the natural action $T^{1/2}$ on the space
$\fg_{-1}=\Pi(\Lambda(\zeta_1,\zeta_2,\zeta_3)\otimes \vvol^{1/2})$
of semidensities with the inverted parity, see formula \eqref{T12}. The Cartan prolong of the pair $(\fg_{-1},\fg_0)_*$ forms the Lie
superalgebra $\fb_{1/2}(3;3)$.

2) Let us describe an embedding of the Lie superalgebra
$\fm\fb(4|5)$ we wish to define into $\fm(4)$. We work in terms of
generating functions and sometimes, by abuse of notation, write $f$
instead of $M_f$.

The exceptional simple Lie superalgebra $\fm\fb(4|5)$ is the
complete prolong in $\fm(4)$ of the trivial central extension
of $\fb_{1/2}(3;3)$ considered as the Lie subsuperalgebra
$\fc(\fb_{1/2}(3;3))\subset \fm(4)$, the
central element $\tau\in\fm$ being added to the 0th component of $\fb_{1/2}(3;3)$.

Formally speaking, $\fm\fb(4|5)$ is defined. To explicitly describe
its elements, let us work a~bit more. To avoid confusion, we assume
that the elements of $\fc(\fb_{1/2}(3;3))$ are generated by functions
of $\eta_1, \eta_2,\eta_3, u_1, u_2, u_3$, whereas the elements of
$\fm(4)$ are generated by functions of $\xi_i$ and $q_i$ for $i=0,
1, 2, 3$, and $\tau$.

\textbf{Recapitulation from \cite{Shch}}. Let $\fn$ be the nilpotent
Lie superalgebra; we wish realize $\fn$ by vector fields, \textit{i.e.},
embed $\fn$ into $\fvect(\dim\fn)$ so that the image of $\fn$ is
spanned by all partial derivatives modulo the vector fields that vanish
at the origin. Such an embedding determines a~nonstandard grading of
certain depth $d$ on $\fvect(\dim\fn)$. Denote $\fvect(\dim\fn)$
endowed with this nonstandard grading by $\fv=
\mathop{\oplus}_{-d\leq k}^{\infty}\fv_k$. Let $\fg_-$ be the
image of $\fn$ in $\fv$. Let $\fg_-$ be spanned by the
left-invariant vector fields $X_i$ on the supergroup $N$ whose Lie
superalgebra is $\fn$. Let the $Y_i$ be the right-invariant vector
fields on the supergroup $N$ that also span $\fg_-$ and such that
$
X_i(e)=Y_i(e),
$
where $e$ is the unit of $N$. We have
\[
 M_{\xi_i}=\pder{q_i}+\xi_i\pder{\tau}; \quad
M_{q_i}=-\pder{\xi_i}+q_i\pder{\tau},
\]
and the operators $Y$ are (as follows from \cite{Shch}) of the form:
\[
 Y_{q_i}=\pder{q_i}-\xi_i\pder{\tau}; \quad
Y_{\xi_i}=\pder{\xi_i}+q_i\pder{\tau}.
\]
According to \cite{Shch}, the equations for generating functions
should be expressed in terms of the operators $Y_{q_i}$ and
$\widehat Y_{\xi_i}:=Y_{\xi_i}\circ \Pty$, \textit{i.e.}, $\widehat
Y_{\xi_i}(f)=(-1)^{p(f)} Y_{\xi_i}(f)$. We have
\[
\fm\fb_-=\fc(\fb_{1/2}(3;3))_-\simeq \fm(4)_-.
\]
Here is the explicit correspondence between elements of degree $-1$:

\begin{equation}
\begin{tabular}{|c|c|}
\hline $\fm\fb_{-1}$ & $\fm_{-1}(4)$\\
\hline $\eta_1\eta_2\eta_3$ & $\xi_0$\\
\hline $\eta_2\eta_3, \eta_3\eta_1, \eta_1\eta_2$ & $q_1, q_2, q_3$\\
\hline $\eta_1, \eta_2, \eta_3$ & $\xi_1, \xi_2,\xi_3$\\
\hline $1$ & $q_0$\\
\hline
\end{tabular}
\end{equation}

The 0th component of $\fc(\fb_{1/2}(3;3))$ is isomorphic to
$\fvect(0|3)$ and its action on the $(-1)$st component is the same
as on the space of semidensities. It is easy to see that the
supertrace of any such operator is equal to $0$. Thus, we get
embeddings $(\fc(\fb_{1/2}(3;3)))_0\subset \fspe(4)$ and
$\fm\fb\subset \fb_{1,0}(4)$.

Comparing the actions of various elements of $\fb_{1/2}(3;3)_0\simeq 
\fvect(0|3)$ and $\fspe(4)\subset \fm(4)_0$ we get the following
identifications (here $(i,j,k)\in A_3$, the group of even permutations):
\[
\renewcommand{\arraystretch}{1.4}
\begin{array}{ll}
\begin{array}{l}
\partial_i\lra -q_0q_i+\xi_j\xi_k,\\
 \eta_i\partial_i\lra \frac
12(-q_0\xi_0-q_i\xi_i+q_j\xi_j+q_k\xi_k),\\
\eta_i\partial_j\lra -q_j\xi_i\text{~~for $i\ne j$}, \\
\eta_j\eta_k\partial_i\lra -\frac 12 q_i^2,\\
\end{array}\begin{array}{l}
\eta_i(\eta_j\partial_j-\eta_k\partial_k)\lra -\frac 12
q_i^2-q_jq_k,\\
\eta_i(\mathop{\sum}\limits_j \eta_j\partial_j)\lra \xi_0\xi_i, \\
\eta_1\eta_2\eta_3\partial_i\lra -\frac 12 \xi_0q_i.
\end{array}\end{array}
\]
Since $\dim \fm(4)_0=17|16$ and $\dim \fm\fb(4|5)_0=13|12$, we need
$4|4$ equations to single out $\fm\fb_0$ in~ $\fm_0$.

Note that we have a~natural $\Zee$-grading of $\fm_{-1}$ induced
by the identification $\fm_{-1}\simeq \Pi(\Lambda(\eta))$; this
grading induces a~$\Zee$-grading on $\fm_0$ as well. Let $\deg$ be
the degree of a~given operator with respect to this $\Zee$-grading;
let $N$ be the number of even$|$odd equations; let $i$ and $j$ run
over the set $\{1,2,3\}$; let the triple $(ijk)$ be a~permutation of
$(123)$. We get
\begin{equation}
\label{tabb} {\footnotesize
\renewcommand{\arraystretch}{1.4}
\begin{tabular}{|c|c|c|c|c|}
\hline
$\deg$ & $\fm_0$ & $\fm\fb_0$ & $N$ & equations\\
\hline

$-3 $ & $q_0^2$ & --- &$0|1$ & $Y_{q_0}^2(f)=0$\\
\hline $-2$ & $q_0\xi_i, $ & ---
&$3|0$ &$Y_{q_0}\widehat Y_{\xi_i}(f)=0$\\

\hline $-1$ & $q_0q_i, \; \xi_i\xi_j$ & $-q_0q_i+\xi_j\xi_k$&$0|3$
&$\left(Y_{q_0}Y_{q_i}+
\widehat Y_{\xi_j}\widehat Y_{\xi_k}\right)(f)=0$\\
\hline

$0$ & $q_i\xi_j$ & $q_i\xi_j,\;i\ne j$ & $1|0$ & \\

 &$q_0\xi_0$ & $-q_0\xi_0-q_i\xi_i+q_j\xi_j+q_k\xi_k$ &
&$\mathop{\sum}\limits_{0\leq i\leq 3}\left(Y_{q_i}\widehat Y_{\xi_i}+
\widehat Y_{\xi_i}Y_{q_i}\right)(f)=0$ \\

\hline $1$ & $\xi_i\xi_0, q_iq_j$ & $\xi_i\xi_0,
q_iq_j$ & $0$ & \\

\hline $2$ & $\xi_0q_i$ & $\xi_0q_i$ & $0$
& \\

\hline $3$ & --- & --- & $0$ & \\
\hline

\end{tabular}}
\end{equation}

The equations are in such a~form that we can discard the hats over
the operators $Y_{\xi_i}$. After that, we see that
the functions that generate elements of the subalgebra $\fm\fb$ are
solutions of the following system of $4|4$ equations (here $i=1,2,3$
and $(i,j,k)\in A_3$):
\[
\renewcommand{\arraystretch}{1.4}
\begin{array}{ll}
Y_{q_0}^2(f)=0, & Y_{q_0}Y_{\xi_i}(f)=0, \\
\left(Y_{q_0}Y_{q_i}- Y_{\xi_j}Y_{\xi_k}\right)(f)=0, &
\mathop{\sum}\limits_{0\leq i\leq 3}\left(Y_{q_i}Y_{\xi_i}+
Y_{\xi_i}Y_{q_i}\right)(f)=0,\end{array}
\]
or, more explicitly,
\begin{equation}
\label{eq1} \frac{\partial^2 f}{\partial
q_0^2}-2\xi_0\pder{q_0}\pderf f\tau=0,
\end{equation}
\begin{equation}
\label{eq2} \frac{\partial^2 f}{\partial q_0\partial
\xi_i}+\left(q_i\pder{q_0}+\xi_0\pder{\xi_i}\right)\pderf f\tau=0,
\end{equation}
\begin{equation}
\label{eq3} \frac{\partial^2 f}{\partial q_0\partial
q_i}-\frac{\partial^2 f}{\partial \xi_j\partial
\xi_k}-\left(\xi_0\pder{q_i}+\xi_i\pder{q_0}-q_j\pder{\xi_k}+q_k\pder{\xi_j}\right)\pderf
f\tau=0,
\end{equation}
\begin{equation}
\label{eq4} \Delta f+E\pderf f\tau=0.
\end{equation}

The last equation means precisely that $\fm\fb$ lies inside of
$\fb_{1,0}(4)$.

Set $q=(q_1,q_2,q_3)$, $\xi=(\xi_1, \xi_2,\xi_3)$. Each of the
equations ~\eqref{eq1}--\eqref{eq4} is of the form
\[
D(f)=(D_0+D_1\partial_\tau)(f)=0,
\]
where the operators $D_0$ and $D_1$ do not depend on $\tau$. On the
space of functions of the form $f=f_0+\tau f_1$, where $f_0, f_1\in
\Cee[q,q_0,\xi,\xi_0]$, each of these equations is equivalent to the
system
\[
\renewcommand{\arraystretch}{1.4}
\begin{cases}
D_0f_1=0, \\
D_0f_0+D_1f_1=0. \end{cases}
\]
Therefore, the solutions of the system ~\eqref{eq1}--\eqref{eq4}
independent of $\tau$ (they generate the image of
$\fc(\fb_{1/2}(3;3))$ in $\fm(4)$), as well as the coefficients of $\tau$
in all the solutions of the system ~\eqref{eq1}--\eqref{eq4}, are
solutions of the following system (here $i=1,2,3$ and $(i,j,k)\in
A_3$):
\begin{equation}
\label{beztau}
\renewcommand{\arraystretch}{1.20}
\begin{array}{ll}
 \frac{\partial^2 f}{\partial q_0^2}=0, \qquad
 \frac{\partial^2 f}{\partial q_0\partial\xi_i}=0, \\
 \frac{\partial^2 f}{\partial q_0\partial
q_i}-\frac{\partial^2 f}{\partial \xi_j\partial \xi_k}=0, \qquad
 \Delta f=0.
 \end{array}
 \end{equation}
The first two groups of equations in~\eqref{beztau} imply that
\[
f=\vf(q,\xi,\xi_0)+q_0\vf_1(q, \xi_0).
\]
Differentiating the equations of the third group in
~\eqref{beztau} with respect to $\xi_i$, and taking into account the
second group we get
\[
 \frac{\partial^3 f}{\partial\xi_1\partial\xi_2\partial\xi_3}=0.
\]
Moreover, the third group of equations ~\eqref{beztau} implies that
\[
 \frac{\partial \vf_1}{\partial q_i}=\frac{\partial^2 \vf}{\partial
\xi_j\partial \xi_k},
\]
\textit{i.e.},
\begin{equation}
 \label{fbeztau}
f=\vf_0(q,\xi,\xi_0)-\Delta(\xi_1\xi_2\xi_3\vf_1)+q_0\vf_1(q,
\xi_0),\quad\text{where $\deg_\xi\vf_0\le 1$.}
 \end{equation}
Conversely, any function of the form ~\eqref{fbeztau} such that
$\Delta f=\Delta\vf_0+\pderf{\vf_1}{\xi_0}=0$ is a~solution of
eq.~\eqref{beztau}, and hence of the system
~\eqref{eq1}--\eqref{eq4}.

Explicit solutions are:
\begin{equation}
 \label{jav1}
f(q), \quad f(q)\xi_0,
\end{equation}
\begin{equation}
 \label{jav2}
F_f=f-q_0\xi_0\Delta f +\xi_0\Delta(\xi_1\xi_2\xi_3\Delta f),\;
\text{ where }\; f=\mathop{\sum}\limits_{1\leq i\leq 3} f_i(q)\xi_i,
\end{equation}
\begin{equation}
 \label{jav3}
 \xi_0\mathop{\sum}\limits_{1\leq i\leq 3}f_i(q)\xi_i, \;\text{ where }\;
\mathop{\sum}\limits_{1\leq i\leq 3}\pderf{f_i}{q_i}=0,
\end{equation}
\begin{equation}
 \label{jav4}
 q_0f(q)-\Delta(\xi_1\xi_2\xi_3 f).
\end{equation}

Now let us find the solutions for which the coefficient of $\tau$
is of the form ~\eqref{jav1}--\eqref{jav4}. Set
\begin{equation}
 \label{Phi}
 \Phi=\mathop{\sum}\limits_{1\leq i\leq 3}
q_i\xi_i.
\end{equation}
The solutions of ~\eqref{jav1} are obviously given by
\[
 \xi_0\left(\mathop{\sum}\limits_{1\leq i\leq 3}f_i(q)\xi_i+\frac
13(\tau-\Phi)\Delta f\right), \text{ where }
f=\mathop{\sum}\limits_{1\leq i\leq 3}f_i(q)\xi_i,
\]
and
\[
kq_0\xi_0f(q)-(k+1)\xi_0\Delta(f\xi_1\xi_2\xi_3)-\tau f, \text{
where } k=\deg f.
\]

Combining this last solution with $F_f$ from ~\eqref{jav2} we get a
solution of the form
\[
 \widetilde F_f=f+ q_0\xi_0\Delta f+\frac
12(\tau-\Phi-q_0\xi_0)\Delta f,\text{ where }f=
\mathop{\sum}\limits_{1\leq i\leq 3} f_i(q)\xi_i.
\]

The functions $F_f$ and $\widetilde F_f$ determine two embeddings of
$\fvect(3|0)$ into $\fm\fb$. From general considerations the
embedded algebras should be obtained from one another by regradings
of the ambient algebra.

Now, let us find solutions with coefficient of $\tau$ of the form
~\eqref{jav2} with $\Delta f=0$. Then, $f$ can be represented as $f=\Delta\widetilde f$ for a~function $\widetilde f$ of the
form
\[
\widetilde f=\sum f^i(q)\xi_j\xi_k, \qquad\text{where $(i,j,k)\in
A_3$}.\] Take
\[\widehat f=q_0\xi_0\Delta \widetilde f +\tau\Delta\widetilde f.\]
Then, $\widehat f$ satisfies the system ~\eqref{eq1}--\eqref{eq2}, and
~\eqref{eq3} implies:
\[
\renewcommand{\arraystretch}{1.4}
\begin{array}{l}
 \frac{\partial^2 \widehat f}{\partial q_0\partial
q_i}-\frac{\partial^2 \widehat f}{\partial \xi_j\partial
\xi_k}-\left(\xi_0\pder{q_i}+\xi_i\pder{q_0}-q_j\pder{\xi_k}+q_k\pder{\xi_j}\right)\pderf
{\widehat f}\tau=\\
 \xi_0\pderf{\Delta\widetilde
f}{q_i}-\xi_0\pderf{\Delta \widetilde f}{q_i}+q_j\pderf{\Delta
\widetilde f}{\xi_k}-q_k\pderf{\Delta \widetilde
f}{\xi_j}=\\
 q_j\left(\pderf{f^i}{q_j}-\pderf{f^j}{q_i}\right)-
q_k\left(\pderf{f^k}{q_i}-\pderf{f^i}{q_k}\right)=\\
 q_i\pderf{f^i}{q_i}+q_j\pderf{f^i}{q_j}+q_k\pderf{f^i}{q_k}-
q_i\pderf{f^i}{q_i}-q_j\pderf{f^j}{q_i}-q_k\pderf{f^k}{q_i}=\\
 E(f^i)-\frac{\partial}{\partial
q_i}(q_if^i+q_jf^j+q_kf^k)+f^i=\\
 E(f^i)- \pder{q_i}\left(\left(\int_{\xi}\Phi\right)\cdot
\widetilde f\right)+f^i.
\end{array}
\]
Therefore, it is easy to fix up $\widehat f$ to get a~solution:
\[
 F=\widetilde f-E(\widetilde f)+q_0\left(\int_{\xi}\Phi \widetilde
f\right)+q_0\xi_0\Delta \widetilde f+\tau\Delta \widetilde f.
\]
It is subject to a~direct verification that $F$ satisfies also
~\eqref{eq4}.

We similarly see that, starting with $\widehat
f=\tau\xi_0\Delta f$, where $f=\sum f^i(q)\xi_j\xi_k$, we arrive at
the function
\[
 F_1=-q_0\xi_0\left(\int_{\xi}\Phi
f\right)+\xi_0(E-1)(f)+\tau\xi_0\Delta f,
\]
which satisfies the system ~\eqref{eq1}--\eqref{eq3}. However,
\[
 \Delta F_1+E\pderf{F_1}{\tau}=-\int_{\xi}(\Phi f).
\]
Adding to $F_1$ the function $\vf=\sum \vf_i\xi_i$ such that $\Delta
\vf=\int_{(\xi)}(\Phi f)$ we get a~solution of the system
~\eqref{eq1}--\eqref{eq4} which, by using the obtained solutions, can be simplified to
\[
 F=\xi_0\left(-\Delta\bigg(\xi_1\xi_2\xi_3\left(\int_{\xi}\Phi
f\right)\bigg)+(E-1)(f)-\tau\Delta f\right). \] So we have obtained
a solution in which the coefficient of $\tau$ is of the form ~\eqref{jav3}.

Similarly, if we want to get a~solution in which the coefficient of $\tau$
is of the form ~\eqref{jav4}, \textit{i.e.}, equal to
\[
q_0\vf(q)-\Delta(\xi_1\xi_2\xi_3\vf),
\] 
we take
\[
F_1=q_0^2\xi_0\vf+\tau\left(q_0\vf(q)-\Delta(\xi_1\xi_2\xi_3\vf)\right)
\]
to satisfy ~\eqref{eq1}, then pass to
\[
F_2=F_1-q_0\Phi\vf-q_0\xi_0\Delta(\xi_1\xi_2\xi_3\vf)
\]
to satisfy ~\eqref{eq2}, and, finally, we pass to
\[
F=q_0(q_0\xi_0-\Phi)\vf-q_0\xi_0\Delta(\xi_1\xi_2\xi_3\vf)+\xi_1\xi_2\xi_3(E+2)\vf+
\tau\left(q_0\vf(q)-\Delta(\xi_1\xi_2\xi_3\vf)\right)
\]
to satisfy ~\eqref{eq3}. Now, one can verify that condition~\eqref{eq4} is
satisfied.

Finally, let us obtain a~solution in which the coefficient of $\tau$ is an
arbitrary function of the form ~\eqref{jav4}. We start with the
function
\[
 F_0=\tau(f-q_0\xi_0\Delta f+\xi_0\Delta(\xi_1\xi_2\xi_3\Delta f)),
\text{ where } f=\mathop{\sum}\limits_{1\leq i\leq 3} f_i(q)\xi_i.
\]
This $F_0$ satisfies ~\eqref{eq1}, and
\[
 \left(q_i\pder{q_0}+\xi_0\pder{\xi_i}\right)\pderf{F_0}{\tau}=-q_i\xi_0\Delta
f+\xi_0f_i.
\]
To satisfy ~\eqref{eq2} we pass to
\[
F_1=q_0\xi_0(f-\Phi\Delta f)+F_0=q_0\xi_0(f-\Phi\Delta f)+
\tau(f-q_0\xi_0\Delta f+\xi_0\Delta(\xi_1\xi_2\xi_3\Delta f)).
\]
Next, to satisfy ~\eqref{eq3}, we pass to
\[
F_2=\xi_0\xi_1\xi_2\xi_3(E+2)(\Delta f)-\Phi f+F_1,
\]
and check that ~\eqref{eq4} is also satisfied. Finally, we get the
solution
\[
F=F_2=\xi_0\xi_1\xi_2\xi_3(E+2)(\Delta f)-\Phi f+
q_0\xi_0(f-\Phi\Delta f)+ \tau(f-q_0\xi_0\Delta
f+\xi_0\Delta(\xi_1\xi_2\xi_3\Delta f)).
\]

\ssec{The W-gradings of $\fm\fb(3|8)$} We will again assume that the
elements of the Lie superalgebra $\fc(\fb_{1/2}(3;3))$ that we have
started with in order to define the Lie superalgebra $\fm\fb(3|8)$
are given by generating functions in $\zeta_1, \zeta_2,\zeta_3, p_1,
p_2,p_3$, whereas the basis element of the added center is $z$. The
regrading of $\fm\fb(3|8)$ is obtained by setting
\[
\deg p_i=2, \quad \deg \zeta_i=1 \text{ for } i=1,2,3, \quad
\deg_{Lie}(f)=\deg(f)-3.
\]

This regrading is naturally extended to the whole Lie superalgebra
$\fm\fb(3|8)$ and is of depth 3. The component $\fg_{-3}$ is spanned
by $1\in \fc(\fb_{1/2}(3;3))$ and the added center $z$.

In this realization, the Lie superalgebra $\fg=\fm\fb(3|8)$ is the
complete prolong of its negative part $\fg_-=\fg_{-3}\oplus
\fg_{-2}\oplus \fg_{-1}$, where
\[
\fg_{-3}=\Span(1,z), \ \fg_{-2}=\Span(\zeta_i\mid i=1,2,3),\;
\fg_{-1}=\Span(p_i, \zeta_j\zeta_k\mid i,j,k=1,2,3),
\] 
and
therefore
\[ \dim\fg_{-3}=0|2, \quad \dim\fg_{-2}=3|0, \quad
\dim\fg_{-1}=0|6.
\]
The non-zero brackets are as follows (here $i=1,2,3$,\; $(i,j,k)\in
A_3)$:
\[
\{p_i,\zeta_i\}=-1, \;\; \; \{\zeta_i\zeta_j,\zeta_k\}=-z, \;\; \;
\{p_i,\zeta_i\zeta_j\}=-\zeta_j, \;\; \;
\{p_i,\zeta_k\zeta_i\}=\zeta_k.
\]
In particular, the center of $\fg_-$ coincides with $\fg_{-3}$.

We use this Lie superalgebra to demonstrate how to work with
generating functions.

First of all, we embed $\fg_-$ into the Lie superalgebra
\[
\fv=\fvect(3|8)=\fder~\,\Cee[u_1,u_2,u_3 \mid \eta_1,\eta_2,\eta_3,\zeta_1,
\zeta_2,\zeta_3, \chi_1, \chi_2]
\]
considered with the grading (for any $i,j$)
\[
\deg \eta_i=\deg \zeta_i=1, \quad \deg u_i=2, \quad \deg \chi_j=3.
\]
In accordance with \cite{Shch}, we find in $\fv_-$ two mutually
commuting families of elements: $X$-vectors and $Y$-vectors. The
table of correspondences reads
\begin{equation}
\renewcommand{\arraystretch}{1.4}
\begin{tabular}{|c|c|c|c|}
\hline $k$ & $\fm\fb_{-k}$ & $X$&$Y$\\
\hline $-1$ & $q_1, q_2, q_3$ & $X_{\eta_i}$&$Y_{\eta_i}$\\
 & $\xi_2\xi_3, \xi_3\xi_1, \xi_1\xi_2$ & $X_{\zeta_i}$&$Y_{\zeta_i}$\\
\hline $-2$ & $\xi_1, \xi_2,\xi_3$ & $X_{u_i}$&$Y_{u_i}$\\
\hline $-3$ & $1, \widehat 1$ & $X_{\chi_j}$&$Y_{\chi_j}$\\
\hline
\end{tabular}
\end{equation}
The non-zero commutation relations for the $X$-vectors, where
$(i,j,k)\in A_3$, are
\[
\renewcommand{\arraystretch}{1.4}
 \begin{array}{ll}
 {}[X_{\eta_i},X_{u_i}]=-X_{\chi_1}, &
[X_{\zeta_i},X_{u_i}]=-X_{\chi_2},\\
{}[X_{\eta_i},X_{\zeta_k}]=-X_{u_j}, &
[X_{\eta_i},X_{\zeta_j}]=X_{u_k}. \end{array}
\]
The non-zero commutation relations for the $Y$-vectors correspond to
the negative of the above structure constants:
\[
\renewcommand{\arraystretch}{1.4}
 \begin{array}{ll}
 {}[Y_{\eta_i},Y_{u_i}]=Y_{\chi_1}, &
[Y_{\zeta_i},Y_{u_i}]=Y_{\chi_2},\\
{}[Y_{\eta_i},Y_{\zeta_k}]=Y_{u_j}, &
[Y_{\eta_i},Y_{\zeta_j}]=-Y_{u_k}.
\end{array}
\]
Let us represent an arbitrary vector field $X\in\fvect(3|8)$ in the
form
\begin{equation}
 \label{fldD}
X=F_1Y_{\chi_1}+ F_2Y_{\chi_2}+\mathop{\sum}\limits_{1\leq i\leq 3}\left(
f_{\xi_i}Y_{\xi_i}+f_{\eta_i}Y_{\eta_i}+f_{u_i}Y_{u_i}\right).
\end{equation}

As it was shown in \cite{Shch}, any $X\in\fm\fb(3|8)$ is completely
determined by a~pair of functions $F_1, F_2$ via the equations, where
$i=1,2,3$ and $(i,j,k)\in A_3$,
\begin{equation}
 \label{coor1}
Y_{\xi_i}(F_1)=0,\quad Y_{\eta_i}(F_1)=-(-1)^{p(f_{u_i})}f_{u_i}=
Y_{\xi_i}(F_2),\quad Y_{\eta_i}(F_2)=0,
\end{equation}
\begin{equation}
 \label{coor2}
Y_{\xi_i}(f_{u_i})= Y_{\eta_i}(f_{u_i})=0,
\end{equation}
\begin{equation}
 \label{coor3}
Y_{\xi_i}(f_{u_j})= (-1)^{p(f_{\eta_k})}f_{\eta_k}, \quad
Y_{\eta_i}(f_{u_j})=-(-1)^{p(f_{\xi_k})}f_{\xi_k},
\end{equation}
\begin{equation}
 \label{coor4}
Y_{\xi_i}(f_{u_k})= -(-1)^{p(f_{\eta_j})}f_{\eta_j}, \quad
Y_{\eta_i}(f_{u_k})=(-1)^{p(f_{\xi_j})}f_{\xi_j}.
\end{equation}

Therefore, the functions $F_1, F_2$ must satisfy the following three
groups of equations:
\begin{equation}
 \label{eqF}
 Y_{\xi_i}(F_1)=0,\quad Y_{\eta_i}(F_1)=
Y_{\xi_i}(F_2),\quad Y_{\eta_i}(F_2)=0~~~ \text{for $i=1,2,3$}.
\end{equation}

The relations ~\eqref{coor1}, ~\eqref{coor3} determine the remaining
coordinates, while the relations ~\eqref{coor2}, ~\eqref{coor4} follow
from ~\eqref{coor1}, ~\eqref{coor3} and the commutation relations that $Y$
obey. Indeed,
\[
Y_{\xi_i}(f_{u_i})=-(-1)^{p(f_{u_i})}Y_{\xi_i}Y_{\xi_i}(F_2)=0
\]
and
\[
Y_{\xi_i}(f_{u_j})=-(-1)^{p(f_{u_j})}Y_{\xi_i}Y_{\xi_j}(F_2)=
(-1)^{p(f_{u_i})}Y_{\xi_j}Y_{\xi_i}(F_2)= -Y_{\xi_j}(f_{u_i}),
\]
since $p(f_{u_j})=p(f_{u_i})=p(X)$. Moreover,
\[
f_{\xi_k}=-(-1)^{p(f_{\xi_k})}Y_{\eta_i}(f_{u_j})=
-Y_{\eta_i}Y_{\xi_j}(F_2)=Y_{\xi_j}Y_{\eta_i}(F_2)+
Y_{u_k}(F_2)=Y_{u_k}(F_2).
\]
We similarly get the expressions for the remaining coordinates:
\[f_{\eta_k}=Y_{u_k}(F_1).\]

Therefore, an arbitrary element $X\in\fm\fb(3|8)$ is of the form
\begin{equation}
 \label{XF}
 X=X^{F}=F_{1}Y_{\chi_1}+F_{2}Y_{\chi_2}+\mathop{\sum}\limits_{1\leq i\leq 3}\left(Y_{u_i}(F_2)Y_{\xi_i}
+Y_{u_i}(F_1)Y_{\eta_i}-(-1)^{p(X)}Y_{\xi_i}(F_2)Y_{u_i} \right),
\end{equation}
where the pair of functions $F=\{F_1, F_2\}$ satisfies the system of
equations ~\eqref{eqF}.

We select the $Y$-vectors so that the equations ~\eqref{eqF} that the
functions $F_1, F_2$ should satisfy are as simple as possible. For
example, take the following $Y$-vectors, where $i=1,2,3$, $s=1,2$, and 
$(i,j,k)\in A_3$:
\[
\renewcommand{\arraystretch}{1.4}
\begin{array}{ll}
Y_{\eta_i}=\partial_{\eta_i}+\xi_k\partial_{u_j}-\xi_j\partial_{u_k}+(\xi_k\eta_j-\xi_j\eta_k)
\partial_{\chi_1}-\xi_j\xi_k\partial_{\chi_2},&
Y_{\xi_i}=\partial_{\xi_i}, \\
Y_{u_i}=\partial_{u_i}+\eta_i\partial_{\chi_1}+\xi_i\partial_{\chi_2},
&
Y_{\chi_s}=\partial_{\chi_s}. \end{array}
\]
Then, the $X$-vectors are of the form
\begin{equation}
\label{Xvect}
\renewcommand{\arraystretch}{1.4}
\begin{array}{ll}
X_{\xi_i}=\partial_{\xi_i}-\eta_j\partial_{u_k}+
\eta_k\partial_{u_j}-\eta_j\eta_k\partial_{\chi_1}+
u_i\partial_{\chi_2},&
X_{\eta_i}=\partial_{\eta_i}+u_i\partial_{\chi_1},\\
 
X_{\chi_1}=\partial_{\chi_1},\ \ \ 
X_{\chi_2}=\partial_{\chi_2}&X_{u_i}=\partial_{u_i}.\end{array}
\end{equation}

The Lie superalgebra $\fm\fb(3|8)$ consists of the vector fields that
preserve the distribution determined by the following equations
for the vector field $D$ of the form ~\eqref{fldD}:
\begin{equation}
\label{distr1} f_{u_1}=f_{u_2}=f_{u_3}=F_1=F_2=0.
\end{equation}

Let us express the coordinates $f$ of the field $D$ in the $Y$-basis
in terms of the standard coordinates in the basis of partial
derivatives:
\[
D=g_{\chi_1}\partial_{\chi_1} +g_{\chi_2}\partial_{\chi_1}+
\mathop{\sum}_{1\leq i\leq 3}\ \left(g_{\xi_i}\partial_{\xi_i}+
g_{\eta_i}\partial_{\eta_i} +g_{u_i}\partial_{u_i}\right).
\]
We get
\[
\renewcommand{\arraystretch}{1.4}
\begin{array}{l}
f_{u_i}=g_{u_i}+g_{\eta_j}\xi_k-g_{\eta_k}\xi_j, \\
F_1=g_{\chi_1}-\mathop{\sum}_{1\leq i\leq 3}\ g_{u_i}\eta_i, \\
F_2=g_{\chi_2}-\mathop{\sum}_{1\leq i\leq 3}\
g_{u_i}\xi_i-\mathop{\sum}_{1\leq i\leq 3,\ (i,j,k)\in A_3}
g_{\eta_i}\xi_j\xi_k.\end{array}
\]
Therefore, in the standard coordinates, the
distribution~\eqref{distr1} is of the form ($i=1,2,3$):
\begin{equation}
 \label{distr2}
\renewcommand{\arraystretch}{1.4}
\begin{array}{l} g_{u_i}+g_{\eta_j}\xi_k-g_{\eta_k}\xi_j=0, \\
g_{\chi_1}-\mathop{\sum}_{1\leq i\leq 3} g_{u_i}\ \eta_i=0, \\
g_{\chi_2}-\mathop{\sum}_{1\leq i\leq 3}
g_{u_i}\xi_i-\mathop{\sum}\limits_{1\leq i\leq 3,\ (i,j,k)\in A_3}
g_{\eta_i}\xi_j\xi_k=0.\end{array}
 \end{equation}
Using the first two relations we rewrite the third one
as
\[
g_{\chi_2}+\mathop{\sum}\limits_{1\leq i\leq 3,\ (i,j,k)\in A_3}
g_{\eta_i}\xi_j\xi_k=0.
\]
Assuming that the pairing between vectors and covectors is given by the formula
\[
\langle f\partial_\xi,gd\xi\rangle =(-1)^{p(g)}fg,
\] we see that the
distribution is specified by the Pfaff
equations given by the 1-forms
\begin{equation}
 \label{distr3}
\renewcommand{\arraystretch}{1.4}
\begin{array}{l}
d u_i+\xi_jd\eta_k-\xi_kd\eta_j, \\ d\chi_1-\mathop{\sum}
\eta_id u_i,
\\
d\chi_2+\mathop{\sum}\limits_{1\leq i\leq 3,\ (i,j,k)\in A_3}
\xi_j\xi_kd\eta_i.\end{array}
 \end{equation}

Let us now solve the system ~\eqref{eqF}.

Since $Y_{\xi_i}=\partial_{\xi_i}$, the equality $Y_{\xi_i}(F_1)=0$
implies that $F_1=F_1(u, \eta,\chi)$. Further, the equality
\[
Y_{\xi_i}(F_2)=Y_{\eta_i}(F_1)=0
\]
takes the form
\[
 \pderf{F_2}{\xi_i}=\pderf{F_1}{\eta_i}+\left(\xi_k\pderf{F_1}{u_j}-
\xi_j\pderf{F_1}{u_k}\right) +\left(\xi_k\eta_j-
\xi_j\eta_k\right)\pderf{F_1}{\chi_1}-\xi_j\xi_k\pderf{F_1}{\chi_2},
\]
whence (since $F_1$ does not depend on $\xi$) 
\begin{equation}
 \label{F2}
 F_2=\mathop{\sum}\limits_{1\leq i\leq 3}\xi_i\pderf{F_1}{\eta_i} -
\mathop{\sum}\limits_{1\leq i\leq 3,~~ (i,j,k)\in
A_3}\xi_i\xi_j\left(\pderf{F_1}{u_k}+\eta_k\pderf{F_1}{\chi_1}\right)
- \xi_1\xi_2\xi_3\pderf{F_1}{\chi_2} +\alpha_2,
\end{equation}
where $\alpha_2=\alpha_2(u,\eta,\chi)$, \textit{i.e.}, does not depend on
$\xi$.

Let us expand the last group of equations 
\[
Y_{\eta_i}(F_2)=0
\] 
in power series in $\xi$. Observe that the coefficient of
$\xi_1\xi_2\xi_3$ vanishes automatically.

The terms of degree-$0$ in $\xi$ are of the form 
\[
\pderf{\alpha_2}{\eta_i}=0 \Longrightarrow
\alpha_2=\alpha_2(u,\chi).
\]
The terms of degree-$1$ in $\xi$ are of the form 
\[
 \xi_k\pderf{\alpha_2}{u_j}- \xi_j\pderf{\alpha_2}{u_k}
+\left(\xi_k\eta_j- \xi_j\eta_k\right)\pderf{\alpha_2}{\chi_1}=
\xi_j\frac{\partial^2F_1}{\partial \eta_i\partial\eta_j}+
\xi_k\frac{\partial^2F_1}{\partial \eta_i\partial\eta_k},
\]
implying that
\begin{equation}
 \label{F1}
 F_1=\mathop{\sum}\limits_{1\leq i\leq 3, \ (i,j,k)\in
A_3}\eta_i\eta_j\pderf{\alpha_2}{u_k} +
\eta_1\eta_2\eta_3\pderf{\alpha_2}{\chi_1} +
\alpha_1(u,\chi)+\mathop{\sum}\limits_{1\leq i\leq 3} f_i(u,\chi)\eta_i.
\end{equation}
Therefore, the generating functions $F_1$ and $F_2$ are completely
determined by the 5 functions 
\[
\text{$\alpha_1, \alpha_2, f_1,f_2,f_3$ which
depend only on $u,\chi$.}
\]

The terms of degree-$2$ in $\xi$ lead to the equation (the
coefficient of $\xi_j\xi_k$):
\begin{equation}
 \label{dgr2xi}
 \mathop{\sum}\limits_{1\leq s\leq 3}\left( \frac{\partial^2F_1}{\partial
u_s\partial\eta_s}- \eta_s\pder{\eta_s}\pderf{F_1}{\chi_1}\right)+
\pderf{F_1}{\chi_1} + \pderf{\alpha_2}{\chi_2} =0.
\end{equation}

Consider this equation by degrees in $\eta$. The degree $0$ yields
the equation
\begin{equation}
 \label{eqfF}
 \mathop{\sum}\limits_{1\leq i\leq 3} (-1)^{p(f_i)}\pderf{f_i}{u_i}+
\pderf{\alpha_1}{\chi_1}+\pderf{\alpha_2}{\chi_2}=0.
\end{equation}
In degrees $1,2,3$ in $\eta$ the equation ~\eqref{dgr2xi} is
automatically satisfied.

Let us express eq.~\eqref{eqfF} in the following more transparent form. Denote
\[
\renewcommand{\arraystretch}{1.4}
\begin{array}{l}f_i=f_i^0+f_i^1\chi_1+f_i^2\chi_2+f_i^{12}\chi_1\chi_2,
\\
\alpha_s=\alpha_s^0+\alpha_s^1\chi_1+\alpha_s^2\chi_2+\alpha_s^{12}\chi_1\chi_2.
\end{array}
\]
Then, equation ~\eqref{eqfF} is equivalent to the system 
\begin{equation}
 \label{eqfFcr}
 \mathop{\sum}\limits_{1\leq i\leq 3} \pderf{f_i^{12}}{u_i}=0,\quad
\alpha_1^{12}-\mathop{\sum}\limits_{1\leq i\leq 3}
\pderf{f_i^2}{u_i}=0,\quad
\alpha_2^{12}+\mathop{\sum}\limits_{1\leq i\leq 3}
\pderf{f_i^1}{u_i}=0,\quad
\alpha_1^1+\alpha_2^2+\mathop{\sum}\limits_{1\leq i\leq 3}
\pderf{f_i^0}{u_i}=0.
\end{equation}

To describe the commutation relations in $\fm\fb(3|8)$ more
explicitly, let us express the vector field ~\eqref{XF} as
\begin{equation}
 \label{XF1}
 X^{F}=x^F+ \mathop{\sum}\limits_{1\leq i\leq 3}\left(f_{\xi_i}Y_{\xi_i}
+f_{\eta_i}Y_{\eta_i}-(-1)^{p(X)}f_{u_i}Y_{u_i} \right),
\text{~where $x^F= F_{1}\partial_{\chi_1}+F_{2}\partial_{\chi_2}$,}
\end{equation}
and observe that, in view of relations ~\eqref{coor1} and ~\eqref{coor3}, we have
\[[X^F, X^G]=X^H,\]
where
\begin{equation}
 \label{H12}
\renewcommand{\arraystretch}{1.4}
\begin{array}{l}
 H_1=[x^F,x^G]_1+\mathop{\sum}\limits_{1\leq i\leq 3}\left(f_{u_i}g_{\eta_i}-(-1)^{p(X^G)}f_{\eta_i}g_{u_i}
 \right),\\
 H_2=[x^F,x^G]_2+\mathop{\sum}\limits_{1\leq i\leq 3}\left(f_{u_i}g_{\xi_i}-(-1)^{p(X^G)}f_{\xi_i}g_{u_i}
 \right).\end{array}
\end{equation}
Observe that it suffices to compute only the \textbf{defining}
components of $H_1$ and $H_2$, \textit{i.e.}, when the pair $F$ is determined
by the set $\alpha_s$, where $s=1,2$, and $f_i$, where $i=1,2,3$,
the pair $G$ is determined by the set $\beta_s$, where $s=1,2$, and
$g_i$, where $i=1,2,3$, and the pair $H$ is determined by the set
$\gamma_s$, where $s=1,2$, and $h_i$, where $i=1,2,3$. Then, with $i=1,2,3$, and $(i,j,k)\in A_3$, we get
\begin{equation}
 \label{gamma1}
\renewcommand{\arraystretch}{1.4}
\begin{array}{ll}
\gamma_1& =\mathop{\sum}\limits_{1\leq i\leq 3}\left(-f_i\pderf{\beta_1}{u_i}+
(-1)^{p(X^G)}\pderf{\alpha_1}{u_i}g_i
 \right)+\\
& 
+\left(\alpha_1\pderf{\beta_1}{\chi_1}+\alpha_2\pderf{\beta_1}{\chi_2}\right)
 -
 (-1)^{p(X^F)p(X^G)}\left(\beta_1\pderf{\alpha_1}{\chi_1}+
 \beta_2\pderf{\alpha_1}{\chi_2}\right);
\end{array}\end{equation}
\begin{equation}
 \label{gamma2}
\renewcommand{\arraystretch}{1.4}
\begin{array}{ll} \gamma_2& =\mathop{\sum}\limits_{1\leq i\leq 3}\left(-f_i\pderf{\beta_2}{u_i}+
(-1)^{p(X^G)}\pderf{\alpha_2}{u_i}g_i
 \right)+\\
& 
+\left(\alpha_1\pderf{\beta_2}{\chi_1}+\alpha_2\pderf{\beta_2}{\chi_2}\right)
 -
 (-1)^{p(X^F)p(X^G)}\left(\beta_1\pderf{\alpha_2}{\chi_1}+
 \beta_2\pderf{\alpha_2}{\chi_2}\right);
\end{array}\end{equation}
\begin{equation}
 \label{hi}
\renewcommand{\arraystretch}{1.4}
\begin{array}{ll} h_i& =-\mathop{\sum}\limits_{1\leq r\leq 3} f_r\pderf{g_i}{u_r}+
\mathop{\sum}\limits_{1\leq r\leq 3}
 \pderf{f_i}{u_r}g_r- \\
& (-1)^{p(X^G)}
 \left(\pderf{\alpha_2}{u_j}\pderf{\beta_1}{u_k}-\pderf{\alpha_2}{u_k}
 \pderf{\beta_1}{u_j}
 -\pderf{\alpha_1}{u_j}\pderf{\beta_2}{u_k} +
 \pderf{\alpha_1}{u_k}\pderf{\beta_2}{u_j}\right)+\\
& 
 +\mathop{\sum}\limits_{s=1, 2} \alpha_s \pderf{g_i}{\chi_s} -
 (-1)^{p(X^F)p(X^G)} \mathop{\sum}\limits_{s=1, 2} \beta_s \pderf{f_i}{\chi_s}.
\end{array}
\end{equation}

In what follows, we identify the vector field $X^F$ with the
collection
\[
\{\alpha_s, f_i\mid s=1,2,\quad i=1,2,3\}.
\]
The bracket of vector fields corresponds to the bracket of such
collections given by eqs. ~\eqref{gamma1}--\eqref{hi}.

Consider now the even part $\fm\fb(3|8)_\ev$ of our algebra. Since
$p(F_1)=p(F_2)=\od$, it follows that $p(\alpha_s)=\od$,
$p(f_i)=\ev$. The component $\fm\fb(3|8)_\ev$ admits a~decomposition into three
subspaces described below
\[
\fm\fb(3|8)_\ev=V_1\oplus V_2\oplus V_3.
\]

The subspace $V_1$ is determined by the collection
\[
 \{\alpha_1=\alpha_2=0,\;\;f_i=f_i(u)\chi_1\chi_2\mid \mathop{\sum}
\pderf{f_i}{u_i}=0\}.
\] 
Eqs. ~\eqref{gamma1}--\eqref{hi} imply that
the vector fields generated by such functions form a~commutative
ideal in $\fm\fb(3|8)_\ev$; we will identify this ideal with
$d\Omega^1(3)$ via the mapping
\[
 \{0,0,f_i\mid i=1,2,3\} \longmapsto -\mathop{\sum}\limits_i
f_idu_jdu_k.\]

Next, the subspace $V_2$ is determined by the collection $\{f_i=0\mid
i=1,2,3\}$. We will identify this space with $\Omega^0(3)\otimes
\fsl(2)$ via the mapping
\[
 \{\alpha(u)(a\chi_1+b\chi_2),\alpha(u)(c\chi_1-a\chi_2),0,0,0\}\longmapsto 
\alpha(u)\otimes \begin{pmatrix} a& c\\b&d \end{pmatrix}.
\]
Eqs. ~\eqref{gamma1}--\eqref{hi} imply that the elements of the subspaces $V_1$ and
$V_2$ commute, whereas the bracket of two
collections from $V_2$ is, in our notation, of the form
\[
[f\otimes A, g\otimes B]=fg\otimes [A,B]+df\wedge dg \cdot \tr AB.
\]

Finally, for $V_3$ we have the following three natural possibilities: we take $f_i=f_i(u)$ for any $i=1,2,3$, while for the
$\alpha$ we make the following choices:
\begin{description}
\item[(a)] $ \alpha_1=-\sum\pderf{f_i}{u_i}\chi_1$,\quad $\alpha_2=0$,
\item[(b)] $ \alpha_1=0$,\quad $ 
\alpha_2=-\sum\pderf{f_i}{u_i}\chi_2$,
\item[(c)] $ \alpha_1=-\frac 12 \sum\pderf{f_i}{u_i}\chi_1$,
\quad $ \alpha_2=-\frac 12
\sum\pderf{f_i}{u_i}\chi_2$.
\end{description}
In the cases (a) and (b), we get two embeddings: $\fvect(3)\subset
\fm\fb(3|8)_\ev$.

The case (c) is, however, more convenient to simplify the
brackets. Thus, we identify $V_3$ with $\fvect(3)$ via
\[
 \left\{-\frac 12 \sum\pderf{f_i}{u_i}\chi_1,\ -\frac 12
\sum\pderf{f_i}{u_i}\chi_2, \ f_1(u),\ f_2(u),\ f_3(u)\right\}
\longmapsto
D_f=-\mathop{\sum}\limits_i f_i(u)\partial_{u_i}. 
\]
The action of $D_f$ on the subspaces $V_1$ (as on the space of
$2$-forms) and $V_2$ (as on the space of vector-valued functions
$\cF\otimes\fsl(2)$) is natural. The bracket of the operators $D_f$
is, however, manifestly different from the usual bracket thanks to the presence of an
extra term:
\[
 [D_f,D_g]=D_fD_g-D_gD_f-\nfrac 12 d(\Div D_f)\wedge d(\Div D_f).
\]

Consider now the odd part $\fm\fb(3|8)_\od$. We have
$p(F_1)=p(F_2)=\ev$, and hence
\[
p(\alpha_s)=\ev,\;\; p(f_i)=\od.
\]

Let $V_4$ consist of the collections with $f_i=0$. We identify $V_4$
with $\Omega^0(3)\vvol^{-1/2}\otimes \Cee^2$ via
\[
 \{\alpha(u)w_1, \ \alpha(u)w_2, 0,0,0\}\longmapsto
\alpha(u)\vvol^{-1/2}\otimes \begin{pmatrix} w_2 \\
-w_1\end{pmatrix}.
\]

Let $V_5$ consist of the collections of the form
\begin{equation}
 \label{V5}
\renewcommand{\arraystretch}{1.4}
\begin{array}{l}
f_i=f_i(u)(v_1\chi_1+v_2\chi_2)~~ \text{for $i=1,2,3$},\quad
 \alpha_1=v_2\sum\pderf{f_i}{u_i}\chi_1\chi_2,\quad
 \alpha_2=-v_1\sum\pderf{f_i}{u_i}\chi_1\chi_2.\end{array}
\end{equation}
We identify $V_5$ with $\Omega^2(3)\vvol^{-1/2}\otimes \Cee^2$, by
assigning to the collection ~\eqref{V5} the element
\[
 \omega\vvol^{-1/2}\otimes \begin{pmatrix} v_1\\v_2 \end{pmatrix},
\text{ ~where~ } \omega=-\sum f_idu_j\wedge du_k.
\]
Having identified $\fm\fb(3|8)_\od=V_4\oplus V_5$, let us summarize. Table \eqref{tablsoo} contains a~description of the spaces $V_i$ and their elements.
\begin{table}[ht]\centering
{\footnotesize%
\begin{equation}
\label{tablsoo}
\begin{tabular}{|c|c|c|c|c|}
 \hline
The space & $\alpha_1$ &$\alpha_2$ & $f_i$ & \\
 \hline
 $V_1\simeq d\Omega^1(3)$ & $0$ & $0$ &
 $f_i(u)\chi_1\chi_2,\; $ & $\omega=\sum
 f_idu_j\wedge du_k, \; $\\
 & & & $\sum\pderf{f_i}{u_i}=0$ &$d\omega=0$\\
 \hline
 $\begin{matrix}V_2\simeq \\
 \Omega^0(3)\otimes \fsl(2)\end{matrix}$ &
$ \alpha(u)(a\chi_1+b\chi_2)$ & $ \alpha(u)(c\chi_1-a\chi_2)$ & 0
& $\alpha(u)\otimes \begin{pmatrix} a&c\\b&-a \end{pmatrix}$\\
 \hline
 $V_3\simeq $ & $-\frac 12 f(u)\chi_1$ & $-\frac 12 f(u)\chi_2$
 & $f_i(u)$ & $D=-\sum f_i(u)\partial_{u_i}$ \\
$ \fvect(3)$ & & & $f(u)=\sum\pderf{f_i}{u_i}$ & $\Div
D=-f(u)$\\
 \hline
 \hline
 $\begin{matrix}V_4\simeq \\
 \Omega^0\vvol^{-1/2}\otimes \Cee^2\end{matrix}$& $ \alpha(u)w_1$ & $ \alpha(u)w_2$ & $0$ &
 $\alpha(u)\vvol^{-1/2}\otimes \begin{pmatrix} w_2 \\ -w_1
 \end{pmatrix}$\\
 \hline
 $V_5\simeq $ &$ v_2 f(u)\chi_1\chi_2$ & $ -v_1 f(u)\chi_1\chi_2$
 &$f_i(u)(v_1\chi_1+v_2\chi_2)$ & $\omega \vvol^{-1/2}\otimes \begin{pmatrix} v_1 \\
 v_2 \end{pmatrix}$\\
$ \Omega^2\vvol^{-1/2}\otimes \Cee^2$ & &
&$f(u)=\sum\pderf{f_i}{u_i}$ & $\omega=\sum
 f_idu_j\wedge du_k$\\
 \hline\hline
\end{tabular}
 \end{equation}\vspace{-4mm}
 }
\end{table}

Having explicitly computed the brackets using expressions
 ~\eqref{gamma1}--\eqref{hi} and presenting the result by means of correspondences
~\eqref{tablsoo}, we obtain formulas almost identical to those
provided in \cite{CCK}.

We have already given the brackets of the even elements. The
brackets of elements of $\fm\fb_\ev$ and $\fm\fb_\od$ are as follows:
\[\footnotesize
\renewcommand{\arraystretch}{1.4}
\begin{array}{ll}
{}[V_1,V_4]: & [\omega, \ \alpha\vvol^{-1/2}\otimes v]=\alpha\cdot
\omega\vvol^{-1/2}\otimes v \in V_5;\\
{} [V_2,V_4]:& [f\otimes A,\ \alpha\vvol^{-1/2}\otimes
v]=(f\alpha-df\wedge
d\alpha)\vvol^{-1/2}\otimes Av \in V_4\oplus V_5;\\
{}[V_3,V_4]:& [D,\ \alpha\vvol^{-1/2}\otimes v]=(D(\alpha)-\frac 12
\Div D\cdot \alpha +\frac 12 d(\Div D)\wedge
d\alpha)\vvol^{-1/2}\otimes
v \in V_4\oplus V_5;\\
{}[V_1,V_5] &=0; \\
{}[V_2,V_5]: & [f\otimes A,\ \omega \vvol^{-1/2}\otimes v]= f\omega
\vvol^{-1/2}\otimes Av\in V_5;\\
{}[V_3,V_5]: & [D,\ \omega \vvol^{-1/2}\otimes v]= (L_D\omega-\frac 12
\Div D\cdot \omega) \vvol^{-1/2}\otimes v\in V_5.
\end{array}
\]
To describe in these terms the bracket of two odd elements, we resort to 
the following natural identifications:
\[\footnotesize
\renewcommand{\arraystretch}{1.4}
\begin{array}{l}
\frac{\Omega^2(3)}{\vvol}\simeq \fvect(3):\quad \frac{\omega}{\vvol}
\lra D_\omega:\\
i_{D_\omega}(\vvol)=\omega, \text{ ~\textit{i.e.},~ } \sum f_idx_j\wedge dx_k
\lra \sum f_i\partial_{i_i};\\
\Lambda^2\Cee^2\simeq \Cee:\; v\wedge w\lra \det
\begin{pmatrix}v_1 & w_1\\v_2&w_2 \end{pmatrix} ; \\
S^2(\Cee^2)\simeq \fsl(2): \; v\cdot w \lra \begin{pmatrix}
-v_1w_2-v_2w_1 & 2v_1w_1\\-2v_2w_2 & v_1w_2+v_2w_1 \end{pmatrix} .
\end{array}\]
The bracket of two odd elements is of the form:
\begin{equation}
 \label{lazha}
\renewcommand{\arraystretch}{1.4}
\begin{array}{ll}
{}[V_4,V_4]: & \left[\frac{f}{\sqrt{\vvol}}\otimes v, \frac
{g}{\sqrt{\vvol}}\otimes w\right]=\frac{df\wedge dg \otimes v\wedge
w}{\vvol}\in\fvect(3)=V_3;\\
{}[V_5,V_5]: &\left[\frac{\omega_1}{\sqrt{\vvol}}\otimes v,
\frac{\omega_2}{\sqrt{\vvol}}\otimes w\right]=
(D_{\omega_1}(\omega_2)-\Div
D_{\omega_1}\cdot \omega_2)v\wedge w\in V_1;\\
{}[V_4,V_5]: &\left[\frac{f}{\sqrt{\vvol}}\otimes
v,\frac{\omega}{\sqrt{\vvol}}\otimes w\right]=\\
&\frac{f\omega}{\vvol}\otimes v\wedge w - \frac
12(fd\omega-\omega\wedge df)\otimes v\cdot w +df\wedge d(\Div
D_\omega)\otimes v\wedge w.\end{array}
\end{equation}
In the last line above, the first summand lies in $V_3$, the second
one in $V_2$, and the third one in~ $V_1$. \textbf{The difference as
compared with} \cite{CCK}: the coefficient of the third summand in
the last line on ~\eqref{lazha} should be $1$, whereas in \cite{CCK}
it is equal to $\frac12$.

To verify, compute the Jacobi identity (it holds for 1 and not for
$\frac12$) for the triple
\[
 u_3du_2\wedge du_3\in V_1\text{~and~ $\nfrac{u_1}{\sqrt{\vvol}}\otimes e_1$,
 $\nfrac{u_2}{\sqrt{\vvol}}\otimes e_2\in V_4$, where $e_1, e_2$ span
 $\Cee^2$.}\hskip 2cm \qed
 \]

\section{Miscellenea}\label{S:9}

\ssec{Odd roots of simple Lie superalgebras of non-vectorial type} For this description for simple Lie superalgebras and their relatives, see table \eqref{8.1}.
{\footnotesize
\begin{equation}\label{8.1}
\renewcommand{\arraystretch}{1.4}
\begin{tabular}{|c|c|c|c|}
\hline $\fg$&$\fg_{\ev}$&$R_{\od}$&$\fg_{\ev}: \fg_{\od}$\cr \hline

$\fsl(n|1)$&$\fsl(n)\oplus\Cee h$&$\pm (\eps_{i}-\delta)$&$\left
(\id\boxtimes\chi_{0}\right)\oplus\left
(\id\boxtimes\chi_{0}\right)^{*}$\cr \hline

$\fsl(n|m)$&$\fsl(n)\oplus\fsl(m)\oplus\Cee h$&$\pm
(\eps_{i}-\delta_{j})$&$\left
(\id_{\fsl(n)}\boxtimes\id_{\fsl(m)}^*\boxtimes
\chi_{0}\right)\oplus$\cr

$n\neq m$&&&$\oplus\left
(\id_{\fsl(n)}\boxtimes\id_{\fsl(m)}^*\boxtimes\chi_{0}\right)^{*}$\cr 
\hline $\fpsl(n|n)$&$\fsl(n)\oplus\fsl(n)$&$\pm
(\eps_{i}-\delta_{j})$&$\left
(\id_{\fsl(n)}\boxtimes\id_{\fsl(n)}^*\right)\oplus$\cr 

$n>2$&&&$\oplus\left
(\id_{\fsl(n)}\boxtimes\id_{\fsl(n)}^*\right)^{*}$\cr
\hline 

$\fosp(1|2n)$&$\fsp(2n)$&$\pm \delta_{j}$&$\id_{\fsp(2n)}$\cr
\hline 

$\fosp(2|2n)$&$\fsp(2n)\oplus\Cee h$&$\pm \eps\pm\delta_{j}$&
$\left(\id_{\fsp(2n)}\boxtimes\chi_{1}\right)\oplus\left
(\id_{\fsp(n)}^*\boxtimes\chi_{-1}\right)$\cr \hline

$\fosp(m|2n)$&$\fsp(2n)\oplus\fo(m)$&$\pm \eps_{i}\pm\delta_{j}$ for
$m$ even&$\id_{\fsp(2n)}\boxtimes\id_{\fo(m)}$\cr 

$m>2$&&$\pm
\eps_{i}\pm\delta_{j}$ and $\pm \delta_{j}$ for $m$ odd&\cr \hline

$\fosp_a(4|2)$&$\fsl_{1}(2)\oplus\fsl_{2}(2)\oplus\fsl_{3}(2)$&$\pm
\eps_{1}\pm \eps_{2}\pm
\eps_{3}$&$\id_{1}\boxtimes\id_{2}\boxtimes\id_{3}$\cr \hline

$\fpsq(n)$&$\fsl(n)$&$\pm (\eps_{i}-\eps_{j})$ and $0$&$\ad$\cr

\hline $\fa\fg_{2}$&$\fsl(2)\oplus\fg_{2}$&$\pm \eps_{i}\pm \delta;
\pm \delta$&$\id_{\fsl(2)}\boxtimes R(\pi_{1})$\cr &&$i=1, 2, 3;\;
\sum \eps_{i}=0$&\cr \hline

$\fa\fb_{3}$&$\fsl(2)\oplus\fo_{7}$&$\frac{1}{2}(\pm \eps_{1}\pm
\eps_{2}\pm \eps_{3}\pm \delta)$&$\id_{\fsl(2)}\boxtimes \spin_{7}$\cr
\hline
\end{tabular}
 \end{equation}\vspace{-4mm}
 }


\ssbegin[Central extensions: $\fg$ simple finite-dimensional]{Theorem}[Central extensions of simple finite-dimensional Lie superalgebras (\cite{L3})]\label{thcentext} The
following list exhausts the central extensions of simple finite-dimensional Lie superalgebras:
\be\label{FDcoc}
\footnotesize
\begin{tabular}{|c|c|c|}
\hline
$\fg$&\emph{Cocycle}&\emph{Name of}\\
&&\emph{extended superalgebra}\\
\hline
$\begin{matrix}\fpsl(2|2)\\
\simeq \fh'(0|4)\end{matrix}$&$\mat{A&B\\ C& D},
\mat{ A'&B'\\ C'& D'}
\longmapsto \left\{\begin{matrix}\tr BB'\\
\tr (BC' +CB' )\\
\tr CC'\end{matrix}\right.$&$\begin{matrix}\fpo'(0|4)\\
\fsl (2|2)\\
\fpo'(0|4)\end{matrix}$\\
\hline 
$\begin{matrix}\fpsl(n|n)\\
\emph{for $n>2$}
\end{matrix}$&$\mat{ A&B\\ C& D}, \mat{ A'&B'\\ C'&
D'}
\longmapsto \tr (BC' +CB')$&$\fsl (n|n)$\\
\hline
$\fpsq (n)$&$\mat{ A&B\\ B& A},
\mat{ A'&B'\\ B'& A'}
\longmapsto \tr BB'$&$\fsq (n)$\\
\hline
$\fspe (4)$&$\mat{ A&B\\ C&
 {-A}^t},
\mat{ A'&B'\\ C'& {-A'}^{t}}
\longmapsto\tr CC'$& $\fas$\\
\hline $\fh'(0|n), ~~n>4$&$H_{f}, H_{f' }\longmapsto \sum {\left.
(\pderf{f}{ \theta
_{j}} \pderf{f'}{\theta_{j}})\right|_{\theta =0}}$, ~~\emph{see }\eqref{2.3.2}&$ \fpo'(0|n)$\\
\hline
\end{tabular}
\ee
\end{Theorem}

For, sometimes unexpected, analogs of these results over fields of positive characteristic, see \cite{BGL1, BGLL1}.

\begin{proof}[Proof \nopoint] (Union of strewn results and ideas). The above table and some of the results from \cite{FuLe}, backed up by Gruson, take care of the superalgebras $\fg$
with reductive $\fg_\ev$; for the corrected previous results concerning non-vectorial Lie superalgebras, see \cite{BoKN}. 

For vectorial Lie superalgebras, the Hochshield-Serre spectral sequence with respect to the
subalgebra $\fg_0$ in the standard grading does the trick. Grozman used his
\textit{SuperLie} code package to confirm these calculations for small
values of the rank, where extra cocycles might have occurred. 
See also \cite{BLS}, where the case of the ground field of characteristic $p>0$ is mainly considered and the notion of \textit{double extension}\index{Extension, double} is thoroughly studied.
\end{proof}

\ssbegin[Central extensions: $\fg$ simple infinite-dimensional]{Theorem}[Central extensions of simple infinite-dimensional vectorial Lie superalgebras]\label{bezN2} The following list exhausts the central
extensions of simple infinite-dimensional vectorial superalgebras:
\[
\footnotesize
\begin{tabular}{|c|c|c|}
\hline
\textup{Algebra}&\textup{Cocycle}&\textup{Name of the extended algebra}\\
\hline $\begin{array}{l}\fh(2n| m)\\
\text{for $(n,m)\neq (1,2)$}
\end{array}$&$ f, g\longmapsto\left.{\{f,
g\}_{\text{P.B.}}}\right| _{p=q=\theta=0}$&
$\fpo(2n|m)$\\
\hline
$\left.\begin{matrix}\fle (n)\\
\fsle'(n)\end{matrix}\right\}$&$\left.{\{f,
g\}_{\text{B.B.}}}\right|
_{q=\xi=0}$&$\left \{
\begin{matrix}\fb (2n| m)\\
\fsb'(n)\end{matrix} \right.$\\
$\fsle'(3;3)$&\textup{2 non-cohomologous cocycles, see} \cite{ShP}&$\fsb'(3;3)$\\
\hline
\end{tabular}
\]
\end{Theorem}


\begin{proof}[Proof \nopoint] (the ideas of). 
The Lie superalgebra $\fsle'(3)$ has two non-cohomologous (odd) central extensions producing
isomorphic Lie superalgebras, see \cite{ShP}.

If $\fg$ has a~grading operator $z$, then it suffices to
consider the space of $z$-invariant cocycles, \textit{i.e.}, elements of 
$(\fg_{-1}\wedge \fg_1)^*$ which constitute, moreover, a~trivial
$\fg_0$-module. This case covers the series $\fg=\fvect$, $\fk$, $\fm$
and the exceptional algebras $\fv\fle(4|3; K)$, $\fk\fas$,
$\fm\fb(4|5; K)$. 

From the explicit description of $\fg_0$-modules $\fg_i$ for $i\leq
1$ we see that such elements only exist for $\fvect$ and $\fk$ (when
$\fg_1$ contains a~direct summand isomorphic to $\fg_{-1}^*$), but
they do not contribute to cohomology.
\end{proof}

\ssec{The extensions $\fe\fg$, and algebras $\fder~\fg$ for simple ``symmetric\rq\rq\ $\fg$ (\cite{BGLL1})}

{\footnotesize%
\begin{equation}\label{8.2}
\renewcommand{\arraystretch}{1.4}
\begin{tabular}{|c|c|c|c|}
\hline $\fg$&$\fpsl(n|n)$, $n>2$&$\fpsl(2|2)$&$\fp\fs\fq(n)$\cr
\hline
$\fe\fg$&$\fsl(n|n)$&$\left\{\begin{matrix}\fsl(2|2)\\
\fpo'(0|4)\\ \fosp_{0}(4|2)\end{matrix}\right .$&$\fs\fq(n)$\cr \hline
$\fder~\fe\fg$&$\fgl(n|n)$&$\left\{\begin{matrix}\fgl(2|2)\\
\fpo(0|4)\ltimes\Cee E\\
\end{matrix}\right .$&$\fq(n)$\cr \hline

$\fder~\fg$&$\fp\fgl(2|2)$&$\fosp_{-1}(4|2)$&$\fp\fq(n)$\\
\hline
\end{tabular}
 \end{equation} 
 }

One of the Cartan
matrices of $\fosp_a(4|2)$ over $\Cee$ is $\footnotesize{\begin{pmatrix} 2 & -1 & 0 \\ -1 & 0
& -a \\ 0 & -1 & 2
\end{pmatrix}}$, as we know (\textit{e.g.}, from \cite{BGL1}), so $\fosp_1(4|2)=\fosp(4|2)$, and $\fosp_a(4|2)$ is simple
except for the cases where $a=0$ and $-1$, where we have exact sequences
\begin{equation}\label{exSeqOsp4,2,a}
\begin{array}{l}
0\tto \fpsl(2|2)\tto \fosp_{-1}(4|2)\tto \fsl(2)\tto 0,\\

0\tto \fsl(2)\tto \fosp_0(4|2)\tto \fpsl(2|2)\tto 0.\\
\end{array}
\end{equation}
In other words, $\fder~\ \fpsl(2|2)\simeq \fosp_{-1}(4|2)$ and on the
space of three nontrivial central extensions of $\fpsl(2|2)$, as well on the space of outer derivations of $\fpsl(2|2)$, one can define a~
structure of the Lie algebra $\fsl(2)$.

\ssec{Finite-dimensional irreducibles over simple vectorial
Lie superalgebras}\label{SS:9.1} 
Let $\cL$ be either $\fvect(0|m)$, or its
subsuperalgebra $\fsvect(0|m)$, or $\fh'(0|m)$ with their
standard filtrations.

Let us prepare a~supply of finite-dimensional irreducible
$\cL$-modules. Take an irreducible $\cL_{0}/\cL_1$-module $M^\chi$
with highest weight $\chi$; assume that the highest weight vector
(as well as the whole $M^\chi$) is even. Construct the induced
$\cL$-module by setting
\[
I(\chi):=\ind^\cL_{\cL_{0}}(M^\chi)=U(\cL)\otimes_{U(\cL_{0})}M^\chi.
\]
It turns out that, for $\chi$ generic, the modules $I(\chi)$ are
irreducible; the exceptional --- reducible --- modules line up into
a~sequence which is exact for $\fvect(0|m)$ and
$\widetilde{\fsvect}(0|m)$ and ``almost exact'' for $\fsvect(0|m)$,
see \cite{BL3}.

To formulate the answer graphically, let us pass from the induced
modules $I(\chi)$ to their duals, the coinduced modules
$T(\chi):=I(\chi^*\otimes(-\tr))$, where $\tr$ is the
$1$-dimensional $\fgl(m)=\cL_{0}$-module determined by the trace also denoted by $\tr$,
and $(-\tr)(X)=-\tr(X)$ for any $X\in \fgl(m)$.

\sssbegin[Irreducibles over $\fvect(0|m)$, $\fsvect(0|m)$ and $\widetilde{\fsvect}(0|m)$]{Theorem}[Irreducibles over $\fvect(0|m)$, $\fsvect(0|m)$ and $\widetilde{\fsvect}(0|m)$, after \cite{BL3}]\label{thBL}
Let $\cL=\fvect(0|m)$, $\fsvect(0|m)$, or $\widetilde{\fsvect}(0|m)$.
Then, the $\cL$-modules $T(\chi)$ are irreducible, except for the
cases where $M^\chi$ is isomorphic to
$S^k(\id)$ or $S^k(\id^*)\otimes (-\tr)$, where $\id$ is the
tautological representation of $\cL_{0}/\cL_1$. These exceptional
modules, denoted by $\Omega^k$ and $\Sigma_{-k}$, respectively, are
interpreted as the spaces of $k$th degree differential forms and,
respectively, $(-k)$th degree integral\footnote{Having discovered
them, Bernstein and Leites called them \textit{integraBle} since
they can be integrated. P.~Deligne, however, calls them
\textit{integral}, cf. \cite[Ch.3, \S3.12]{Del}; perhaps to rhyme with differential or after the incorrect (since the authors' interpretation was lost) English translation of the Russian version of \cite{BL1}.} forms on the $(0|m)$-dimensional
supermanifold. The exterior differential and the Berezin integral
unite the exceptional modules in one --- infinite in both directions ---
exact sequence: \[ 
\cdots\stackrel{d}{\tto}
\Sigma_{-k}\stackrel{d}{\tto} \cdots \stackrel{d}{\tto}
\Sigma_{-1}\stackrel{d}{\tto} \Sigma_{0}\stackrel{\int}{\tto}
\Omega^0\stackrel{d}{\tto} \Omega^1\stackrel{d}{\tto} \cdots
\stackrel{d}{\tto} \Omega^k\stackrel{d}{\tto} \cdots , 
\] 
and
therefore the kernels and the images of the above maps are the only
irreducible submodules of the exceptional modules if $\cL$ is of
series $\fvect$ or $\widetilde{\fsvect}$.

For $\cL=\fsvect(0|m)$, since $\Sigma_{0}\simeq \Omega^0$, the module
of functions, $\Omega^0$, has the additional submodule $\Vol_{0}$ consisting of the functions with integral zero.
\end{Theorem}

We see that, except for $m=2$ (since
$\fvect(0|2)\simeq\fsl(1|2)\simeq\fosp(2|2)$, this case was already
considered), a~rank-1 operator can occur only in the modules
$T(\chi)$ coinduced from the 1-dimensional module, or in a~quotient
of such module, or in the modules obtained from such module by reversing their parity. The corresponding cases are listed in Table
~\eqref{rank1T}. For example, the corresponding rank-1 operators from
$\fvect(0|m)$ and $\widetilde \fsvect(0|m)$ are 
\[
\renewcommand{\arraystretch}{1.4}
 \xi_{1}\dots \xi_{m}\partial_{\xi_{i}}\; \text{ for any $i$ realized
 in}\begin{cases}\text{$\cF:=\Lambda(\xi)$ and $\Vol^\lambda(\xi)$},&\text{for $m$
even and $\lambda\in\Cee$},\\
\text{$\Pi(\cF):=\Pi(\Lambda(\xi))$ or $\Pi(\Vol^\lambda(\xi))$},&\text{for $m$
odd and $\lambda\in\Cee$}.\end{cases} 
\]
In the case of $\fsvect(0|m)$, the operators $(\prod_{j\neq i, t}\xi_{j})\pder{\xi_{i}}$ acting on $\cF:=\Lambda(\xi)$ have rank 1.

Several of our exceptional vectorial superalgebras exist due to the
fact that adding the trivial center to the above-considered $\cL$ we
get, in certain small dimensions, new nontrivial Cartan prolongs. For details and such new (in 1996) examples, see \cite{Sh14}.

Let $\cL=\fh(0|m)$ preserve the symplectic 2-form $\omega$. Then, the
exceptional modules are still realized in the spaces of 
differential forms (since $\Sigma_{-i}\simeq \Omega ^i$ thanks to
nondegeneracy of $\omega$). The complete description of irreducible
modules over $\cL=\fh'(0|m)$ and all its relatives --- central
extensions (\textit{i.e.}, finite-dimensional Poisson Lie superalgebras) and
superalgebras of derivations is due to Shapovalov
\cite{Sha1,Sha2}. Except for the cases of small dimension where an
occasional isomorphism occurs (such as $\fh'(0|4)\simeq \fp\fsl(2|2)$),
and which, therefore, are already considered, there are no operators
of rank 1.

\ssec{On primitive Lie superalgebras}\label{SS:10.2} Recall that
a~filtered Lie (super)algebra $\fg$ (or an associated graded
algebra) is called \textit{primitive}\index{Lie (super)algebra, primitive} if it has a~maximal subalgebra of finite codimension which does not contain a~ homogeneous ideal of
$\fg $.

Let
\begin{equation}
\label{pr-1} \fr\subset\fvect(\xi)\text{~~ be a~subalgebra such that
$\pr_{-1}(\fr)=\fvect_{-1}(\xi)$,} \end{equation} where $\pr_{-1}$ is
the projection of $\fvect(\xi)$ to
$\fvect_{-1}(\xi)=\Span(\pder{\xi_{1}}, \dots , \pder{\xi_{n}})$.

Let $\fh=\mathop{\oplus}_{i\geq -d}\ \fh_i$ be
a~\textit{primitive} Lie superalgebra (with
$\mathop{\oplus}_{i\geq 0}\ \fh_i$ being its maximal
subalgebra) such that $\fh_{-1}$ is an irreducible $\fh_{0}$-module
and $\fh_{-}$ is generated by $\fh_{-1}$. Set:
\begin{equation}
\label{(*)}
\renewcommand{\arraystretch}{1.4}
\begin{array}{l}
\fg_{min}=\fh\otimes \Lambda(n)\ltimes\fr, \ \ \ 
\fg_{max}=(\fh_{-1}\otimes \Lambda(n),\ \fh_{0}\otimes
\Lambda(n)\ltimes \fr)_{*}, 
\end{array}
\end{equation}
where $(V, \fg)_{*}$ is the Cartan prolong of $(V, \fg)$. It is not
difficult to see that any Lie superalgebra $\fp$ such that
$\fg_{min}\subseteq \fp\subseteq\fg_{max}$, is primitive.

\ssbegin[On primitive Lie superalgebras]{Conjecture}[On primitive Lie superalgebras] Every primitive Lie superalgebra $\fp$ is of the form $\fg_{min}$ for
a~suitable $\fh$. \end{Conjecture}\index{Conjecture} 

Regardless of whether this Conjecture is true or not, ``describe primitive Lie
superalgebras" is a~wild problem because there are indescribably many Lie superalgebras
satisfying condition ~\eqref{pr-1}.

\ssec{On the Lie superalgebra $E(2, 2)$ from \cite{K3}}\label{SS:10.3}
While this paper was being written and \cite{Sh5, Sh14} refereed, there
appeared a~proof of our classification in the case of compatible grading, see ~\cite{Klc}. In it,
V.~Kac states, among other things, that there exist six
exceptional simple vectorial Lie
superalgebras. In addition to the five (considered as abstract Lie
superalgebras) exceptions already described in \cite{Sh5,Sh14}, Kac
distinguished one more ``exception'', $E(2,2)$. Actually, it belongs to a series.

\sssbegin[$E(2, 2)\simeq\fb'_{\infty}(2; 2)$]{Lemma}\label{lE22} $E(2, 2)\simeq\fb'_{\infty}(2; 2)$. \end{Lemma}

\begin{proof} Indeed, by definition, $E(2, 2)$ is the Cartan prolong
$\fg=(\fg_{-1}, \fg_0)_*$, where ``the $\fg_0$-module $\fg_{-1}$
gives rise to the linear Lie superalgebra $\text{spin}_4$". This
means that $\fg_{0}$ is the semidirect sum
$\fg_{0}=\fhei(0|4)\ltimes \fsl(2)$ and $\fg_{-1}$ is the space of the 
unique irreducible representation of $\fhei(0|4)$ for which the
central element corresponds to the identity operator. Hence, we may
identify $\fg_{-1}$ either with $\Lambda(2)=\Lambda(\xi_1,\xi_2)$ or
(due to uniqueness) with
$\Pi(\Lambda(2))=\Pi(\Lambda(\xi_1,\xi_2))$. The second scenario is
more convenient for our purposes.

The action of $\fhei(0|4)$ is given by the operators $(1, \xi_1, \xi_2,
\partial_{\xi_1}, \partial_{\xi_2})$, while $\fsl(2)$ acts via the
tautological representation on the even component of $\fg_{-1}$ and
via the trivial representation on the odd component of $\fg_{-1}$.
Now, what is $\fb'_{\infty}(2; 2)$ or,
in Kac's notation, ``$SKO(2, 3; 1)$ with subprincipal grading''? It is the Cartan prolong
$\fh=(\fh_{-1},\fh_0)_*$, where
\[
\begin{array}{l}
\fh_{-1}=\Pi(\Lambda(2))=\Pi(\Lambda(\eta_1,\eta_2));\\
\fh_0=\Span(1,\eta_1,\eta_2)\ltimes \fsvect(0|2) \simeq \left
(\Span(1,\eta_1,\eta_2)\ltimes \Span(\partial_{\eta_1},
\partial_{\eta_2}) \right )\ltimes \fsl(2)\simeq \fg_0.
\end{array}
\]
On $\Pi(\Span(\eta_1,\eta_2))$, \textit{i.e.}, on $(\fh_{-1})_\ev$, the subalgebra $\fsl(2)$ acts via
the tautological representation, and on 
$\Pi(\Span(1, \eta_1\eta_2))$, \textit{i.e.}, on 
$(\fg_{-1})_\od$, it acts via the trivial representation.

We get an isomorphism $\fg_{-1}\oplus\fg_0\simeq 
\fh_{-1}\oplus\fh_0$, and hence an isomorphism $\fg\simeq \fh$. \end{proof} 

In the paper \cite{CK2} which corrects several statements of
\cite{K3}, the result of this lemma is briefly stated, but the
description of $\fg_0$ in \cite{CK2} should still be corrected. The
paper \cite{CCK}, in its turn, corrects \cite{CK2}, but not completely: see eq.~\eqref{lazha}.


\section{On history}\label{SS:1.11}

Lie superalgebras first appeared in topology, in the 1940s, over finite fields, in (co)homology theory. The collection of elements of homotopy groups with the Whitehead multiplication is an example of the Lie superring. These examples drew so little attention of researchers that the difference between ``Lie superalgebras" and ``graded Lie algebras" was clearly formulated only after \cite{WZ} introduced super terminology.


Since the 1960s, Berezin tried to describe Bose and Fermi particles in a~parallel manner. He conjectured that to do so one has to consider anticommuting indeterminates on equal footing with commuting coordinates. To see that the classical differential and algebraic geometries are particular cases of a~ wider area of mathematics, one had to define the object on which these indeterminates play the role of functions. This problem (proof of Berezin's conjecture) was solved in \cite{L0}, where what is now called \textit{superscheme} was defined. For further mathematical development in smooth and analytic settings, see \cite{B3, LSoS, MaG, Ld}. The vital inspirational role of Berezin in several areas of mathematics and physics is described in \cite{Shi, N, KNV, KMLT}. 

Meanwhile physicists arrived at the new fundamental philosophical paradigm transcending the limits of pure physics: ``\textbf{We live on a~supermanifold}". This idea is still being developed; it became contagious thanks to a~ work by Wess and Zumino \cite{WZ} where the term ``SUSY" was introduced in the currently accepted sense, different from the now abandoned Wheeler's usage of it, see \cite{Whe}. In \cite{WZ}, perspectives of \textit{supersymmetry} were outlined more clearly than in several pioneering works on the subject, see introductions in \cite{DSB}, and also \cite{Del, Shi} written decades after \cite{WZ} was published. Note that applications of SUSY to high energy physics are still conjectural, whereas results in solid state physics are already handable, see \cite{Ef}.

\textbf{Simple finite-dimensional Lie superalgebras over $\Cee$}. After \cite{WZ} was published, the perspectives in theoretical physics and manifest mathematical beauty elucidating certain classical facts encouraged several researchers to tackle the classification of simpleLie superalgebras. Various
\textit{particular} types (with non-degenerate even
form, with reductive even part) were classified by Djokovi\'c (with
Hochschild) \cite{DH, Dj1, Dj2, Dj3}, Scheunert, Nahm, and Rittenberg \cite{SNR},
I.~Kaplansky (\cite{FrK}, \cite{Kap, Kapp}). Note that there was no understanding of Lie superalgebras as algebras of operators until mid-1970s when $\fgl(m|n)$ was defined (\cite{Lber, Kapp, SNR}) and the space of superderivations of the Grassmann algebra was interpreted as a~(vectorial) Lie superalgebra. 
In the talk \cite{K0}, see Kac's preface to \cite{DSB}, V.~Kac gave the first examples of (relatives of) simple finite-dimensional vectorial Lie
superalgebras. 
Kac called them ``some algebras"; we now
identify them as $\fvect(0|m)$, $\fsvect(0|m)$ and $\fh(0|m)$. 
Several years later, Kac considered \textit{all} types of simple finite-dimensional Lie superalgebras. 
Kaplansky's preprints-1975 \cite{Kapp} helped Kac to repair gaps in his 
classification (\cite{K1}, \cite{K1C}). Various types of simple Lie
superalgebras have to be treated separately, so it is
instructive to compare different techniques of \cite{SNR} and
\cite{FrK,Kapp, Kap}, with \cite{K2}, where earlier partial
results, later disregarded, are cited. For a~2-page-long summary of Kac's ideas (\cite{K1.5, K2}), see \cite[\S8]{Ld}.
None of these researchers considered deformations with odd
parameters classified, at last, in~ \cite{Ld}. 

\textbf{Lie superalgebras with Cartan matrices}. Kac was the first to attribute analogs of Cartan matrices to certain Lie
superalgebras over $\Cee$, and to observe that one Lie superalgebra can
have several inequivalent Cartan matrices. The definition of Cartan matrices and analogs of Dynkin diagrams was later improved in \cite{CCLL}. For the classification of finite-dimensional Lie superalgebras with indecomposable Cartan matrix in characteristic $p>0$, see \cite{BGL}.

Over $\Cee$, the classification of all Cartan matrices of a~given Lie superalgebra 
is due to Serganova, see \cite{Se2, HS}. Lebedev extended Sergnova's method to the case of characteristic $p>0$, see~ \cite{Leb}, where this place in \cite{BGL} is corrected. For the classification of Lie superalgebras of GLAPG type with indecomposable Cartan matrix, see \cite{Se2, vdL, HS}. 
For the classification of hyperbolic Lie algebras and their superizations, see \cite{CCLL}.

\textbf{Simple loop Lie superalgebras}. For  their conjectural list, see \cite{LSS}. To prove the completeness of this list \`a la \cite{M} is an \textbf{Open problem};\index{Problem, open} for a~partial result, see~ \cite{HS}; for Cartan matrices of loop superalgebras, see \cite{BLS}. 

\textbf{Stringy Lie superalgebras}. For the classification of simple stringy Lie superalgebras (with Laurent coefficients, based on the classification of vectorial, i.e. with polynomial coefficients), see \cite{GLS1}; for the classification of their central extensions, see \cite{KvdL}.
Note that there is a~mysterious difference between presentations of $\fk\fas$
generated by polynomials and those of $\fk\fas^L$ generated by Laurent
polynomials as the superscript indicates, see \cite{GLS4}: we \textbf{conjecture}\index{Conjectuer} that there are infinitely
many relations between ``natural" (from a~certain point of view) generators of $\fk\fas$, whereas $\fk\fas^L$ is finitely presented.

\textbf{Vectorial Lie superalgebras}. It took a~while to distinguish a~reasonable classification problem.
A classification of \textbf{primitive} Lie superalgebras was claimed in
\cite{K1.5, K2}, but was obviously wrong; Theorem~ 10 and even Conjecture~ 1 in \cite{K2} dealing with \textbf{simple} Lie superalgebras, were also wrong. After \cite{L1} appeared, it became clear that the
list of simple vectorial Lie superalgebras ``doubles" 
Cartan's list at least twice: by superdimension and by the ``odd'' counterparts such as $\fle$ of $\fh$, and $\fm$ of
$\fk$, together with $\fsle$ and $\fsm$. Leites dubbed the new algebras by the next letters after $\fk$ ($\fl$ and $\fm$); soon $\fl$ was changed to $\fle$ to eliminate confusions of $\fsle$ with $\fsl$. 
Grozman found $\fb_{\lambda}(n)$ in 1978, see \cite{G}, but was not understood; $\fb_{\lambda}(n)$ was rediscovered in \cite{ALSh}.

In \cite{ALSh}, we claimed that \textbf{one vectorial Lie superalgebra has
finitely many distinct \textit{non-isomorphic}
realizations as a~\textit{W-filtered or W-graded} Lie
superalgebra}; for a~ proof, see Subsection~\ref{SS:7.2}. This phenomenon was totally missed in \cite{K2}, and in
Kac's later papers it was stated
that to describe all filtered vectorial superalgebras is ``probably impossible''. Over fields of positive characteristic, the non-standard gradings were appreciated at approximately the same time as in \cite{ALSh}, but not as much as they deserve, see \cite{KLS}.

Observe that, as early as in \cite{ALSh}, we wrote that the problem ``classify
primitive Lie superalgebras'' is wild, cf. \cite{K2}; for details, see Subsection~
\ref{SS:10.2}. In 1988, we had almost all examples of simple W-graded vectorial Lie
superalgebras over $\Cee$ (see \cite{L1, ALSh, Sh3, LSh1}). The adjective ``exceptional'', first used in
\cite{Sh3}, implies a~claim for classification of ``serial''
examples. Shchepochkina discovered all the exceptional examples
\cite{Sh3, Sh5, Sh14}. 

In {1996--97}, we have explicitly announced the
classification of W-graded simple complex vectorial Lie
superalgebras at various conferences, \textit{e.g.}, at a~seminar of E.~Ivanov, JINR, Dubna (July, 1996), Voronezh Winter
School (Jan. 12--18, 1997), and at the conference in honor of D.~Buchsbaum,
Northeastern Univ., Boston (October 1997), see \cite{LSh2$'$}, as an aside remark
to another, \textbf{still open}, problem\index{Problem, open} ``how to superize representations of quivers?", cf.~ \cite{LSh2, BZ}. The proof
of classification of simple W-graded Lie superalgebras, 
already outlined 
in \cite{LSh1}, was completed by that time, \textbf{except for W-gradings compatible with parity}. 
Getting our confirmation that the classification is announced (November 1996),
Kac, sometimes together with Cheng, considered the case of compatible
W-gradings (\cite{CK}) complementing our proof, and published  partial results
on filtered deformations \cite{CK1, CK1a, CK2, Klc, K2, Kcf,
Kga, Kbj}.
Cantarini helped Cheng and Kac with corrections, cf. \cite{CCK} with \S\ref{Smb} and Subsection~\ref{SS:10.3}. 

The papers \cite{Sh14} and \cite{Sh5} were being refereed for about two years each.
When the details of the talk \cite{Klc} were published in \cite{K3} (accepted 42 days after submission) we saw, to
our amazement, that Kac's method is very similar to ours, published and 
lectured about. However, Kac 
considered \textit{the case of compatible gradings} we left out. The crucial contribution of \cite{K3} to our approach is the proof of existence and description of a~certain remarkable grading of every simple vectorial Lie superalgebra. These particular gradings help to prove the completeness of our list of W-graded simple vectorial Lie superalgebras (conjectural for compatible gradings).

In this article, we classified deformations of simple vectorial Lie superalgebras.

\textbf{Related results} 1) The description of the irreducible continuous representations of simple
vectorial Lie superalgebras was started in \cite{BL3} for the
general series $\fvect(m|n)$ in its standard realization following
the pioneering works by A.~Rudakov (\cite{R1}, \cite{R2}). For
a~review of results on representations for other series of vectorial
Lie superalgebras obtained by 2000, see \cite{GLS3}. Some of these
results can explain and simplify proof of the results of
\cite{BKL1}.

For a~ description of irreducible representations of 11 of the 15
exceptional simple vectorial Lie superalgebras, and the 
related description of invariant differential operators, see
\cite{GLS3}. In \cite{KR1,KR2, KR3}, the representations of 3 of the 15
exceptions are further studied; to explain the discrepancies between answers in \cite{GLS3} and in \cite{CaCK, R3, CaCK1} is an \textbf{open problem}.\index{Problem, open}

2) For the classification of simple vectorial Lie superalgebras over
$\Ree$, see the original and unexpected result \cite{CaKa1}. Serganova classified real and quaternionic forms of finite-dimensional simple Lie superalgebras, see \cite{Se2}.

3) For a~conjectural list of simple stringy (resp., loop) Lie superalgebras over
$\Cee$, see \cite{GLS1} (resp., \cite{LSS}). Serganova classified real forms of loop and stringy superalgebras, see \cite{Se3}; to classify their quaternionic forms remains an \textbf{open problem}.\index{Problem, open}

4) For the classification of simple finite-dimensional Jordan superalgebras over
$\Cee$, see \cite{Kj, HK}. For the same approach in the infinite-dimensional case, see \cite{Lint} and references in it.

\ssec{Comparison of notation} Strangely, not all of the
results from \cite{ALSh} or \cite{Ko1, LSh3} are reproduced in \cite{K3}.
Cheng and Kac considered the parameter $\vec r$ only for exceptional Lie
superalgebras. In this way, they missed 34 types of W-gradings, and most of the occasional isomorphisms, cf. Subsection~\ref{Occasion}; confusing notation obscured ``drop-outs"
listed in Subsubsection~\ref{over}. Kac and his co-authors often changed, probably for non-mathematical reasons, the notation of simple Lie superalgebras previously discovered by other researchers; some of Kac's notation confused finite-dimensional Lie (super)algebras with their root systems. 

We concede that our notation for exceptional simple W-filtered vectorial
superalgebras~ $\cL$ is rather long, though it reflects the way they are constructed and the
geometry preserved ($\fk\fle$ 
is a~subalgebra of
$\fk$ related to $\fs\fle$, and so on). But to write
$\fe(\sdim (\cL/\cL_0)^{*})$, as we once suggested, is to create confusion: $\sdim (\cL/\cL_0)^{*}$ can coincide for non-isomorphic superalgebras $\cL$.

The half-page-long review Zbl 0929.17026 of \cite{K3} is succinct and (unlike MR99m:17006)
void of emotions: ``The main result of the paper is a~ complete classification of simple infinite-dimensional linearly compact Lie superalgebras.\index{Lie (super)algebra, linearly compact} The list of these superalgebras contains 8 series of completed $\Zee$-graded Lie superalgebras, 2 series of filtered deformations, and 6 exceptional Lie superalgebras." However, the correct answer is quite different, see statements in Subsections~\ref{ssres} and \ref{SS:10.3}. 
 In their papers, Cantarini, Cheng, and Kac considered only
a~\textit{part} of problems A to C. 
{\footnotesize
\begin{equation}
\label{1.5.1}
\renewcommand{\arraystretch}{1.4}
\begin{array}{l}
\begin{tabular}{|l|l|}
\hline 
our notation&Kac's notation\cr 
\hline 
$\fvect(m|n;
r)$ &$W(m|n)$\cr
 \hline 
 $\fsvect(m|n; r)$&$S(m|n)$\cr $\fsvect'(1|n;
r)$ &$S(1|n)$\cr 
$\widetilde{\fsvect}(0|m)$, any $m\geq 3$&$\widetilde {S}(0|n)$ for $n\geq 4$ even\cr
\hline
\hline $\fk(2m+1|n; r)$&$K(2m+1|n)$\cr \hline $\fh(2m|n; r)$
&$H(2m|n)$\cr 
$\fh_{\lambda}(2|2)\simeq \fb_{\lambda}(2;
2)$& not singled out\cr
$\fh'(0|n)$&$H(0|n)$\cr \hline
$\widetilde{\fs\fb}_{\mu}(2^{n-1}-1|2^{n-1})$& appears as a~filtered\\
for $n$ even& deform; not recognized \\
&as $\Zee$-graded \cr \hline
\end{tabular}\quad
\renewcommand{\arraystretch}{1.4}
\begin{tabular}{|l|l|}
\hline our notation&Kac's notation\cr
\hline $\fm(n; r)$&$KO(n|n+1)$\cr 
\hline $\fsm(n;
r)$&$SKO(n|n+1)$\cr 
\hline $\fb_{\lambda}(n; r)=\fb_{a,b}(n; r)$&$SKO(n|n+1;
\beta)$ for\cr
for $\lambda:=\frac{2a}{n(a-b)}\in\Cee\Pee^1$&$\beta(=\frac{b}{a}=1-\frac{2}{n\lambda})\in\Cee$\cr
\hline $\fb_{1}'(n; r)$&$SKO(n|n+1; \frac{n-2}n)$ \cr 
\hline
$\fb_{\infty}'(n; r)$& $SKO(n|n+1; 1)$ \cr 
\hline $\fle(n;
r)$&$HO(n|n)$\cr 
\hline $\fsle'(n; r)$&$SHO(n|n)$\cr 
\hline
$\widetilde{\fs\fb}_{\mu}(2^{n-1}|2^{n-1}-1)$& not considered\cr 
for $n$ odd&\\
\hline
singular deformations&not considered\\
of $\fb_{\mu}(n)$& \cr 
\hline
\end{tabular}
\end{array}\end{equation}
}

\normalsize

\subsection*{Acknowledgements} We are thankful to D.~Alekseevsky, J.~Bernstein, P.~Grozman,
I.~Kantor, Yu.~Kotchetkov, A.~Krutov, A.~Lebedev, G.~Post and A.~Rudakov for help.
Thanks are
due to M.~Kontsevich and A.~Lebedev for shrewd comments.

Financial support of D.L. by NFR (1987--96), and of I.Shch.
by INTAS grant 96-0538 is thankfully acknowledged. D.L. also thanks the
Mittag--Leffler Institute and Harvard University (1986-87), IHES (a~
month each summer during 1987--2003); IAS (1989), and the Hebrew
University of Jerusalem and Tel-Aviv University (1990-91), where the 
main results and ideas were presented (1986--2003) and 
MPIM-Bonn, where the first (\cite{LSh1}, 1987) and final, but without classification of deformations, draft (\cite{LSh5}, 2003) were typed, for hospitality.


\begin{thebibliography}{9999}

\bibitem[ALSh]{ALSh}
Alekseevsky D., Leites D., Shchepochkina I., Examples of simple Lie
superalgebras of vector fields. C. R. Acad. Bulg. Sci., v.~34,
no.~9, (1980), 1187--1190 (in Russian)



\bibitem[AM]{AM} 
Avetisyan M. Y., Mkrtchyan R. L., Uniqueness of universal dimensions and configurations of points and lines. Proceedings of Science (PoS) Regio2020 (2021) 005; \url{http://arxiv.org/pdf/2101.10860}

\bibitem[Bat]{Bat}
Batchelor M., The structure of supermanifolds, Trans. Amer. Math. Soc., {\bf 253} (1979), 329--338.

\bibitem[BE]{BE}
Benkart, G., Elduque, A., The Tits construction and the exceptional
simple classical Lie superalgebras. Q. J. Math., v.~54 (2003),
no.~2, 123--137

\bibitem[BGP]{BGP}
Benkart G., Gregory Th., Premet A., The recognition theorem for
graded Lie algebras in prime characteristic. American Mathematical
Society, 2009, 145 pp.



\bibitem[B3]{B3}
Berezin F., \textit{Introduction to superanalysis}. 2nd edition,
revised and expanded. Edited by D.~Leites. With Appendices by
D.~Leites, V.~Shander, I.~Shchepochkina. MCCME, Moscow,
2013. 432 pp. (in Russian) \url{https://staff.math.su.se/mleites/books/berezin-2013-vvedenie-2nd-ed.pdf}

\bibitem[BL1]{BL1}
Bernstein J.N., Leites D.A., 
Integral forms and Stokes's formula on
supermanifolds. Funct. Anal. Appl.,
\textbf{11} (1977), no.~1, 55--56

\bibitem[BL2]{BL2}
Bernstein J.N., Leites D.A., How to integrate differential forms
on supermanifolds. Functional Anal. Appl.,
v.~11 (1977), 219--221.

\bibitem[BL3]{BL3}
Bernstein J.N., Leites D.A., Irreducible representations of Lie
superalgebras of type $W$ and $S$, C.~ R.~ de l'Acad. de Sci. Bulgare,
v.~ 32 (1979), 277--278. (in Russian; English translation in:
Selecta Math. Soviet., v.~3 (1983/84), no.~1, 63--68)

\bibitem[BL4]{BL4}
Bernstein J., Leites D., Invariant differential operators and
irreducible representations of Lie superalgebras of vector fields.
Sel. Math. Sov., v.~1, no.~2, (1981), 143--160

\bibitem[BL5]{BL5}
Bernstein J.N., Leites D.A., Irreducible representations of type $Q$, odd trace and odd
determinant. C.~R.~ Acad. Bulg. Sci.,
\textbf{35} (1982), no.~3, 285--286



\bibitem[BoKN]{BoKN}
Boe B. D., Kujawa J. R., Nakano D. N.,
Cohomology and support varieties for Lie superalgebras. Trans. Amer. Math. Soc. 362 (2010), 6551--6590 

\bibitem[Bbk]{Bbk}
Bourbaki N., {\em Lie groups and Lie algebras. Chapters $4-6$.}
Translated from the 1968 French original by A.~ Pressley.
Elements of Mathematics. 
Springer, Berlin, 2002. xii+300
pp.




\bibitem[BGLL]{BGLL}
Bouarroudj S., Grozman P., Lebedev A., Leites D., Divided power
(co)homology. Presentations of simple finite-dimensional modular Lie
superalgebras with Cartan matrix. Homology, Homotopy and
Applications, v.~12 (2010), no.~1, 237--278;
\url{http://arxiv.org/pdf/0911.0243}

\bibitem[BGLL1]{BGLL1}
Bouarroudj S., Grozman P., Lebedev A., Leites D. (with Appendix by A.~Krutov), Derivations and
central extensions of simple modular Lie algebras and superalgebras. Symmetry, Integrability and Geometry: Methods and Applications (SIGMA) 19 (2023), 59 pages;
\url{http://arxiv.org/pdf/1307.1858}

\bibitem[BGLLS]{BGLLS}
Bouarroudj S., Grozman P., Lebedev A., Leites D., Shchepochkina I.,
Simple vectorial Lie algebras in characteristic~ $2$ and their
superizations. Symmetry, Integrability and Geometry: Methods and Applications (SIGMA) 16 (2020), 089, 101 pages; \url{http://arxiv.org/pdf/1510.07255}

\bibitem[BGL]{BGL}
Bouarroudj S., Grozman P., Leites D., Classification of finite-dimensional modular Lie superalgebras with indecomposable Cartan
matrix. Symmetry, Integrability and Geometry: Methods and
Applications (SIGMA), v.~5 (2009), 060, 63 pages;
\url{http://arxiv.org/pdf/math.RT/0710.5149}

\bibitem[BGL1]{BGL1}
Bouarroudj S., Grozman P., Leites D., 
Deformations of symmetric simple modular Lie superalgebras. Symmetry, Integrability and Geometry: Methods and Applications (SIGMA) 19 (2023), 63 pages; \url{http://arxiv.org/pdf/0807.3054} 


\bibitem[BGLS]{BGLS}
Bouarroudj S., Grozman P., Leites D., Shchepochkina I., Minkowski
superspaces and superstrings as almost real-complex supermanifolds.
Theor. and Mathem. Physics, v.~173(3): 1687--1708 (2012);
\url{http://arxiv.org/pdf/1010.4480}

\bibitem[BKLS]{BKLS}
Bouarroudj S., Krutov A., Leites D., Shchepochkina I., Non-degenerate invariant (super)symmetric bilinear forms on 
simple Lie (super)algebras. Algebras and Repr. Theory, 21(5) (2018), 897--941; \url{http://arxiv.org/pdf/1806.05505}

\bibitem[BLW]{BLW}
Bouarroudj S., Lebedev A., Wagemann F., Deformations of the Lie
algebra $\fo(5)$ in characteristics~$3$ and~$2$. Math. Notes,
(2011), v.~89, no.~5--6, 777--791; \url{http://arxiv.org/pdf/0909.3572}

\bibitem[BLLS]{BLLS}
Bouarroudj S., Lebedev A., Leites D., Shchepochkina I., Lie algebra
deformations in characteristic~2. Math. Research Letters, v.~22
(2015) no.~2, 353--402; \url{http://arxiv.org/pdf/1301.2781}

\bibitem[BLLS1]{BLLS1}
Bouarroudj S., Lebedev A., Leites D., Shchepochkina I.,
Classifications of simple Lie superalgebras in characteristic~ $2$. Internat. Math. Res. Not. (2021) Iss. 1, (2023), 54--94; 
\url{http://arxiv.org/pdf/1407.1695}

\bibitem[BoLe]{BoLe}
Bouarroudj S., Leites D., Invariant differential operators in positive characteristic. J.~Algebra. V. 499, (2018), 281--297:
\url{http://arxiv.org/pdf/1605.09500}

\bibitem[BLS]{BLS}
Bouarroudj S., Leites D., Shang J., Computer-aided study of double extensions of restricted Lie superalgebras preserving the non-degenerate closed 2-forms in characteristic~ 2. Experimental Math. {\bf 31} (2019) no. 2, 676--688 (2022);
\url{http://arxiv.org/pdf/1904.09579}


\bibitem[BZ]{BZ}
Bovdi V. A., Zubkov A.N.:
Super-representations of quivers and related polynomial semi-invariants. Int. J. Algebra Comput. 30(04): (2020) 883--902; \url{http://arxiv.org/pdf/1912.00627}

\bibitem[BKL1]{BKL1}
Boyallian, C., Kac, V., Liberati, J., Irreducible modules over
finite simple Lie conformal superalgebras of type $K$. J. Math.
Phys., v.~51 (2010) 063507, 37 pp.; \url{http://arxiv.org/pdf/1003.4420}

\bibitem[Bu]{Bu}
Buttin C., The\'eorie des op\'erateurs differentiels gradu\'e sur les
formes diff\'erentielles. Bull. Soc. Math. France, v. 102, 1974, no.
1, 49--73



\bibitem[Ca1]{Ca1}
Cantarini N., $\Zee$-graded Lie superalgebras of infinite depth and finite growth. Ann. Sc. Norm. Super. Pisa Cl. Sci. (5) 1 (2002), no. 3, 545--568

\bibitem[Ca2]{Ca2}
Cantarini N., A classification of $\Zee$-graded Lie superalgebras of infinite depth. J. Algebra Appl. 1 (2002), no. 4, 425--449




\bibitem[CCK]{CCK}
Cantarini, N.; Cheng, S.-J.; Kac, V. Errata to: ``Structure of some
$\Zee$-graded Lie superalgebras of vector fields" [Transform.
Groups, v.~4 (1999), no.~2-3, 219--272; MR1712863 (2001b:17037) by
Cheng and Kac]. Transform. Groups, v.~9 (2004), no.~4, 399--400.

\bibitem[CaCK]{CaCK}
Cantarini, N., Caselli, F., Kac, V., Classification of degenerate Verma modules for E(5, 10). Comm. Math. Phys. 385 (2021), no. 2, 963--1005 

\bibitem[CaCK1]{CaCK1}
Cantarini N., Caselli F., Kac V., A Lie conformal superalgebra and duality of representations for $E(4,4)$; \url{http://arxiv.org/pdf/2303.17450}

\bibitem[CaKa1]{CaKa1}
Cantarini N., Kac V. G., Automorphisms and forms of simple infinite
dimensional linearly compact Lie superalgebras, J. of Geom. Methods
in Phys., v.~3, nos. 5 and 6 (2006), 1--23.
\url{http://arxiv.org/pdf/math.QA/0601292}

\bibitem[CaKa2]{CaKa2}
Cantarini N., Kac V. G., Infinite-dimensional primitive linearly
compact Lie superalgebras. Adv. in Math., v.~207, no.~1, (2006)
328--419; \url{http://arxiv.org/pdf/math.QA/0511424}


\bibitem[CaKa3]{CaKa3}
Cantarini N., Kac V. G., Classification of linearly compact simple
rigid superalgebras. IMRN 17 (2010), 3341--3393;
\url{http://arxiv.org/pdf/0909.3100} 

\bibitem[Car]{Car}
Carles R., Un example d'algebre de Lie resolubles rigides, au deuxieme groupe de cohomologie non nul et pour lesquelles l'application quadratique de D.~S.~Rim est
injective. Comptes Rendus Acad. Sci. Paris 300 (1985), no.~14, 467--469

\bibitem[Cart]{Cart}
Cartan \'E., {\OE}uvres complètes. Partie I. (French) [Complete works. Part I] Groups de Lie. [Lie groups] Second edition. \'Editions du Centre National de la Recherche Scientifique (CNRS), Paris, 1984. xxx+1356 pp. 

\bibitem[CCLL]{CCLL}
Chapovalov~D., Chapovalov~M., Lebedev~A., Leites~D. The
classification of almost affine (hyperbolic) Lie superalgebras.
v.~17, Special issue in memory of F.~Berezin, (2010) 103--161;
\url{http://arxiv.org/pdf/0906.1860}

\bibitem[Che]{Che}
Cheng Sh.-J. Differentiable simple Lie superalgebras and
representations of semi-simple Lie superalgebras. J. Algebra, v.~
 173, (1995) no\ 1. 1--43

\bibitem[CK]{CK}
Cheng Shun-Jen, Kac V., A new $N = 6$ superconformal algebra. Comm.
Math. Phys., v.~186 (1997), 219--231


\bibitem[CK1]{CK1}
Cheng Sh.-J., Kac V., Generalized Spencer cohomology and filtered
deformations of $\Zee$-graded Lie superalgebras. Adv. Theor. Math.
Phys., v.~2 (1998), no.~5, 1141--1182;
\url{http://arxiv.org/pdf/math.RT/9805039}

\bibitem[CK1a]{CK1a}
Cheng, Shun-Jen; Kac, V., Addendum: ``Generalized Spencer cohomology
and filtered deformations of ${\Zee}$-graded Lie superalgebras"
[Adv. Theor. Math. Phys., v.~2 (1998), no.~5, 1141--1182; MR1688484
(2000d:17025)]. Adv. Theor. Math. Phys., v.~8 (2004), no.~4,
697--709

\bibitem[CK2]{CK2}
Cheng Sh.-J, Kac V., Structure of some $\Zee$-graded Lie
superalgebras of vector fields, Transformation groups, v.~4, (1999),
219--272

 
\bibitem[Co]{Co}
Coleman A. J., The greatest mathematical paper of all time. Math. Intelligencer 11 (1989), no. 3, 29--38; \url{https://doi.org/10.1007/BF03025189}

\bibitem[CSS]{CSS}
Cushing D. Stewart D., Stagg G., A Prolog assisted search for new simple Lie algebras; \url{http://arxiv.org/pdf/2207.01094}


\bibitem[Del]{Del}
Deligne P., Etingof P.,
Freed D., Jeffrey L., Kazhdan D., Morgan J., Morrison D., Witten E.,
(eds.). \textit{Quantum fields and strings: a~course for
mathematicians}. 
American Mathematical Society,
Providence, RI; Institute for Advanced Study (IAS), Princeton, NJ,
1999. Vol. 1: xxii+723,  Vol. 2: xxiv+727--1501.

\bibitem[DLZ]{DLZ}
Deligne P., Lehrer G. I., Zhang R. B.,
The first fundamental theorem of invariant theory for the orthosymplectic super group. 
Adv. Math. 327 (2018), 4--24. 

\bibitem[DWN]{DWN}
DeWitt B., van Nieuwenhuizen P., Explicit construction of the
exceptional superalgebras $F(4)$ and $G(3)$. J. Math. Phys., v.~23
(1982), no.~10, 1953--1963

\bibitem[Di]{Di}
Dixmier J. \textit{Enveloping algebras}. Revised reprint of the 1977
translation. Graduate Studies in Mathematics, v.~11. American
Mathematical Society, Providence, RI, 1996. xx+379 pp.


\bibitem[Dj1]{Dj1} 
Djokovi\'c, D. \v Z., Semisimplicity of $2$-graded
Lie algebras. II. Illinois J. Math., v.~20 (1976), no.~1, 134--143.

\bibitem[Dj2]{Dj2} 
Djokovi\'c, D. \v Z., Classification of some $2$-graded Lie algebras.
J. Pure Appl. Algebra, v.~7 (1976), no.~2, 217--230.

\bibitem[Dj3]{Dj3} 
Djokovi\'c, D. \v Z., Isomorphism of some simple $2$-graded Lie
algebras. Canad. J. Math., v.~29 (1977), no.~2, 289--294.

\bibitem[DH]{DH} 
Djokovi\'c, D. \v Z., Hochschild, G., Semisimplicity of 2-graded Lie algebras, Illinois J. Math., vol. 20 (1976),
107--123



\bibitem[DGS]{DGS}
Donin J., Gurevich D., Shnider S., Double quantization on some
orbits in the coadjoint representations of simple Lie groups. Comm.
Math. Phys., v.~204 (1999), no.~1, 39--60

\bibitem[DSB]{DSB}
Duplij S., Siegel W., Bagger J. (eds.) \textit{Concise Encyclopedia
of Supersymmetry and Noncommutative Structures in Mathematics and
Physics}, 2nd edition, Berlin, New York: Springer, 2005, 561pp.

\bibitem[Ef]{Ef}
Efetov K.B., \textit{Supersymmetry in Disorder and Chaos}, Cambridge University Press. (1999) xiii + 441.

\bibitem[Eg]{Eg}
Egorov G., \textit{How to superize $\mathfrak{gl}(\infty)$}. In:
Mickelsson J., Pekonen O. (eds.) \textit{Topological and geometrical
methods in field theory}. Proceedings of~the Second International
Symposium held in Turku, May 26--June 1, 1991, World Sci.
Publishing, River Edge, NJ, 1992, 135--146; \url{http://arxiv.org/pdf/2205.07046}

\bibitem[El3]{El3}
Elduque, A. Models of some simple modular Lie superalgebras, Pacific
J. Math., v.~ \textbf{240} (2009), no.~1, 49--83;
\url{http://arxiv.org/pdf/0805.1304}

\bibitem[FaKa]{FaKa}
Fattori D., Kac V. Classification of finite simple Lie conformal
superalgebras. Special issue in celebration of Claudio Procesi's
60th birthday. J. Algebra, v.~258 (2002), no.~1, 23--59.


\bibitem[FF]{FF}
Feigin B., Fuchs D., Cohomologies of Lie groups and Lie algebras.
In: \textit{Lie groups and Lie algebras. II. Discrete subgroups of
Lie groups and cohomologies of Lie groups and Lie algebras}.
Encyclopaedia of Mathematical Sciences, 21. Springer-Verlag, Berlin,
(2000) 125--223

\bibitem[FL]{FL}
Feigin B. L., Leites D. A., New Lie superalgebras of string theories In: M. A. Markov, V. I. Man'ko and A. E. Shabad (eds.) \textit{Group
theoretical methods in physics}, Vol. 1--3 (Zvenigorod, 1982), 623--629 (in Russian), Harwood Academic Publ., Chur, (1985) 


\bibitem[FLS]{FLS}
Feigin B. L., Leites D. A., Serganova V. V., Kac-Moody
superalgebras. In: \textit{Group theoretical methods in physics}, vv.~1--3
(Zvenigorod, 1982), Harwood Academic Publ., Chur, 1985, 631--637

\bibitem[FiFu]{FiFu}
Fialowski, A., Fuchs, D. Singular deformations of Lie algebras.
Example: deformations of the Lie algebra $L_1$. Topics in
singularity theory, Amer. Math. Soc. Transl. Ser. 2, v.~180, Amer.
Math. Soc., Providence, RI, 1997, 77--92;
\url{http://arxiv.org/pdf/q-alg/9706027}

\bibitem[FSS]{Dic}
Frappat L., Sciarrino, A., Sorba, P., \textit{Dictionary on Lie
algebras and superalgebras.} Academic Press, Inc., San Diego, CA,
2000. xxii+410 pp.



\bibitem[FPS]{FPS}
Frenkel I., Penkov I., Serganova V., A categorification of the
boson-fermion correspondence via representation theory of $\fsl(8)$. Comm. Math. Phys. 341 (2016), no. 3, 911--931; \url{http://arxiv.org/pdf/1405.7553}


\bibitem[FrK]{FrK}
Freund P., Kaplansky I., Simple supersymmetries. J. Math. Phys.,
v.~17 (1976), no.~2, 228--231

\bibitem[FuLe]{FuLe}
Fuchs D., Leites D., Cohomology of Lie superalgebras. C. R. Acad.
Bulg. Sci., v.~37 (1984), no.~12, 1595--1596

\bibitem[Fu]{Fu}
Fuks [Fuchs] D., \textit{Cohomology of infinite-dimensional Lie
algebras}, Consultants Bureau, NY, 1986

\bibitem[FLSf]{FLSf}
Fulp R., Lada T., Stasheff J., Noether's variational theorem II and
the BV formalism, Proceedings of the 22nd Winter School ``Geometry
and Physics" (Srn\' i, 2002). Rend. Circ. Mat. Palermo (2) Suppl.
no.~71 (2003), 115--126; \url{http://arxiv.org/pdf/math.QA/0204079}

\bibitem[Ga]{Ga}
Gaw\c{e}dzki K., Supersymmetries --- mathematics of supergeometry,
Ann. Inst. H. Poincar\'e Sect. a~(N.S.), {\bf 27} (1977), 335--366.


\bibitem[GoHo1]{GoHo1}
Golod P., A deformation of the affine Lie algebra $A^{(1)}_1$ and
Hamiltonian systems on the orbits of its subalgebras, In:
\textit{Group-theoretical methods in physics}. Markov M. et al (eds.) Proc.
of the 3rd seminar, Yurmala, 1985, v.1, Moscow, Nauka, 1986,
368--376 (in Russian; English translation by VNU Sci. Press,
Utrecht, (1986).)

\bibitem[GPS]{GPS}
Gomis J., Paris J., Samuel S., Antibracket, Antifields and
Gauge-Theory Quantization, Phys. Rept., v.~ 259 (1995) 1--191

\bibitem[GOV]{GOV}
Gorbatsevich, V. V., Onishchik, A. L., Vinberg, E. B.,
\textit{Foundations of Lie theory and Lie transformation groups}.
Encyclopaedia Math. Sci., v.~20, Springer, Berlin, 1993; MR
95f:22001]. Springer-Verlag, Berlin, 1997. vi+235 pp.

\bibitem[GSW]{GSW}
Green M., Schwarz J., Witten E. \textit{Superstring theory}, vv.1,
2. Second edition. Cambridge Monographs on Mathematical Physics.
Cambridge University Press, Cambridge, 1988. x+470 pp., xii+596 pp.

\bibitem[GN]{GN}
Gritsenko V., Nikulin V., On classification of Lorentsian Kac-Moody
algebras. Russian Mathematical Surveys, 57 (2002) no.5, 921--979;
\url{http://arxiv.org/pdf/math/0201162}


\bibitem[G]{G}
Grozman P., Classification of bilinear invariant operators on tensor
fields. Functional Anal. Appl., v.~14 (1980), no.~2, 127--128; for
proofs, see Grozman P., Invariant bilinear differential operators. Communications in Mathematics, vol 30 no 3 (2022), 129--188; \url{http://arxiv.org/pdf/math/0509562}

\bibitem[Gr]{Gr}
Grozman P., (2013) \url{http://www.equaonline.com/math/SuperLie}

\bibitem[GL1]{GL1}
Grozman P., Leites D., Defining relations for Lie superalgebras with
Cartan matrix. Czech. J. Phys., v.~51, (2001), no.~1, 1--22;
\url{http://arxiv.org/pdf/hep-th/9702073}

\bibitem[GL2]{GL2}
Grozman P., Leites D., Lie superalgebras of supermatrices of complex
size. Their generalizations and related integrable systems. In:
E.~Ram\'irez de Arellano, M.~V.~Shapiro, L.~M.~Tovar and
N.~L.~Vasilevski (eds.) \textit{Proc. Internatnl. Symp. Complex
Analysis and related topics}, Mexico, 1996, Birkh\"auser Verlag,
1999, 73--105; \url{http://arxiv.org/pdf/math.RT/0202177}

\bibitem[GL3]{GL3}
Grozman P., Leites D., Link invariants and Lie superalgebras. J.
Nonlin. Math. Phys., v.~12 (2005), Suppl. 1, 372--379; \url{https://doi.org/10.2991/jnmp.2005 .12.s1.30}

\bibitem[GL4]{GL4}
Grozman P., Leites D., From supergravity to ballbearings. In:
J.~Wess, E.~Ivanov (eds.), \textit{Supersymmetries and quantum
symmetries}, (SQS'97, 22--26 July, 1997), Lecture Notes in Phys.,
v.~ 524, 1999, 58--67.



\bibitem[GL5]{GL5}
Grozman P., Leites D., Structures of $G(2)$ type and nonholonomic distributions in
characteristic $p$. Lett. Math. Phys. \textbf{74}
(2005), no.~3, 229--262; \url{http://arxiv.org/pdf/math.RT/0509400}

\bibitem[GLP]{GLP}
Grozman P., Leites D., Poletaeva E., Defining relations for simple
Lie superalgebras. Lie superalgebras without Cartan matrix. In:
E.~Ivanov et. al. (eds.) \textit{Supersymmetries and Quantum
Symmetries} (SQS'99, 27--31 July, 1999), Dubna, JINR, 2000,
387--396; Homology,
Homotopy and Applications, vol 4 (2), 2002, 259--275;; \url{http://arxiv.org/pdf/math.RT/0202152}


\bibitem[GLS1]{GLS1}
Grozman P., Leites D., Shchepochkina I., Lie superalgebras of string
theories. Acta Mathematica Vietnamica, v.~26, 2001, no.~1, 27--63;
\url{http://arxiv.org/pdf/hep-th/9702120}

\bibitem[GLS2]{GLS2}
Grozman P., Leites D., Shchepochkina I., The analogs of the Riemann
tensor for exceptional structures on supermanifolds. In:
S.~K.~Lando, O.~K.~Sheinman (eds.) \textit{Proc. International
conference ``Fundamental Mathematics Today" (December 26--29, 2001)
in honor of the 10th Anniversary of the Independent University of
Moscow}, IUM, MCCME 2003, 89--109; \url{http://arxiv.org/pdf/math.RT/0509525}

\bibitem[GLS3]{GLS3}
Grozman P., Leites D., Shchepochkina I., Invariant differential
operators on supermanifolds and The Standard Model. In: Olshanetsky
M., Vainshtein A., (eds.) \textit{M.~Marinov memorial volume}, World
Sci., (2002), 508--555; \url{http://arxiv.org/pdf/math.RT/0202193}


\bibitem[GLS4]{GLS4}
Grozman P., Leites D., Shchepochkina I., Defining relations for the
exceptional Lie superalgebras of vector fields pertaining to The Standard Model. In: C.~Duval,
L.~Guieu and V.~Ovsienko (eds.) \textit{The orbit method in geometry
and physics} (A.~A.~Kirillov Festschrift), Progress in Mathematics,
Progress in Mathematics, Birkh\"auser, 2003, 101--146;
\url{http://arxiv.org/pdf/math-ph/0202025}

\bibitem[Gr1]{Gr1}
Gruson C., Finitude de l'homologie de certains modules de dimension
finie sur une super alg\'ebre de Lie. (French) [Finiteness of the
homology of certain finite-dimensional modules over a~Lie
superalgebra] Ann. Inst. Fourier (Grenoble), v.~47 (1997), no.~2,
531--553

\bibitem[Gr2]{Gr2}
Gruson C., Sur la cohomologie des super alg\`ebres de Lie
\'etranges. (French) [Cohomology of strange Lie superalgebras]
Transform. Groups, v.~5 (2000), no.~1, 73--84

\bibitem[Gr3]{Gr3}
Gruson C., Sur l'id\'eal du cone autocommutant des super alg\`ebres
de Lie basiques classiques et \'etranges. (French) [On the ideal of
the self-commuting cone of basic classical and strange Lie
superalgebras] Ann. Inst. Fourier (Grenoble), v.~50 (2000), no.~3,
807--831


\bibitem[GQS]{GQS}
Guillemin V., Quillen D., Sternberg S., The classification of the
complex primitive infinite pseudogroups. Proc. Nat. Acad. Sci.
U.S.A., v.~55, (1966), 687--690;

id., The classification of the irreducible complex algebras of
infinite type. J. Analyse Math., v.~18, (1967), 107--112

\bibitem[Gu]{Gu}
Guillemin, V., Infinite-dimensional primitive Lie algebras. J.
Differential Geometry, v.~4, (1970), 257--282


\bibitem[GoHo2]{GoHo2}
Holod P., Hidden symmetry of the Landau-Lifshitz equation, hierarchy
of higher equations, and the dual equation for an asymmetric chiral
field. Theor. and Mathem. Phys., v.~70, no.~1, (1987) 11--19

\bibitem[HP]{HP} 
Hoyt C., Penkov I., \textit{Classical Lie algebras at infinity}, Springer Monographs in Mathematics, 2022, ISBN 978-3-030-89659-1 

\bibitem[HS]{HS}
Hoyt, C.; Serganova, V. Classification of finite-growth general
Kac-Moody superalgebras. Comm. Algebra, v.~35 (2007), no.~3,
851--874

\bibitem[Hm]{Hm}
Humphreys J, E. \textit{Reflection Groups and Coxeter Groups}.
Cambridge University Press, (1992) 204 pp.


\bibitem[HK]{HK}
H\"ogben L., Kac V., Erratum: ``Classification of simple $Z$-graded
Lie superalgebras and simple Jordan superalgebras" [Comm. Algebra,
v.~5 (1977), no.~13, 1375--1400; MR 58 \#16806] by Kac V. Comm.
Algebra, v.~11 (1983), no.~10, 1155--1156

\bibitem[IoMa]{IoMa}
Iohara K., Mathieu O., A global version of Grozman's theorem, Math. Z. 274 (3) (2013) 955--992; \url{http://arxiv.org/pdf/1204.0695}

\bibitem[J]{J}
Joseph A., Infinite-dimensional Lie algebras in mathematics and
physics. \textit{Groups theoretical methods in Physics} (Proc. Third
Internat. Colloq., Centre Phys. Th\'eor., Marceille, 1974), v.~2,
582--665



\bibitem[K]{K}
Kac, V. G. Simple irreducible graded Lie algebras of finite growth.
Math. USSR-Izv., v.~2 (1968), 1271--1311.

\bibitem[K0]{K0}
Kac V. G., Some algebras related to the quantum field theory, XI-th
All-Union [National USSR] Algebr. Coll., Kishinev, 1971, 140--141 (in Russian)

\bibitem[KfiD]{KfiD}
Kac V.G., Description of filtered Lie algebras with which graded lie
algebras of Cartan type are associated. Mathematics of the
USSR-Izvestiya, v.~ 8(4) (1974), 801--835

\bibitem[K1]{K1}
Kac V. G., Classification of simple Lie superalgebras. Functional Anal. Appl., v.~ 9 (1975), no.~3,
263--265

\bibitem[K1C]{K1C}
Kac V. G., Letter to the editors: ``Classification of simple Lie
superalgebras" (Funkcional. Anal. i Prilo\v zen., v.~ 9 (1975),
no.~3, 91--92). (in Russian) Funct. Anal. Appl., 10:2 (1976), 163

\bibitem[K1.5]{K1.5}
Kac V. G., A sketch of Lie superalgebra theory. Comm. Math. Phys.,
v.~ 53 (1977), no.~1, 31--64

\bibitem[K2]{K2}
Kac V. G., Lie superalgebras. Adv. Math., v.~26, 1977, 8--96 

\bibitem[Kj]{Kj}
Kac V., Classification of simple $Z$-graded Lie superalgebras and
simple Jordan superalgebras, Comm. Algebra, v.~5 (1977), no.~13,
1375--1400; For a~correction, see \cite{HK}.

\bibitem[Kch]{Kch}
Kac V., Characters of typical representations of classical Lie
superalgebras. Comm. Algebra, v.~5 (1977), no.~8, 889--897

\bibitem[Kcl]{Kcl}
Kac V., Representations of classical Lie superalgebras. In: \textit{Differential
geometrical methods in mathematical physics}, II (Proc. Conf., Univ.
Bonn, Bonn, 1977), pp. 597--626, Lecture Notes in Math., 676,
Springer, Berlin, 1978

 \bibitem[Kpr]{Kpr}
Kac V., Some remarks on Lie superalgebra theory (Preprint ca. 1978)



\bibitem[Kb3]{Kb3}
Kac V. G., \textit{Infinite Dimensional Lie Algebras}. 3rd ed.
Cambridge Univ. Press, Cambridge, 1992

\bibitem[Kcf]{Kcf}
Kac V., Superconformal algebras and transitive group actions on
quadrics. Comm. Math. Phys., v.~186 (1997), no.~1, 233--252; For a~correction, see \cite{Kcfe}.



\bibitem[Klc]{Klc}
Kac V., Classification of infinite-dimensional simple linearly
compact Lie superalgebras. Talk at a~Conference dedicated to
S.~Novikov's 60-th birthday, Steklov Math. Inst., Moscow, May, 1998;
ESI-preprint 605 (withdrawn, see \cite{K3})

\bibitem[K3]{K3}
Kac V., Classification of infinite-dimensional simple linearly
compact Lie superalgebras. Adv. Math., v.~139 (1998), no.~1, 1--55;
see MR99m:17006



\bibitem[Kga]{Kga}
Kac V., Classification of infinite-dimensional simple groups of
supersymmetry and quantum field theory. GAFA 2000 (Tel Aviv, 1999).
Geom. Funct. Anal. 2000, Special Volume, Part I, 162--183;
\url{http://arxiv.org/pdf/math.QA/9912235}

\bibitem[Kcfe]{Kcfe}
Kac V., Erratum [Superconformal algebras and transitive group actions on
quadrics]: Comm. Math. Phys., v.~217 (2001), no.~3, 697--698


\bibitem[Kbj]{Kbj}
Kac V. Classification of supersymmetries. Proceedings of the
International Congress of Mathematicians, v.~I (Beijing, 2002),
Higher Ed. Press, Beijing, 2002, 319--344;
\url{http://arxiv.org/pdf/math-ph/0302016}



\bibitem[KvdL]{KvdL}
Kac V. G., van de Leur J. W. On classification of superconformal
algebras. In: S. J. Gates, Jr., C.~R.\ Preitschopf and W. Siegel
(eds.) \textit{Strings '88} Proceedings of the workshop held at the
University of Maryland, College Park, Maryland, May 24--28, 1988.
World Scientific Publishing Co., Inc., Teaneck, NJ, 1989, 77--106

\bibitem[KR1]{KR1}
Kac, V. G., Rudakov A., Representations of the exceptional Lie
superalgebra E(3,6): I. Degeneracy conditions. Transform. Groups,
v.~7 (2002), no.~1, 67--86; \url{http://arxiv.org/pdf/math-ph/0012049};

\bibitem[KR2]{KR2}
Kac, V. G., Rudakov A., Representations of the exceptional Lie superalgebra E(3,6): II.
Four series of degenerate modules. Comm. Math. Phys., v.~222 (2001),
no.~3, 611--661; \url{http://arxiv.org/pdf/math-ph/0012050};

\bibitem[KR3]{KR3}
Kac, V. G., Rudakov A., Complexes of modules over exceptional Lie superalgebras
$E(3,8)$ and $E(5,10)$. Int. Math. Res. Not. 2002, no.~19,
1007--1025


\bibitem[KW]{KW}
Kac V., Wakimoto M., Integrable highest weight modules over affine
superalgebras and number theory. In: \textit{Lie theory and geometry}. Progr.
Math., v.~123, Birkh\"auser Boston, Boston, MA, 1994, 415--456

\bibitem[Kapp]{Kapp}
Kaplansky I., Newsletters: Graded Lie algebras I, II, preprints, Univ. Chicago,
Chicago, Ill., 1975, see
\url{https://web.osu.cz/~Zusmanovich/links/files/kaplansky/}

\bibitem[Kap]{Kap}
Kaplansky I., Superalgebras. Pacific J. Math., v.~86 (1980), no.~1,
93--98

\bibitem[KNV]{KNV}
Karabegov A., Neretin Yu., Voronov Th., Felix Alexandrovich Berezin and his works; 
\url{http://arxiv.org/pdf/1202.3930}

\bibitem[KMLT]{KMLT}
Karpel E., Minlos R. (compilors), Leites D., Tyutin I. (eds.) \textit{Recollections on Felix Alexandrovich
Berezin --- the founder of supermathematics}. MCCME,
Moscow, 382 pp. (in Russian) 

\bibitem[Ki]{Ki}
Killing W., Die Zusammensetzung der stetigen endlichen Transformations-gruppen, Math. Ann. 31 (1888) 252--290, Math. Ann. 33 (1888) 1--48, Math. Ann. 34 (1889) 57--122, Math. Ann. 36 (1890) 161--189

\bibitem[Kir]{Kir}
Kirillov S.A., Sandwich algebras in simple finite-dimensional Lie algebras. Ph.D. thesis, Nizhny Novgorod, (1992) (in Russian); \url{https://web.osu.cz/~Zusmanovich/files}


\bibitem[KT]{KT}
Konstein S. E., Tyutin I., General form of the deformation of the
Poisson superbracket on a~$(2, n)$-dimensional superspace, Theor.
and Mathem. Physics, 155(2): (2008) 734--753;
\url{http://arxiv.org/pdf/hep-th/0610308}

\bibitem[KT1]{KT1}
Konstein S. E., Tyutin I.V., The number of independent traces and
supertraces on symplectic reflection algebras, J. of Nonlinear
Mathematical Physics, 21:3, (2014) 308--335;
\url{http://arxiv.org/pdf/1308.3190}

\bibitem[KV]{KV}
Konstein S. E., Vasiliev M. A. Supertraces on the algebras of
observables of the rational Calogero model with harmonic potential.
J. Math. Phys., v.~37 (1996), no.~6, 2872--2891;

\bibitem[Kon]{Kon}
Kontsevich M., Deformation quantization of Poisson manifolds, Lett.
Math. Phys., v.~66 (2003), 157--216.

\bibitem[Kor]{Kor}
Kornyak V. V., A method for splitting cochain complexes in computing
the cohomology of Lie (super) algebras. Program. Comput. Software,
v.~28 (2002), no.~2, 76--80


\bibitem[KD]{KD}
Kostrikin, A. I., Dzhumadildaev A. S., Modular Lie algebras: new
trends In: Yu. Bahturin (ed.), \textit{Algebra. Proc. of the
International Algebraic Conference on the Occasion of the 90th
Birthday of A.G. Kurosh (May, 1998. Moscow)}, de Gruyter, Berlin
(2000) 181--203


\bibitem[Kz]{Kz}
Koszul J.-L., Connections and splittings of supermanifolds.
Differential Geom. Appl., v.~ 4 (1994), no.~2, 151--161

\bibitem[Ko1]{Ko1}
Kotchetkoff [Kochetkov] Yu., D\'eformations des superalgèbres de Buttin et quantification. (French) [Deformations of Buttin superalgebras and quantization] C. R. Acad. Sci. Paris Sér. I Math. 299 (1984), no. 14, 643--645.


\bibitem[KN]{KN}
Krichever I.M., Novikov S.P., Algebras of Virasoro type, Riemann
surfaces and structures of the theory of solitons . Funktional Anal.
i. Prilozhen. 21, No.2 (1987), 46--63.

\bibitem[KrLe]{KrLe}
Krutov A., Lebedev A., On Gradings Modulo 2 of Simple Lie Algebras
in Characteristic 2. Symmetry, Integrability and Geometry: Methods and Applications SIGMA 14 (2018), 130, 27 pp. \url{http://arxiv.org/pdf/1711.00638}

\bibitem[KLLS]{KLLS}
Krutov A., Lebedev A., Leites D., Shchepochkina I.,
Non-degenerate invariant symmetric bilinear forms on simple Lie superalgebras in characteristic~ $2$. Oberwolfach preprint OWP 2020-02 \url{http://publications.mfo.de/handle/mfo/3697} 


\bibitem[KLLS1]{KLLS1}
Krutov A., Lebedev A., Leites D., Shchepochkina I., Non-degenerate invariant symmetric bilinear forms on simple Lie superalgebras in characteristic $2$. Linear Algebra and Its Appl. (2022) \url{http://arxiv.org/pdf/2102.11653}

\bibitem[KLS]{KLS}
Krutov A., Leites D., Shchepochkina I., Non-integrable distributions with an infinite-dimensional Lie (super)algebra 
of symmetries; \url{http://arxiv.org/pdf/2309.16370} 

\bibitem[Leb]{Leb}
Lebedev A., On reflections of roots in positive characteristic; \url{http://arxiv.org/pdf/2311.02710}



\bibitem[LR]{LR}
Lecomte P.A.B., Roger C., Rigidity of current Lie algebras of
complex simple Lie type, J. London Math. Soc., v.~\textbf{37}
(1988), 232--240

\bibitem[L0]{L0}
Leites D., Spectra of graded commutative rings. Uspehi Matem. Nauk,
v. 30, 1974, no. 3, 209--210 (in Russian)

\bibitem[Lber]{Lber}
Leites D., A certain analog of the determinant. Russian Math.
Surveys, \textbf{35} (1975), no.~3, 156--157

\bibitem[L1]{L1}
Leites D., New Lie superalgebras and mechanics. Soviet Math.
Doklady, v.~18, no. 5, 1977, 1277--1280

\bibitem[L2]{L2}
Leites D., Lie superalgebras. In: \textit{Modern Problems of
Mathematics. Recent developments}, v.~25, VINITI, Moscow, 1984,
3--49 (in Russian; English translation in: J. Soviet Math., v.~30 (6),
1985, 2481--2512)

\bibitem[L3]{L3}
Leites D., \textit{Supermanifold theory}, Karelia Branch of the USSR
Acad. Sci., Petrozavodsk, 1983, 200 pp. (in Russian; expanded in
\cite{L5} and \cite{LSoS})


\bibitem[L5]{L5}
Leites D. (ed.) \textit{Seminar on Supermanifolds}, Reports of
Stockholm University, nos. 1--34 (1987--89); \url{https://staff.math.su.se/mleites/sos.html}

\bibitem[L6]{L6}
Leites D., \textit{Quantization. Supplement $3$}. In: F.~Berezin,
M.~Shubin. \textit{Schr\" odinger equation}, Kluwer, Dordrecht,
1991, 483--522


\bibitem[L7]{L7}
Leites D., Selected problems of supermanifold theory, Duke Math. J.,
v.~54, (1987) no.~2, 649--656

\bibitem[L8]{L8}
Leites D., The Riemann tensor for nonholonomic manifolds. Homology,
Homotopy and Applications, v.~4 (2), (2002), 397--407;
\url{http://arxiv.org/pdf/math.RT/0202213}

\bibitem[LSoS]{LSoS}
Leites D. (ed.) \textit{Seminar on supersymmetry v. $1$. Algebra and
Calculus: Main chapters}, (J.~Bernstein, D.~Leites, V.~Molotkov,
V.~Shander), MCCME, Moscow, 2012, 410 pp (in Russian) 


\bibitem[Lq]{Lq}
Leites D. New simple Lie superalgebras as queerified associative algebras. Adv. Theor. and Mathem. Physics. V. 26, (2022) No. 9, 3189--3206; \url{http://arxiv.org/pdf/2203.06917} 

\bibitem[Ld]{Ld}
Leites D. On odd parameters in geometry. J. Lie theory, Volume 33, no. 4 (2023) 965--1004; \url{http://arxiv.org/pdf/2210.17096} 

\bibitem[Lint]{Lint}
Leites D. On unconventional integration on supermanifolds and cross ratio on classical superspaces. J. of Lie Theory. Volume 33, no. 2 (2023) 527--546; \url{http://arxiv.org/pdf/math.RT/0202194}



\bibitem[LPS]{LPS}
Leites D., Poletaeva E., Serganova V., On Einstein equations on
manifolds and supermanifolds, J. Nonlinear Math. Physics, v.~9,
(2002), no.~4, 394--425; \url{http://arxiv.org/pdf/math.DG/0306209}

\bibitem[LSS]{LSS}
Leites D., Saveliev M., Serganova V., Embeddings of $\mathfrak{osp} (N|2)$ and completely integrable
systems. In: M.~Markov,
V.~Man'ko (eds.) \textit{Proc. International Conf. Group-theoretical
Methods in Physics}, Yurmala, May, 1985. Nauka, Moscow, 1986,
377--394 MR 89h:17042 (English translation: VNU Sci Press, 1987,
255--297).

\bibitem[LS]{LS}
Leites D., Serganova V., Metasymmetry and Volichenko algebras. Phys.
Lett. B 252 (1990), no.~1, 91--96. 

id., Symmetries wider than supersymmetries. In: S.~Duplij and
J.~Wess (eds.) \textit{Noncommutative structures in mathematics and
physics}, Proc. NATO Advanced Research Workshop, Kiev, 2000. Kluwer,
13--30 

\bibitem[LS1]{LS1}
Leites D., Serganova V., Defining relations for classical Lie
superalgebras. I. Superalgebras with Cartan matrix or Dynkin-type
diagram. In: J.~Mickelsson, O.~Pekonen (eds.) \textit{Topological and geometrical methods in field theory} (Turku,
1991), World Sci. Publishing, River Edge, NJ, 1992, 194--201 

\bibitem[LSe]{LSe}
Leites D., Sergeev A., Orthogonal polynomials of discrete variable
and Lie algebras of complex size matrices. In: \textit{Proceedings of
M.~Saveliev memorial conference}, MPI, Bonn, February, 1999; 
Theoret. and Math. Phys., v.~123 (2000), no.~2, 582--608;
MPI-1999-36, 49--70 (\url{http://www.mpim-bonn.mpg.de})

\bibitem[LSh1]{LSh1}
Leites D., Shchepochkina I., Toward classification of simple
vectorial Lie superalgebras. In: \cite{L5}, 31/1988-14, 235--278;

Leites D., Toward classification of classical Lie superalgebras. In:
Nahm W., Chau L. (eds.) \textit{Differential geometric methods in
theoretical physics} (Davis, CA, 1988), NATO Adv. Sci. Inst. Ser. B
Phys., v.~245, Plenum, New York, (1990), 633--651

\bibitem[LSh2]{LSh2}
Leites D., Shchepochkina I., Quivers and Lie superalgebras, Czech.
J. Phys. , v.~47, no.~12, (1997), 1221--1229;

\bibitem[LSh2$'$]{LSh2$'$}
Leites D., Shchepochkina I., Quivers and Lie superalgebras. In: Eisenbud D., Martsinkovsky
A., Weyman J., (eds.) \textit{Commutative Algebra, Representation
theory and Combinatorics. Conference in honor of D.~Buchsbaum},
Northeastern U., Boston, October 18--20, 1997, 67

\bibitem[LSh3]{LSh3}
Leites D., Shchepochkina I., How should the antibracket be
quantized? Theoret. and Math. Phys., v.~126 (2001),
no.~3, 281--306;
\url{http://arxiv.org/pdf/math-ph/0510048}

\bibitem[LSh4]{LSh4}
Leites D., Shchepochkina I., The Howe duality and Lie superalgebras. In: S.~Duplij and J.~Wess (eds.) \textit{Noncommutative structures
in mathematics and physics}, Proc. NATO Advanced Research Workshop,
Kiev, 2000. Kluwer, 2001, 93--112; \url{http://arxiv.org/pdf/math.RT/0202181}.


\bibitem[LSh5]{LSh5}
Leites D., Shchepochkina I., Classification of simple Lie
superalgebras of vector fields, \texttt{preprint MPIM-Bonn 2003-28
(\url{http://www.mpim-bonn.mpg.de/preblob/2178})}

\bibitem[LZ]{LZ}
Lopatkin V., Zusmanovich P., Commutative Lie algebras and commutative cohomology
in characteristic~ $2$. Commun. Contemp. Math. 23 (2021), no. 5, Paper No. 2050046, 20 pp.; \url{http://arxiv.org/pdf/1907.03690}


\bibitem[MaG]{MaG}
Manin Yu., \textit{Gauge field theory and complex geometry}. Second
edition. Springer-Verlag, Berlin, (1997) xii+346 pp.

\bibitem[MaAG]{MaAG}
Manin Yu., \textit{Introduction to the Theory of Schemes}. Springer-Verlag, Berlin, (2018) 221 pp.

\bibitem[M]{M}
Mathieu O., Class of simple graded Lie algebras of finite growth, Inv.
Math., 108, 1992, 455--519

\bibitem[My]{My}
Montgomery S.,
Constructing simple Lie superalgebras from associative graded algebras. J.~Algebra 195 (1997), no. 2, 558--579.

\bibitem[N]{N}
Neretin Yu., Roger C. (traducreur), ``La m\'ethode de la second quantification" de F.A.Berezin. Regards quarante ans plus tard; \url{https://hal.science/hal-00478476} 

\bibitem[O]{O}
Ochiai T., Classification of the finite nonlinear primitive Lie
algebras. Trans. Amer. Math. Soc., v.~124, 1966, 313--322

\bibitem[OV]{OV}
Onishchik A., Vinberg E., \textit{Lie groups and algebraic groups.}
Springer-Verlag, Berlin, 1990. xx+328 pp.

\bibitem[PS]{PS}
Penkov I., Serganova V., Categories of integrable $\fsl(\infty)$-,
$\fo(\infty)$-, $\fsp(\infty)$-modules, In: J.~Adams, B.~ Lian,, and
S.~ Sahi (eds.), ``Representation Theory and Mathematical Physics",
Contemporary Mathematics 557, AMS (2011), 335--357

\bibitem[PU]{PU}
Pinczon G., Ushirobira R., Supertrace and
superquadratic Lie structure on the Weyl algebra, and applications
to formal inverse Weyl transform. Lett. Math. Phys., v.~74(3),
(2005), 263--291.

\bibitem[Po1]{Po1}
Poletaeva E., Structure functions on the usual and exotic symplectic
and periplectic supermanifolds. In: \textit{Differential geometric
methods in theoretical physics} (Rapallo, 1990), Lecture
Notes in Phys., v.~375, Springer, Berlin, 1991, 390--395


\bibitem[Po2]{Po2}
Poletaeva E., Semi-infinite cohomology and superconformal algebras. Ann. Inst. Fourier (Grenoble) 51 (2001), no.3, 745--768; \url{http://arxiv.org/pdf/math/0012195}

\bibitem[Po3]{Po3}
Poletaeva E., Analogs of Riemann and Penrose tensors on supermanifolds;
\url{http://arxiv.org/pdf/math.RT/0510165}

\bibitem[Ri]{Ri}
Richardson R. W., On the rigidity of semi-direct products of Lie algebras.
Pacific J.~Math. 22 (1967), no.~2, 339--344

\bibitem[Ru]{Ru}
Rudakov A. N., Deformations of simple Lie algebras. Mathematics of the USSR-Izvestiya
(1971), 5(5), 1120--1126


\bibitem[R1]{R1}
Rudakov A.N., Irreducible representations of infinite-dimensional
Lie algebras of the Cartan type, Mathematics of the USSR-Izvestia,
v.~8 (1974), 835--866.

\bibitem[R2]{R2}
Rudakov A.N., Irreducible representations of infinite-dimensional
Lie algebras of type $S$ and $H$, Mathematics of the USSR-Izvestia,
v.~9 (1975), 465--480.


\bibitem[R3]{R3}
Rudakov A.N., Morphisms of Verma modules over exceptional Lie superalgebra $E(5,10)$;
\url{http://arxiv.org/pdf/1003.1369}

\bibitem[SV]{SV}
Saveliev M. V., Vershik A. M., Continuum analogues of
contragredient Lie algebras (Lie algebras with a~Cartan operator and
nonlinear dynamical systems). Comm. Math. Phys., v.~126 (1989),
no.~2, 367--378


\bibitem[SNR]{SNR}
Scheunert M., Nahm W., Rittenberg V., Classification of all simple
graded Lie algebras whose Lie algebra is reductive. I, II.
Construction of the exceptional algebras. J. Math. Phys.,
v.~17 (1976), no.~9, 1626--1639, 1640--1644

\bibitem[Sch]{Sch}
Scheunert M., \textit{The theory of Lie superalgebras. An
introduction}. Lecture Notes in Mathematics, 716. Springer, Berlin,
1979. x+271 pp.



\bibitem[Se1]{Se1}
Serganova V., Classification of simple real Lie superalgebras and
symmetric superspaces. Functional
Anal. Appl., v.~17 (1983), no.~3, 200--207

\bibitem[Se2]{Se2}
Serganova V., Automorphisms of simple Lie superalgebras. Izv. Akad.
Nauk SSSR Ser. Mat., v.~48 (1984), no.~3, 585--598; (in Russian)
English translation: Math. USSR-Izv., v.~24 (1985), no.~3, 539--551


\bibitem[Se3]{Se3}
Serganova V., Automorphisms of the Lie superalgebras of string theory, Funct. Anal. Appl., 19:3 (1985), 226--228

\bibitem[Se4]{Se4}
Serganova V., On generalizations of root systems, Comm. in Algebra 24: 13 (1996) 4281--4299

\bibitem[Se5]{Se5}
Serganova V., Characters of irreducible representations of simple
Lie superalgebras. \textit{Proceedings of the International Congress of
Mathematicians}, Vol. II (Berlin, 1998). Doc. Math. 1998, Extra Vol.
II, 583--593 (electronic)



\bibitem[S0]{S0}
Sergeev A. N., Irreducible representations of solvable Lie
superalgebras, In: \cite{L5}, 22/1988, 1--12 and Represent. Theory,
v.~3 (1999), 435--443; \url{http://arxiv.org/pdf/math.RT/9810109}

\bibitem[S1]{S1}
Sergeev~A. The invariant polynomials on simple Lie superalgebras.
Represent. Theory, v.~3 (1999) 250--280.

\bibitem[SeVe]{SeVe}
Sergeev A.\,N., Veselov A.\,P. Grothendieck rings of basic classical
Lie superalgebras. Ann. Math, v.~173, (2011) no.\,2. 663--703;
\url{http://arxiv.org/pdf/0704.2250}

\bibitem[Shan]{Shan}
Shander V.N. Analogues of the Frobenius and Darboux theorems for supermanifolds. Comptes rendus de l'Academie bulgare des sciences. 36 (1983) 309--311

\bibitem[Sha1]{Sha1}
Shapovalov A.V., Finite-dimensional irreducible representations of
Hamiltonian Lie superalgebras, Mathematics of the USSR-Sbornik,
v.~35, (1979), 541--554

\bibitem[Sha2]{Sha2}
Shapovalov A. V., Invariant differential operators and irreducible
representations of finite-dimensional Hamiltonian and Poisson Lie
superalgebras. Serdica Math. J., v.~7, (1981), no.~4, 337--342 (in Russian)


\bibitem[Sh3]{Sh3}
Shchepochkina I., Exceptional simple infinite-dimensional Lie
superalgebras. C. R. Acad. Bulg. Sci., v.~36 (1983),
no.~3, 313--314 (in Russian)

\bibitem[Sh5]{Sh5}
Shchepochkina I., The five simple exceptional Lie superalgebras of
vector fields. Funct. Anal. Appl., v.~33
(1999), no.~3, 208--219 (a preliminary version of [Sh14]);
\url{http://arxiv.org/pdf/hep-th/9702121}

\bibitem[Sh14]{Sh14}
Shchepochkina I., Five exceptional simple Lie superalgebras of
vector fields and their fourteen regradings. Represent. Theory, v.~
3, (1999), 373--415

\bibitem[ShM]{ShM}
Shchepochkina I., Maximal subalgebras of classical Lie
superalgebras. 
In: C.~Duval,
L.~Guieu and V.~Ovsienko (eds.) \textit{The orbit method in geometry
and physics} (A.~A.~Kirillov Festschrift), Progress in Mathematics,
Birkh\"auser, (2003) 445--472 (a~shorter version:
\url{http://arxiv.org/pdf/hep-th/9702122})


\bibitem[Shch]{Shch}
Shchepochkina I., How to realize Lie algebras by vector fields.
Theoret. and Math. Phys., v.~147 (2006) no.~3, 821--838;
\url{http://arxiv.org/pdf/math.RT/0509472}

\bibitem[ShP]{ShP}
Shchepochkina I., Post G., Explicit bracket in an exceptional simple
Lie superalgebra. Internat. J. Algebra Comput., v.~8 (1998), no.~4,
479--495; \url{http://arxiv.org/pdf/physics/9703022}

 \bibitem[ShS]{ShS}
Sherman A., Silberberg L., A queer Kac-Moody construction; \url{http://arxiv.org/pdf/2309.09559}


\bibitem[Shi]{Shi}
Shifman M. (ed.), \textit{Felix Berezin. The Life and Death of the Mastermind of Supermathematics}. World Scientific, (2007) 256 pp.

\bibitem[Shei]{Shei}
Sheinman O., \textit{Current algebras on Riemann surfaces}. De Gruyter
Expositions in Mathematics, 58, Walter de Gruyter GmbH \& Co,
Berlin--Boston, (2012), 150 pp.

\bibitem[Sm2]{Sm2}
Shmelev Yu. The contact algebras are rigid. (in Russian) C. R. Acad.
Bulg. Sci., v.~36 (1983), no.~5, 569--570.

\bibitem[Sk]{Sk}
Skryabin S., The normal shapes of symplectic and contact forms over
algebras of divided powers. (1985) VINITI deposition 8504-B86 (in Russian; for an edited translation into English, see \url{http://arxiv.org/pdf/1906.11496})

\bibitem[Sk91]{Sk91}
Skryabin S. M., Modular Lie algebras of Cartan type over algebraically non-closed fields.I. Communications in Algebra, V. 19 (1991) 6, 1629--1741 

\bibitem[Sk95]{Sk95}
Skryabin S. M., Modular Lie algebras of Cartan type over algebraically non-closed fields.II. Communications in Algebra, V. 23 (1995) 4, 1403--1453


\bibitem[St]{St}
Sternberg S., \textit{Lectures on differential geometry}, Chelsey,
2nd edition, (1985)

\bibitem[S]{S}
Strade H., Simple {L}ie algebras over fields of positive
characteristic. I--III. Structure theory, \textit{de Gruyter Expositions in
Mathematics}, Vol.~38, Walter de Gruyter \& Co., Berlin (2004), viii+540
pages, (2009) vi+385~pages, (2012) x+239~pages.


\bibitem[Sud]{Sud}
Sudbery A., Octonionic description of exceptional Lie superalgebras.
J. Math. Phys., v.~24 (1983), no.~8, 1986--1988

\bibitem[T]{T}
Tanaka N., On differential systems, graded Lie algebras and
pseudogroups. J. Math. Kyoto Univ., v.~10, (1970), 1--82;


id., On the equivalence problems associated with simple graded Lie
algebras. Hokkaido Math. J., v.~8 (1979), no.~1, 23--84

\bibitem[Tyut]{Tyut}
Tyutin I. V., The general form of the $*$-commutator on the Grassmann algebra. Theor.Math.Phys. 128 (2001) 1271--1292; \url{http://arxiv.org/pdf/hep-th/0101068}



\bibitem[vdL]{vdL}
van de Leur J., A classification of contragredient Lie superalgebras
of finite growth. Comm. Algebra, v.~17 (1989), no.~8, 1815--1841

\bibitem[V]{V}
Vey J., D\'eformations du crochet de Poisson d'une vari\'et\'e
symplectique. Comm. Math. Helv., v.~50, (1975), 421--454

\bibitem[Vo]{Vo}
Voronov F.F., Quantization on supermanifolds and an analytic proof of the Atiyah--Singer index theorem, J. Soviet Math., 64:4 (1993), 993--1069



\bibitem[W]{W}
Weisfeiler B., Infinite-dimensional filtered Lie algebras and their
connection with the graded Lie algebras. Funct. Anal. Appl., 2:1 (1968), 88--89

\bibitem[WZ]{WZ}
Wess, J., Zumino, B., Supergauge transformations in four dimensions. Nuclear Physics B. 70 (1) (1974) 39--50

\bibitem[Whe]{Whe}
Wheeler J. A., \emph{Einstein’s Vision}, Springer, Berlin et al (1968), 108 pp.

\bibitem[Y]{Y}
Yamaguchi K., Differential systems associated with simple graded Lie
algebras. \textit{Progress in differential geometry}, Adv. Stud.
Pure Math., v.~22, Math. Soc. Japan, Tokyo, (1993), 413--494



\end{thebibliography}
\end{document}